\documentclass[a4paper,11pt]{amsart}

\usepackage{amscd}
\usepackage{amsmath}
\usepackage{amsthm}
\usepackage{graphicx}
\usepackage{txfonts}
\usepackage{dsfont}
\usepackage{tikz-cd}
\usepackage[titletoc]{appendix}

\usepackage{geometry}
\usepackage{hyperref}
\usepackage{tikz}

\theoremstyle{definition}
\newtheorem{Defn}{Definition}[section]

\theoremstyle{plain}
\newtheorem{Lemma}[Defn]{Lemma}
\newtheorem{prop}[Defn]{Proposition}
\newtheorem{theorem}[Defn]{Theorem}
\newtheorem{Corollary}[Defn]{Corollary}
\newtheorem{main}{Theorem}

\theoremstyle{remark}
\newtheorem{remark}[Defn]{Remark}
\newtheorem{Example}[Defn]{Example}

\newtheorem{notation}[Defn]{Notation}

\DeclareMathOperator{\E}{\Upsilon}
\DeclareMathOperator{\K}{K}
\DeclareMathOperator{\T}{T}
\DeclareMathOperator{\G}{\mathbf{F}}
\DeclareMathOperator{\IL}{J}

\DeclareMathOperator{\alg}{A}
\DeclareMathOperator{\comalg}{\hat{A}}
\DeclareMathOperator{\Lie}{L}
\DeclareMathOperator{\mola}{\mathcal{E}}

\definecolor{dPurple}{rgb}{0.7,0,0.1}

\title{PROPs associated to Lawvere theories and their relation to polynomial functors}
\author{Minkyu Kim}
\date{}
\address{KIAS, Seoul, South Korea}
\email{kimminq@kias.re.kr}
\urladdr{https://minq92.github.io/Gaeul-Autumn/}

\begin{document}

\maketitle

\begin{abstract}
    Several adjunctions between functor categories have been studied and applied previously.
    These include Powell’s adjunction between functor categories on free groups and on the linear PROP associated with the Lie operad, as well as those implicit in the equivalence of Pirashvili between functors on projective modules and modules over wreath products.
    In this paper, for a Lawvere theory $\mathcal{C}$ with a zero object, we construct a natural linear PROP $\tilde{\Phi}_{\mathcal{C}}$ together with a canonical adjunction between $\mathcal{C}$-modules and $\tilde{\Phi}_{\mathcal{C}}$-modules.
    We show that, under mild structural conditions on $\mathcal{C}$, the adjunction is compatible with polynomial degree, inducing a correspondence between polynomial $\mathcal{C}$-modules and truncated $\tilde{\Phi}_{\mathcal{C}}$-modules.
    This framework recovers the constructions of Powell and Pirashvili, while also yielding new examples, including adjunctions for functor categories on modules over a ring, as well as on free nilpotent groups of class $\leq c$ and, more generally, on free $\mathcal{R}$-semisimple groups, where $\mathcal{R}$ is a radical functor for groups.
\end{abstract}

\section{Introduction}

Let $\mathds{k}$ be a unital commutative ring.
For a category $\mathcal{C}$, a $\mathcal{C}$-module is a functor $\mathcal{C} \to \mathds{k}\mbox{-}\mathsf{Mod}$, where $\mathds{k}\mbox{-}\mathsf{Mod}$ denotes the category of $\mathds{k}$-modules.
For a unital $\mathds{k}$-algebra $A$, let $\mathbf{P}_{A}$ denote the category of finitely generated projective left $A$-modules.
If $\mathds{k}$ is a field of characteristic zero, then the category of homogeneous polynomial $\mathbf{P}_A$-modules of degree $d$ is equivalent to the category of right $(A \thicksim \mathfrak{S}_d)$-modules \cite{wreathMac,Mac1995}.
Here, $\mathfrak{S}_d$ is the $d$-th symmetric group and $A \thicksim \mathfrak{S}_d$ denotes the wreath product.
For $\mathds{k} = \mathds{Z}$, an analogous statement holds for the category of $\mathbf{P}_A$-modules of polynomial degree at most $d$, modulo those of degree at most $d-1$ \cite{PolyPira}.

Let $\G^{\mathsf{o}}$ denote the opposite category of finitely generated free groups.
In \cite{powell2024analytic}, Powell constructed, in characteristic zero, an adjunction between the category of $\G^{\mathsf{o}}$-modules and the category of modules over the $\mathds{k}$-linear PROP associated to the Lie operad, where a $\mathds{k}$-linear PROP is a $\mathds{k}$-linear permutative category with object set $\mathds{N}$ whose monoidal structure on objects is addition.
When restricted to {\it analytic} $\G^{\mathsf{o}}$-modules, the adjunction induces an equivalence of categories, compatible with the polynomial degree of $\G^{\mathsf{o}}$-modules.
This result has been applied to the study of various $\G^{\mathsf{o}}$-modules such as outer functors \cite{powell2024outer}, and those arising from the Habiro-Massuyeau category \cite{habiro2021kontsevich,vespa2022functors}.

In \cite{kim2024analytic}, Powell's adjunction is refined over a general ground ring $\mathds{k}$ by introducing a distinguished subclass of analytic $\G^{\mathsf{o}}$-modules, called {\it primitive} $\G^{\mathsf{o}}$-modules.

The purpose of this paper is to develop a systematic method for constructing adjunctions that incorporate the refinement of Powell’s adjunction and the adjunction implicit in the equivalences of Pirashvili and Macdonald.
In addition to recovering these known examples, the present work also yields new adjunctions involving functor categories over Lawvere theories, compatible with polynomial functors.
Recall that a Lawvere theory \cite{Lawvere1963} is a category with finite products, whose objects are $\mathds{N}$, with the products on objects given by addition.

The motivation behind our systematic approach is to formulate the previous results within a unified and general framework.
A long-term goal of this line of research is to obtain analogous results for modules
over the Habiro-Massuyeau category \cite{habiro2021kontsevich}, which will be addressed in subsequent work.

\subsection{Statement of the results} \label{202604101313}

In this paper, we establish a unified description of the above adjunctions in terms of an eigenring-type construction.
The classical eigenring construction producing a subquotient ring was introduced in \cite{Ore1932} to study a formal aspect of differential operators.
In the present work, we develop a generalization of this construction to monads in bicategories.
For our applications, we work within a particular bicategory $\mathsf{Mat}_{\mathds{k}}$ \cite{betti1983variation}, where monads are equivalent to small $\mathds{k}$-linear categories (see Section \ref{202601051252} for a brief review).
For a small category $\mathcal{C}$, let $\mathtt{L}_{\mathcal{C}}$ denote the monad in $\mathsf{Mat}_{\mathds{k}}$ corresponding to the $\mathds{k}$-linearization of $\mathcal{C}$.
The category of $\mathcal{C}$-modules is equivalent to that of $\mathtt{L}_{\mathcal{C}}$-modules.

Let $\mathtt{T}$ be a monad in $\mathsf{Mat}_\mathds{k}$ and $\mathtt{J}$ a left ideal of $\mathtt{T}$.
In this paper, we introduce the following:
\begin{itemize}
    \item a subquotient monad of $\mathtt{T}$, denoted by $\mathrm{E}_{\mathtt{T}}(\mathtt{J})$ and called the {\it eigenmonad};
    \item a canonical $(\mathtt{T}, \mathrm{E}_{\mathtt{T}}(\mathtt{J}))$-bimodule structure on $\mathtt{T}/\mathtt{J}$.
\end{itemize}
These constructions give rise to an adjunction between $\mathtt{T}\mbox{-}\mathsf{Mod}$ and $\mathrm{E}_{\mathtt{T}}(\mathtt{J})\mbox{-}\mathsf{Mod}$, called the {\it eigenmonad adjunction}, where $\mathtt{T}\mbox{-}\mathsf{Mod}$ denotes the category of left $\mathtt{T}$-modules.
There is an analogous adjunction for right modules.

We specialize these constructions to monads arising from Lawvere theories.
Let $\mathcal{C}$ be a Lawvere theory with a zero object.
In this paper, we introduce a distinguished left ideal of the monad $\mathtt{L}_{\mathcal{C}}$, which is denoted by $\mathtt{I}^{\mathsf{pr}}_{\mathcal{C}}$ and called the {\it primitivity ideal} of $\mathcal{C}$.
We write $\Phi_{\mathcal{C}}$ for the eigenmonad $\mathrm{E}_{\mathtt{L}_{\mathcal{C}}} (\mathtt{I}^{\mathsf{pr}}_{\mathcal{C}})$.

We will consider a condition on $\mathcal{C}$, denoted by (ZM*), which imposes a structural constraint on the object $1 \in \mathcal{C}$. 
This is satisfied by all examples considered in this paper.

We now state the first main result of the paper, which associates a canonical PROP to any Lawvere theory with a zero object.
\begin{main}[see Theorem \ref{202604071720} and Corollary \ref{202512090929}] \label{202604091210}
    Let $\mathcal{C}$ be a Lawvere theory with a zero object.
    Then the $\mathds{k}$-linear category $\tilde{\Phi}_\mathcal{C}$ corresponding to the monad $\Phi_{\mathcal{C}}$ carries a natural structure of a $\mathds{k}$-linear PROP.
    Furthermore, if $\mathcal{C}$ satisfies (ZM*), then $\tilde{\Phi}_\mathcal{C} (n,m) = 0$ for $n,m\in\mathds{N}$ such that $n < m$.
\end{main}

Applying the general construction above, we obtain an adjunction relating $\mathtt{L}_{\mathcal{C}}$-modules and $\Phi_\mathcal{C}$-modules:
\begin{main}[see Corollary \ref{202512081855} and Theorem \ref{202509201658}] \label{202512111518}
    Let $\mathcal{C}$ be a Lawvere theory with a zero object.
    The eigenmonad adjunction specializes as follows:
    $$
    \begin{tikzcd}
             \Phi_{\mathcal{C}}\mbox{-}\mathsf{Mod} \arrow[r, shift right=1ex, ""{name=G}] & \mathtt{L}_{\mathcal{C}}\mbox{-}\mathsf{Mod} \arrow[l, shift right=1ex, ""{name=F}]
            \arrow[phantom, from=G, to=F, , "\scriptscriptstyle\boldsymbol{\top}"].
    \end{tikzcd}
    $$
    Furthermore, if $\mathcal{C}$ satisfies (ZM*), then the associated adjoint functors respect the polynomial degree of $\mathcal{C}$-modules, in the sense that for $d\in\mathds{N}$, the adjunction sends polynomial $\mathcal{C}$-modules of degree at most $d$ to $d$-truncated $\Phi_{\mathcal{C}}$-modules, and conversely.
\end{main}

The proofs of the statements under the assumption (ZM*) in Theorems \ref{202604091210} and \ref{202512111518}
rely on a structural interaction between the primitivity ideal and polynomial functor theory established in Section \ref{202512141933}.

The adjoint functors in Theorem \ref{202512111518} are induced by the canonical 
$(\mathtt{L}_{\mathcal{C}}, \Phi_{\mathcal{C}})$-bimodule 
$\mathtt{L}_{\mathcal{C}}/\mathtt{I}^{\mathsf{pr}}_{\mathcal{C}}$.
Central to the study of this bimodule and of the monad $\Phi_\mathcal{C}$ are the following constructions:
\begin{itemize}
    \item The operad $\mathfrak{O}_{\mathcal{C}}$ defined by restricting the PROP $\tilde{\Phi}_\mathcal{C}$ to operations with a single output;
    \item The right $\mathtt{L}_{\mathfrak{S}}$-module $\Psi_\mathcal{C}$, obtained from the canonical $(\mathtt{L}_{\mathcal{C}}, \Phi_{\mathcal{C}})$-bimodule $\mathtt{L}_{\mathcal{C}}/\mathtt{I}^{\mathsf{pr}}_{\mathcal{C}}$, where $\mathfrak{S}$ denotes the category of finite sets and bijections.
\end{itemize}
We will regard $\Psi_{\mathcal{C}}$ as a $\mathds{k}$-linear species (cf. \cite{Loday2012,aguiar2010monoidal}) via the identification
of right $\mathtt{L}_{\mathfrak{S}}$-modules with $\mathds{k}$-linear species.
It turns out that $\Psi_\mathcal{C}$ carries a natural coaugmented comonoid structure with respect to the Day convolution product, which allows one to consider its primitive elements.
In fact, the underlying species of the operad $\mathfrak{O}_\mathcal{C}$ naturally identifies with the primitive elements of $\Psi_\mathcal{C}$; see Theorem \ref{202603221137}.

The construction of $\Psi_\mathcal{C}$ comes equipped with a canonical $(\mathtt{L}_{\mathfrak{S}}, \mathtt{L}_{\mathfrak{S}})$-bimodule map $\alpha_\mathcal{C} : \mu\Psi_\mathcal{C} \to \mathtt{L}_{\mathcal{C}}/\mathtt{I}^{\mathsf{pr}}_{\mathcal{C}}$, where $\mu$ denotes the canonical construction assigning an $(\mathtt{L}_{\mathfrak{S}}, \mathtt{L}_{\mathfrak{S}})$-bimodule to a right $\mathtt{L}_{\mathfrak{S}}$-module.
We show that $\Psi_\mathcal{C}$ reconstructs $\mathtt{L}_{\mathcal{C}}/\mathtt{I}^{\mathsf{pr}}_{\mathcal{C}}$ in the sense that $\alpha_\mathcal{C}$ is an isomorphism.
Likewise, applying $\mu$ to the operad $\mathfrak{O}_\mathcal{C}$, we obtain a canonical monad map $\mola_\mathcal{C} : \mu\mathfrak{O}_\mathcal{C} \to \Phi_\mathcal{C}$.
In all examples considered in this paper, this map is an isomorphism; moreover, it is an isomorphism whenever $\mathds{k}$ is a field.
We do not know whether it is always an isomorphism.
These constructions fit into the following commutative diagram:
$$
\begin{tikzcd}
    \Phi_{\mathcal{C}} \ar[r, hookrightarrow] & \mathtt{L}_{\mathcal{C}}/\mathtt{I}^{\mathsf{pr}}_{\mathcal{C}}\\
    \mu \mathfrak{O}_{\mathcal{C}} \ar[r] \ar[u, "\mola_\mathcal{C}"] & \mu \Psi_{\mathcal{C}}.  \ar[u, "\alpha_\mathcal{C}"]
\end{tikzcd}
$$

\vspace{2mm}
As a universal example, we consider the category $\mathbf{W}$ of finitely generated free monoids. 
The opposite category $\mathbf{W}^{\mathsf{o}}$ is a Lawvere theory satisfying (ZM*).

\begin{main} \label{202604082143}
    The monad map $\mola_{\mathbf{W}^{\mathsf{o}}} : \mu \mathfrak{O}_{\mathbf{W}^{\mathsf{o}}} \to \Phi_{\mathbf{W}^{\mathsf{o}}}$ is an isomorphism.
    Moreover, the operad $\mathfrak{O}_{\mathbf{W}^{\mathsf{o}}}$ is the Lie operad, and the species $\Psi_{\mathbf{W}^{\mathsf{o}}}$ is the linear-order species.
\end{main}

The category $\G^{\mathsf{o}}$ is also a Lawvere theory satisfying (ZM*).
The same conclusions hold for
$\G^{\mathsf o}$, which can be refined as follows:
\begin{main} \label{202604071759}
    The groupification functor $\mathbf{W}^{\mathsf{o}} \to \G^{\mathsf{o}}$ induces an isomorphism of monads $\Phi_{\mathbf{W}^{\mathsf{o}}} \to \Phi_{\G^{\mathsf{o}}}$ and an isomorphism of the associated bimodules $\mathtt{L}_{\mathbf{W}^{\mathsf{o}}} / \mathtt{I}^{\mathsf{pr}}_{\mathbf{W}^{\mathsf{o}}} \to \mathtt{L}_{\G^{\mathsf{o}}} / \mathtt{I}^{\mathsf{pr}}_{\G^{\mathsf{o}}}$.
\end{main}

The adjunction in Theorem \ref{202512111518}, applied to $\G^{\mathsf{o}}$, recovers Powell's adjunction \cite{powell2024analytic,kim2024analytic}.
Theorems \ref{202604082143} and \ref{202604071759} follow from Corollary \ref{202603251841} together with Theorem \ref{202509091931}.

In this paper, we go further by considering several variants of $\G^{\mathsf{o}}$.
For $c \in \mathds{N}^\ast$, let $\mathbf{N}_{c}$ denote the category of finitely generated free nilpotent groups of class $\leq c$.
The opposite category $\mathbf{N}_{c}^{\mathsf{o}}$ is a Lawvere theory satisfying (ZM*).

\begin{main}[see Theorem \ref{202604021256}] \label{202512111051}
    The monad map $\mola_{\mathbf{N}_{c}^{\mathsf{o}}} : \mu\mathfrak{O}_{\mathbf{N}_{c}^{\mathsf{o}}} \to \Phi_{\mathbf{N}_{c}^{\mathsf{o}}}$ is an isomorphism.
    Moreover, the operad $\mathfrak{O}_{\mathbf{N}_{c}^{\mathsf{o}}}$ is the operad governing nilpotent Lie algebras of class $\leq c$, while the species $\Psi_{\mathbf{N}_{c}^{\mathsf{o}}}$ can be viewed as a nilpotent truncation of the linear-order species.
\end{main}

In view of the above theorem, the following result is a refinement of Powell's equivalence \cite{powell2024analytic}:
\begin{main}[see Corollary \ref{202604051720}] \label{202604111753}
    If the ground ring $\mathds{k}$ is a field of characteristic zero, then the adjunction in Theorem \ref{202512111518} applied to $\mathcal{C} = \mathbf{N}_{c}^{\mathsf{o}}$ induces an equivalence of the categories of $\mu\mathfrak{O}_{\mathbf{N}_{c}^{\mathsf{o}}}$-modules and of analytic $\mathbf{N}_{c}^{\mathsf{o}}$-modules.
\end{main}

As another application, we consider the category of free right $R$-modules of finite rank, which we denote by $\mathbf{M}_{R}$, where $R$ is a unital ring.
This is a Lawvere theory satisfying (ZM*).
\begin{main}[see Theorem \ref{202603241057} and Corollary \ref{202604051802}] \label{202512111638}
    The monad map $\mola_{\mathbf{M}_{R}} : \mu\mathfrak{O}_{\mathbf{M}_{R}} \to \Phi_{\mathbf{M}_{R}}$ is an isomorphism.
    Moreover, the operad $\mathfrak{O}_{\mathbf{M}_{R}}$ is concentrated in arity $1$ with value $R_\mathds{k} {:=} \mathds{k} \otimes_{\mathds{Z}} R$, and the species $\Psi_{\mathbf{M}_{R}}$ is given by $\Psi_{\mathbf{M}_{R}}(X) = R_\mathds{k}^{\otimes X}$.
\end{main}

Using the identifications in Theorem \ref{202512111638}, the adjunction in Theorem \ref{202512111518}, applied to $\mathbf{M}_{R}$, induces the equivalence given by Pirashvili \cite{PolyPira} when restricted to polynomial $\mathbf{M}_{R}$-modules.

Theorem \ref{202512111051} applied to $c=1$ agrees with Theorem \ref{202512111638} applied to $R = \mathds{Z}$, under the isomorphism $\mathbf{M}_{\mathds{Z}} \cong \mathbf{N}_1^{\mathsf{o}}$.

The results in Theorems \ref{202604082143}, \ref{202604071759}, \ref{202512111051} and \ref{202512111638} are summarized in the following table.
\begin{table}[ht]
\centering
\begin{tabular}{c|c|c}
\hline
$\mathcal C$
&
the operad $\mathfrak{O}_{\mathcal{C}}$
&
the species $\Psi_{\mathcal{C}}$
\\
\hline
$\mathbf{W}^{\mathsf{o}}, \G^{\mathsf{o}}$
&
Lie operad
&
linear-order species $\mathsf{Lin}$
\\
\hline
$\mathbf{N}_c^{\mathsf{o}}$
&
$c$-nilpotent Lie operad
&
$c$-nilpotent truncation of $\mathsf{Lin}$
\\
\hline
$\mathbf{M}_R$
&
$R_\mathds{k}$ in arity $1$
&
tensor-power species $R_\mathds{k}^{\otimes}$
\\
\hline
\end{tabular}
\end{table}

\subsection{Position of this paper in the series} \label{202601081338}

This is the second paper in our series, following \cite{kim2025poly}.
Except for Section \ref{202512141933}, all results and proofs in this paper are independent of that work.
In the following, we summarize the role of the present paper in this series.
\begin{itemize}
    \item The general framework introduced in this paper--namely the eigenmonad adjunction--extends the framework in the first paper \cite{kim2025poly}.
    \item One main construction of \cite{kim2025poly} is the polynomiality ideal. In Section \ref{202512141933}, we investigate its relationship with the primitivity ideal introduced in the present paper.
    \item This paper also provides part of the foundation for a subsequent paper where we will study functors from the Habiro-Massuyeau category.
    One closely related result, announced earlier as noted in \cite{katada2025modules}, was obtained independently in \cite{katada2025modules} using a different approach.
\end{itemize}

\begin{remark}
    An earlier version of this work was posted on arXiv.
    Some material from that version has since been developed into a separate
    paper \cite{kim2025poly}, while the present paper substantially revises and strengthens the remaining part.
\end{remark}

\subsection{Outline}
This paper is divided into three parts.
In Part \ref{202512151906}, we present a general framework for formulating the main results of this paper.
In Part \ref{202512151908}, we introduce the primitivity ideal of a Lawvere theory and study the associated eigenmonad adjunction.
In Part \ref{202604131423}, applying the previous results, we study the eigenmonad adjunctions associated with the category of free nilpotent groups of class $\leq c$ and the category of free $R$-modules.

Part \ref{202512151906} is organized as follows.
In Section \ref{202601051252}, we give an overview of the bicategory $\mathsf{Mat}_{\mathds{k}}$, in particular its monad theory, formulating an analogue of the tensor-hom adjunction.
In Section \ref{202410161744}, we introduce the general framework involving eigenmonads and the eigenmonad adjunction.

Part \ref{202512151908} is outlined as follows.
In Section \ref{202509031755}, we review some preliminaries on Lawvere theories and introduce the condition (ZM*).
In Section \ref{202512162136}, we introduce the notion of a primitivity ideal $\mathtt{I}^{\mathsf{pr}}_{\mathcal{C}}$ for a Lawvere theory $\mathcal{C}$ with a zero object, and prove some fundamental properties for the associated eigenmonad.
In Section \ref{202512181624}, we give some preliminaries for studying right $\mathtt{L}_{\mathfrak{S}}$-modules giving an overview of linear species.
In Section \ref{202606121546}, we introduce a canonical $(\mathtt{L}_\mathfrak{S},\mathtt{L}_\mathfrak{S})$-bimodule associated to a right $\mathtt{L}_\mathfrak{S}$-module.
In Section \ref{202604131411}, we give a structural theorem for the quotient $\mathtt{L}_{\mathcal{C}}/\mathtt{I}^{\mathsf{pr}}_{\mathcal{C}}$.
In Section \ref{202603211037}, we study the PROP corresponding to the monad $\Phi_\mathcal{C}$ by using a natural comonoid species induced by the quotient $\mathtt{L}_{\mathcal{C}}/\mathtt{I}^{\mathsf{pr}}_{\mathcal{C}}$.
As a consequence, we show that the eigenmonad adjunction for $\G^{\mathsf{o}}$ recovers Powell's adjunction \cite{powell2024analytic}.
In Section \ref{202509201216}, we construct a canonical map with codomain $\mathtt{L}_{\mathcal{C}}/\mathtt{I}^{\mathsf{pr}}_{\mathcal{C}}$, and prove that, under the condition (ZM*), it is an epimorphism.
In Section \ref{202512141933}, we give several key interactions between the primitivity ideal and the polynomiality ideal introduced in \cite{kim2025poly}.

Part \ref{202604131423} is structured as follows.
In Section \ref{202512171132}, we study the Lawvere theory arising from free $\mathcal{R}$-semisimple groups, where $\mathcal{R}$ is a radical functor for groups, within the framework developed in Part \ref{202512151908}.
In Section \ref{202410161800}, we study the category $\mathbf{M}_{R}$ within the same framework.

In appendices, we give several deferred proofs.
In Appendix \ref{202604101548}, we give proofs of Lemmas \ref{202604011630} and \ref{202603321744}.
In Appendix \ref{202512031342}, we establish Lemmas \ref{202604011516} and \ref{202604021122}.

\section*{Notation}

\begin{itemize}
    \item $\mathds{N}$ and $\mathds{N}^\ast$ respectively denote the set of nonnegative integers and that of positive integers.
    \item $\mathds{Z}$ denotes the ring of integers.
    \item ${\bf n}$ denotes the set of positive integers less than or equal to $n$.
    \item For a category $\mathcal{C}$, $\mathcal{C}^{\mathsf{o}}$ denotes the opposite category of $\mathcal{C}$.
\end{itemize}

{\Large \part{General theory} \label{202512151906}}

In this part, we present a general framework for formulating the main results of this paper.
We introduce and investigate the eigenmonad construction in the bicategory of matrices which will be recalled in Section \ref{202512141926}.
We show that this construction naturally induces an adjunction of module categories, which we call the eigenmonad adjunction.
This leads to the explicit statement given in the beginning of Section \ref{202604101313}.

\vspace{3mm}
\section{Basic constructions} \label{202601051252}

In this section, we present a framework that will be used to formulate the main results.

\subsection{The bicategory of matrices} \label{202512141926}

Throughout this paper, we fix a unital commutative ring $\mathds{k}$, unless specified otherwise.
We briefly recall a bicategory $\mathsf{Mat}_{\mathds{k}}$ and fix our notation as follows.
This is the bicategory of matrices introduced in \cite{betti1983variation}; for a general treatment, we refer the reader to that paper.
\begin{itemize}
    \item Objects: sets $\mathcal{X}, \mathcal{Y}, \mathcal{Z}, \cdots$.
    
    \item The category of 1-morphisms from $\mathcal{X}$ to $\mathcal{Y}$ is defined as follows:
    \begin{itemize}
        \item 1-morphisms: $(\mathcal{Y} \times \mathcal{X})$-indexed $\mathds{k}$-modules 
        $\mathtt{F} = \{ \mathtt{F}(Y,X) \}_{Y \in \mathcal{Y}, X \in \mathcal{X}}$, 
        viewed as matrices of the $\mathds{k}$-modules with rows indexed by $\mathcal{Y}$ and columns indexed by $\mathcal{X}$.
        
        \item 2-morphisms: maps between such indexed modules 
        that preserve the indices, denoted by 
        $\hom_{\mathcal{Y}, \mathcal{X}}(\mathtt{F}, \mathtt{F}')$.
    \end{itemize}
\end{itemize}
We use the notation $\mathtt{F}:\mathcal{X} \rightsquigarrow \mathcal{Y}$ to denote a 1-morphism from $\mathcal{X}$ to $\mathcal{Y}$.
For 1-morphisms $\mathtt{F}: \mathcal{Y} \rightsquigarrow \mathcal{Z}$ and $\mathtt{G} : \mathcal{X} \rightsquigarrow \mathcal{Y}$, we denote by $\mathtt{F} \otimes \mathtt{G} : \mathcal{X} \rightsquigarrow \mathcal{Z}$ the composition that is defined by using the matrix product formula:
$$
(\mathtt{F} \otimes \mathtt{G}) (Z,X) {:=} \bigoplus_{Y \in \mathcal{Y}} \mathtt{F} (Z,Y) \otimes \mathtt{G} (Y,X) .
$$
Let $\mathds{I}_{\mathcal{X}} : \mathcal{X} \rightsquigarrow \mathcal{X}$ denote the 1-identity on $\mathcal{X}$ given by the diagonal matrix consisting of $\mathds{k}$, i.e. $\mathds{I}_{\mathcal{X}} (X,X) {:=} \mathds{k}$ and $\mathds{I}_{\mathcal{X}} (Y,X) {:=} 0$ for $Y \neq X$.

The category of 1-morphisms from $\mathcal{X}$ to $\mathcal{Y}$ forms an abelian category.
Consequently, notions such as kernels, cokernels, monomorphisms and epimorphisms have their usual meanings.

\begin{Defn}
    For a $1$-morphism $\mathtt{F}:\mathcal{X} \rightsquigarrow \mathcal{Y}$, a {\it submodule} $\mathtt{G} \subset \mathtt{F}$ is a $(\mathcal{Y} \times \mathcal{X})$-indexed $\mathds{k}$-submodule; that is, $\mathtt{G} (Y,X) \subset \mathtt{F}(Y,X)$ for all $X$ and $Y$.
    Given such a submodule, we define $\mathtt{F}/\mathtt{G} :\mathcal{X} \rightsquigarrow \mathcal{Y}$ to be the quotient module determined by $\mathtt{F}/\mathtt{G} (Y,X) {:=} \mathtt{F}(Y,X)/\mathtt{G}(Y,X)$.
\end{Defn}

In this paper, we mainly study 1-endomorphisms on $\mathds{N}$.
As in the case of matrices, we introduce the following terminologies:
\begin{Defn} \label{202512071707}
    We say that a 1-morphism $\mathtt{F} : \mathds{N} \rightsquigarrow \mathds{N}$ is
\begin{itemize}
    \item {\it upper triangular} if $\mathtt{F}(m,n) \cong 0$ whenever $m > n$, and
    \item {\it diagonal} if $\mathtt{F}(m,n) \cong 0$ whenever $m \neq n$.
\end{itemize}
\end{Defn}

The composition of $1$-morphisms in $\mathsf{Mat}_\mathds{k}$ from any side is {\it closed}.
In other words, the composition on either side admits a right adjoint in the following sense:
\begin{Defn}[Internal Hom] \label{202410121637}
    For a $1$-morphism $\mathtt{E} : \mathcal{Y} \rightsquigarrow \mathcal{X}$, we define $\mathrm{Hom^{L}}_{\mathcal{X}} ( \mathtt{E} , -)$ to be a right adjoint to $\mathtt{E} \otimes (-)$: for an object $\mathcal{X}$, $1$-morphisms $\mathtt{G} : \mathcal{Z} \rightsquigarrow \mathcal{X}$ and $\mathtt{F} : \mathcal{Z} \rightsquigarrow \mathcal{Y}$, the $1$-morphism $\mathrm{Hom^{L}}_{\mathcal{Z}} ( \mathtt{E} , \mathtt{G}) : \mathcal{Z} \rightsquigarrow \mathcal{Y}$ satisfies a natural isomorphism:
    \begin{align*}
        \hom_{\mathcal{X},\mathcal{Z}} ( \mathtt{E} \otimes \mathtt{F} , \mathtt{G} ) \cong  \hom_{\mathcal{Y},\mathcal{Z}} ( \mathtt{F} , \mathrm{Hom^{L}}_{\mathcal{X}} ( \mathtt{E}, \mathtt{G} ) ).
    \end{align*}
    Analogously, we introduce $\mathrm{Hom^{R}}_{\mathcal{Z}} ( \mathtt{F} , -)$ as a right adjoint to $(-) \otimes \mathtt{F}$.
\end{Defn}

\begin{remark}
    Note that we have an explicit description of $\mathrm{Hom^{L}}_{\mathcal{X}} ( \mathtt{E}, \mathtt{G} )$.
    Indeed, for $Y \in \mathcal{Y}, Z \in \mathcal{Z}$, the $\mathds{k}$-module $\left( \mathrm{Hom^{L}}_{\mathcal{X}} ( \mathtt{E}, \mathtt{G} ) \right) (Y,Z)$ consists of index-preserving maps $\mathtt{E}(-,Y) \to \mathtt{G} (-,Z)$.
    Likewise, $\left( \mathrm{Hom^{R}}_{\mathcal{Z}} ( \mathtt{F}, \mathtt{G} ) \right) (X,Y)$ is a $\mathds{k}$-module of index-preserving maps from $\mathtt{F} (Y, -)$ to $\mathtt{G} ( X,-)$.
\end{remark}

For $1$-morphisms $\mathtt{F}_{i} : \mathcal{Y}_{i} \rightsquigarrow \mathcal{X}, ~i \in \{1,2,3\}$ in $\mathsf{Mat}_{\mathds{k}}$, the composition of maps induce the following associative pairing:
\begin{align} \label{202410121142}
    \mathrm{Hom^{L}}_{\mathcal{X}} ( \mathtt{F}_1, \mathtt{F}_2 ) \otimes \mathrm{Hom^{L}}_{\mathcal{X}} ( \mathtt{F}_2, \mathtt{F}_3) \to \mathrm{Hom^{L}}_{\mathcal{X}} ( \mathtt{F}_1, \mathtt{F}_3 ) .
\end{align}

\subsection{Monads in $\mathsf{Mat}_{\mathds{k}}$} \label{202606121310}

In this section, we briefly recall the notion of a monad in the bicategory $\mathsf{Mat}_{\mathds{k}}$.

Let $\mathcal{X}$ be an object in $\mathsf{Mat}_{\mathds{k}}$.
A monad in $\mathsf{Mat}_{\mathds{k}}$ on $\mathcal{X}$ consists of a 1-morphism $\mathtt{T} : \mathcal{X} \rightsquigarrow\mathcal{X}$ equipped with 2-morphisms $\nabla : \mathtt{T} \otimes \mathtt{T} \to \mathtt{T}$, the {\it multiplication}, and $\eta : \mathds{I}_{\mathcal{X}} \to \mathtt{T}$, the {\it unit}, satisfying the axioms of a monoid object.
In this paper, a monad in $\mathsf{Mat}_{\mathds{k}}$ will be simply called a monad.

\begin{notation}
    We denote by $1_X \in \mathtt{T}(X,X)$ the unit evaluated at $X \in \mathcal{X}$, and by $f \circ g \in \mathtt{T}(X,Z)$ the multiplication of $f \in \mathtt{T}(X,Y), g\in\mathtt{T}(Y,Z)$.
\end{notation}

There is an equivalence between the category of monads on $\mathcal{X}$ and the category of $\mathds{k}$-linear categories with object set $\mathcal{X}$.
Given a monad $\mathtt{T}$, one obtains a category $\tilde{\mathtt{T}}$ by defining the morphism set $$\tilde{\mathtt{T}} (X,Y) {:=} \mathtt{T} ( Y,X), \quad X,Y\in\mathcal{X} .$$
The following is of particular interest:
\begin{Example} \label{202604091040}
    For a category $\mathcal{C}$ with object set $\mathcal{X}$, let $\mathtt{L}_{\mathcal{C}}$ denote the monad on $\mathcal{X}$, induced by the $\mathds{k}$-linearization of $\mathcal{C}$.
    In particular, $\mathtt{L}_{\mathcal{C}}(Y,X) = \mathds{k} [\mathcal{C}(X,Y)]$.
    For instance, let $\mathfrak{S}$ be the category whose objects are $\mathds{N}$, with $\mathfrak{S}(n,m) = \mathfrak{S}_n$ the $n$-th symmetric group if $n = m$ and $\mathfrak{S}(n,m) = \emptyset$ if $n \neq m$.
    The induced monad $\mathtt{L}_{\mathfrak{S}}$ is diagonal.
    Further examples arise from Lawvere theories which will be recalled in Section \ref{202509031755}.
\end{Example}

\subsection{Modules in $\mathsf{Mat}_{\mathds{k}}$}

In this section, we recall the notion of a module in $\mathsf{Mat}_{\mathds{k}}$ and present an analogue of tensor-hom adjunction for 1-morphisms in $\mathsf{Mat}_{\mathds{k}}$.

Let $\mathcal{X},\mathcal{Y}$ be objects of $\mathsf{Mat}_{\mathds{k}}$.
Let $\mathtt{T}$ and $\mathtt{S}$ be monads on $\mathcal{X}$ and $\mathcal{Y}$ respectively.
Let $\mathtt{M} : \mathcal{Y} \rightsquigarrow \mathcal{X}$ be a 1-morphism.
A {\it left $\mathtt{T}$-action} on $\mathtt{M}$ consists of a 2-morphism $\phi : \mathtt{T} \otimes \mathtt{M} \to \mathtt{M}$ satisfying the axioms of a left action, i.e. $\phi (t \otimes \phi (t^\prime \otimes v)) = \phi ((t \circ t^\prime) \otimes v)$ and $\phi ( 1_{X_1} \otimes v) =v$, where $t \in \mathtt{T} (X_3,X_2), t^\prime \in \mathtt{T}(X_2,X_1)$ and $v \in \mathtt{M} (X_1,Y)$.
Throughout this paper, we use the notation $$t \rhd v = \phi (t \otimes v).$$
One analogously defines a {\it right $\mathtt{S}$-action} and the notation $v \lhd s$.

A {\it left $\mathtt{T}$-module with domain $\mathcal{Y}$} is a 1-morphism $\mathtt{M} : \mathcal{Y} \rightsquigarrow \mathcal{X}$ with a left $\mathtt{T}$-action.
In this paper, a {\it left $\mathtt{T}$-module} is a left $\mathtt{T}$-module with domain $\ast$, a singleton.
For a left $\mathtt{T}$-module $\mathtt{M}$, we write $\mathtt{M}(X) := \mathtt{M}(X, \ast)$.
One defines {\it right $\mathtt{S}$-modules} ({\it with codomain $\mathcal{X}$}) and adopts the same conventions analogously.

\begin{Example} \label{202604062227}
    Let $V$ be a $\mathds{k}$-module.
    We denote by $V^{\otimes}$ the right $\mathtt{L}_{\mathfrak{S}}$-module such that $V^{\otimes} (0) {:=} \mathds{k}$ and $V^{\otimes} (n) {:=} V^{\otimes n}$ equipped with the place permutation $\mathfrak{S}_n$-action.
\end{Example}
\begin{notation}
    We denote by $\mathtt{T}\mbox{-}\mathsf{Mod}$ the category of left $\mathtt{T}$-modules and $\mathtt{T}$-homomorphisms.
    Likewise, $\mathsf{Mod}\mbox{-}\mathtt{T}$ denotes the category of right $\mathtt{T}$-modules.
\end{notation}

A {\it $(\mathtt{T},\mathtt{S})$-bimodule} consists of a 1-morphism $\mathtt{M} : \mathcal{Y} \rightsquigarrow \mathcal{X}$ equipped with a left $\mathtt{T}$-action and a right $\mathtt{S}$-action which commute with each other.

There is an equivalence between the category of $(\mathtt{T},\mathtt{S})$-bimodules and the category of $\mathds{k}$-bilinear functors $\tilde{\mathtt{T}} \times \tilde{\mathtt{S}}^{\mathsf{o}} \to \mathds{k}\mbox{-}\mathsf{Mod}$.
Via this equivalence, most of the results of this paper can be translated into corresponding results for functors with values in $\mathds{k}\mbox{-}\mathsf{Mod}$.

\begin{Defn}[Balanced Tensor product]
\label{202409261711}
    Consider a right $\mathtt{S}$-module $\mathtt{M}$ with codomain $\mathcal{X}$ and a left $\mathtt{S}$-module $\mathtt{N}$ with domain $\mathcal{Z}$.
    We define the morphism $\mathtt{M} \otimes_{\mathtt{S}} \mathtt{N} : \mathcal{Z} \rightsquigarrow \mathcal{X}$ in $\mathsf{Mat}_{\mathds{k}}$ to be the coequalizer, in $(\mathcal{X}\times\mathcal{Z})$-indexed modules, of
    $$
    \begin{tikzcd}
        \mathtt{M} \otimes \mathtt{S} \otimes \mathtt{N} \ar[r, shift right, "\mathrm{id}_{\mathtt{M}} \otimes \lhd"'] \ar[r, shift left, "\rhd \otimes \mathrm{id}_{\mathtt{N}}"] & \mathtt{M} \otimes \mathtt{N}
    \end{tikzcd}
    $$
    Here, $\lhd : \mathtt{M} \otimes \mathtt{S} \to \mathtt{M}$ is the right $\mathtt{S}$-action on $\mathtt{M}$ and $\rhd : \mathtt{S}\otimes\mathtt{N} \to \mathtt{N}$ is the left $\mathtt{S}$-action on $\mathtt{N}$.
\end{Defn}
Let $\mathtt{T},\mathtt{R}$ be monads in $\mathsf{Mat}_{\mathds{k}}$ on objects $\mathcal{X}$ and $\mathcal{Z}$ respectively.
We note that, for a $(\mathtt{T},\mathtt{S})$-bimodule $\mathtt{M}$ and an $(\mathtt{S},\mathtt{R})$-bimodule $\mathtt{N}$, $\mathtt{M} \otimes_{\mathtt{E}} \mathtt{N}$ admits the $(\mathtt{T},\mathtt{R})$-bimodule structure inherited from $\mathtt{M}$ and $\mathtt{N}$.
    
Extending the internal hom (see Definition \ref{202410121637}) to general modules, we introduce the following:
\begin{Defn}
\label{202409261651}
    Let $\mathtt{E}$ and $\mathtt{G}$ be left $\mathtt{T}$-modules with domains $\mathcal{Y}$ and $\mathcal{Z}$ respectively.
    We define the morphism $\mathrm{Hom^{L}}_{\mathtt{T}} ( \mathtt{E}, \mathtt{G} ) : \mathcal{Z} \rightsquigarrow \mathcal{Y}$ in $\mathsf{Mat}_{\mathds{k}}$ to be the $(\mathcal{Y} \times \mathcal{Z})$-indexed submodule of $\mathrm{Hom^{L}}_{\mathcal{X}} ( \mathtt{E}, \mathtt{G} )$ consisting of maps preserving the $\mathtt{T}$-actions. 
\end{Defn}

Let $\mathtt{E}$ be a $(\mathtt{T},\mathtt{S})$-bimodule and $\mathtt{G}$ be a $(\mathtt{T},\mathtt{R})$-bimodule.
Then $\mathrm{Hom^{L}}_{\mathtt{T}} ( \mathtt{E}, \mathtt{G} )$ naturally admits a $(\mathtt{S},\mathtt{R})$-bimodule structure.

\begin{prop} \label{202407291625}
    We have a natural isomorphism for $\mathtt{E},\mathtt{F}$ and $\mathtt{G}$:
    \begin{align*} 
        \hom_{\mathtt{T},\mathtt{R}} ( \mathtt{E} \otimes_{\mathtt{S}} \mathtt{F}, \mathtt{G} ) \cong \hom_{\mathtt{S},\mathtt{R}} ( \mathtt{F}, \mathrm{Hom^{L}}_{\mathtt{T}} ( \mathtt{E}, \mathtt{G} ) ) ,
    \end{align*}
    where $\hom_{\mathtt{T},\mathtt{R}}$ denotes the set of $(\mathtt{T},\mathtt{R})$-homomorphisms between bimodules.
    In particular, the functor $\mathrm{Hom^{L}}_{\mathtt{T}} ( \mathtt{E} , -)$ on the category of $(\mathtt{T},\mathtt{R})$-bimodules gives a right adjoint to the functor $\mathtt{E} \otimes_{\mathtt{S}} (-)$.
    Analogously, the functor $\mathrm{Hom^{R}}_{\mathtt{R}} ( \mathtt{F}, -)$ gives a right adjoint to the functor $(-)\otimes_{\mathtt{S}} \mathtt{F}$.
\end{prop}
\begin{proof}
    The natural isomorphism is induced by the adjunction considered in Definition \ref{202410121637}.
\end{proof}

Let $\mathtt{T}$ be a monad in $\mathsf{Mat}_{\mathds{k}}$.
Consider left $\mathtt{T}$-modules $\mathtt{M}_i$ with domain $\mathcal{Y}_i$ for $i \in \{1,2,3\}$.
By restriction, the pairing (\ref{202410121142}) induces the following associative pairing:
\begin{align} \label{202410051553}
    \mathrm{Hom^{L}}_{\mathtt{T}} ( \mathtt{M}_1, \mathtt{M}_2 ) \otimes \mathrm{Hom^{L}}_{\mathtt{T}} ( \mathtt{M}_2, \mathtt{M}_3) \to \mathrm{Hom^{L}}_{\mathtt{T}} ( \mathtt{M}_1, \mathtt{M}_3 ) .
\end{align}

\begin{Defn}
    For a left $\mathtt{T}$-module $\mathtt{M}$ with domain $\mathcal{Y}$, we define the monad $\mathrm{End^{L}}_{\mathtt{T}} (\mathtt{M} )$ on $\mathcal{Y}$ to be $ \mathrm{Hom^{L}}_{\mathtt{T}} ( \mathtt{M}, \mathtt{M} )$.
    We call this the {\it $\mathtt{T}$-endomorphism monad of $\mathtt{M}$}.
    Considering $\mathtt{M}_1=\mathtt{M}_2=\mathtt{M}_3 = \mathtt{M}$ in the pairing (\ref{202410051553}), we obtain the monad operation.
    The unit $\eta : \mathds{I}_{\mathcal{Y}} \to \mathrm{End^{L}}_{\mathtt{T}} (\mathtt{M} )$ is given by the assignment of the identities on $\mathtt{M} (Y)$ for each $Y \in \mathcal{Y}$.
\end{Defn}

\begin{remark}
    Observe that the monad operation on $\mathrm{End^{L}}_{\mathtt{T}} (\mathtt{M} )$ is given by the opposite of the usual composition of endomorphisms on $\mathtt{M}$.
\end{remark}

\begin{notation} \label{202410121644}
    We are mainly interested in {\it left} $\mathtt{T}$-modules, so that we introduce the following notation:
    \begin{align*}
        \mathrm{Hom}_{\mathtt{T}} ( \mathtt{M}_1, \mathtt{M}_2 ) {:=} \mathrm{Hom^{L}}_{\mathtt{T}} ( \mathtt{M}_1, \mathtt{M}_2 ), \quad \mathrm{End}_{\mathtt{T}} (\mathtt{M} ) {:=} \mathrm{End^{L}}_{\mathtt{T}} (\mathtt{M} ).
    \end{align*}
\end{notation}

\subsection{Upper triangular monads} 

In this section, we investigate monads which are upper triangular in the sense of Definition \ref{202512071707}.

\begin{Defn} \label{202512090927}
    Let $\mathtt{T}$ be an upper triangular monad in $\mathsf{Mat}_{\mathds{k}}$ on $\mathds{N}$.
    A left $\mathtt{T}$-module $\mathtt{M}$ is {\it $d$-truncated} if $\mathtt{M} (n) = 0$ for $n > d$.
    We define $\mathtt{T}\mbox{-}\mathsf{Mod}^{d-\mathsf{trun}}$ to be the full subcategory of $\mathtt{T}\mbox{-}\mathsf{Mod}$ consisting of $d$-truncated $\mathtt{T}$-modules.
\end{Defn}

\begin{Defn}
    For a left $\mathtt{T}$-module $\mathtt{M}$, we define an $\mathds{N}$-indexed module $\tau_{d}(\mathtt{M})$ by $\tau_{d}(\mathtt{M}) (n) {:=} \mathtt{M} (n)$ if $n \leq d$, and $\tau_{d}(\mathtt{M}) (n) {:=} 0$ otherwise.
\end{Defn}

This gives a well-defined $\mathtt{T}$-submodule $\tau_{d}(\mathtt{M})$ of $\mathtt{M}$, since $\mathtt{T}$ is upper triangular.
Furthermore, there is a filtration of left $\mathtt{T}$-modules:
$$
\tau_{0}(\mathtt{M}) \subset \cdots \subset \tau_{d}(\mathtt{M}) \subset \tau_{d+1}(\mathtt{M}) \subset \cdots \subset \mathtt{M} .
$$

It is clear that the subcategory $\mathtt{T}\mbox{-}\mathsf{Mod}^{d-\mathsf{trun}}$ is thick in $\mathtt{T}\mbox{-}\mathsf{Mod}^{d-\mathsf{trun}}$.
Recall from \cite{Gabriel1962} the notion of Serre quotient category:
\begin{prop} \label{202512081028}
    For an upper triangular monad $\mathtt{T}$ and $d \in \mathds{N}$, the assignment of $\tau_{d}(\mathtt{M}) / \tau_{d-1}(\mathtt{M})$ to $\mathtt{M}$ induces an equivalence of categories:
    $$
     \mathtt{T}\mbox{-}\mathsf{Mod}^{d-\mathsf{trun}}/\mathtt{T}\mbox{-}\mathsf{Mod}^{(d-1)-\mathsf{trun}} \simeq  \mathtt{T}(d,d)\mbox{-}\mathsf{Mod} ,
    $$
    where the category of the left hand side denotes the Serre quotient category.
\end{prop}
\begin{proof}
    For a $\mathtt{T}(d,d)$-module $V$, we construct a left $\mathtt{T}$-module $\mathtt{M}$ by $\mathtt{M} (n) = V$ if $n =d$ and $\mathtt{M}(n) = 0$ otherwise.
    Conversely, any $\mathtt{T}$-module $\mathtt{M}$ induces a $\mathtt{T}(d,d)$-module $\mathtt{M}(d)$.
    These constructions give functors between $\mathtt{T}(d,d)\mbox{-}\mathsf{Mod} $ and $\mathtt{T}\mbox{-}\mathsf{Mod}^{d-\mathsf{trun}}$.
    Moreover, the latter factors through the Serre quotient category in the statement by definition, and these functors induce an equivalence of categories.
\end{proof}

\section{The eigenmonad construction}
\label{202410161744}

This section is devoted to the main construction of this paper: given a monad and a multiplicative submodule, we present a natural subquotient, which we call the {\it eigenmonad}.
We also show that there are canonical adjunctions between module categories associated with the eigenmonad.
Throughout this section, we fix an object $\mathcal{X} \in \mathsf{Mat}_{\mathds{k}}$ and a monad $\mathtt{T}$ in $\mathsf{Mat}_{\mathds{k}}$ on $\mathcal{X}$.

\subsection{Definitions}
\label{202512151221}

In this section, we introduce the notions of idealizer monad and eigenmonad.

\begin{Defn} \label{202410051557}
    A submodule $\mathtt{J} \subset \mathtt{T}$ is {\it multiplicative} if $\mathtt{J} \circ \mathtt{J} \subset \mathtt{J}$.
    In other words, for $X,Y,Z\in\mathcal{X}$, $\mathtt{J}(X,Y) \subset \mathtt{T}(X,Y)$ is a $\mathds{k}$-submodule such that $f \circ g \in \mathtt{J}(X,Z)$ for $f \in \mathtt{J}(X,Y)$ and $g \in \mathtt{J}(Y,Z)$.
\end{Defn}

\begin{Defn}[Idealizer Monad] \label{202410011802}
    Let $\mathtt{J}$ be a multiplicative submodule of $\mathtt{T}$.
    We define the {\it idealizer monad of $\mathtt{J}$} to be the submodule $\mathrm{Iz}_{\mathtt{T}} (\mathtt{J})$ of $\mathtt{T}$ such that
    \begin{align} \label{202409261426}
    \left( \mathrm{Iz}_{\mathtt{T}} (\mathtt{J}) \right) (Y,X) {:=} \{ f \in \mathtt{T} (Y,X) |~ \mathtt{J} \circ f \subset \mathtt{J} (-,X), ~ f \circ \mathtt{J} \subset \mathtt{J} (Y,-) \} , \quad X,Y \in \mathcal{X}.
    \end{align}
    If the monad $\mathtt{T}$ is clear from the context, then we write $\mathrm{Iz} (\mathtt{J})$ for $\mathrm{Iz}_{\mathtt{T}} (\mathtt{J})$.
\end{Defn}

\begin{remark}
    We have $\mathtt{J} \subset \mathrm{Iz} (\mathtt{J})$, since $\mathtt{J}$ is assumed to be multiplicative.
    It turns out that $\mathrm{Iz} (\mathtt{J})$ is the {\it maximal} submonad of $\mathtt{T}$ containing $\mathtt{J}$ as a two-sided ideal.
\end{remark}

\begin{Defn}[Eigenmonad] \label{202512162200}
    The {\it eigenmonad of $\mathtt{T}$ by $\mathtt{J}$} is defined to be the quotient monad,
    $$\mathrm{E}_{\mathtt{T}} ( \mathtt{J}) {:=} \mathrm{Iz}_{\mathtt{T}} (\mathtt{J}) /\mathtt{J} .$$
    We write $\mathrm{E} (\mathtt{J})$ for $ \mathrm{E}_{\mathtt{T}} ( \mathtt{J})$ if $\mathtt{T}$ is clear from the context.
\end{Defn}

\begin{remark}
    To the best of our knowledge, the origin of this definition goes back to Ore \cite{Ore1932}.
    Ore introduced the notion of an eigenring to study a formal property of differential operators.
    See \cite[Section 0.7]{Cohn1985} for the modern approach.
    Our definition extends it to monads in the bicategory $\mathsf{Mat}_{\mathds{k}}$ on $\mathcal{X}$, in that, for $\mathcal{X} = \ast$ the singleton, it recovers the classical notion.
\end{remark}

\begin{remark}
    Although the above definition and the related constructions in this part are developed using the specific bicategory $\mathsf{Mat}_{\mathds{k}}$, our construction considered here extends to any bicategory that is closed and whose composition is right exact.
\end{remark}

In this paper, we mainly study eigenmonads by a particular multiplicative submodule, as follows.

\begin{Defn}
    A submodule $\mathtt{J} \subset \mathtt{T}$ is a {\it left ideal} of $\mathtt{T}$ if $\mathtt{T} \circ \mathtt{J} \subset \mathtt{J}$.
    In a symmetric fashion, a submodule $\mathtt{J} \subset \mathtt{T}$ is a {\it right ideal} of $\mathtt{T}$ if we have $\mathtt{J} \circ \mathtt{T} \subset \mathtt{J}$.
    A {\it two-sided ideal} of $\mathtt{T}$ is a left ideal which is also a right ideal of $\mathtt{T}$.
\end{Defn}

\begin{remark} \label{202410161125}
    If $\mathtt{J}$ is a left or right ideal of $\mathtt{T}$, the definition of the eigenmonad $\mathrm{E}_{\mathtt{T}}(\mathtt{J})$ is simplified.
    Indeed, the second or the first condition of (\ref{202409261426}) is redundant, respectively.
\end{remark}

\begin{remark}
    Based on the equivalence between linear categories and monads in $\mathsf{Mat}_{\mathds{k}}$, this notion corresponds to that of a (left or right) ideal of an additive category introduced in \cite[Section 3]{Mitchell1972}.
\end{remark}

The following is immediate from the definitions.
\begin{prop} \label{202410041610}
    Let $\mathtt{J} \subset \mathtt{T}$ be a multiplicative submodule.
    The following conditions are equivalent to each other:
    \begin{enumerate}   
        \item $\mathtt{J}$ is a two-sided ideal.
        \item $\mathrm{Iz} ( \mathtt{J} ) = \mathtt{T}$.
    \end{enumerate}
    Under either condition, $\mathrm{E} ( \mathtt{J}) \cong \mathtt{T}/\mathtt{J}$.
\end{prop}

\begin{prop}[Naturality] \label{202512071648}
    Let $\pi: \mathtt{T}\to \mathtt{Q}$ be a monad epimorphism.
    \begin{enumerate}
        \item For a left ideal $\mathtt{J} \subset \mathtt{T}$, the image $\mathtt{I} {:=} \pi ( \mathtt{J})$ is a left ideal of $\mathtt{Q}$.
        \item $\pi ( \mathrm{Iz}_{\mathtt{T}}(\mathtt{J})) \subset \mathrm{Iz}_{\mathtt{Q}}(\mathtt{I})$.
        \item $\pi$ induces a monad map $\mathrm{E}_{\mathtt{T}}(\mathtt{J}) \to \mathrm{E}_{\mathtt{Q}}(\mathtt{I})$, which is functorial with respect to monad epimorphisms $\pi: (\mathtt{T},\mathtt{J}) \to (\mathtt{Q},\mathtt{I})$ of such pairs.
    \end{enumerate}
\end{prop}
\begin{proof}
    The first part follows from the assumption that $\pi$ is an epimorphism.
    For $f \in \mathrm{Iz}_{\mathtt{T}}(\mathtt{J}) (Y,X)$, we have $\mathtt{I} \circ \pi (f )= \pi ( \mathtt{J}) \circ \pi (f ) \subset  \pi ( \mathtt{J} \circ f) \subset \pi (\mathtt{J}) = \mathtt{I}$, so we obtain the second statement.
    The third one follows from the second one.
    The functoriality is clear from the construction of the map $\mathrm{E}_{\mathtt{T}}(\mathtt{J}) \to \mathrm{E}_{\mathtt{Q}}(\mathtt{I})$.
\end{proof}

\subsection{Associated adjunctions}
\label{202404111312}

In this section, for a left ideal $\mathtt{J} \subset \mathtt{T}$, we introduce an adjunction between $\mathtt{T}$-modules and $\mathrm{E}_{\mathtt{T}}( \mathtt{J})$-modules.

\begin{Lemma} \label{202401241122}
    The monad structure on $\mathtt{T}$ induces a $(\mathtt{T}, \mathrm{E}_{\mathtt{T}} ( \mathtt{J}))$-bimodule structure on $\mathtt{T}/\mathtt{J}$.
\end{Lemma}
\begin{proof}
    The regular left $\mathtt{T}$-action on $\mathtt{T}$ induces that on $\mathtt{T}/\mathtt{J}$.
    The regular right $\mathtt{T}$-action on $\mathtt{T}$ itself induces the right $\mathrm{E}_{\mathtt{T}} ( \mathtt{J})$-action on $\mathtt{T}/\mathtt{J}$.
    For $f \in \left( \mathtt{I}_{\mathtt{T}}(\mathtt{J}) \right) (Y,X)$, we have $\mathtt{J} (Z,Y) \circ f \subset \mathtt{J} (Z,X)$, so that it gives a well-defined map 
    $$\left( \mathtt{T}/\mathtt{J} \right) (Z,Y) \otimes \left( \mathtt{I}_{\mathtt{T}}(\mathtt{J}) \right) (Y,X) \to \left( \mathtt{T}/\mathtt{J} \right) (Z,X); \quad (v \mod{\mathtt{J}}) \otimes f \mapsto (v \circ f \mod{\mathtt{J}}) .$$
    Furthermore, for $f \in \mathtt{J} (Y,X)$, we have $v \circ f \mod{\mathtt{J}} = 0$.
    Since this map is induced by the regular right $\mathtt{T}$-action on $\mathtt{T}$, this construction yields a right $\mathrm{E}_{\mathtt{T}}(\mathtt{J})$-action on $\mathtt{T}/\mathtt{J}$, which commutes with the left $\mathtt{T}$-action.
\end{proof}

\begin{Defn} \label{202403061037}
    The quotient $\mathtt{T}/\mathtt{J}$ endowed with the above $(\mathtt{T}, \mathrm{E}_{\mathtt{T}} ( \mathtt{J}))$-bimodule structure is referred to as the {\it canonical $(\mathtt{T}, \mathrm{E}_{\mathtt{T}} ( \mathtt{J}))$-bimodule associated with $(\mathtt{T} , \mathtt{J} )$}, or, if $\mathtt{T}$ is specified, the {\it canonical $(\mathtt{T}, \mathrm{E}_{\mathtt{T}} ( \mathtt{J}))$-bimodule associated with $\mathtt{J}$}.
\end{Defn}

This bimodule structure of $\mathtt{T}/\mathtt{J}$ is of fundamental importance for our main results.
By applying the constructions in Definitions \ref{202409261711} and \ref{202409261651}, we obtain the following adjunctions:
\begin{prop}[Eigenmonad adjunctions] \label{202409261707}
    \begin{enumerate}
        \item We have an adjunction between the category of left $\mathrm{E}_{\mathtt{T}} ( \mathtt{J})$-modules and that of left $\mathtt{T}$-modules:
        $$
        \begin{tikzcd}
                \mathtt{T}/\mathtt{J} \otimes_{\mathrm{E}_{\mathtt{T}} ( \mathtt{J})} (-) : \mathrm{E}_{\mathtt{T}} ( \mathtt{J})\mbox{-}\mathsf{Mod} \arrow[r, shift right=1ex, ""{name=G}] & \mathtt{T}\mbox{-}\mathsf{Mod} : \mathrm{Hom^{L}}_{\mathtt{T}} (\mathtt{T}/\mathtt{J} , -) \arrow[l, shift right=1ex, ""{name=F}]
                \arrow[phantom, from=G, to=F, , "\scriptscriptstyle\boldsymbol{\top}"].
        \end{tikzcd}
        $$
        \item We have an adjunction between the category of right $\mathrm{E}_{\mathtt{T}} ( \mathtt{J})$-modules and that of right $\mathtt{T}$-modules:
        $$
        \begin{tikzcd}
                  (-) \otimes_{\mathtt{T}} \mathtt{T}/\mathtt{J} : \mathsf{Mod}\mbox{-}\mathtt{T} \arrow[r, shift right=1ex, ""{name=G}] & \mathsf{Mod}\mbox{-}\mathrm{E}_{\mathtt{T}} ( \mathtt{J}) : \mathrm{Hom^{R}}_{\mathrm{E}_{\mathtt{T}} ( \mathtt{J})} (\mathtt{T}/\mathtt{J} , -) \arrow[l, shift right=1ex, ""{name=F}]
                \arrow[phantom, from=G, to=F, , "\scriptscriptstyle\boldsymbol{\top}"].
        \end{tikzcd}
        $$
    \end{enumerate}
\end{prop}
\begin{proof}
    Both of the statements follow from the tensor-hom adjunction given in Proposition \ref{202407291625}.
\end{proof}

We refer to these adjunctions as the {\it eigenmonad adjunctions associated with $(\mathtt{T}, \mathtt{J})$}.

\begin{remark} \label{202512151851}
    These adjunctions generalize those given in our previous paper \cite{kim2025poly} where we only consider a two-sided ideal $\mathtt{J}$.
\end{remark}

\subsection{Vanishing module}

In this section, we give a useful description of the functor $\mathrm{Hom}_{\mathtt{T}} (\mathtt{T}/\mathtt{J} , -)$.
Subsequently, we obtain three equivalent descriptions of eigenmonads.
The results presented here are applied to concrete examples throughout this paper.

\begin{Defn} \label{202311281540}
    Let $\mathtt{M}$ be a left $\mathtt{T}$-module with arbitrary domain $\mathcal{Y}$.
    For a left ideal $\mathtt{J}$ of $\mathtt{T}$, we define the {\it $\mathtt{J}$-vanishing module} $\mathrm{V} ( \mathtt{M} ; \mathtt{J})$ to be the submodule of $\mathtt{M}$ such that
    $$\left(\mathrm{V} ( \mathtt{M} ;  \mathtt{J})\right) (X, Y) {:=} \{ f \in \mathtt{M} (X,Y) |~ \mathtt{J} \rhd f \cong 0 \}, \quad (X,Y) \in \mathcal{X} \times \mathcal{Y} .$$
    We endow the $\mathtt{J}$-vanishing module with a left $\mathrm{E}_{\mathtt{T}} ( \mathtt{J})$-module structure as follows.
\end{Defn}

\begin{Lemma} 
    The left $\mathtt{T}$-module structure on $\mathtt{M}$ induces a left $\mathrm{E}_{\mathtt{T}} ( \mathtt{J})$-module structure on $\mathrm{V}( \mathtt{M} ;  \mathtt{J})$.
\end{Lemma}
\begin{proof}
    It is sufficient to show that, for any $X_1,X_2 \in \mathcal{X}$ and $Y \in \mathcal{Y}$, the assignment of $f \rhd v \in \left( \mathrm{V}  ( \mathtt{M} ;  \mathtt{J}) \right) (X_2,Y)$ to $(f \mod{\mathtt{J}}) \in \left( \mathrm{E}_{\mathtt{T}} ( \mathtt{J}) \right) (X_2,X_1)$ and $v \in \left( \mathrm{V} (\mathtt{M} ;  \mathtt{J}) \right) (X_1,Y)$ is well-defined.
    In fact, for any $X \in \mathcal{X}$ and $g \in \mathtt{J}(X,X_2)$, we have $g \rhd ( f \rhd v) =( g \circ f) \rhd v = 0$, since $f$ lies in the idealizer monad.
    Hence, we obtain the result.
\end{proof}

As in Definition \ref{202409261651}, for a left $\mathtt{T}$-module $\mathtt{M}$ (with any domain $\mathcal{Y}$), the hom module $\mathrm{Hom}_{\mathtt{T}} ( \mathtt{T} , \mathtt{M} )$ naturally carries a left $\mathtt{T}$-module structure.
This is naturally isomorphic to $\mathtt{M}$ by the Yoneda lemma.
This observation leads to an equivalent description of the vanishing module:

\begin{Lemma} \label{202405241347}
    The evaluation at the unit of the monad $\mathtt{T}$ induces a natural isomorphism of left $\mathrm{E}_{\mathtt{T}} ( \mathtt{J})$-modules.
    $$
    \mathrm{Hom}_\mathtt{T} ( \mathtt{T}/ \mathtt{J} , \mathtt{M} ) \stackrel{\cong}{\longrightarrow} \mathrm{V} ( \mathtt{M} ;  \mathtt{J}) .
    $$
\end{Lemma}
\begin{proof}
    Applying the functor $\mathrm{Hom}_\mathtt{T} (- , \mathtt{M})$ to the short exact sequence of left $\mathtt{T}$-modules $0 \to \mathtt{J} \to \mathtt{T} \to \mathtt{T}/\mathtt{J} \to 0$, we obtain the exact sequence $0\to \mathrm{Hom}_\mathtt{T} (\mathtt{T}/\mathtt{J} , \mathtt{M}) \to \mathrm{Hom}_\mathtt{T} (\mathtt{T} , \mathtt{M}) \to \mathrm{Hom}_\mathtt{T} (\mathtt{J} , \mathtt{M})$.
    The exactness and the Yoneda isomorphism $\mathrm{Hom}_\mathtt{T} ( \mathtt{T} , \mathtt{M} ) \cong \mathtt{M}$ lead to the isomorphism in the statement.
\end{proof}

Recall the endomorphism monad from Notation \ref{202410121644}.
\begin{prop} \label{202401241330}
    The eigenmonad $\mathrm{E}_{\mathtt{T}} ( \mathtt{J})$ has equivalent descriptions as follows:
    \begin{enumerate}
        \item The quotient monad $\mathrm{Iz}_{\mathtt{T}} (\mathtt{J}) / \mathtt{J}$ of the idealizer monad by $\mathtt{J}$.
        \item The $\mathtt{J}$-vanishing module $\mathrm{V} ( \mathtt{T}/ \mathtt{J} ;  \mathtt{J})$.
        \item The endomorphism monad $\mathrm{End}_{\mathtt{T}} ( \mathtt{T} / \mathtt{J})$ of $\mathtt{T}/\mathtt{J}$.
    \end{enumerate}
\end{prop}
\begin{proof}
    The inclusion $\mathrm{Iz} ( \mathtt{J} ) \to \mathtt{T}$ induces a monomorphism $\mathrm{E} ( \mathtt{J} ) \to \mathtt{T} / \mathtt{J}$.
    The image is contained in $\mathrm{V} ( \mathtt{T}/ \mathtt{J} ;  \mathtt{J})$, since $\mathtt{J} \rhd \mathrm{E} ( \mathtt{J} ) = \mathtt{J} \rhd ( \mathrm{Iz}(\mathtt{J})/\mathtt{J}) = 0$ by the definition of the idealizer monad.
    Hence, we obtain (2).
    (3) follows from (2) and the isomorphism in Lemma \ref{202405241347}.
\end{proof}

We highlight several features of the above descriptions:
\begin{remark}
    \begin{enumerate}
        \item When studying specific eigenmonads, it is useful to apply the second form of Proposition \ref{202401241330}.
        \item The third part establishes that the maximal right action that commutes with the left $\mathtt{T}$-action on $\mathtt{T}/\mathtt{J}$ is exactly the right $\mathrm{E}_{\mathtt{T}} ( \mathtt{J})$-action on $\mathtt{T}/\mathtt{J}$.
        \item
        Moreover, the third formulation of the eigenmonad applies to any bicategory whose hom-categories are abelian, and whose composition is closed and right exact.
    \end{enumerate}
\end{remark}

Lemma \ref{202405241347} implies that the vanishing module construction gives an alternative description of the right adjoint of the adjunction (1) of Proposition \ref{202409261707}.
In a symmetric fashion, we treat the second adjunction by introducing the following:
\begin{Defn}
    Let $\mathtt{N}$ be a right $\mathtt{T}$-module with arbitrary codomain $\mathcal{Y}$.
    For a left ideal $\mathtt{J}$ of $\mathtt{T}$, we define the {\it $\mathtt{J}$-covanishing module} $\Lambda ( \mathtt{N} ; \mathtt{J})$ to be the quotient module $\mathtt{N} / (\mathtt{N} \lhd \mathtt{J} )$.
\end{Defn}

One may verify that the $\mathtt{J}$-covanishing module has a right $\mathrm{E}_{\mathtt{T}} ( \mathtt{J})$-module structure.
Moreover, we have a natural isomorphism of right $\mathrm{E}_{\mathtt{T}} ( \mathtt{J})$-modules:
\begin{align*}
    \mathtt{N} \otimes_{\mathtt{T}} \mathtt{T}/\mathtt{J} \stackrel{\cong}{\longrightarrow} \Lambda ( \mathtt{N} ; \mathtt{J})   .
\end{align*}

\begin{remark}
    We remark that, under this natural isomorphism, the second adjunction of Proposition \ref{202409261707} generalizes the adjunction given in \cite{kim2025poly} where $\mathtt{J}$ is assumed to be a two-sided ideal.
\end{remark}

\subsection{Left ideals, adjunctions and properties}

In this section, we present a natural restriction of the first adjunction in Proposition \ref{202409261707}, which is referred to as the {\it refined eigenmonad adjunction}.
To this end, we begin with the following general definition.
\begin{Defn} \label{202402011603}
    Let $\mathcal{Y}$ be an object of $\mathsf{Mat}_\mathds{k}$.
    \begin{itemize}
        \item A left $\mathtt{T}$-module $\mathtt{M}$ with domain $\mathcal{Y}$ is {\it $\mathtt{J}$-vanishingly generated} if the composition $$\mathtt{T} \otimes \mathrm{V} ( \mathtt{M} ; \mathtt{J}) \stackrel{\mathrm{id}_{\mathtt{T}}\otimes \iota}{\to} \mathtt{T}\otimes \mathtt{M} \stackrel{\rhd}{\to} \mathtt{M}$$ where $\iota : \mathrm{V} ( \mathtt{M} ; \mathtt{J}) \to \mathtt{M}$ denotes the inclusion and $\rhd$ is the $\mathtt{T}$-action, is an epimorphism.
        \item A left $\mathrm{E}_{\mathtt{T}} (\mathtt{J})$-module $\mathtt{N}$ with domain $\mathcal{Y}$ is {\it $\mathtt{J}$-vanishingly extensible} if the composition $$\mathtt{N} \cong \mathds{I}_{\mathcal{X}} \otimes \mathtt{N} \stackrel{\eta_{\mathtt{T}}\otimes\mathrm{id}_{\mathtt{N}}}{\to} \mathtt{T} \otimes \mathtt{N} \twoheadrightarrow \mathtt{T}/\mathtt{J} \otimes_{\mathrm{E}_{\mathtt{T}}(\mathtt{J} )} \mathtt{N}$$ where $\eta_{\mathtt{T}}$ is the unit of $\mathtt{T}$, gives a monomorphism.
    \end{itemize}
\end{Defn}

\begin{Example}
    The regular left $\mathrm{E}_{\mathtt{T}} (\mathtt{J})$-module $\mathrm{E}_{\mathtt{T}} (\mathtt{J})$ is $\mathtt{J}$-vanishingly extensible, since $\mathtt{T}/\mathtt{J} \otimes_{\mathrm{E}_{\mathtt{T}}(\mathtt{J} )} \mathrm{E}_{\mathtt{T}}(\mathtt{J} ) \cong \mathtt{T}/\mathtt{J}$ and $\mathrm{E}_{\mathtt{T}}(\mathtt{J} )\hookrightarrow \mathtt{T}/\mathtt{J}$ by (2) of Proposition \ref{202401241330}.
\end{Example}

\begin{Example}
    The left $\mathtt{T}$-module $\mathtt{T}/\mathtt{J}$ is $\mathtt{J}$-vanishingly generated; however the regular left $\mathtt{T}$-module $\mathtt{T}$ itself may not be $\mathtt{J}$-vanishingly generated in general.
\end{Example}

\begin{remark}
    Recall the eigenmonad adjunction from Proposition \ref{202409261707}.
    It is clear that a left $\mathtt{T}$-module $\mathtt{M}$ is $\mathtt{J}$-vanishingly generated if and only if the associated counit $\mathtt{T} / \mathtt{J} \otimes_{\mathrm{E}_{\mathtt{T}}(\mathtt{J})} \mathrm{V} ( \mathtt{M} ;  \mathtt{J}) \to \mathtt{M}$ is an epimorphism.
    Analogously, a left $\mathrm{E}_{\mathtt{T}} (\mathtt{J})$-module $\mathtt{N}$ is $\mathtt{J}$-vanishingly extensible if and only if the associated unit $\mathtt{N} \to \mathrm{V} \left( \mathtt{T}/\mathtt{J} \otimes_{\mathrm{E}_{\mathtt{T}}(\mathtt{J} )} \mathtt{N} ;  \mathtt{J} \right)$ is a monomorphism.
\end{remark}

\begin{Lemma} \label{202512081840}
    If $\mathrm{E}_{\mathtt{T}}(\mathtt{J}) \hookrightarrow \mathtt{T}/\mathtt{J}$ admits a retract in right $\mathrm{E}_{\mathtt{T}} ( \mathtt{J})$-modules, then every left $\mathrm{E}_{\mathtt{T}} ( \mathtt{J})$-module is $\mathtt{J}$-vanishingly extensible.
\end{Lemma}
\begin{proof}
    For a left $\mathrm{E}_{\mathtt{T}}(\mathtt{J})$-module $\mathtt{N}$, the retract $r$ tensored with the identity on $\mathtt{N}$ induces a retract of $\mathtt{N} \to \mathtt{T}/\mathtt{J} \otimes_{\mathrm{E}_{\mathtt{T}}(\mathtt{J} )} \mathtt{N}$.
\end{proof}

Using the above concepts, we give some relationships among left ideals, adjunctions of module categories and isomorphism-invariant properties of modules.
For this purpose, we introduce the following:

\begin{Defn}
    Let $\mathtt{T}\mbox{-}\mathcal{M}\mathsf{od}$ denote the class of isomorphism classes of left $\mathtt{T}$-modules.
    We denote by $\mathcal{V}(\mathtt{J})$ the subclass of $\mathtt{T}\mbox{-}\mathcal{M}\mathsf{od}$ consisting of isomorphism classes of $\mathtt{J}$-vanishingly generated left $\mathtt{T}$-modules.
\end{Defn}

The map $\mathcal{V}$ is contravariant in the following sense:
\begin{Lemma} \label{202405291701}
    Let $\mathtt{J}, \mathtt{J}^\prime$ be left ideals of $\mathtt{T}$.
    Consider the following conditions:
    \begin{enumerate}
        \item $\mathtt{J} \subset \mathtt{J}^\prime$.
        \item $\mathcal{V} (\mathtt{J}^\prime) \subset \mathcal{V} (\mathtt{J})$.
    \end{enumerate}
    The condition (1) implies (2).
    Furthermore, if $\mathtt{J}$ is a two-sided ideal, then (2) implies (1).
\end{Lemma}
\begin{proof}
    We prove the first statement.
    Let $\mathtt{M}$ be a left $\mathtt{T}$-module.
    By definition, we have $\mathrm{V} (\mathtt{M} ; \mathtt{J}^\prime) \subset \mathrm{V} ( \mathtt{M} ; \mathtt{J})$.
    Hence, if $\mathtt{M}$ is $\mathtt{J}^\prime$-vanishingly generated, then it is $\mathtt{J}$-vanishingly generated.

    We now prove the second statement.
    Note that $\mathtt{T}/ \mathtt{J}^\prime$ is $\mathtt{J}^\prime$-vanishingly generated, so, by the assumption, this is $\mathtt{J}$-vanishingly generated.
    Hence, we have $\mathrm{V} ( \mathtt{T} / \mathtt{J}^\prime ; \mathtt{J}) = \mathtt{T} / \mathtt{J}^\prime$ which implies $\mathtt{J} \subset \mathtt{J}^\prime$.
\end{proof}

In our earlier paper in this series \cite{kim2025poly}, we presented some relationships among {\it two-sided} ideals of a monad, adjunctions of module categories and isomorphism-invariant properties of modules.
We extend these relationships to left ideals of a monad; a brief remark is given below.

\begin{remark}
    The map $\mathcal{V}$ can be understood using the following schematic diagram.
    $$
    \mathcal{V} : ~\{ \mathrm{Left~ideals~of~}\mathtt{T}\} \to \{\mathrm{Adjunctions~involving~}\mathtt{T}\mbox{-}\mathsf{Mod}\} \to \{\mathrm{Subclasses~of~}\mathtt{T}\mbox{-}\mathcal{M}\mathsf{od}\} .
    $$
    Given a left ideal of a monad $\mathtt{T}$, we have the eigenmonad adjunctions, specifically the first part of Proposition \ref{202409261707}.
    This is represented in the first map in the diagram.
    The second map assigns to an adjunction the class of $\mathtt{T}$-modules for which the counit is an epimorphism, using some general arguments on adjunctions \cite{kim2025poly}.
    In particular, it involves a refinement of adjunctions, which yields the following proposition.
\end{remark}

We denote by $\mathtt{T}\mbox{-}\mathsf{Mod}^{\mathtt{J}\mathrm{-vg}}$ the full subcategory of $\mathtt{T}\mbox{-}\mathsf{Mod}$ of $\mathtt{J}$-vanishingly generated $\mathtt{T}$-modules.
We denote by $\mathrm{E}_{\mathtt{T}} ( \mathtt{J})\mbox{-}\mathsf{Mod}^{\mathtt{J}\mathrm{-ve}}$ the full subcategory of $\mathrm{E}_{\mathtt{T}} ( \mathtt{J})\mbox{-}\mathsf{Mod}$ of $\mathtt{J}$-vanishingly extensible $\mathrm{E}_{\mathtt{T}} ( \mathtt{J})$-modules.
\begin{prop}[Refined eigenmonad adjunction]
\label{202402011440}    
    The restriction of the eigenmonad adjunction associated with $(\mathtt{T},\mathtt{J})$ induces the following adjunction:
    $$
    \begin{tikzcd}
            \mathtt{T}/\mathtt{J} \otimes_{\mathrm{E}_{\mathtt{T}} ( \mathtt{J})} (-) : \mathrm{E}_{\mathtt{T}} ( \mathtt{J})\mbox{-}\mathsf{Mod}^{\mathtt{J}\mathrm{-ve}} \arrow[r, shift right=1ex, ""{name=G}] & \mathtt{T}\mbox{-}\mathsf{Mod}^{\mathtt{J}\mathrm{-vg}} : \mathrm{Hom}_{\mathtt{T}} (\mathtt{T}/\mathtt{J} , -) \arrow[l, shift right=1ex, ""{name=F}]
            \arrow[phantom, from=G, to=F, , "\scriptscriptstyle\boldsymbol{\top}"].
    \end{tikzcd}
    $$
    Furthermore, both of these functors are faithful.
\end{prop}

\begin{remark}
    An analogous construction can be carried out for the second one of Proposition \ref{202409261707}.
\end{remark}

{\Large \part{The primitivity ideal} \label{202512151908}}

We now specialize the preceding constructions to the monad and ideal arising from a Lawvere theory.
Recall that a Lawvere theory \cite{Lawvere1963} is a category with finite products, whose objects are $\mathds{N}$, with the products on objects given by addition.
For a Lawvere theory $\mathcal{C}$ with a zero object, we introduce the primitivity ideal and study the associated eigenmonad.
In particular, we establish Theorems \ref{202604091210} and \ref{202512111518}.
As universal examples, we provide explicit computations for the opposite categories of free monoids and free groups.
As a consequence, we obtain Theorems \ref{202604082143} and \ref{202604071759}.

\vspace{3mm}

\section{Lawvere theories}
\label{202509031755}

In this section, we review Lawvere theories and present several examples.
We also introduce some structural conditions, which we call (ZM*), on Lawvere theories.

\subsection{Definitions and basic properties}

In this section, we give some basic properties of Lawvere theories and fix our notation.
Let $\mathcal{C}$ be a Lawvere theory.

\begin{Defn} \label{202606081632}
    Let $\Delta_n \in \mathcal{C}(n,2n)$ denote the diagonal map of $n \in \mathds{N}$.
    We denote by $\Delta^{(0)}_n \in \mathcal{C}(n,0)$ the unique morphism.
    We also set $\Delta^{(1)}_n {:=} \mathrm{id}_n$ and $\Delta^{(2)}_n {:=} \Delta_n$.
    For $m \geq 3$, we recursively define the {\it $m$-fold diagonal map} $\Delta^{(m)}_n \in \mathcal{C}(n,nm)$ by $\Delta^{(m)}_n {:=} (\Delta^{(m-1)}_n \times \mathrm{id}_n)\circ \Delta_n$.
    If $n = 1$, then we write simply $\Delta^{(m)}$.
\end{Defn}

\begin{Defn} \label{202604091226}
    Let $p_{n,k} \in \mathcal{C} (n,1)$ denote the $k$-th projection from $n$ to $1$.
\end{Defn}

\begin{notation}
    Let $\mathcal{C}_n {:=} \mathcal{C}(n,1)$.
\end{notation}

Throughout this paper, we identify the set $\mathcal{C}(n,m)$ with $\mathcal{C}_n^{\times m}$ by using the universality of products:
\begin{align} \label{202606181948}
    \mathcal{C}(n,m) \stackrel{\cong}{\to} \mathcal{C}_n^{\times m} ; \quad f \mapsto ( p_{m,1} \circ f, \cdots, p_{m,m} \circ f ) .
\end{align}
The inverse is given by
$$
(g_1, \cdots, g_m) \mapsto \left( \prod^m_{k=1} g_k \right) \circ \Delta^{(m)}_{n}.
$$
For instance, under the identification (\ref{202606181948}), we have
\begin{align} 
    (p_{n,1}, p_{n,2} , \cdots, p_{n,n}) &= \mathrm{id}_n \in \mathcal{C}(n,n) , \label{202606190854} \\
    (p_{n,1},  \cdots, p_{n,k-1}, p_{n,k},p_{n,k}, p_{n,k+1}, \cdots, p_{n,n}) &= \mathrm{id}_{k-1} \times \Delta \times \mathrm{id}_{n-k} \in \mathcal{C}(n,n+1) . \label{202606201739}
\end{align}

\begin{Example} \label{202603121546}
    The opposite categories of the following are Lawvere theories:
    \begin{itemize}
        \item $\mathbf{W}$, the category of finitely generated free monoids;
        \item $\G$, the category of finitely generated free groups.
    \end{itemize}
    In the Lawvere theories $\mathbf{W}^{\mathsf{o}}$ and $\G^{\mathsf{o}}$, finite products are given by free products of monoids and groups, respectively.
    Throughout this paper, we fix a countable set $\mathfrak{X} = \{ x_1, x_2, \cdots, x_n, \cdots \}$.
    We identify $n \in \mathds{N}$ with the free monoid $\mathsf{W}_n$ generated by the $n$ elements $x_1,x_2,\cdots, x_n$.
    Then, for $\mathcal{C} = \mathbf{W}^{\mathsf{o}}$, $\mathcal{C}_n = \mathcal{C}(n,1)$ is identified with $\mathsf{W}_n$, where $p_{n,k}$ corresponds to $x_k$.
    Similarly, for $\mathcal{C}= \G^{\mathsf{o}}$, the set $\mathcal{C}_n$ is identified with the free group $\mathsf{F}_n$ on the generators $x_1,x_2,\cdots, x_n$.
\end{Example}

For $g \in \mathcal{C}(l,n)$, we denote by $g^\ast : \mathcal{C} (n,m) \to \mathcal{C} (l,m)$ the map $f \mapsto f \circ g$.
\begin{Lemma} \label{202509042125}
    Let $g\in \mathcal{C}(l,n)$ and $f_1,\cdots, f_m \in \mathcal{C}_n$.
    We have $g^\ast (f_1, \cdots, f_m) = ( g^\ast (f_1), \cdots, g^\ast (f_m)) \in \mathcal{C}(l,m)$.
\end{Lemma}
\begin{proof}
    \begin{align*}
        &g^\ast (f_1, \cdots, f_m)  = (f_1, \cdots, f_m) \circ g = (f_1 \times \cdots \times  f_m) \circ \Delta^{(m)}_n\circ g  , \\
        =& (f_1 \times \cdots \times  f_m) \circ (g^{\times m}) \circ \Delta^{(m)}_l  = (g^\ast (f_1) \times \cdots \times g^\ast (f_m)) \circ \Delta^{(m)}_l = (g^\ast (f_1), \cdots, g^\ast (f_m)) .
    \end{align*}
\end{proof}

It is useful to understand how products in a Lawvere theory are described under the bijection (\ref{202606181948}).

\begin{Defn} \label{202604121701}
    An {\it $m$-partition of $n$} is an $m$-tuple $(n_1,\cdots,n_m) \in \mathds{N}^{m}$ such that $n_1+\cdots + n_m = n$.
    The {\it map associated with the partition}, $\rho : {\bf n} \to {\bf m}$, is defined by $\rho (i) {:=} l$ if $\sum^{l-1}_{j=1} n_j < i \leq \sum^{l}_{j=1} n_j$.
\end{Defn}

\begin{Defn} \label{202606181952}
    For $S \subset {\bf n}$ with $|S| =r$, we define $P_{n,S} \in \mathcal{C} (n,r)$ to be the product whose $i$-th component is $\mathrm{id}_1$ for $i \in S$ and the unique morphism otherwise. Clearly, $P_{n,\{k\}} = p_{n,k}$.
\end{Defn}
Let $(n_1 , \cdots ,n_m )$ be an $m$-partition of $n$ and $\rho_0$ be the map associated with the partition.
Let $g_k \in \mathcal{C}(n_k,1)$.
Then, by the universality of products in $\mathcal{C}$, we have
\begin{align} \label{202606182102}
    g_1 \times \cdots \times g_m = (f_1, \cdots, f_m) ,
\end{align}
where $f_k = g_k \circ P_{n,\rho_0^{-1} (k)}$.

\subsection{The functor $\mathfrak{S} \to \mathcal{C}$}

In this section, we investigate the functor $\mathfrak{S} \to \mathcal{C}$ induced by permutations of the object $1 \in \mathcal{C}$.
More precisely, it is defined as follows.

\begin{Defn} \label{202603261255}
    For $n \in \mathds{N}$ and $\sigma \in \mathfrak{S}_n$, we define $$E_\sigma {:=} (p_{n,\sigma^{-1}(1)}, \cdots, p_{n,\sigma^{-1}(n)})  \in \mathcal{C} (n,n)$$ via the identification (\ref{202606181948}). 
    These morphisms determine a functor $E : \mathfrak{S} \to \mathcal{C}$ which is the identity on objects and which sends $\sigma \in \mathfrak{S}_n$ to $E_{\sigma}$.
\end{Defn}

Via this functor, we regard $\mathcal{C}(n,m)$ as a $(\mathfrak{S}_m, \mathfrak{S}_n)$-set.
Hence, $\mathtt{L}_{\mathcal{C}}(m,n)$ is a $(\mathfrak{S}_m,\mathfrak{S}_n)$-bimodule.

Since a category with products is canonically a symmetric monoidal category, the above construction can be described in terms of the symmetry isomorphisms.

\begin{Defn} \label{202603281205}
    Let $(n_1,\cdots,n_m)$ be an $m$-partition of $n$, and $\tau \in \mathfrak{S}_m$.
    We denote by $\tau_{n_1,\cdots,n_m} \in \mathfrak{S}_n$ the permutation obtained by permuting the consecutive blocks of lengths $n_1, \cdots, n_m$, according to $\tau$, while preserving the order within each block. 
\end{Defn}

The following remark is intended to clarify the meaning of the above definition.
\begin{remark}
    Let $\rho$ be the map associated with the partition $(n_1,\cdots,n_m)$ and $\rho^\prime$ the map associated with the permuted partition $(n_{\tau^{-1} (1)}, \cdots, n_{\tau^{-1}(m)})$.
    Then $\tau_{n_1,\cdots,n_m} \in \mathfrak{S}_n$ is the unique bijection such that for $k \in \mathbf{m}$, the restriction $\tau_{n_1,\cdots,n_m} \mid_{\rho^{-1}(k)} : \rho^{-1}(k) \to {\bf n}$ factors into an order-preserving map $\rho^{-1}(k) \to (\rho^\prime)^{-1} (\tau (k))$.
\end{remark}

The following is immediate:
\begin{Lemma}
    For $\sigma,\tau\in\mathfrak{S}_m$, we have $\sigma_{n_{\tau(1)},\cdots,n_{\tau(m)}}  \circ \tau_{n_1,\cdots,n_m} = (\sigma \tau)_{n_1,\cdots,n_m}$.
\end{Lemma}

\begin{Lemma} \label{202606182054}
    Let $(n_1,\cdots,n_m)$ be an $m$-partition of $n$.
    For $f_k \in \mathcal{C}(1,n_k), ~k\in {\bf m}$ and $\tau \in \mathfrak{S}_m$, we have
    $$
    (f_1 \times \cdots \times f_m) \circ E_\tau = E_{\tau_{n_1,\cdots,n_m}} \circ ( f_{\tau^{-1}(1)} \times \cdots \times f_{\tau^{-1}(m)}) .
    $$
    An analogous statement holds for $g_k \in \mathcal{C} (n_k,1), ~ k \in {\bf m}$.
\end{Lemma}
\begin{proof}
    By the above lemma, it suffices to prove the statement when $\tau$ is a transposition.
    Since products endow $\mathcal{C}$ with a symmetric monoidal category structure, the desired equality follows from the naturality of the symmetry.
\end{proof}

\subsection{The condition (ZM*)}
\label{202606181932}

In this section, we introduce some structural conditions on Lawvere theories:
(Z) the category $\mathcal{C}$ has a zero object; (M) the functor assigning $\mathcal{C}_n = \mathcal{C}(n,1)$ to $n \in \mathrm{Obj}(\mathcal{C}) = \mathds{N}$ is endowed with a specified factorization through the forgetful functor $U: \mathsf{Mon} \to \mathsf{Set}$ where $\mathsf{Mon}$ denotes the category of monoids and $\mathsf{Set}$ denotes the category of sets:
$$
\begin{tikzcd}
    & \mathsf{Mon} \ar[d,"U"] \\
    \mathcal{C}^{\mathsf{o}} \ar[r, "\mathcal{C}(-\mbox{,}1)"'] \ar[ur, "\exists"] & \mathsf{Set}
\end{tikzcd}
$$
and (M*) for each $n$, the monoid $\mathcal{C}_n$ is generated by the union $\bigcup^{n}_{k=1} p_{n,k}^\ast ( \mathcal{C}_1 )$.
The condition (M) is equivalent to the object $1$ being endowed with a specified monoid object structure.

\begin{notation} \label{202512162232}
    When $\mathcal{C}$ satisfies the condition (M), we denote by $f \star g$ the monoid product of $f,g \in \mathcal{C}_n$.
    Moreover, we use the same notation for the induced product of the monoid algebra $\mathds{k}[\mathcal{C}_n^{\times m}]$.
    When $\mathcal{C}$ satisfies the condition (Z), we denote by $e \in \mathcal{C}_n$ the zero morphism.
\end{notation}

A Lawvere theory $\mathcal{C}$ satisfies the condition (ZM*) if all the above conditions hold.
More strongly, we say that $\mathcal{C}$ satisfies the condition (ZCM*) if it satisfies (ZM*) and each monoid $\mathcal{C}_n$ is a commutative monoid.

\begin{Example}
    The Lawvere theories $\mathbf{W}^{\mathsf{o}}$ and $\G^{\mathsf{o}}$ satisfy (ZM*).
\end{Example}

The category $\mathbf{W}^{\mathsf{o}}$ satisfies a universal property among Lawvere theories $\mathcal{C}$ subject to the condition (M).
More precisely, since the Lawvere theory $\mathbf{W}^{\mathsf{o}}$ is generated by the monoid object $1$, we have a unique functor $\mathbf{W}^{\mathsf{o}} \to \mathcal{C}$ which preserves finite products and the monoid object $1$.

\begin{Example} \label{202606091733}
    The category of free right $R$-modules of finite rank, denoted by $\mathbf{M}_{R}$, further satisfies the condition (ZCM*).
    More precisely, we identify the object $R^n \in \mathbf{M}_{R}$ freely generated by $n$ elements with $n$, so that we regard $\mathbf{M}_{R}$ as a category with objects $\mathds{N}$.
    Moreover, we identify the morphism set $\mathbf{M}_{R} (n,m)$ with    $\mathrm{M}_{m,n} (R)$ the set of matrices with $m$ rows and $n$ columns, consisting of elements of $R$.
\end{Example}

\begin{Example}
    There exists an example that satisfies (Z) but not (M).
    Let $\mathsf{Fin}_\ast$ be the category of pointed finite sets.
    The opposite category $\mathsf{Fin}_\ast^{\mathsf{o}}$ is a Lawvere theory satisfying (Z) while the object $1$ of $\mathsf{Fin}_\ast^{\mathsf{o}}$ does not admit any monoid object structure.
\end{Example}

There are some variants of the category $\G$ as follows.
\begin{Defn} \label{202512311144}
    A {\it radical functor for groups} is a functor that assigns a normal subgroup $\mathcal{R} (G) \subset G$ to each group $G$ such that $\mathcal{R} (G/\mathcal{R}(G)) \cong 1$.
    We define $|\mathcal{R}| \in\mathds{N}$ to be the nonnegative integer $r$ such that $\mathcal{R}(\mathds{Z}) = r \mathds{Z}$.
\end{Defn}

Examples of $\mathcal{R}$ include $\gamma_{c}$ the $c$-th component of the lower central series (see \cite{kim2025poly} for further examples).
Obviously, $|\gamma_{c}| = 0$.

\begin{Defn}
    Let $\mathcal{R}$ be a radical functor for groups.
    A group $G$ is {\it $\mathcal{R}$-semisimple} if $\mathcal{R}(G) \cong 1$.
    We define $\G_{\mathcal{R}}$ to be the category of finitely generated free $\mathcal{R}$-semisimple groups.
\end{Defn}
The opposite category $\G_{\mathcal{R}}^{\mathsf{o}}$ is a Lawvere theory satisfying the condition (ZM*).
For instance, if $\mathcal{R} \equiv 1$, then $\G_{\mathcal{R}} = \G$.
If $\mathcal{R} = \gamma_{c+1}$, then $\G_{\mathcal{R}}$ is the category of finitely generated free nilpotent groups of class $\leq c$.
We write $\mathbf{N}_{c}$ for $\G_{\gamma_{c+1}}$.

\section{The primitivity ideal of Lawvere theories}
\label{202512162136}

In this section, we introduce the notion of the primitivity ideal of a Lawvere theory $\mathcal{C}$.
Using the results of Part \ref{202512151906}, we present the associated eigenmonad and the eigenmonad adjunction.
We also introduce the notion of a primitively generated $\mathtt{L}_{\mathcal{C}}$-module.

\subsection{Definition and basic properties} \label{202512141943}

In this section, we introduce the primitivity ideal and prove several basic properties.
Let $\mathcal{C}$ be a Lawvere theory having a zero object.
\begin{Defn} \label{202408081349}
    Let $\Delta \in \mathcal{C}(1,2)$ be the diagonal morphism and $\eta \in \mathcal{C}(0,1)$ be the zero morphism.
    For $n \in \mathds{N}^\ast$ and $k \in {\bf n}$, we define $\bar{\Delta}^{n,k} \in \mathtt{L}_{\mathcal{C}} (n+1,n)$ by
    $$
    \bar{\Delta}^{n,k} {:=} \mathrm{id}_{k-1} \times \Delta \times \mathrm{id}_{n-k} - \mathrm{id}_{k} \times \eta \times \mathrm{id}_{n-k} - \mathrm{id}_{k-1} \times \eta \times \mathrm{id}_{n-k+1}  .
    $$
    We write $\bar{\Delta}$ for $\bar{\Delta}^{1,1}$.
    
    We define the {\it primitivity ideal of $\mathcal{C}$} to be the left ideal of $\mathtt{L}_{\mathcal{C}}$ generated by $\bar{\Delta}^{n,k} $ for $n\in \mathds{N}^\ast$ and $k \in {\bf n}$.
    We denote this by $\mathtt{I}^{\mathsf{pr}}_{\mathcal{C}}$.
    In other words, for $n,m \in \mathds{N}$, the $\mathds{k}$-submodule $\mathtt{I}^{\mathsf{pr}}_{\mathcal{C}} (m,n) \subset \mathtt{L}_{\mathcal{C}} (m,n)$ is generated by $f \circ \bar{\Delta}^{n,k} $ for $f \in \mathtt{L}_{\mathcal{C}}(m,n+1)$ and $k \in {\bf n}$.
\end{Defn}

\begin{Example} \label{202604041521}
    \begin{enumerate}
        \item For $m \in \mathds{N}$, $\mathtt{I}^{\mathsf{pr}}_{\mathcal{C}} (m,0) = 0$.
        \item For $n \in \mathds{N}^\ast$, $\mathtt{I}^{\mathsf{pr}}_{\mathcal{C}} (0,n) = \mathtt{L}_{\mathcal{C}} (0,n) \cong \mathds{k}$.
        Indeed, $\mathtt{L}_{\mathcal{C}} (0,n) = \mathds{k} [\mathcal{C}(n,0)]$ has one generator $e_{n} \in \mathcal{C} (n,0)$, the zero morphism, and we have $- e_{n+1} \circ \bar{\Delta}^{n,k} = e_n$.
    \end{enumerate}
\end{Example}

\begin{Lemma} \label{202601051527}
    The ideal $\mathtt{I}^{\mathsf{pr}}_{\mathcal{C}} \subset \mathtt{L}_{\mathcal{C}}$ is a monoidal ideal in the sense that $$f_1 \times u \times f_2 \in \mathtt{I}^{\mathsf{pr}}_{\mathcal{C}} (m_1+m+m_2,n_1+n+n_2)$$ for $f_i \in \mathtt{L}_{\mathcal{C}} (m_i,n_i),~i\in\{1,2\}$ and $u \in \mathtt{I}^{\mathsf{pr}}_{\mathcal{C}} (m,n)$.
\end{Lemma}
\begin{proof}
    We may put $u = v \circ \bar{\Delta}^{n,k}$ for some $k \in {\bf n}$ and $v$.
    By definition, we have $f_1 \times u \times f_2 = (f_1 \times v \times f_2) \circ \bar{\Delta}^{n_1+n+n_2,k+n_1}$, which lies in the ideal $\mathtt{I}^{\mathsf{pr}}_{\mathcal{C}}$.
\end{proof}

Recall that $\mathtt{L}_{\mathcal{C}}(m,n)$ is a $(\mathfrak{S}_m,\mathfrak{S}_n)$-bimodule by the construction of $E_{\sigma}$ in Definition \ref{202603261255}.
\begin{Lemma} \label{202509071205}
    The submodule $\mathtt{I}^{\mathsf{pr}}_{\mathcal{C}} (m,n) \subset \mathtt{L}_{\mathcal{C}} (m,n)$ is a sub-$(\mathfrak{S}_m,\mathfrak{S}_n)$-bimodule.
\end{Lemma}
\begin{proof}
    Since $\mathtt{I}^{\mathsf{pr}}_{\mathcal{C}}  \subset \mathtt{L}_{\mathcal{C}}$ is a left ideal, $\mathtt{I}^{\mathsf{pr}}_{\mathcal{C}} (m,n) \subset \mathtt{L}_{\mathcal{C}} (m,n)$ is closed under the left $\mathfrak{S}_m$-action.
    Let $\sigma \in \mathfrak{S}_n$.
    Lemma \ref{202606182054} implies that, for $k\in {\bf n}$, $\bar{\Delta}^{n,k} \circ E_{\sigma} = E_{\sigma^\prime} \circ \bar{\Delta}^{n,\sigma^{-1}(k)}$ where $\sigma^\prime = \sigma_{1,\cdots,2,\cdots,1} \in \mathfrak{S}_{n+1}$.
    Here, $(1,\cdots,1,2,1,\cdots,1)$ is the $n$-partition of $n+1$ whose $\sigma^{-1} (k)$-th component is $2$.
    Thus, $\bar{\Delta}^{n,k} \circ E_{\sigma}$ is contained in the primitivity ideal.
    This proves that $\mathtt{I}^{\mathsf{pr}}_{\mathcal{C}} (m,n) \subset \mathtt{L}_{\mathcal{C}} (m,n)$ is closed under the right $\mathfrak{S}_n$-action.
\end{proof}

In what follows, we exhibit some elements of the ideal $\mathtt{I}^{\mathsf{pr}}_\mathcal{C}$.
Recall Definitions \ref{202604091226} and \ref{202606181952}.
It is obvious that
$$
P_{n,{\bf n}\backslash \{k\}} = (p_{n,1} , p_{n,2} , \cdots, p_{n,k-1} , p_{n,k+1} , \cdots  , p_{n,n} ) .
$$
\begin{Lemma} \label{202604121607}
    For $k\in{\bf n}$, we have $P_{n,{\bf n}\backslash \{k\}} \in \mathtt{I}^{\mathsf{pr}}_{\mathcal{C}}(n-1,n)$.
\end{Lemma}
\begin{proof}
    Using this identity (\ref{202606201739}) and Lemma \ref{202509042125}, we obtain
    \begin{align*}
        &(\mathrm{id}_{k-1} \times \Delta \times \mathrm{id}_{n-k})^\ast (p_{n+1,1} , p_{n+1,2} , \cdots, p_{n+1,k-1} , p_{n+1,k+2} , \cdots  , p_{n+1,n+1} ) \\
        =& (p_{n,1} , p_{n,2} , \cdots, p_{n,k-1} , p_{n,k+1} , \cdots  , p_{n,n} ) = P_{n,{\bf n}\backslash \{k\}} .
    \end{align*}
    Via the identification in (\ref{202606181948}), we have
    $$
    \mathrm{id}_{l} \times \eta \times \mathrm{id}_{n-l} = (p_{n,1},\cdots,p_{n,l}, e, p_{n,l+1}, \cdots, p_{n,n})
    $$
    where $e \in \mathcal{C}_n$ is the zero morphism.
    Applying this to $l\in\{k-1,k\}$, we obtain
    \begin{align*}
        &\{ (\mathrm{id}_{k-1} \times \eta \times \mathrm{id}_{n-k+1})^\ast +  (\mathrm{id}_{k} \times \eta \times \mathrm{id}_{n-k})^\ast \}  (p_{n+1,1} , p_{n+1,2} , \cdots, p_{n+1,k-1} , p_{n+1,k+2} , \cdots  , p_{n+1,n+1} )  \\
        =& (p_{n,1} , p_{n,2} , \cdots, p_{n,k-1} , p_{n,k+1} , \cdots  , p_{n,n} ) + (p_{n,1} , p_{n,2} , \cdots, p_{n,k-1} , p_{n,k+1} , \cdots  , p_{n,n} ) = 2 P_{n,{\bf n}\backslash \{k\}} .
    \end{align*}
    By the definition of $\mathtt{I}^{\mathsf{pr}}_{\mathcal{C}}(n-1,n)$, the above results imply that $P_{n,{\bf n}\backslash \{k\}} \equiv 0 \mod{\mathtt{I}^{\mathsf{pr}}_{\mathcal{C}}(n-1,n)}$.
\end{proof}

\begin{Defn} \label{202512090939}
    Let $\mathcal{C}_n^{(k)}$ be the image of the map
    $$
    P_{n,{\bf n}\backslash \{k\}}^\ast : \mathcal{C}_{n-1} \to \mathcal{C}_n .
    $$
\end{Defn}

\begin{Example}
    Recall the notation in Example \ref{202603121546}.
    For $\mathcal{C} = \mathbf{W}^{\mathsf{o}}$, $\mathcal{C}_n^{(k)}$ is the submonoid of $\mathsf{W}_n$ generated by $x_1, \cdots, x_{k-1},x_{k+1},\cdots,x_n$.
    For $\mathcal{C} = \G^{\mathsf{o}}$, $\mathcal{C}_n^{(k)}$ is the subgroup of $\mathsf{F}_n$ generated by $x_1, \cdots, x_{k-1},x_{k+1},\cdots,x_n$.
\end{Example}

\begin{Example} \label{202606191434}
    For distinct $i,j \in {\bf n}$, we have $p_{n,j} \in \mathcal{C}_n^{(i)}$.
    In fact, the map in Definition \ref{202512090939} sends $p_{n-1,l}$ to $p_{n,l}$ if $l < i$, and to $p_{n,l+1}$ if $l \geq i$.
\end{Example}

\begin{prop} \label{202509042136}
    Let $\mathcal{C}$ be a Lawvere theory with a zero object.
    For $(w_1, \cdots, w_m) \in \mathcal{C} (n,m)$, if there exists $k \in {\bf n}$ such that $\{ w_1, \cdots, w_m \} \subset \mathcal{C}_n^{(k)}$, then $(w_1, \cdots, w_m) \in \mathtt{I}^{\mathsf{pr}}_{\mathcal{C}} (m,n)$.
\end{prop}
\begin{proof}
    Since $\{ w_1, \cdots, w_m \} \subset \mathcal{C}_n^{(k)}$, there exists $\rho^\prime \in \mathcal{C} (n-1,m)$ such that $\rho = \rho^\prime \circ P_{n,{\bf n}\backslash \{k\}}$.
    By Lemma \ref{202604121607}, we obtain $\rho = \rho^\prime \circ P_{n,{\bf n}\backslash \{k\}} \in \mathtt{I}^{\mathsf{pr}}_{\mathcal{C}} (m,n)$, since $\mathtt{I}^{\mathsf{pr}}_{\mathcal{C}}$ is a left ideal of $\mathtt{L}_{\mathcal{C}}$.
\end{proof}

\subsection{The eigenmonad of the primitivity ideal} \label{202606091353}

Let $\mathcal{C}$ be a Lawvere theory with a zero object, so that the primitivity ideal $\mathtt{I}^{\mathsf{pr}}_{\mathcal{C}}$ is defined.
In this section, we record several general properties of the eigenmonad associated with $(\mathtt{L}_\mathcal{C}, \mathtt{I}^{\mathsf{pr}}_{\mathcal{C}})$, obtained by specializing the results of Section \ref{202410161744} to $\mathtt{I}^{\mathsf{pr}}_{\mathcal{C}}$.

\begin{notation} \label{202511281314}
    Let $\Phi_{\mathcal{C}}$ denote the eigenmonad $\mathrm{E}_{\mathtt{L}_{\mathcal{C}}} ( \mathtt{I}^{\mathsf{pr}}_{\mathcal{C}})$, or equivalently the endomorphism monad $\mathrm{End}_{\mathtt{L}_{\mathcal{C}}} ( \mathtt{L}_{\mathcal{C}} / \mathtt{I}^{\mathsf{pr}}_{\mathcal{C}})$ (see Proposition \ref{202401241330}).
\end{notation}

As in Definition \ref{202403061037}, a monad and a left ideal determine a canonical bimodule.
Specializing to the monad $\mathtt{L}_\mathcal{C}$ with the ideal $\mathtt{I}^{\mathsf{pr}}_{\mathcal{C}}$, we obtain {\it the canonical $(\mathtt{L}_{\mathcal{C}}, \Phi_\mathcal{C})$-bimodule} $\mathtt{L}_{\mathcal{C}}/ \mathtt{I}^{\mathsf{pr}}_{\mathcal{C}}$.

Both $\mathtt{L}_{\mathcal{C}}$ and $\Phi_{\mathcal{C}}$ admit canonical monad maps from $\mathtt{L}_{\mathfrak{S}}$.
The map $\mathfrak{S}_n \to \mathcal{C}(n,n);  \sigma \mapsto E_\sigma$ (see Definition \ref{202603261255}) induces a monad map $\mathtt{L}_{\mathfrak{S}} \to \mathtt{L}_{\mathcal{C}}$.
The following lemma shows that this further induces a monad map $\mathtt{L}_{\mathfrak{S}} \to \Phi_{\mathcal{C}}$:
\begin{Lemma} \label{202603211051}
    The map $\mathtt{L}_{\mathfrak{S}} \to \mathtt{L}_{\mathcal{C}}$ factors through the idealizer monad $\mathrm{Iz} (\mathtt{I}^{\mathsf{pr}}_{\mathcal{C}})  \subset \mathtt{L}_{\mathcal{C}}$.
\end{Lemma}
\begin{proof}
    It suffices to show that $\mathtt{I}^{\mathsf{pr}}_{\mathcal{C}} (-,n) \circ E_{\sigma} \subset \mathtt{I}^{\mathsf{pr}}_{\mathcal{C}} (-,n)$ for $\sigma \in \mathfrak{S}_n$.
    This follows from Lemma \ref{202509071205}.
\end{proof}
Via the maps obtained above, we often regard $\mathtt{L}_{\mathcal{C}}/ \mathtt{I}^{\mathsf{pr}}_{\mathcal{C}}$ as $(\mathtt{L}_{\mathfrak{S}}, \mathtt{L}_{\mathfrak{S}})$-bimodule, which will be studied in detail in Section \ref{202604131411}.

The construction $\mathcal{C} \mapsto \Phi_{\mathcal{C}}$ is natural in the following sense.

\begin{prop} \label{202512071739}
    Let $\mathcal{C},\mathcal{D}$ be Lawvere theories having zero objects.
    Let $\pi : \mathcal{C} \to \mathcal{D}$ be a full functor which is the identity on objects and preserves finite products.
    \begin{enumerate}
        \item We have a natural monad map $\pi : \Phi_{\mathcal{C}} \to \Phi_{\mathcal{D}}$.
        \item  The map $\mathtt{L}_{\mathcal{C}}/ \mathtt{I}^{\mathsf{pr}}_{\mathcal{C}} \to \mathtt{L}_{\mathcal{D}}/ \mathtt{I}^{\mathsf{pr}}_{\mathcal{D}}$ is an $(\mathtt{L}_{\mathcal{C}}, \Phi_\mathcal{C})$-bimodule map where $\mathtt{L}_{\mathcal{D}}/ \mathtt{I}^{\mathsf{pr}}_{\mathcal{D}}$ is regarded as an $(\mathtt{L}_{\mathcal{C}}, \Phi_\mathcal{C})$-bimodule via the monad maps $\pi: \mathtt{L}_{\mathcal{C}} \to \mathtt{L}_{\mathcal{D}}$ and $\pi : \Phi_{\mathcal{C}} \to \Phi_{\mathcal{D}}$.
    \end{enumerate}
\end{prop}
\begin{proof}
    By the hypothesis on $\pi$, the induced monad map $\pi: \mathtt{L}_{\mathcal{C}} \to \mathtt{L}_{\mathcal{D}}$ is an epimorphism which assigns $\bar{\Delta}^{n,k}$ in $\mathtt{L}_{\mathcal{D}} (n+1,n)$ to $\bar{\Delta}^{n,k}$ in $\mathtt{L}_{\mathcal{C}} (n+1,n)$.
    Thus, $\pi ( \mathtt{I}^{\mathsf{pr}}_{\mathcal{C}} ) = \mathtt{I}^{\mathsf{pr}}_{\mathcal{D}}$.
    By Proposition \ref{202512071648}, we obtain the natural monad map in the statement (1).
    The second assertion is clear from the definition of the canonical bimodules.
\end{proof}

We obtain the following corollary from Proposition \ref{202409261707}.
This adjunction plays a fundamental role in establishing the main results of the present paper.

\begin{Corollary} \label{202512081855}
    The eigenmonad adjunction associated with $(\mathtt{L}_{\mathcal{C}} , \mathtt{I}^{\mathsf{pr}}_{\mathcal{C}})$ yields 
    $$
    \begin{tikzcd}
             \mathbb{L} : \Phi_{\mathcal{C}}\mbox{-}\mathsf{Mod} \arrow[r, shift right=1ex, ""{name=G}] & \mathtt{L}_{\mathcal{C}}\mbox{-}\mathsf{Mod} :  \mathbb{R} \arrow[l, shift right=1ex, ""{name=F}]
            \arrow[phantom, from=G, to=F, , "\scriptscriptstyle\boldsymbol{\top}"].
    \end{tikzcd}
    $$
    These adjoint functors satisfy the following properties:
    \begin{enumerate}
        \item $\mathbb{L} = \mathtt{L}_{\mathcal{C}}/\mathtt{I}^{\mathsf{pr}}_{\mathcal{C}} \otimes_{\Phi_{\mathcal{C}}} (-)$
        \item $\mathbb{R} = \mathrm{Hom}_{\mathtt{L}_{\mathcal{C}}} (\mathtt{L}_{\mathcal{C}}/\mathtt{I}^{\mathsf{pr}}_{\mathcal{C}} , -) \cong \mathrm{V} (- ; \mathtt{I}^{\mathsf{pr}}_{\mathcal{C}})$ and $\mathbb{R}((\mathtt{L}_{\mathcal{C}}/\mathtt{I}^{\mathsf{pr}}_{\mathcal{C}}) (-,n)) = \Phi_{\mathcal{C}} (-,n)$.
    \end{enumerate}
\end{Corollary}
\begin{proof}
    The adjunction is immediate from (1) of Proposition \ref{202409261707}.
    The equivalent descriptions of the right adjoint functor follows from Lemma \ref{202405241347}.
    The last claim $\mathbb{R}((\mathtt{L}_{\mathcal{C}}/\mathtt{I}^{\mathsf{pr}}_{\mathcal{C}}) (-,n)) = \Phi_{\mathcal{C}} (-,n)$ follows from Proposition \ref{202401241330}.
\end{proof}

\begin{Defn}
    A left $\mathtt{L}_{\mathcal{C}}$-module is {\it primitively generated}, or {\it primitive}, if it is $\mathtt{I}^{\mathsf{pr}}_{\mathcal{C}}$-vanishingly generated in the sense of Definition \ref{202402011603}.
    In other words, a left $\mathtt{L}_{\mathcal{C}}$-module $\mathtt{M}$ is primitively generated if the counit associated with the adjunction in Corollary \ref{202512081855} $\mathtt{L}_{\mathcal{C}}/\mathtt{I}^{\mathsf{pr}}_{\mathcal{C}} \otimes_{\Phi_{\mathcal{C}}} \mathrm{V} (\mathtt{M}; \mathtt{I}^{\mathsf{pr}}_{\mathcal{C}}) \to \mathtt{M}$ is an epimorphism.
    We denote by $\mathtt{L}_{\mathcal{C}} \mbox{-}\mathsf{Mod}^{\mathsf{prim}}$ the full subcategory of primitively generated left $\mathtt{L}_{\mathcal{C}}$-modules.
\end{Defn}

\begin{remark}
    This concept generalizes the notion of the same name in \cite{kim2024analytic}.
\end{remark}

\subsection{Some components of the eigenmonad}

As fundamental examples, we compute several components of $\Phi_{\mathcal{C}}$ .

\begin{Lemma}
    For $m,n\in\mathds{N}$ such that $mn=0$, we have
    \begin{align*}
        \Phi_\mathcal{C} (m,n) \cong 
        \begin{cases}
            \mathds{k} & \mathrm{if~}m=n=0,\\
            0 & \mathrm{otherwise}
        \end{cases}
    \end{align*}
\end{Lemma}
\begin{proof}
    Since $0$ is a zero object of $\mathcal{C}$, $\mathtt{L}_{\mathcal{C}}(m,0) \cong \mathds{k}$ is generated by the morphism $\eta^{\times m}$.
    Recall from Example \ref{202604041521} that $\mathtt{I}^{\mathsf{pr}}(r,0) \cong 0$.
    These imply that $\mathrm{Iz} (\mathtt{I}^{\mathsf{pr}}_{\mathcal{C}} ) (0,0) = \mathds{k}$, so $\Phi_\mathcal{C} (0,0) \cong \mathds{k}$.
    For $m \in \mathds{N}^\ast$ and $k \in {\bf m}$, $\bar{\Delta}^{m,k} \circ \eta^{\times m} = - \eta^{\times m+1}$, so $\left( \mathrm{Iz} (\mathtt{I}^{\mathsf{pr}}_{\mathcal{C}}) \right) (m,0) = 0$ which proves $\Phi_\mathcal{C} (m,0) \cong 0$.
    By Example \ref{202604041521}, we have $\mathtt{L}_\mathcal{C}/\mathtt{I}^{\mathsf{pr}}_{\mathcal{C}}(0,n) \cong 0$ for $n \in \mathds{N}^\ast$.
    Therefore, $\Phi_\mathcal{C} (0,n) \cong 0$.
\end{proof}

\begin{Lemma} \label{202601081632}
    $\Phi_{\mathcal{C}}(1,1) = \mathtt{L}_{\mathcal{C}} (1,1) / \mathtt{I}^{\mathsf{pr}}_{\mathcal{C}} (1,1)$.
\end{Lemma}
\begin{proof}
    By the universality of products and the hypothesis that $0 \in \mathcal{C}$ is a zero object, we deduce $\bar{\Delta} \circ f = (f \times f ) \circ \bar{\Delta}$ for $f\in \mathcal{C}(1,1)$.
    This leads to $\left( \mathrm{Iz} (\mathtt{I}^{\mathsf{pr}}_{\mathcal{C}}) \right) (1,1) = \mathtt{L}_{\mathcal{C}}(1,1)$, so the statement follows.
\end{proof}

\begin{prop}\label{202509091604}
    For $\mathcal{C} \in \{ \G^{\mathsf{o}}, \mathbf{W}^\mathsf{o} \}$, we have a $\mathds{k}$-algebra isomorphism $\mathds{k} \cong \Phi_{\mathcal{C}}(1,1)$.
\end{prop}
\begin{proof}
    Note that the $\mathds{k}$-module $\mathtt{L}_{\G^{\mathsf{o}}} (1,1) = \mathds{k}[\mathsf{F}_1]$ has a basis consisting of $x_1^{a}$ for $a \in \mathds{Z}$.
    By Lemma \ref{202601081632}, we have $\Phi_{\G^{\mathsf{o}}}(1,1) \cong \mathtt{L}_{\G^{\mathsf{o}}} (1,1) / \mathtt{I}^{\mathsf{pr}}_{\G^{\mathsf{o}}} (1,1)$, so it suffices to investigate $\mathtt{I}^{\mathsf{pr}}_{\G^{\mathsf{o}}} (1,1)$.
    By definition, this is generated by $g \circ \bar{\Delta}$ for $g \in \G^{\mathsf{o}}(2,1)$.
    By $\G^{\mathsf{o}}(2,1) \cong \mathsf{F}_2$, the morphism $g$ is represented by $x_1^{a_1}x_2^{b_1}x_1^{a_2}x_2^{b_2}\cdots x_1^{a_k}x_2^{b_k}$ for a finite set $\{ a_1,b_1, \cdots, a_k,b_k\} \subset \mathds{Z}$.
    Since we have
    $$
    g \circ \bar{\Delta} = x_1^{\sum a_j+\sum b_j} - x_1^{\sum a_j} - x_1^{\sum b_j} ,
    $$
    the $\mathds{k}$-module $\mathtt{I}^{\mathsf{pr}}_{\G^{\mathsf{o}}} (1,1)$ is generated by $x_1^{n+m} - x_1^{n} - x_1^{m}$ for $n,m \in \mathds{Z}$.
    Thus, $\mathds{k} \to \mathtt{L}_{\G^{\mathsf{o}}} (1,1) / \mathtt{I}^{\mathsf{pr}}_{\G^{\mathsf{o}}} (1,1); 1 \mapsto x_1$ gives a $\mathds{k}$-module isomorphism.
    This map is a $\mathds{k}$-algebra isomorphism, since $x_1 \in \mathsf{F}_1 \cong \G^{\mathsf{o}}(1,1)$ gives the identity.

    The above proof also works for $\mathcal{C} = \mathbf{W}^{\mathsf{o}}$ except that the exponents of $x_i$'s (the generators of the free monoids in this case) are required to be non-negative.
\end{proof}

\section{Linear species} \label{202512181624}

When studying the canonical $(\mathtt{L}_{\mathcal{C}}, \Phi_\mathcal{C})$-bimodule associated with a Lawvere theory $\mathcal{C}$ with a zero object, right $\mathtt{L}_{\mathfrak{S}}$-modules naturally arise.
To study them, it is convenient to work with linear species via their equivalence to right $\mathtt{L}_{\mathfrak{S}}$-modules.
In this section, we collect several classical constructions related to linear species.
We also discuss comonoid species and give several technical results that will be used in later sections.

\subsection{$\mathds{k}$-linear species and right $\mathtt{L}_{\mathfrak{S}}$-modules} \label{202604091335}

In this section, we briefly recall the notions of linear species.
In this paper, a {\it $\mathds{k}$-linear species} is a functor
$F : \mathrm{Bij}^{\mathsf{o}} \to \mathds{k}\mbox{-}\mathsf{Mod}$, where $\mathrm{Bij}$ denotes the category of finite sets and bijections.
Our definition differs slightly from the standard convention found in references such as \cite{aguiar2010monoidal} where a species is a functor from $\mathrm{Bij}$. This choice is made for convenience; since $\mathrm{Bij} \cong \mathrm{Bij}^{\mathsf{o}}$, there is no essential difference between the two approaches.

\begin{Defn}
    Let $\mathsf{Sp}_\mathds{k}$ denote the category of $\mathds{k}$-linear species and natural transformations.
\end{Defn}

The category $\mathsf{Sp}_\mathds{k}$ is an abelian category in which kernels and cokernels are computed objectwise.

Recall from Example \ref{202604091040} the category $\mathfrak{S}$ and the associated monad $\mathtt{L}_{\mathfrak{S}}$ on $\mathds{N}$.
Note that $\mathfrak{S}$ can be viewed as a full subcategory of $\mathsf{Bij}$.
A $\mathds{k}$-linear species naturally gives rise to a right $\mathtt{L}_{\mathfrak{S}}$-module.
Conversely, a right $\mathtt{L}_{\mathfrak{S}}$-module $\mathtt{M}$ naturally extends to a $\mathds{k}$-linear species via the construction $F (X) {:=} \mathtt{M} (n) \otimes_{\mathfrak{S}_n} \mathds{k} \mathrm{Bij} ( X, {\bf n})$ where $n$ is the cardinality of $X$.
These constructions yield an equivalence of categories:
\begin{align} \label{202603251315}
    \mathsf{Sp}_{\mathds{k}}  \simeq \mathsf{Mod}\mbox{-}\mathtt{L}_{\mathfrak{S}}.
\end{align}
In this paper, we do not distinguish right $\mathtt{L}_{\mathfrak{S}}$-modules from $\mathds{k}$-linear species based on this equivalence.
\begin{notation}
    For a $\mathds{k}$-linear species $F$, we use the same notation for the corresponding right $\mathtt{L}_{\mathfrak{S}}$-module, writing $F (n) $ for $F( {\bf n})$, where ${\bf n} {:=} \{ 1,2,\cdots, n\}$.
\end{notation}

\begin{Example} \label{202603251354}
    For a $\mathds{k}$-module $V$, let $I_V$ be the $\mathds{k}$-linear species defined by $I_V (X) {:=} V$ if $|X|=1$ and $I_V(X) {:=} 0$ otherwise.
\end{Example}

\begin{Example} \label{202606112022}
    The $\mathds{k}$-linear species corresponding to the right $\mathtt{L}_{\mathfrak{S}}$-module $\mathds{k}^{\otimes}$ introduced in Example  \ref{202604062227} is the {\it exponential} species
   which we denote by $\mathsf{Exp}$.
    We shall freely use both notations.
    More generally, the linear species induced by $V^{\otimes}$, if we use the same notation, is nothing but the assignment of $V^{\otimes} (X) = V^{\otimes X}$ to any finite set $X$.
\end{Example}

\begin{Example} 
    Let $\mathsf{Lin}$ denote the {\it linear-order} species: $\mathsf{Lin} (X)$ is the free $\mathds{k}$-module generated by the set of linear-orders on $X$.  
    By convention, $\mathsf{Lin}(\emptyset) = \mathds{k}$.
    There is a canonical map $\mathsf{Lin} \to \mathsf{Exp}$ which sends any orders to the unit of $\mathds{k}$.
\end{Example}

We will make use of the following notation when we restrict $\mathsf{Lin}$ to subsets of $\mathds{N}$:
\begin{notation} \label{202606111507}
    Throughout this paper, we fix a countable set $\Theta = \{ \theta_1, \theta_2, \cdots, \theta_n, \cdots \}$.
    For $X \subset \mathds{N}^\ast$, we write $\theta_{i_1} \cdots \theta_{i_n}$ for the generator of $\mathsf{Lin} (X)$ corresponding to the order $i_1 \prec \cdots \prec i_n$ where $X = \{ i_1,\cdots,i_n \}$.
    Furthermore, we also use the same notation for the right $\mathtt{L}_{\mathfrak{S}}$-module corresponding to $\mathsf{Lin}$.
\end{notation}
For instance, $\mathsf{Lin} (1)$ is generated by the single element $\theta_1$, whereas $\mathsf{Lin} (2)$ has two generators $\theta_1\theta_2$ and $\theta_2\theta_1$ with $\mathfrak{S}_2$ acting by place permutations.

We now turn to a class of examples that admit a $\mathds{k}$-linear symmetric operad structure.
We will simply refer to a $\mathds{k}$-linear symmetric operad as an operad.
The reader is referred to \cite{Loday2012} for the concepts related to operads.
An operad $\mathfrak{O}$ is a sequence of right $\mathfrak{S}_{n}$-modules $\mathfrak{O}(n),~n\in\mathds{N}$ equipped with an operad product.
For an operad $\mathfrak{O}$, we use the same notation for the $\mathds{k}$-linear species corresponding to $\mathfrak{O}$.

\begin{remark}
    The exponential species and the linear-order species admit well-known operad structures: these are operads for unital commutative algebras and unital associative algebras respectively, while, in this paper, we do not consider those structures.
\end{remark}

\begin{Example} \label{202606121312}
    Let $B$ be a unital $\mathds{k}$-algebra.
    The algebra structure naturally induces an operad structure on $I_B$.
    For instance, $I_\mathds{k}$ is the initial operad.
\end{Example}

\begin{Example}
    Let $\mathfrak{Lie}$ denote the {\it Lie} operad; thus, $\mathfrak{Lie} (X)$ is the multilinear component of the free Lie algebra on the set $X$.
    Via the classical embedding of the Lie operad into the linear-order species, we regard $\mathfrak{Lie}$ as a subspecies of $\mathsf{Lin}$.
    As notational convention, for $X \subset \mathds{N}^\ast$, we use $\theta_i,~i\in X$, as the generators of $\mathfrak{Lie}(X)$, instead of $i \in X$.
    For instance, $\mathfrak{Lie} (2)$ is generated by the single element $[\theta_1,\theta_2]=\theta_1\theta_2 - \theta_2\theta_1$.
\end{Example}

Using the Lie operad $\mathfrak{Lie}$, we also introduce the following operad.
\begin{Defn}
    For $c \in \mathds{N}^\ast$, we denote by $\mathfrak{Lie}_{c}$ the operad for nilpotent Lie algebras of class $\leq c$.
    More precisely, it is given by the quotient operad of $\mathfrak{Lie}$ such that
    \begin{align*}
        \mathfrak{Lie}_{c} (n) {:=}
        \begin{cases}
            \mathfrak{Lie}(n) & n \leq c, \\
            0 & n > c .
        \end{cases}
    \end{align*}
\end{Defn}

\begin{Example} \label{202606121319}
    For instance, $\mathfrak{Lie}_1 = I_\mathds{k}$.
\end{Example}

\subsection{Day convolution of $\mathds{k}$-linear species}

In this section, we review the symmetric monoidal structure on $\mathsf{Sp}_\mathds{k}$ arising from the {\it Day convolution}:
$$
(F \odot G)(X) {:=} \bigoplus_{X_1 \amalg X_2 = X} F(X_1) \otimes G(X_2) ,
$$
where $F,G$ are $\mathds{k}$-linear species.
The unit object is given by $\mathds{k}$, which we regard as the species concentrated on the empty set $\emptyset$.
The associator and the symmetry are obtained from those for the tensor products.

We recursively define a $\mathds{k}$-linear species $F^{\odot m}, ~ m\in\mathds{N}$ by
\begin{align*}
    F^{\odot m} {:=} F^{\odot m-1} \odot F, \quad \mathrm{and~} F^{\odot 0} {:=} \mathds{k} .
\end{align*}
    
A symmetric monoidal category structure allows one to formulate monoids, comonoids, and Hopf monoids (for instance, see \cite[Chapter 1]{aguiar2010monoidal}).
We briefly recall these notions here in order to fix conventions for the remainder of the paper.
A comonoid in $\mathsf{Sp}_{\mathds{k}}$, or a {\it comonoid species}, is a $\mathds{k}$-linear species $\mathtt{C}$ equipped with a comultiplication $\Delta : \mathtt{C} \to \mathtt{C} \odot \mathtt{C}$ and a counit $\epsilon : \mathtt{C} \to \mathds{k}$ satisfying the coassociativity and counit axioms.
Several additional properties of comonoids will play an important role in what follows.
A comonoid species $\mathtt{C}$ is said to be {\it cocommutative} if the composition $\mathtt{C} \stackrel{\Delta}{\to} \mathtt{C} \odot \mathtt{C} \stackrel{\sigma}{\to} \mathtt{C} \odot \mathtt{C}$ coincides with $\Delta$, where $\sigma$ denotes the symmetry in $\mathsf{Sp}_{\mathds{k}}$. 
A comonoid species $\mathtt{C}$ is said to be {\it connected} if $\epsilon : \mathtt{C}(\emptyset) \stackrel{\cong}{\to} \mathds{k}$.
A {\it coaugmented comonoid species} is a comonoid $\mathtt{C}$ endowed with a comonoid morphism $\eta : \mathds{k} \to \mathtt{C}$.

Dually, a {\it monoid species} is a $\mathds{k}$-linear species equipped with a multiplication $\nabla: \mathtt{C}\odot\mathtt{C} \to \mathtt{C}$ and a unit $\eta: \mathds{k} \to \mathtt{C}$.

A {\it bimonoid species} is a $\mathds{k}$-linear species endowed with both a comonoid structure and a monoid structure such that the multiplication and the unit are comonoid morphisms.
A bimonoid species is said to be {\it (co)commutative} if its underlying (co)monoid is (co)commutative.
A bimonoid species is {\it bicommutative} if it is both commutative and cocommutative.

A {\it Hopf monoid species} is a bimonoid species admitting an antipode.

\begin{Example} \label{202604091327}
    Let $V$ be a $\mathds{k}$-module.
    The species $V^{\otimes}$ in Example \ref{202606112022} admits a bicommutative Hopf monoid species structure as follows.
    Note that there is a natural isomorphism $$(V^{\otimes})^{\odot m} (X) \cong \bigoplus \bigotimes^{m}_{i=1} V^{\otimes X_i}  \cong \bigoplus_{f:X \to {\bf m}} V^{\otimes X} ,$$
    where the first direct sum is taken over decompositions $X = \coprod^{m}_{j=1} X_j$ and the second is taken over maps $f:X \to {\bf m}$.
    The first isomorphism follows from the definition of the Day convolution.
    Using this isomorphism, the comultiplication $\Delta : V^{\otimes} \to (V^{\otimes})^{\odot 2}$ and the multiplication $\nabla : (V^{\otimes})^{\odot 2} \to V^{\otimes}$ 
    are described as the diagonal map and the folding map: 
    $$
    V^{\otimes X}  
    \stackrel{\Delta}{\to} \bigoplus_{f:X \to {\bf 2}} V^{\otimes X} , \quad  \bigoplus_{f:X \to {\bf 2}} V^{\otimes X}  \stackrel{\nabla}{\to} V^{\otimes X} .
    $$
    The unit and the counit are obvious.
    The antipode $S : V^{\otimes}\to V^{\otimes}$ is given by $(-1)^{|X|} : V^{\otimes X} \to V^{\otimes X}$.
\end{Example}

\begin{Example} \label{202604101536}
    Recall Notation \ref{202606111507}.
    The linear-order species $\mathsf{Lin}$ has a cocommutative Hopf monoid structure as follows:
    \begin{enumerate}
        \item 
        The multiplication $\nabla : \mathsf{Lin}^{\odot 2} \to \mathsf{Lin}$ is determined by concatenation:
        $$
        \nabla ( \theta_{1} \cdots \theta_{k} \otimes \theta_{k+1} \cdots \theta_{n}) = \theta_{1} \cdots \theta_{n} .
        $$
        \item 
        The comultiplication $\Delta : \mathsf{Lin} \to \mathsf{Lin}^{\odot 2}$ is given by
        $$
        \Delta ( \theta_{1} \cdots \theta_{n} ) = \sum \theta_{\sigma(1)} \cdots \theta_{\sigma(r)} \otimes \theta_{\sigma(r+1)} \cdots \theta_{\sigma(n)} ,
        $$
        where the sum is taken over $0 \leq r \leq n$ and $(r,n-r)$-shuffles $\sigma$, i.e. $\sigma \in \mathfrak{S}_n$ such that $\sigma(1) < \cdots < \sigma(r)$ and $\sigma(r+1) < \cdots < \sigma(n)$.
        \item The counit $\mathsf{Lin} \to \mathds{k}$ and the unit $\mathds{k} \to \mathsf{Lin}$ are clear.
        \item 
        The antipode $S : \mathsf{Lin} \to \mathsf{Lin}$ is given by $S( \theta_{1} \cdots \theta_{n} ) = (-1)^{n} \theta_{n} \cdots \theta_{1}$.
    \end{enumerate}
\end{Example}

\subsection{Primitive elements in comonoid species}

Let $\mathtt{C}$ be a coaugmented comonoid species.
In this section, we study the kernel of the {\it reduced comultiplication} $\bar{\Delta} : \mathtt{C} \to \mathtt{C}^{\odot 2}$ and its generalization, where the reduced comultiplication is defined by
$$\bar{\Delta} {:=} (\Delta - \eta \odot \mathrm{id}_{\mathtt{C}} - \mathrm{id}_{\mathtt{C}} \odot \eta)  .$$
We give two key applications (Propositions \ref{202603252001} and \ref{202604021602}), which are crucial for calculating the eigenmonad $\Phi_\mathcal{C}$ for a Lawvere theory $\mathcal{C}$.

\begin{Defn} \label{202603241659}
    Let $\mathtt{C}$ be a coaugmented comonoid species.
    For $m \in \mathds{N}^\ast$, we define $\mathrm{Pr}^m (\mathtt{C}) \subset \mathtt{C}^{\odot m}$ to be the joint kernel of the maps $\mathrm{id}_{\mathtt{C}}^{\odot k-1} \odot \bar{\Delta} \odot \mathrm{id}_{\mathtt{C}}^{\odot m-k}$ for $k \in {\bf m}$ corresponding to the exact sequence:
    \begin{align*}
        0 \to \mathrm{Pr}^m (\mathtt{C}) \to \mathtt{C}^{\odot m} \to \bigoplus^{m}_{k=1} \mathtt{C} \odot \cdots \odot \overbrace{(\mathtt{C} \odot \mathtt{C})}^{k} \odot \cdots \odot \mathtt{C} .
    \end{align*}
    We write $\mathrm{Pr}(\mathtt{C})$ for $\mathrm{Pr}^1 (\mathtt{C}) = \mathrm{Ker} (\bar{\Delta})$.
    We also define $\mathrm{Pr}^0 ( \mathtt{C}) {:=} \mathds{k}$.
\end{Defn}

\begin{remark}
    For a coaugmented coalgebra in the usual sense, an element of the kernel of the reduced comultiplication is called primitive.
    The above definition is an analogue of this notion.
\end{remark}

We denote by $\bar{\mathtt{C}}$ the kernel of the counit $\epsilon : \mathtt{C} \to \mathds{k}$.
It is straightforward to verify that the reduced comultiplication $\bar{\Delta}$ restricts to a map $\bar{\mathtt{C}} \to \bar{\mathtt{C}}^{\odot 2}$.

\begin{Lemma} \label{202603221101}
    $\mathrm{Pr}^m (\mathtt{C}) \subset \bar{\mathtt{C}}^{\odot m}$.
\end{Lemma}
\begin{proof}
    Consider the following diagram:
    $$
    \begin{tikzcd}[row sep=small]
        & & \bar{\mathtt{C}}^{\odot m} \ar[d, "h"] &\\
        0 \ar[r] & \mathrm{Pr}^m (\mathtt{C}) \ar[ur, "\exists"] \ar[r] & \mathtt{C}^{\odot m} \ar[d, "f"'] \ar[r] & \bigoplus^m_{k=1} \mathtt{C} \odot \cdots \odot \overbrace{(\mathtt{C} \odot \mathtt{C})}^{k} \odot \cdots \mathtt{C} \ar[dl, "g"] \\
        & & \bigoplus^m_{k=1} \mathtt{C} \odot \cdots \odot \underbrace{\mathds{k}}_{k} \odot \cdots \mathtt{C}
    \end{tikzcd}
    $$
    Here, the map $f$ is the sum of $\mathrm{id}_{\mathtt{C}}^{\odot k-1} \odot \epsilon \odot \mathrm{id}_{\mathtt{C}}^{\odot m-k},~ k \in {\bf m}$, the map $g$ is the direct sum of $-\mathrm{id}_\mathtt{C}^{\odot k-1} \odot \epsilon^{\odot 2} \odot \mathrm{id}_{\mathtt{C}}^{\odot m-k},~ k \in {\bf m}$, and $h$ is induced by the inclusion $\bar{\mathtt{C}} \hookrightarrow \mathtt{C}$.
    By the isomorphism $\mathtt{C} \cong \bar{\mathtt{C}} \oplus \mathds{k}$, the vertical sequence is left exact.
    Since $\epsilon^{\odot 2} \circ \bar{\Delta} = - \epsilon$, the right lower triangle commutes.
    Hence, the composition of the inclusion $\mathrm{Pr}^m (\mathtt{C}) \to \mathtt{C}^{\odot m}$ with $f$ is trivial.
    By the left-exactness of the vertical sequence, we obtain a map $\mathrm{Pr}^m (\mathtt{C}) \to \bar{\mathtt{C}}^{\odot m}$ which is automatically injective and which fits into the above commutative diagram.
\end{proof}

\begin{Lemma} \label{202603221113}
    If $\mathtt{C}$ is connected, then for $m\in\mathds{N}$, we have $\mathrm{Pr}(\mathtt{C})^{\odot m} (m) \cong \mathrm{Pr}^m (\mathtt{C}) (m) = \bar{\mathtt{C}}^{\odot m} (m)$.
\end{Lemma}
\begin{proof}
    The connectivity of $\mathtt{C}$ and Lemma \ref{202603221101} imply that $\mathrm{Pr} (\mathtt{C})(0) = \bar{\mathtt{C}} (0) = 0$ and $\mathrm{Pr} (\mathtt{C}) (1) = \bar{\mathtt{C}}(1)$.
    Hence,
    \begin{align*}
        \mathrm{Pr} (\mathtt{C})^{\odot m} (m) &= \bigoplus_{\sigma \in \mathfrak{S}_m} \bigotimes^m_{k=1} \mathrm{Pr} (\mathtt{C}) (\{\sigma (k) \} ) = \bigoplus_{\sigma \in \mathfrak{S}_m} \bigotimes^m_{k=1}\bar{\mathtt{C}}(\{\sigma (k) \} ) = \bar{\mathtt{C}}^{\odot m} (m) .
    \end{align*}
    This result proves the claim.
\end{proof}

We defer the proofs of the following two lemmas, which provide conditions under which $\mathrm{Pr}(\mathtt{C})^{\odot m} \to \mathrm{Pr}^m (\mathtt{C})$ is an isomorphism, since they require arguments that are not needed elsewhere in the paper.
See Appendix \ref{202604101548}.

\begin{Lemma} \label{202604011630}
    Let $\mathtt{C}$ be a coaugmented comonoid species.
    If the ground ring $\mathds{k}$ is a field, then $\mathrm{Pr}(\mathtt{C})^{\odot m} \stackrel{\cong}{\to} \mathrm{Pr}^m (\mathtt{C})$.
\end{Lemma}

Even if $\mathds{k}$ is not a field, there is still a sufficient condition.

\begin{Lemma} \label{202603321744}
    Let $\mathtt{C}$ be a coaugmented comonoid species.
    Assume that, for any $n \in \mathds{N}$, there exists a map, $\varphi_n : \bar{\mathtt{C}}^{\odot 2} (n) \to \bar{\mathtt{C}} (n)$ making the following diagram commute:
    $$
    \begin{tikzcd}
        \bar{\mathtt{C}} (n) \ar[r, "\bar{\Delta}"] \ar[d, "\bar{\Delta}" ]& \bar{\mathtt{C}}^{\odot 2} (n) \ar[d, "\varphi_n"] \\
        \bar{\mathtt{C}}^{\odot 2} (n) & \bar{\mathtt{C}} (n) . \ar[l, "\bar{\Delta}"]
    \end{tikzcd}
    $$
    Then $\mathrm{Pr}(\mathtt{C})^{\odot m}  \stackrel{\cong}{\to} \mathrm{Pr}^m (\mathtt{C})$.
    Furthermore, for $n \in \mathds{N}$, the map $\varphi_n \circ \bar{\Delta} :\bar{\mathtt{C}} (n) \to \bar{\mathtt{C}} (n)$ is an idempotent and the image of the idempotent $\mathrm{id}_{\mathtt{C}(n)} - \varphi_n \circ \bar{\Delta}$ coincides with $\mathrm{Pr}(\mathtt{C})(n)$.
\end{Lemma}

This lemma is fundamental for our applications.
We present two examples to which Lemma \ref{202603321744} applies.
Recall $I_V$ from Example \ref{202603251354}.
We also recall the Hopf monoid $V^{\otimes}$ from Example \ref{202604091327}.
\begin{prop} \label{202603252001}
    The inclusion $I_V \hookrightarrow V^{\otimes}$ induces
    $\mathrm{Pr}^{m} ( V^{\otimes}) \cong I_V^{\odot m}$.
\end{prop}
\begin{proof}
    We will introduce the maps $\varphi_n$ as in the previous lemma for $\mathtt{C} = V^{\otimes}$.
    Note that $\bar{\mathtt{C}} (0) = 0$, which implies $\bar{\mathtt{C}}^{\odot 2} (n) = \bigoplus V^{\otimes n}$ where the direct sum is take over surjective maps ${\bf n}\to{\bf 2}$.
    So, $\varphi_n$ for $n \in \{0,1\}$ is trivial.
    We construct $\varphi_n$ for $n \geq 2$.
    For $n\in\mathds{N}^\ast$, let $f_n: {\bf n}\to{\bf 2}$ be the map such that $f_n^{-1} (1) =\{1\}$.
    Let $\varphi_n :\bar{\mathtt{C}}^{\odot 2} (n) = \bigoplus V^{\otimes n} \to V^{\otimes n} = \bar{\mathtt{C}}(n)$ be the projection to the $f_n$-component.
    By the definition, for $x \in V^{\otimes n}$, all components of $\bar{\Delta} (x) \in \bar{\mathtt{C}}^{\odot 2} (n) = \bigoplus V^{\otimes n}$ are $x$.
    Hence, $\bar{\Delta} \circ \varphi_n \circ \bar{\Delta} (x) = \bar{\Delta}(x)$.
    Thus, we obtain $\bar{\Delta} \circ \varphi_n \circ \bar{\Delta} = \bar{\Delta}$.
    By Lemma \ref{202603321744}, $\mathrm{Pr}^m (V^{\otimes} ) \cong \mathrm{Pr} (V^{\otimes})^{\odot m}$ and the inclusion $I_V \hookrightarrow V^{\otimes}$ is a map onto the image of $\mathrm{id}_{V^{\otimes n}} - \varphi_n \circ \bar{\Delta}$, so $I_V(n) \cong \mathrm{Pr}(V^{\otimes}) (n)$.
\end{proof}

\begin{notation}
    For a set $U$, let $\mathcal{P}_f (U)$ denote the set of finite subsets of $U$.
\end{notation}

\begin{Defn} \label{202604021629}
    Let $\mathtt{C}$ be a bimonoid species.
    Consider a family of subsets $\mathfrak{B}_X \subset \bar{\mathtt{C}}(X)$ for $X \in \mathcal{P}_f (\mathds{N}^{\ast})$ together with a total order $\preceq$ on the set $\coprod_{X \in \mathcal{P}_f (\mathds{N}^{\ast})} \mathfrak{B}_X$.
    We say that $( \{ \mathfrak{B}_X \mid X \in \mathcal{P}_f (\mathds{N}^{\ast}) \} , \preceq)$ forms a {\it PBW (Poincar\'e-Birkhoff-Witt) datum} for $\mathtt{C}$ if the following conditions hold for each $X \in \mathcal{P}_f (\mathds{N}^{\ast})$:
    \begin{enumerate}
        \item Any $h \in \mathfrak{B}_X$ is primitive in the sense that $\bar{\Delta} (h) = 0$.
        \item 
        The set $$\{ h_1\cdots h_r \in \bar{\mathtt{C}} (X) \mid r \in \mathds{N}^\ast,~ X = \coprod^r_{j=1} X_j ,~ h_j \in \mathfrak{B}_{X_j}, ~ h_1 \succ \cdots \succ h_r \}$$
        gives a basis of the $\mathds{k}$-module $\bar{\mathtt{C}} (X)$.
        Here, $h_1\cdots h_r$ denotes the element of $\bar{\mathtt{C}} (X)$ obtained by applying the $r$-fold multiplication map $\mathtt{C}^{\odot r} (X) \to \mathtt{C} (X)$ to $h_1 \otimes \cdots \otimes h_r\in \bigotimes^r_{j=1} \mathtt{C} (X_j)$
    \end{enumerate}
\end{Defn}

\begin{remark}
    Note that the above notion is not symmetric: for a bijection $f : Y \to X$, the induced map $\mathtt{C}(X) \to \mathtt{C}(Y)$ does not necessarily send $\mathfrak{B}_X$ to $\mathfrak{B}_Y$.
\end{remark}

The bimonoid species in Example \ref{202604101536} provides an example:
\begin{Lemma}[see Appendix \ref{202512031342} for the proof] \label{202604011516}
    There exists a PBW datum $(\{ \mathfrak{B}_X \mid X \in \mathcal{P}_f (\mathds{N}^{\ast}) \}, \preceq )$ for $\mathsf{Lin}$ such that $\mathfrak{B}_X \subset \mathfrak{Lie}(X)$ and $\mathfrak{Lie} (X) \cong \mathds{k} [\mathfrak{B}_X]$ for $X \in \mathcal{P}_f (\mathds{N}^{\ast})$.
\end{Lemma}

\begin{prop} \label{202604021602}
    Let $\mathtt{C}$ be a bimonoid endowed with a PBW datum as above.
    Then we have $\mathrm{Pr}^m (\mathtt{C}) \cong \mathrm{Pr}(\mathtt{C})^{\odot m}$.
    Furthermore, for $X \in \mathcal{P}_f ( \mathds{N}^\ast)$, $\mathrm{Pr}(\mathtt{C})(X) \cong \mathds{k} [\mathfrak{B}_X]$.
\end{prop}
\begin{proof}
    For $X \in \mathcal{P}_f ( \mathds{N}^\ast)$, we define a map $\varphi_X : \bar{\mathtt{C}}^{\odot 2} (X) \to \bar{\mathtt{C}} (X)$.
    By the condition (2) of Definition \ref{202604021629}, $\bar{\mathtt{C}}^{\odot 2} (X)$ has a basis consisting of $h_1 \cdots h_k \otimes h_1^\prime \cdots h^\prime_s$ where $k,s \in \mathds{N}^\ast$, $X = \coprod^{k}_{j=1} X_{1j} \amalg \coprod^s_{i=1} X_{2i}$, $h_j \in \mathfrak{B}_{X_{1j}}$ such that $h_1 \succ \cdots \succ h_k$; and $h^\prime_i \in \mathfrak{B}_{X_{2i}}$ such that $h^\prime_1 \succ \cdots \succ h^\prime_s$.
    Define
    \begin{align*}
        \varphi_X (h_1 \cdots h_k \otimes h_1^\prime \cdots h^\prime_s) {:=}
        \begin{cases}
            h_1 \cdots h_k h_1^\prime & \mathrm{if~} s= 1 \mathrm{~and~} h_k \succ h_1^\prime ,\\
            0 & \mathrm{otherwise}.
        \end{cases}
    \end{align*}
    This linearly extends to a map $\varphi_X : \bar{\mathtt{C}}^{\odot 2} (X) \to \bar{\mathtt{C}} (X)$.
    We will prove $\bar{\Delta} \circ \varphi_X \circ \bar{\Delta} = \bar{\Delta}$.
    
    Consider a basis of $\bar{\mathtt{C}}(X)$, $h_1\cdots h_k$, where $X = \coprod^k_{i=1} X_i$, $h_j \in \mathfrak{B}_{X_i}$ and $h_1 \succ \cdots \succ h_k$.
    Since $h_i$'s are primitive by (1) of Definition \ref{202604021629}, we have
    \begin{align*}
        \bar{\Delta} ( h_1 \cdots h_k) =  \sum h_{\sigma(1)} \cdots h_{\sigma(r)} \otimes h_{\sigma(r+1)} \cdots h_{\sigma(k)} \in \bar{\mathtt{C}}^{\odot 2} (X) ,
    \end{align*}
    where the sum is taken over $1 \leq r \leq k-1$ and $(r,k-r)$-shuffles $\sigma$.
    In particular, if $k= 1$, then $\bar{\Delta} ( h_1 ) = 0$.
    So, $\bar{\Delta} \circ \varphi_X \circ \bar{\Delta} (h_1) = \bar{\Delta} (h_1)$.
    If $k \geq 2$, then by the definition of $\varphi_X$, 
    \begin{align*}
        \varphi_X \circ \bar{\Delta} ( h_1 \cdots h_k) = \varphi_X \left( \sum h_{\sigma(1)} \cdots h_{\sigma(r)} \otimes h_{\sigma(r+1)} \cdots h_{\sigma(k)} \right) =  h_1 \cdots h_k  .
    \end{align*}
    This leads to $\bar{\Delta} \circ \varphi_X \circ \bar{\Delta} ( h_1 \cdots h_k) = \bar{\Delta} ( h_1 \cdots h_k)$.
    Therefore, $\bar{\Delta} \circ \varphi_X \circ \bar{\Delta} = \bar{\Delta}$.
    From Lemma \ref{202603321744}, we obtain $\mathrm{Pr}^m (\mathtt{C}) \cong \mathrm{Pr}(\mathtt{C})^{\odot m}$.
    Furthermore, the above computations imply that the image of $\mathrm{id}_{\mathtt{C}(X)} - \varphi_X \circ \bar{\Delta}$ equals $\mathds{k} [\mathfrak{B}_X]$.
\end{proof}

\section{A canonical $(\mathtt{L}_{\mathfrak{S}},\mathtt{L}_{\mathfrak{S}})$-bimodule construction $\mu$}
\label{202606121546}

In this section, we give a canonical construction of $(\mathtt{L}_{\mathfrak{S}},\mathtt{L}_{\mathfrak{S}})$-bimodules from right $\mathtt{L}_{\mathfrak{S}}$-modules that will be used throughout the paper. 
Recall the monad $\mathtt{L}_{\mathfrak{S}}$ on $\mathds{N}$ from Example \ref{202604091040}.
Under the natural identification of $\mathtt{L}_{\mathfrak{S}}$-modules with linear species, the constructions below are the standard ones from the theory of linear species.

\subsection{Construction and basic properties} 

In this section, we introduce a construction of a natural $(\mathtt{L}_\mathfrak{S},\mathtt{L}_\mathfrak{S})$-bimodule $\mu \mathtt{M}$ for a right $\mathtt{L}_{\mathfrak{S}}$-module $\mathtt{M}$.

\begin{Defn} \label{202603131024}
    Let $\mathtt{M}$ be a right $\mathtt{L}_{\mathfrak{S}}$-module.
    Let $\mu \mathtt{M}$ be the $(\mathtt{L}_\mathfrak{S},\mathtt{L}_\mathfrak{S})$-bimodule defined by
    \begin{align*}
        \mu \mathtt{M} (m,n) {:=} \bigoplus \bigotimes^{m}_{i=1} \mathtt{M} (n_i) \otimes_{\mathfrak{S}_{n_1} \times \cdots \mathfrak{S}_{n_m}} \mathds{k} [\mathfrak{S}_n] , \quad n,m \in \mathds{N} ,
    \end{align*}
    where the direct sum is taken over $m$-partitions $(n_1,\cdots,n_m)$ of $n$.
    We recall Definition \ref{202603281205}.
    The left $\mathfrak{S}_m$-action and right $\mathfrak{S}_n$-action on $\bigotimes^m_{k=1} [v_k] \otimes \sigma \in \mu\mathtt{M}(m,n)$ are defined by
    \begin{align*}
        \tau \rhd \left( \bigotimes^m_{k=1} [v_k] \otimes \sigma \right) &= \bigotimes^m_{k=1} [v_{\tau^{-1}(k)}] \otimes \tau_{n_1,\cdots,n_m} \sigma , \quad  \tau \in \mathfrak{S}_m ,\quad \mathrm{and;} \\
        \left( \bigotimes^m_{k=1} [v_k] \otimes \sigma \right)\lhd \pi &= \bigotimes^m_{k=1} [v_k] \otimes \sigma\pi , \quad \pi \in \mathfrak{S}_n .
    \end{align*}
\end{Defn}

The following is immediate from the definition:
\begin{Lemma} \label{202603100950}
    We have $\mathtt{M}(0) = 0$ if and only if $\mu\mathtt{M}$ is upper triangular in the sense of Definition \ref{202512071707}.
    Furthermore, $\mathtt{M}(n) = 0$ for $n \neq 1$ if and only if $\mu\mathtt{M}$ is diagonal.
\end{Lemma}

Using the equivalence of categories in (\ref{202603251315}), we give alternative descriptions of $\mu\mathtt{M}$:
\begin{Lemma} \label{202603211811}
    Let $n,m\in \mathds{N}$.
    There are natural isomorphisms for right $\mathtt{L}_{\mathfrak{S}}$-modules $\mathtt{M}$:
    \begin{enumerate}
        \item
        $\mu \mathtt{M} (m,n) \cong \bigoplus_{f} \bigotimes^{m}_{i=1} \mathtt{M} (f^{-1}(i))$ where $f$ ranges over all maps from ${\bf n}$ to ${\bf m}$.
        \item 
        $\mu \mathtt{M} (m,n) \cong \mathtt{M}^{\odot m} (n)$.
    \end{enumerate}
    Furthermore, under the isomorphism in (2), the $(\mathfrak{S}_m, \mathfrak{S}_n)$-bimodule structure obtained from the symmetry in $\mathsf{Sp}_\mathds{k}$ agrees with that in Definition \ref{202603131024}.
\end{Lemma}
\begin{proof}
    To prove (1), we construct a map from $\mu \mathtt{M}(m,n)$ to the right hand side.
    Let $(n_1,\cdots,n_m)$ be a partition of $n$ and $\rho: {\bf n} \to {\bf m}$ the associated map.
    Given $\sigma \in \mathfrak{S}_n$ , we put $f = \rho \circ \sigma$.
    The composition $\sigma |_{f^{-1}(i)} : f^{-1}(i) \to \rho^{-1}(i) \cong {\bf n_i}$ gives an element in $\mathrm{Bij} ( f^{-1}(i), {\bf n_i})$ where $\rho^{-1}(i) \cong {\bf n_i}$ preserves the order.
    So, for $v_i \in \mathtt{M}  (n_i)$, we have $v_i \lhd \sigma |_{f^{-1}(i)} \in \mathtt{M} (f^{-1}(i))$.
    Then $\bigotimes^m_{i=1} v_i \otimes \sigma \in \mu \mathtt{M} (m,n) \mapsto \bigotimes^m_{i=1} v_i \lhd \sigma |_{f^{-1}(i)}$ yields a well-defined map with the obvious inverse.
    (2) is just a reformulation of (1).
    The last claim is clear from the above results.
\end{proof}

\subsection{Monads associated to operads} \label{202512151740} 

As we briefly reviewed in Section \ref{202604091335}, operads may be viewed as right $\mathtt{L}_{\mathfrak{S}}$-modules.
For an operad $\mathfrak{O}$, the operad composition extends to the monad structure on $\mu\mathfrak{O}$ as follows.
For a map $f : X \to Y$ between finite sets, the operad composition induces a map
\begin{align*}
    \gamma_f: \mathfrak{O} (Y) \otimes \bigotimes_{y \in Y} \mathfrak{O}(f^{-1}(y) ) \to \mathfrak{O}(X) .
\end{align*}
Let $g : Y \to Z$ be another map between finite sets.
By tensoring the extended operad compositions, we obtain a map:
{\small
\begin{align*}
     \bigotimes_{z \in Z} \mathfrak{O}(g^{-1}(z)) \otimes \bigotimes_{y \in Y} \mathfrak{O} (f^{-1}(y))\cong \bigotimes_{z \in Z} \left( \mathfrak{O}(g^{-1}(z)) \otimes \bigotimes_{y \in g^{-1}(z)} \mathfrak{O} (f^{-1}(y)) \right) \stackrel{\bigotimes_{z\in Z} \gamma_{f|_{g^{-1}(z)}}}{\longrightarrow} \bigotimes_{z \in Z} \mathfrak{O} ((g \circ f)^{-1}(z)) .
\end{align*}
}
Considering this construction for all the maps $f: {\bf l} \to {\bf n}$ and $g: {\bf n} \to {\bf m}$, we obtain a map
\begin{align*}
    \mu\mathfrak{O}(m,n) \otimes \mu\mathfrak{O}(n,l) \to \mu\mathfrak{O}(m,l) .
\end{align*}
This gives a monad structure $\mu\mathfrak{O} \otimes \mu\mathfrak{O} \to \mu\mathfrak{O}$ together with the unit $\eta : \mathds{I}_{\mathds{N}} \to \mu\mathfrak{O}$ induced by the unit of the operad $\mathfrak{O}$.

\begin{remark} \label{202604051201}
    Under the bijection between monads in $\mathsf{Mat}_{\mathds{k}}$ and small $\mathds{k}$-linear categories (see Section \ref{202606121310}), the monad $\mu\mathfrak{O}$ corresponds to the category associated with $\mathfrak{O}$, denoted by $\mathsf{Cat}\mathfrak{O}$, that appears in the literature (for instance, \cite{powell2024analytic,powell2024outer,kim2024analytic}).
\end{remark}

\begin{Example} \label{202603100930}
    Recall the operad $I_B$ from Example \ref{202606121312}.
    The monad $\mu I_B$ is diagonal.
    The monad $\mu I_\mathds{k}$ is isomorphic to the monad $\mathtt{L}_{\mathfrak{S}}$.
    By Example \ref{202606121319}, we also have $\mu\mathfrak{Lie}_1 \cong \mathtt{L}_{\mathfrak{S}}$.
\end{Example}

We recall the wreath product $B \thicksim \mathfrak{S}_n$ of a $\mathds{k}$-algebra $B$  \cite{wreathMac}.
It underlying $\mathds{k}$-module is defined to be $B^{\otimes n} \otimes \mathds{k} [\mathfrak{S}_n]$.
The multiplication is defined by 
$$
(a_1 \otimes \sigma_1) \cdot (a_2 \otimes \sigma_2) {:=} a_1 (\sigma_1 \rhd a_2) \otimes  \sigma_1 \sigma_2 .
$$
Here, $(-) \rhd (-) : \mathfrak{S}_n  \times B^{\otimes n} \to B^{\otimes n}$ is the left place permutation action.

\begin{prop} \label{202604061428}
    For $n \in \mathds{N}$, we have an algebra isomorphism $\mu I_B (n,n) \cong B \thicksim \mathfrak{S}_n$.
\end{prop}
\begin{proof}
    By the definition of $\mu I_B$, the underlying modules of $\mu I_B (n,n)$ and $B \thicksim \mathfrak{S}_n$ are identified.
    Furthermore, one may directly verify that this preserves the algebra structures.
\end{proof}

\section{A structural theorem for the quotient $\mathtt{L}_{\mathcal{C}}/\mathtt{I}^{\mathsf{pr}}_{\mathcal{C}}$} \label{202604131411}

In this section, we give a structural theorem for the quotient $\mathtt{L}_{\mathcal{C}}/\mathtt{I}^{\mathsf{pr}}_{\mathcal{C}}$ where $\mathcal{C}$ is a Lawvere theory with a zero object.
Recall from Section \ref{202606091353} that $\mathtt{L}_{\mathcal{C}}/\mathtt{I}^{\mathsf{pr}}_{\mathcal{C}}$ admits a canonical $(\mathtt{L}_{\mathcal{C}}, \Phi_{\mathcal{C}})$-bimodule structure.
We show that the underlying $(\mathtt{L}_\mathfrak{S}, \mathtt{L}_\mathfrak{S})$-bimodule can be reconstructed from a natural right $\mathtt{L}_\mathfrak{S}$-module via the construction $\mu$ in Definition \ref{202603131024}.
In the application, we give an explicit computation of the quotient $\mathtt{L}_{\G^{\mathsf{o}}}/ \mathtt{I}^{\mathsf{pr}}_{\G^{\mathsf{o}}}$.

\subsection{Statement of the theorem}

Let $\mathcal{C}$ be a Lawvere theory with a zero object, so that the primitivity ideal $\mathtt{I}^{\mathsf{pr}}_{\mathcal{C}}$ is defined.
Recall that $\mathtt{I}^{\mathsf{pr}}_{\mathcal{C}} (1,n) \subset \mathtt{L}_{\mathcal{C}} (1,n)$ is closed under the $\mathfrak{S}_n$-action (see Lemma \ref{202509071205}).

\begin{Defn} \label{202604041850}
    We define a right $\mathtt{L}_{\mathfrak{S}}$-module $\Psi_{\mathcal{C}}$ by $$\Psi_{\mathcal{C}} (n){:=} (\mathtt{L}_{\mathcal{C}}/\mathtt{I}^{\mathsf{pr}}_{\mathcal{C}} )(1,n), \quad n \in \mathds{N} .$$
\end{Defn}

\begin{Example} \label{202603221114}
    Since $\mathtt{I}^{\mathsf{pr}}_{\mathcal{C}} (1,0) = 0$, we have $\Psi_{\mathcal{C}}(0)=  \mathtt{L}_{\mathcal{C}} (1,0) \cong \mathds{k}$.
\end{Example}

\begin{notation}
    For $f \in \mathtt{L}_{\mathcal{C}} (1,n)$, we denote by $[f] \in \Psi_{\mathcal{C}} (n)$ the class of $f$.
\end{notation}

Recall from Definition \ref{202603131024} the $(\mathtt{L}_\mathfrak{S},\mathtt{L}_\mathfrak{S})$-bimodule $\mu \Psi_{\mathcal{C}}$.
Lemma \ref{202601051527} implies that the following determines a well-defined map:
\begin{Defn} \label{202603131134}
    We define an $(\mathtt{L}_\mathfrak{S},\mathtt{L}_\mathfrak{S})$-bimodule map, denoted by $\alpha_{\mathcal{C}} : \mu \Psi_{\mathcal{C}} \to \mathtt{L}_{\mathcal{C}}/\mathtt{I}^{\mathsf{pr}}_{\mathcal{C}}$.
    Let $m,n\in\mathds{N}$.
    For an $m$-partition $(n_1,\cdots,n_m)$ of $n$, the component
    $$
    \bigotimes^{m}_{j=1} \Psi_{\mathcal{C}} (n_j) \otimes_{\mathfrak{S}_{n_1}\times \cdots \times \mathfrak{S}_{n_m}} \mathds{k} [\mathfrak{S}_n] \to \mathtt{L}_{\mathcal{C}}/\mathtt{I}^{\mathsf{pr}}_{\mathcal{C}} (m,n) ,
    $$
    is given by
    $$
    \bigotimes^{m}_{j=1} [f_j] \otimes \sigma \mapsto \left( \prod^m_{j=1} f_j \right) \circ E_\sigma \mod{\mathtt{I}^{\mathsf{pr}}} (m,n),
    $$
    where $f_j \in \mathtt{L}_{\mathcal{C}} (1,n_j),~j \in {\bf m}$ and $\sigma \in \mathfrak{S}_n$.  
    Here, $E_{\sigma} \in \mathcal{C}(n,n)$ was introduced in Definition \ref{202603261255}.
\end{Defn}

\begin{Lemma}
    The above construction yields a $(\mathtt{L}_\mathfrak{S},\mathtt{L}_\mathfrak{S})$-bimodule map.
\end{Lemma}
\begin{proof}
    It suffices to show that the above map preserves the left $\mathfrak{S}_m$-action and the right $\mathfrak{S}_n$-action.
    Clearly, the right $\mathfrak{S}_n$-action is preserved.
    We now recall $\tau_{n_1,\cdots,n_m} = \tilde{\tau}$ from Definition \ref{202603281205}.
    Then we have
    $$
    \alpha_\mathcal{C} \left(\tau \rhd  ( \bigotimes^{m}_{j=1} [f_j] \otimes \sigma ) \right) = \alpha_\mathcal{C} \left( \bigotimes^{m}_{j=1} [f_{\tau^{-1} (j)}] \otimes \tilde{\tau} \sigma \right) = \left( \prod^m_{j=1} f_{\tau^{-1}(j)} \right) \circ E_{\tilde{\tau}\sigma} .
    $$
    This is computed as follows: 
    $$\left( \prod^m_{j=1} f_{\tau^{-1}(j)} \right) \circ E_{\tilde{\tau}\sigma} = \left( \prod^m_{j=1} f_{\tau^{-1}(j)} \right) \circ E_{\tilde{\tau}} \circ  E_{\sigma} = E_{\tau} \circ \prod^m_{k=1} f_k \circ  E_{\sigma} = \tau \rhd \alpha_\mathcal{C} \left( \bigotimes^{m}_{j=1} [f_j] \otimes \sigma \right),$$
    where $E_{\tau} \circ \prod^m_{k=1} f_k =  \prod^m_{k=1} f_{\tau^{-1} (k)} \circ E_{\tilde{\tau}}$ is obtained from Lemma \ref{202606182054}.
\end{proof}

\begin{theorem} \label{202603161027}
    The map $\alpha_\mathcal{C}: \mu \Psi_{\mathcal{C}} \to \mathtt{L}_{\mathcal{C}}/\mathtt{I}^{\mathsf{pr}}_{\mathcal{C}}$ gives an isomorphism of $(\mathtt{L}_\mathfrak{S}, \mathtt{L}_\mathfrak{S})$-bimodules.
\end{theorem}

The proof is postponed to Section \ref{202604011319}.
The theorem shows that the underlying $(\mathtt{L}_\mathfrak{S}, \mathtt{L}_\mathfrak{S})$-bimodule of $\mathtt{L}_{\mathcal{C}}/\mathtt{I}^{\mathsf{pr}}_{\mathcal{C}}$ is completely determined by its $(1,-)$-components $\Psi_\mathcal{C}$.

Recall from Section \ref{202606091353} that $\mathtt{L}_{\mathcal{C}}/\mathtt{I}^{\mathsf{pr}}_{\mathcal{C}}$ is a $(\mathtt{L}_\mathcal{C},\Phi_{\mathcal{C}})$-bimodule.
Via the isomorphism in Theorem \ref{202603161027}, the $(\mathtt{L}_\mathcal{C},\Phi_{\mathcal{C}})$-bimodule structure on $\mathtt{L}_{\mathcal{C}}/\mathtt{I}^{\mathsf{pr}}_{\mathcal{C}}$ transfers to $\mu \Psi_{\mathcal{C}}$.
\begin{notation} \label{202603201036}
    Let $\tilde{\mu}\Psi_{\mathcal{C}}$ denote the $(\mathtt{L}_\mathcal{C},\Phi_{\mathcal{C}})$-bimodule.
\end{notation}

\subsection{The case $\mathcal{C} = \G^{\mathsf{o}}$}
\label{202509201222}

Recall the linear-order species $\mathsf{Lin}$ together with  Notation \ref{202606111507} and the notation in Example \ref{202603121546}.
There is a unique right $\mathfrak{S}_n$-module map $\mathsf{Lin} (n) \to \mathds{k} [\mathsf{F}_n] \cong \mathtt{L}_{\G^{\mathsf{o}}} (1,n)$ which assigns $x_{1} \cdots x_{n}$ to $\theta_{1} \cdots \theta_{n}$.
The composition with the quotient map $\mathtt{L}_{\G^{\mathsf{o}}} (1,n) \to \Psi_{\G^{\mathsf{o}}} (n)$ induces $\mathcal{F}_n : \mathsf{Lin} (n) \to \Psi_{\G^{\mathsf{o}}} (n)$.
These maps assemble into a map of right $\mathtt{L}_{\mathfrak{S}}$-modules:
\begin{align} \label{202603171311} 
    \mathcal{F} : \mathsf{Lin} \to \Psi_{\G^{\mathsf{o}}} .
\end{align}
The main result of this section is the following theorem:
\begin{theorem} \label{202509091931}
    The map $\mathcal{F}$ gives an isomorphism of right $\mathtt{L}_{\mathfrak{S}}$-modules:
    $$\mathsf{Lin} \cong \Psi_{\G^{\mathsf{o}}} .$$
    Hence, by Theorem \ref{202603161027}, we have an isomorphism of $(\mathtt{L}_\mathfrak{S}, \mathtt{L}_\mathfrak{S})$-bimodules:
    $$
    \mu\mathsf{Lin} \cong \mathtt{L}_{\G^{\mathsf{o}}} / \mathtt{I}^{\mathsf{pr}}_{\G^{\mathsf{o}}} .
    $$
\end{theorem}

\begin{notation} \label{202604101645}
    By the theorem, $\mu\mathsf{Lin}$ admits an $(\mathtt{L}_{\G^{\mathsf{o}}} , \Phi_{\G^{\mathsf{o}}})$-bimodule structure; we write $\tilde{\mu}\mathsf{Lin}$ for $\mu\mathsf{Lin}$ equipped with this structure.
\end{notation}

In the rest of this section, we establish the theorem. 
It is useful to recall the Magnus expansion for free groups and related constructions.
These are closely related to the congruence modulo the primitivity ideal $\mathtt{I}^{\mathsf{pr}}_{\G^{\mathsf{o}}}$.

With Example \ref{202603121546} and Notation \ref{202606111507} in mind, we fix the following:
\begin{notation} \label{202601091343}
    We denote by $\comalg (n)$ the $\mathds{k}$-algebra of power series in noncommutative indeterminates of $\theta_i \in \Theta$ for $i \in {\bf n}$.
\end{notation}
The Magnus expansion is the map $\chi : \mathsf{F}_n \to \comalg (n)$ which preserves the product and satisfies:
\begin{align*}
    \chi (x_i) = 1+\theta_i , \quad \chi (x_i^{-1}) = (1+\theta_i)^{-1} =  1-\theta_i +\theta_i ^2 - \theta_i^3 +  \cdots  .
\end{align*}

\begin{Defn} \label{202512011406}
    For a finite set $X$, let $V_{X}$ be the additive monoid consisting of maps $X \to \mathds{N}$.
    We also define $\mathds{1} = \mathds{1}_{X} \in V_{X}$ to be the constant map with value $1 \in \mathds{N}$.
    When $X = {\bf n}$, we often use $V_n = V_{\bf n}$.
\end{Defn}

\begin{Defn} \label{202601071222}
    For $\delta \in V_{n}$, the $\delta$-component $\alg (n)_{\delta}$ is defined to be the $\mathds{k}$-submodule of $\comalg (n)$ generated by the products of $\theta_i$'s having $\delta(k)$ copies of $\theta_k \in \Theta$ for $k \in {\bf n}$.
\end{Defn}

We then have
$$\comalg (n) \cong \prod_{\delta \in V_n} \alg (n)_{\delta} ,$$
where the product is taken in the category of $\mathds{k}$-modules.
For $\delta \in V_{n}$, we denote by $\pi_{\delta} = \pi_{n,\delta} : \comalg (n) \to \alg (n)_{\delta}$ the projection.

\begin{Defn} \label{202509110944}
    Let $n \subset \mathds{N}$.
    For $\delta \in V_{n}$, we introduce 
    $$\T_{\delta}  = \T_{n,\delta} {:=} \pi_{\delta} \circ \chi : \mathds{k} [\mathsf{F}_n] \to \alg (n)_{\delta} .$$
    For $S \subset {\bf n}$, we also denote by $\T_{S} {:=} \T_{\mathds{1}_S}$ where $\mathds{1}_S \in V_n$ takes values $1$ on $S$ and $0$ otherwise.
    When $S = {\bf n}$, we often use $\T_n {:=} \T_{\bf n}$.
\end{Defn}

\begin{Example}
    Consider $f = x_1^2 x_2^{-1} x_1^3 \in \mathsf{F}_2$.
    Then $$\chi(f) = (1+\theta_1)^2 ( 1+\theta_2)^{-1} (1+\theta_1)^3 =  1 + \cdots + (-2 \theta_1\theta_2 - 3 \theta_2\theta_1) + \cdots .$$
    Hence, $\T_{2} (f) = -2 \theta_1\theta_2 - 3 \theta_2\theta_1$.
\end{Example}

We identify $\alg (n)_{\mathds{1}} = \mathsf{Lin}(n)$.
By Definition \ref{202509110944}, applied to $\delta = \mathds{1} \in V_{n}$, we obtain a map $\mathds{k} [\mathsf{F}_n] \to \alg (n)_{\mathds{1}}$ which can be regarded as
$$\T_{n} : \mathtt{L}_{\G^{\mathsf{o}}}(1,n) \to \mathsf{Lin}(n).$$

\begin{Lemma} \label{202509081651}
    Let $S \subset {\bf n}$.
    For $a,b \in \mathsf{F}_n$, we have $$\T_{S} (ab) = \sum_{S_1 \amalg S_2=S} \T_{S_1} (a) \T_{S_2} (b) . $$
\end{Lemma}
\begin{proof}
    Let $\delta = \mathds{1}_S \in V_n$.
    Then
    \begin{align*}
        \T_{\delta} (ab) = \pi_{\delta} \circ \chi (ab) = \pi_{\delta} (\chi(a)\chi(b)) = \sum_{\delta_1+\delta_2=\delta} \pi_{\delta_1} (\chi(a)) \pi_{\delta_2} (\chi(a))
        = \sum_{\delta_1+\delta_2=\delta} \T_{\delta_1} (a) \T_{\delta_2} (b) .
    \end{align*}
    If we set $S_i$ as the support of $\delta_i$, then we obtain the result.
\end{proof}

\begin{Lemma} \label{202603171154}
    The congruence modulo $\mathtt{I}^{\mathsf{pr}}(1,n)$ on $\mathtt{L}_{\G^{\mathsf{o}}}(1,n)$ is the $\mathds{k}$-linear relation generated by the following:
    for $i \in {\bf n}$ and $w_1, \cdots, w_N$ contained in the subgroup of $\mathsf{F}_n$ generated by $x_1,\cdots,x_{i-1},x_{i+1},\cdots,x_n$,
    \begin{align*}
        w_1  x_i^{a_1+b_1} w_2 x_i^{a_2+b_2}  \cdots w_N x_i^{a_N+b_N} \equiv w_1  x_i^{a_1}  w_2 x_i^{a_2}  \cdots w_N x_i^{a_N} + w_1  x_{i}^{b_1}  w_2 x_{i}^{b_2}  \cdots w_N x_{i}^{b_N}  , \quad a_j,b_j \in \mathds{Z}
    \end{align*}
\end{Lemma}
\begin{proof}
    Since $\mathtt{I}^{\mathsf{pr}}(1,n)$ is generated by $\rho \circ \bar{\Delta}^{n,i}$ for $\rho \in \G^{\mathsf{o}}(n+1,1)$, it suffices to explicitly calculate $\rho \circ \bar{\Delta}^{n,i}$.
    Based on Lemma \ref{202509071205}, we only discuss the case that $i=n$.
    Let $\rho = w_1  x_n^{a_1}  x_{n+1}^{b_1}  w_2 x_n^{a_2} x_{n+1}^{b_2}  \cdots w_N x_n^{a_N} x_{n+1}^{b_N} \in \mathsf{F}_{n+1} \cong \G^{\mathsf{o}}(n+1,1)$ where $w_j,a_j,b_j$ are given as in the statement. 
    Under the bijection $\G^{\mathsf{o}}(n,1) \cong \mathsf{F}_{n}$, we have the following:
    \begin{align*}
        \rho \circ (\mathrm{id}_{n-1} \ast \Delta ) &= w_1  x_n^{a_1+b_1} w_2 x_n^{a_2+b_2}  \cdots w_N x_n^{a_N+b_N} , \\
        \rho \circ (\mathrm{id}_{n} \ast \eta )&= w_1  x_n^{a_1}  w_2 x_n^{a_2}  \cdots w_N x_n^{a_N}  , \\
        \rho \circ (\mathrm{id}_{n-1} \ast \eta \ast \mathrm{id}_{1})&= w_1  x_{n}^{b_1}  w_2 x_{n}^{b_2}  \cdots w_N x_{n}^{b_N} .
    \end{align*}
    Hence, the relation $\rho \circ \bar{\Delta}^{n,n} \equiv 0 \mod{\mathtt{I}^{\mathsf{pr}}}$, and the relations obtained from permutations of variables, are equivalent to that in the statement.
\end{proof}

\begin{Lemma} \label{202509042159}
    Let $f,g \in \mathtt{L}_{\G^{\mathsf{o}}}(1,n)$.
    If $f \equiv g \mod{\mathtt{I}^{\mathsf{pr}}} (1,n)$, then $\T_{n}(f) = \T_{n} (g)$.
    Hence, we obtain a map $\bar{\T}_n : \Psi_{\G^{\mathsf{o}}} (n) \to \mathsf{Lin} (n)$.
\end{Lemma}
\begin{proof}
    It suffices to prove that $\T_n$ respects the relation described in Lemma \ref{202603171154}.
    This is proved by using Lemma \ref{202509081651}.
    For instance, if $N=1$, then we have
    \begin{align*}
        \T_{n} (w_1 x_i^{a_1} x_i^{b_1} ) = \sum_{S_1\amalg S_2 \amalg S_3 ={\bf n}}  \T_{S_1} (w_1)  \T_{S_2} (x_i^{a_1}) \T_{S_3} (x_i^{b_1})  .
    \end{align*}
    If $S_1 = {\bf n}$, then the corresponding term vanishes by the hypothesis on $w_1$.
    Unless either $S_2$ or $S_3$ is empty, the corresponding term also vanishes.
    Hence, 
    \begin{align*}
        \T_{n} (w_1 x_i^{a_1} x_i^{b_1} ) = \sum_{S_1 \amalg S_2 ={\bf n}}  \T_{S_1} (w_1)  \T_{S_2} (x_i^{a_1})  + \sum_{S_1 \amalg S_3 ={\bf n}}  \T_{S_1} (w_1)  \T_{S_3} (x_i^{b_1}) = \T_{n} (w_1 x_i^{a_1} ) + \T_{n} (w_1 x_i^{b_1} ) .
    \end{align*}
    The same proof works for general $N$.
\end{proof}

We now recall the map $\mathcal{F} : \mathsf{Lin} \to \Psi_{\G^{\mathsf{o}}}$ from (\ref{202603171311}).

\begin{Lemma} \label{202509061220}
    For $f \in \mathtt{L}_{\G^{\mathsf{o}}}(1,n)$, we have
    \begin{align*}
        \mathcal{F}_n (\T_{n} (f)) =  [f].
    \end{align*}
\end{Lemma}
\begin{proof}
    It suffices to prove the claim for $f \in \G^{\mathsf{o}}(n,1)$.
    Suppose that $f \in \G^{\mathsf{o}}(n,1) \cong \mathsf{F}_n$ is given by $f=   x_{\sigma(1)}^{a_1} x_{\sigma(2)}^{a_2} \cdots x_{\sigma(n)}^{a_n}$ where $\sigma \in \mathfrak{S}_n$ and $\{ a_1, \cdots, a_n \} \subset \mathds{Z}$.
    One can directly check from definition that 
    $$\T_{n}(f) = (a_1 a_2 \cdots a_n)  \theta_{\sigma(1)} \theta_{\sigma(2)} \cdots \theta_{\sigma(n)} . $$
    We claim that $f$ coincides with $(a_1 a_2 \cdots a_n) f$ modulo $\mathtt{I}^{\mathsf{pr}} (1,n)$.
    Then, in this specific case, we can obtain $\mathcal{F}_n (\T_{n} (f)) =  [f]$.
    
    By iteratively applying Lemma \ref{202603171154}, one may prove that every $f \in \mathtt{L}_{\G^{\mathsf{o}}}(1,n)$ can be identified with a linear combination of the above types, say $f^\prime$, modulo $\mathtt{I}^{\mathsf{pr}}_{\G^{\mathsf{o}}}(1,n)$.
    By Lemma \ref{202509042159}, we have $\mathcal{F}_n (\T_{n} (f)) = \mathcal{F}_n (\T_{n} (f^\prime))$.
    Hence, the previous argument gives $\mathcal{F}_n (\T_{n} (f^\prime)) = [f^\prime]$.
    This proves the desired statement.

    We now prove the above claim.
    If there exists $i$ such that $a_i = 0$, then Proposition \ref{202509042136}, applied to $\mathcal{C} = \G^{\mathsf{o}}$, implies that $f \equiv 0\mod{\mathtt{I}^{\mathsf{pr}}}_{\G^{\mathsf{o}}}$.
    We now assume that $a_i \neq 0$ for all $i$.
    By Lemma \ref{202509071205}, it suffices to prove the assertion for $\sigma = \mathrm{id}$.
    We start with the case that $a_n > 0$.
    By Lemma \ref{202603171154},
    \begin{align*}
        f \equiv x_{1}^{a_1} x_{2}^{a_2} \cdots x_{n-1}^{a_{n-1}}  x_{n}^{a_n-1} + x_{1}^{a_1} x_{2}^{a_2} \cdots  x_{n-1}^{a_{n-1}} x_{n} \mod{\mathtt{I}^{\mathsf{pr}}_{\G^{\mathsf{o}}}(1,n)},
    \end{align*}
    so, by induction, we can obtain
    $$
    f \equiv a_n x_{1}^{a_1} x_{2}^{a_2} \cdots  x_{n-1}^{a_{n-1}} x_{n} \mod{\mathtt{I}^{\mathsf{pr}}_{\G^{\mathsf{o}}}(1,n)} .
    $$
    This is also the case for $a_n < 0$:
    by Lemma \ref{202603171154}, we have modulo $\mathtt{I}^{\mathsf{pr}}_{\G^{\mathsf{o}}} (1,n)$
    \begin{align*}
        x_{1}^{a_1} x_{2}^{a_2} \cdots  x_{n-1}^{a_{n-1}} x_{n}^{a_n}  x_{n}^{-a_n}  \equiv&  x_{1}^{a_1} x_{2}^{a_2} \cdots  x_{n-1}^{a_{n-1}} x_{n}^{a_n}  + x_{1}^{a_1} x_{2}^{a_2} \cdots  x_{n-1}^{a_{n-1}} x_{n}^{-a_n}  , \\
        \equiv& x_{1}^{a_1} x_{2}^{a_2} \cdots  x_{n-1}^{a_{n-1}} x_{n}^{a_n} -a_n x_{1}^{a_1} x_{2}^{a_2} \cdots  x_{n-1}^{a_{n-1}} x_n  ,
    \end{align*}
    while Proposition \ref{202509042136}, applied to $\mathcal{C} = \G^{\mathsf{o}}$, implies that 
    \begin{align*}
        x_{1}^{a_1} x_{2}^{a_2} \cdots  x_{n-1}^{a_{n-1}} x_{n}^{a_n}  x_{n}^{-a_n} = x_{1}^{a_1} x_{2}^{a_2} \cdots  x_{n-1}^{a_{n-1}}  \equiv 0 \mod{\mathtt{I}^{\mathsf{pr}}_{\G^{\mathsf{o}}}(1,n)} .
    \end{align*}
    The same arguments for other $a_i$'s complete the proof.
\end{proof}

\begin{proof}[Proof of Theorem \ref{202509091931}]
    Let $n \in \mathds{N}$.
    We have $\mathcal{F}_n \circ \bar{\T}_n = \mathrm{id}$ by Lemma \ref{202509061220}, so it is sufficient to prove that $\bar{\T}_n \circ \mathcal{F}_n = \mathrm{id}$.
    Let $f = \theta_{\sigma(1)}\theta_{\sigma(2)} \cdots \theta_{\sigma(n)} \in \mathsf{Lin}(n)$ be the generator corresponding to $\sigma \in \mathfrak{S}_n$.
    Then we obtain
    \begin{align*}
        \bar{\T}_n ( \mathcal{F}_n (f)) =& \bar{\T}_n \left(   x_{\sigma(1)}x_{\sigma(2)}  \cdots  x_{\sigma(n)} \right) = (\pi_{\mathds{1}} \circ \chi) \left( x_{\sigma(1)} x_{\sigma(2)} \cdots x_{\sigma(n)} \right)  , \\
        =& \pi_{\mathds{1}} \left( \sum \left(   \theta_{\sigma(1)}^{\epsilon_1} \theta_{\sigma(2)}^{\epsilon_2} \cdots \theta_{\sigma(n)}^{\epsilon_n} \right) \right) =  \theta_{\sigma(1)} \theta_{\sigma(2)}\cdots\theta_{\sigma(n)} =f ,
    \end{align*}
    where the sum is taken all over $\epsilon_j \in \{0,1\}$ for $j \in {\bf n}$.
\end{proof}

\begin{remark}
    It is well known that the linear-order species $\mathsf{Lin}$, or equivalently the associated symmetric module, forms the operad governing associative unital algebras.
    By the constructions in Section \ref{202512151740}, $\mu\mathsf{Lin}$ admits a monad structure.
    We note that the map $P : \mathtt{L}_{\G^{\mathsf{o}}} \twoheadrightarrow  \mathtt{L}_{\G^{\mathsf{o}}}/ \mathtt{I}^{\mathsf{pr}}_{\G^{\mathsf{o}}} \cong \mu\mathsf{Lin}$ obtained from Theorem \ref{202509091931} is not a monad map for nontrivial $\mathds{k}$.
    Indeed, $P ( x_1 , x_1 ) = \theta_1 \otimes e  + e \otimes  \theta_1$ and $P ( x_1 ) = 0$ where we regard $(x_1,x_1) \in \mathsf{F}_1^{\times 2} = \G^{\mathsf{o}}(1,2)$ and $x_1 \in \mathsf{F}_2 = \G^{\mathsf{o}}(2,1)$, so $P ( x_1) \circ P ( x_1 , x_1 ) = 0$ while $P ((x_1 ) \circ (x_1 , x_1 ) ) = P (x_1) = \theta_1$, where $x_1 \in \mathsf{F}_1 = \G^{\mathsf{o}}(1,1)$.
\end{remark}

\begin{remark}
    Recall the Lawvere theory $\mathbf{W}^{\mathsf{o}}$ from Section \ref{202509031755}.
    Using the same method as above, one can establish that the canonical map $\mathsf{Lin} \to \Psi_{\mathbf{W}^{\mathsf{o}}}$, defined analogously to $\mathcal{F}$, is an isomorphism, with inverse given by a partial Magnus expansion. However, rather than pursuing this parallel argument, we instead give a conceptual proof of the isomorphism $\mathsf{Lin} \cong \Psi_{\mathbf{W}^{\mathsf{o}}}$ in a later section (see Corollary \ref{202603251841}).
\end{remark}

\subsection{Proof of Theorem \ref{202603161027}}
\label{202604011319}

In this section, we prove Theorem \ref{202603161027} which states that the map $\alpha_\mathcal{C} : \mu\Psi_{\mathcal{C}} \to \mathtt{L}_{\mathcal{C}}/\mathtt{I}^{\mathsf{pr}}_{\mathcal{C}}$ is an isomorphism.
We construct a map $\delta_\mathcal{C}$ that is shown to be an inverse by Lemmas \ref{202603161242} and \ref{202603161244}.

\subsubsection{Preliminaries for maps}

As the first step, we introduce several constructions of maps between finite sets.

\begin{Defn} \label{202603161239}
    Let $\rho : {\bf n} \to {\bf m}$ be a map.
    For $k \in {\bf n}$, let $S_\rho (k)$ be the set,
    $$
     \{ i \in {\bf n} \mid \rho (i) < \rho (k) \mathrm{~or~} (\rho(i) = \rho (k) \mathrm{~and~} i \leq k)\}.
    $$
    Let $\sigma_\rho (k)$ be the cardinality of $S_\rho (k)$.
    This determines a map $\sigma_\rho : {\bf n} \to {\bf n}$, since $\emptyset \neq S_\rho (k) \subset {\bf n}$.
\end{Defn}

\begin{Lemma}
    Let $\rho_0$ be the map associated with the partition $(n_1 , \cdots ,n_m )$ (see Definition \ref{202604121701}).
    Then $\sigma_{\rho_0}$ is the identity.
\end{Lemma}
\begin{proof}
    From the definition of $\rho_0$, we obtain $S_{\rho_0} (k) = {\bf k}$ for any $k \in {\bf n}$.
\end{proof}

\begin{Lemma}
    For any map $\rho$, $\sigma_\rho$ is a bijection, i.e. $\sigma_\rho \in \mathfrak{S}_n$.
\end{Lemma}
\begin{proof}
    It suffices to show that $\sigma_\rho$ is injective since both the domain and the codomain have the same cardinality.
    Let $k_1,k_2 \in {\bf n}$.
    We assume that $\sigma_\rho (k_1) = \sigma_\rho (k_2)$.
    We first show that $\rho (k_1)=\rho (k_2)$.
    If $\rho (k_1) < \rho (k_2)$, then $S_{\rho}(k_1) \subset S_{\rho}(k_2)$ follows from the definition of $S_\rho$.
    Furthermore, $k_2 \in S_\rho (k_2) \backslash S_{\rho}(k_1)$, since $\rho (k_1) < \rho (k_2)$.
    Hence, $S_{\rho}(k_1) \subsetneq S_{\rho}(k_2)$ so $\sigma_{\rho} ( k_1) < \sigma_{\rho} (k_2)$, which contradicts the hypothesis.
    Likewise, the inequality $\rho (k_2) < \rho (k_1)$ cannot occur, so $\rho (k_1)=\rho (k_2)$.
    
    We now prove $k_1 = k_2$.
    If $k_1 < k_2$, then $S_{\rho}(k_1) \subset S_{\rho}(k_2)$, since $\rho (k_1)=\rho (k_2)$ and $k_1 < k_2$.
    We also obtain $k_2 \in S_\rho (k_2) \backslash S_{\rho}(k_1)$ from the inequality $k_1 < k_2$.
    This implies $\sigma_\rho (k_1) < \sigma_\rho (k_2)$, which contradicts the hypothesis.
    Analogously, the inequality $k_2 < k_1$ is also impossible.
\end{proof}

\begin{Defn} \label{202606101643}
    Let $\rho: {\bf n} \to {\bf m}$ be a map and write $n_k = |\rho^{-1}(k)|,~ k \in {\bf m}$.
    Note that $(n_1,\cdots,n_m)$ is an $m$-parition of $n$.
    Let $\iota_k = \iota_{\rho^{-1}(k)} : {\bf n_k} \to {\bf n}$ be the order-preserving map with image $\rho^{-1}(k)$.
    Let $\tilde{\rho} : {\bf n} \to {\bf nm}$ be the composition of $\coprod^{m}_{k=1} \iota_k : {\bf n} \to {\bf nm}$ and $\sigma_\rho : {\bf n} \to {\bf n}$.
    Clearly, the map $\tilde{\rho}$ is injective.
\end{Defn}

\begin{Lemma} \label{202603170954}
    For $j \in {\bf n}$, 
    we have $\tilde{\rho} (j) = n(\rho (j) -1) + j$.
\end{Lemma}
\begin{proof}
    Let $j \in {\bf n}$ and $s = \rho (j)$.
    Since $j \in \rho^{-1}(s)$, we can take a unique $l \in {\bf n_s}$ such that $j = \iota_s (l)$.
    In other words, $j$ is the $l$-th smallest element of $\rho^{-1}(s)$.
    Hence, by the definition of $\sigma_\rho$, we have
    $$\sigma_\rho (j) = n_1 + \cdots + n_{s-1} + l.$$
    This leads to 
    $$\tilde{\rho}(j) = (\coprod^{m}_{k=1} \iota_k) \circ \sigma_\rho (j) = (\coprod^{m}_{k=1} \iota_k) (n_1 + \cdots + n_{s-1} + l) = n(s-1) + j . $$
\end{proof}

\subsubsection{The map $\delta_\mathcal{C}$}

We now construct the inverse $\delta_{\mathcal{C}}$ of $\alpha_\mathcal{C}$.
To this end, we extend the functor $E$ in Definition \ref{202603261255} to injective maps:
\begin{Defn} \label{202603161044}
    Let $\iota : {\bf n} \to {\bf m}$ be an injective map.
    We define $E_\iota \in \mathcal{C}(n,m)$ as the morphism characterized by
    \begin{align*}
        p_{m,k} \circ E_\iota =
        \begin{cases}
            p_{n,l} & \mathrm{if~} k \in \iota ({\bf n}) \mathrm{~and~}k = \iota (l) , \\
            e & \mathrm{otherwise}.
        \end{cases}
    \end{align*}
    In particular, if $n=0$, then $E_{\iota} = e$, the zero morphism.

    We also make use of the following notation.
    Let $S \subset {\bf n}$ with $|S| = r$.
    We denote by $\iota_S : {\bf r} \to {\bf n}$ the unique order-preserving map whose image is $S$.
    We define $$E_{n,S} {:=} E_{\iota_S} \in \mathcal{C}(r,n).$$
    If $S \subset {\bf n}$ is a singleton, say $S = \{ a\}$, then we also denote by $E_{n,a} {:=} E_{n,\{a\}}$.
\end{Defn}

For instance, $E_{n,{\bf n} \backslash \{k\}} = \mathrm{id}_{k-1} \times \eta \times \mathrm{id}_{n-k}$.

For $n_1,n_2 \in \mathds{N}$, we identify ${\bf n_1} \amalg {\bf n_2} \cong {\bf n_1 +n_2}$ by identifying $k \in {\bf n_1}$ with $k \in {\bf n_1 +n_2}$ and $k \in {\bf n_2}$ with $k+ n_1 \in {\bf n_1 +n_2}$.
Then the following is immediate from definitions:
\begin{Lemma} \label{202606191144}
    \begin{enumerate}
        \item For injective maps $\iota : {\bf n} \to {\bf m}$ and $\iota^\prime : {\bf l} \to {\bf n}$, we also have $E_{\iota} \circ E_{\iota^\prime} = E_{\iota \circ \iota^\prime}$.
        \item For injective maps $\iota_i : {\bf n_i} \to {\bf m_i}$, we have $E_{\iota_1} \times E_{\iota_2} = E_{\iota_1 \amalg \iota_2}$ where $\iota_1 \amalg \iota_2$ is the map ${\bf n_1 +n_2} \cong {\bf n_1} \amalg {\bf n_2} \to {\bf m_1} \amalg {\bf m_2} \cong {\bf m_1 +m_2}$.
    \end{enumerate}
    In particular, the assignments ${\bf n} \mapsto n$ and $\iota \mapsto E_{\iota}$ define a monoidal functor from the category of finite sets and injective maps to $\mathcal{C}$.
\end{Lemma}

\begin{Lemma} \label{202603170910}
    Let $n,m,r\in \mathds{N}$.
    Let $S \subset {\bf n}$ with $|S| = r$, and $S^\prime \subset {\bf m}$ with $|S^\prime| = n$.
    Then we have $E_{m,S^\prime} \circ E_{n,S} = E_{m,\iota_{S^\prime}(S)}$.
\end{Lemma}
\begin{proof}
    The result is obtained from the previous lemma:
    \begin{align*}
        E_{m,S^\prime} \circ E_{n,S} = E_{\iota_{S^\prime}} \circ E_{\iota_S} = E_{\iota_{S^\prime} \circ \iota_S} =  E_{m,\iota_{S^\prime}(S)} .
    \end{align*}
\end{proof}

\begin{Defn}
    Let $\rho : {\bf n} \to {\bf m}$ be a map and write $|\rho^{-1}(k)| = n_k,~k\in {\bf m}$.
    Let $F_{m,n,\rho} = F_\rho : \mathtt{L}_{\mathcal{C}} (m,n) \to \bigotimes^{m}_{k=1} \Psi_{\mathcal{C}} (n_k)$ be the map defined by $\mathds{k}$-linearly extending
    \begin{align*}
        F_\rho (f) {:=} \bigotimes^{m}_{k=1} [f_k \circ E_{n,\rho^{-1}(k)} ] , \quad f= (f_1,\cdots,f_m) \in \mathcal{C}_n^{\times m} = \mathcal{C}(n,m) .
    \end{align*}
\end{Defn}

\begin{Example} \label{202606201827}
    If $m=1$, then $F_\rho$ is the quotient map $\mathtt{L}_{\mathcal{C}} (1,n) \to \Psi_{\mathcal{C}} (n)$.
\end{Example}

In this section, for $S \subset \mathds{N}$, let $$S_k {:=} \{ i \mid i \in S, i \leq k\} \cup \{ i+1 \mid i \in S,  k<i\}.$$
Then the following is immediate.
\begin{Lemma} \label{202603161155}
    If $k \not\in S$, then $S_k = S_{k-1}$.
\end{Lemma}

For $n\in\mathds{N}$ and $k \in {\bf n}$, let $\Delta^{n,k}$ denote the morphism $\mathrm{id}_{k-1} \times \Delta \times \mathrm{id}_{n-k} \in \mathcal{C} (n, n+1)$.

\begin{Lemma} \label{202603141056}
    Let $S \subset {\bf n}$ with $|S| = r$. 
    For each $k \in {\bf n}$, the following identity holds in $\mathtt{L}_{\mathcal{C}}/\mathtt{I}^{\mathsf{pr}}_{\mathcal{C}} (n+1,r)$:
    \begin{align*}
        \Delta^{n,k} \circ E_{n,S} =
        \begin{cases}
            E_{n+1,S_k} + E_{n+1,S_{k-1}} & \mathrm{if~} k \in S , \\
            E_{n+1,S_k} & \mathrm{otherwise}.
        \end{cases}
    \end{align*}
\end{Lemma}
\begin{proof}
    Recall from Definition \ref{202603161044} the notation $\iota_S$.
    We begin with the case where $k \not\in S$.
    It is clear that $\{k,k+1\} \cap S_k = \emptyset$.
    For $i \in \{ k,k+1 \}$, $p_{n+1,i} \circ \Delta^{n,k} \circ E_{n,S} = p_{n,i} \circ E_{n,S} = e$.
    For $i \not\in S_k$, we have $p_{n+1,i} \circ \Delta^{n,k} \circ E_{n,S} = p_{n,i} \circ E_{n,S} = e$ if $i <k$; and $p_{n+1,i} \circ \Delta^{n,k} \circ E_{n,S} = p_{n,i-1} \circ E_{n,S} = e$ if $i > k+1$.
    For $i \in S_k$, we have $p_{n+1,i} \circ \Delta^{n,k} \circ E_{n,S} = p_{n,i} \circ E_{n,S} = p_{r,l}$ if $\iota_S (l)=i < k$; and $p_{n+1,i} \circ \Delta^{n,k} \circ E_{n,S} = p_{n,i-1} \circ E_{n,S} = p_{r,l}$ if $\iota_S (l) = i-1 > k$, where we use $\iota_{S_k} (l) =\iota_S (l) +1$.
    Hence, by the universality of products in $\mathcal{C}$, we obtain $\Delta^{n,k} \circ E_{n,S} = E_{n+1,S_k}$.
    
    We now consider the case that $k \in S$.
    Let $l \in {\bf r}$ such that $k = \iota_S (l)$.
    By the universality of products in $\mathcal{C}$, we have $\Delta^{n,k} \circ E_{n,S} = E_{n+1,S^\prime} \circ \Delta^{r,l}$ where $S^\prime = \{ i \mid i\in S, i \leq k \} \cup \{ i +1 \mid i \in S , i\geq k \}$.
    By the definition of $\mathtt{I}^{\mathsf{pr}}$, we have
    \begin{align*}
        E_{n+1,S^\prime} \circ \Delta^{r,l} \equiv E_{n+1,S^\prime} \circ ( \mathrm{id}_{l} \times \eta \times \mathrm{id}_{r-l} + \mathrm{id}_{l-1} \times \eta \times \mathrm{id}_{r-l+1} ) \mod{\mathtt{I}^{\mathsf{pr}}_{\mathcal{C}} (n+1,r)}.
    \end{align*}
    Since $\mathrm{id}_{l} \times \eta \times \mathrm{id}_{r-l} = E_{r,{\bf r+1} \backslash \{l+1\}}$, we obtain $E_{n+1,S^\prime} \circ (\mathrm{id}_{l} \times \eta \times \mathrm{id}_{r-l} )= E_{n+1,S_k}$ from Lemma \ref{202603170910}.
    Similarly, $E_{n+1,S^\prime} \circ (\mathrm{id}_{l-1} \times \eta \times \mathrm{id}_{r-l+1} )= E_{n+1,S_{k-1}}$.
\end{proof}

\begin{Lemma}
    The map $F_\rho$ vanishes on $\mathtt{I}^{\mathsf{pr}}_{\mathcal{C}} (m,n)$.
    Hence, $F_\rho$ induces a map $$\bar{F}_{\rho} : \mathtt{L}_{\mathcal{C}} / \mathtt{I}^{\mathsf{pr}}_{\mathcal{C}}(m,n) \to \bigotimes^{m}_{i=1} \Psi_{\mathcal{C}} (n_i).$$
\end{Lemma}
\begin{proof}
    We shall prove that $F_\rho ( f \circ \bar{\Delta}^{n,k}) = 0$ for $f \in \mathcal{C} (n+1,m)$ and $k \in {\bf n}$.  
    This is equivalent to showing that $F_\rho ( f \circ \Delta^{n,k}) = F_\rho ( f \circ (\mathrm{id}_k \times \eta \times \mathrm{id}_{n-k})) + F_\rho ( f \circ (\mathrm{id}_{k-1} \times \eta \times \mathrm{id}_{n-k+1}))$.
    Let $f = (f_1, \cdots, f_m) \in \mathcal{C}(n+1,m)$ where $f_k \in \mathcal{C}(n+1,1)$.
    We begin by calculating $F_\rho ( f \circ \Delta^{n,k})$ which equals $F_\rho ( f_1 \circ \Delta^{n,k} , \cdots f_m \circ \Delta^{n,k} )$ by Lemma \ref{202509042125}.
    We set $\rho (k) = i_0$.
    By Lemma \ref{202603141056}, we obtain
    \begin{align*}
        &F_\rho ( f_1 \circ \Delta^{n,k} , \cdots , f_m \circ \Delta^{n,k} )  \\
        =&  \bigotimes^{m}_{i=1} [f_i \circ \Delta^{n,k} \circ E_{n,\rho^{-1}(i)}] = \bigotimes^{i_0-1}_{i=1} [f_i \circ \Delta^{n,k} \circ E_{n,\rho^{-1}(i)}]  \otimes [f_{i_0} \circ \Delta^{n,k} \circ E_{n,\rho^{-1}(i_0)}] \otimes  \bigotimes^{m}_{i=i_0+1} [f_i \circ \Delta^{n,k} \circ E_{n,\rho^{-1}(i)}]   , \\
        =& \bigotimes^{i_0-1}_{i=1} [f_i \circ E_{n+1,\rho^{-1}(i)_{k}}]  \otimes \Bigl\{ [ f_{i_0} \circ E_{n+1,\rho^{-1}(i_0)_{k}} ] + [ f_{i_0} \circ E_{n+1,\rho^{-1}(i_0)_{k-1}} ] \Bigr\}  \otimes \bigotimes^{m}_{i=i_0+1} [f_{i} \circ E_{n+1,\rho^{-1}(i)_k}] .
    \end{align*}
    By Lemma \ref{202603161155}, this result coincides with $\bigotimes^{m}_{i=1} [f_i \circ E_{n+1,\rho^{-1}(i)_{k}}]+ \bigotimes^{m}_{i=1} [f_i \circ E_{n+1,\rho^{-1}(i)_{k-1}}]$.

    We now calculate $F_\rho ( f \circ (\mathrm{id}_c \times \eta \times \mathrm{id}_{n-c}))$, where $c \in {\bf n}$.
    Recall that $\mathrm{id}_c \times \eta \times \mathrm{id}_{n-c} = E_{n+1,{\bf n+1}\backslash \{c+1\}}$.
    Lemma \ref{202603141056} gives
    \begin{align*}
        F_\rho ( f \circ (\mathrm{id}_c \times \eta \times \mathrm{id}_{n-c}))  =\bigotimes^{m}_{i=1} [f_i \circ E_{n+1,{\bf n+1}\backslash \{c+1\}} \circ E_{n,\rho^{-1}(i)} ]   =  \bigotimes^{m}_{i=1} [f_i \circ E_{n+1,\rho^{-1}(i)_{c}}] ,
    \end{align*}
    since $E_{n+1,{\bf n+1}\backslash \{c+1\}} \circ E_{n,\rho^{-1}(i)} = E_{n+1,\rho^{-1}(i)_{c}}$ by Lemma \ref{202603170910}.
    Applying this to $c \in \{ k-1,k\}$, we obtain
    \begin{align*}
        &F_\rho ( f \circ (\mathrm{id}_k \times \eta \times \mathrm{id}_{n-k})) + F_\rho ( f \circ (\mathrm{id}_{k-1} \times \eta \times \mathrm{id}_{n-k+1}))  \\
        =& \bigotimes^{m}_{i=1} [f_i \circ E_{n+1,\rho^{-1}(i)_{k}}]+ \bigotimes^{m}_{i=1} [f_i \circ E_{n+1,\rho^{-1}(i)_{k-1}}] .
    \end{align*}
    Comparing this with the above result, we obtain the desired statement.
\end{proof}

Using $\bar{F}_{\rho}$ and $\sigma_\rho$ defined before, we introduce the following map:
\begin{Defn}
    Let $\delta_{\mathcal{C}}$ be the map $\mathtt{L}_{\mathcal{C}}/ \mathtt{I}^{\mathsf{pr}}_{\mathcal{C}} \to \mu\Psi_{\mathcal{C}}$ with $(m,n)$-component given by 
    $$
    \delta_{\mathcal{C}} ( u) {:=} \sum_{\rho} \bar{F}_{\rho} (u) \otimes \sigma_\rho \in \bigoplus \bigotimes^{m}_{i=1} \Psi_{\mathcal{C}} (n_i) \otimes_{\mathfrak{S}_{n_1} \times \cdots \mathfrak{S}_{n_m}} \mathds{k}[\mathfrak{S}_n] ,
    $$
    where $\rho$ runs over all maps $\rho : {\bf n} \to {\bf m}$ and the direct sum is taken over $m$-partitions $(n_1,\cdots,n_m)$ of $n$.
\end{Defn}

\begin{Example}
    By Example \ref{202606201827}, the $(1,-)$-component of $\delta_\mathcal{C}$ is the identity on $\Psi_\mathcal{C}$.
\end{Example}

\subsubsection{Proof of $\delta_\mathcal{C} \circ \alpha_\mathcal{C} = \mathrm{id}$}

As part of the proof of Theorem \ref{202603161027}, we prove that the composition $\delta_\mathcal{C} \circ \alpha_\mathcal{C}$ is the identity on $\mu\Psi_\mathcal{C}$.
Recall the projection $P_{n,S}$ from Definition \ref{202606181952}.

\begin{Lemma} \label{202603141952}
    For $S,S^\prime \subset {\bf n}$ such that $|S| = r$ and $|S^\prime|=r^\prime$, we have
    \begin{align*}
        P_{n,S^\prime} \circ E_{n,S} \mod{\mathtt{I}^{\mathsf{pr}}_{\mathcal{C}}(r^\prime, r)}  =
        \begin{cases}
            0 & \mathrm{if~} S\backslash S^\prime \neq \emptyset , \\
            \mathrm{id}_{|S|} & \mathrm{if~}  S = S^\prime .
        \end{cases}
    \end{align*}
    In particular, $P_{n,S^\prime} \circ E_{n,S} =0\mod{\mathtt{I}^{\mathsf{pr}}}_{\mathcal{C}} (r^\prime, r)$ if $r^\prime < r$; or $r= r^\prime$ and $S \neq S^\prime$.
\end{Lemma}
\begin{proof}
    Observe that the morphism $E_{n,S}$ coincides with the product $\prod^{n}_{i=1} h_i$ where $h_i = \mathrm{id}_1$ if $i \in S$, and $h_i = \eta$ otherwise.
    Likewise, $P_{n,S^\prime}$ is the product $\prod^{m}_{i=1} v_i$ where $v_i = \mathrm{id}_1$ if $i \in S^\prime$ and $v_i = \epsilon$ otherwise.
    Then $P_{n,S^\prime} \circ E_{n,S} = \prod^{n}_{i=1} v_i \circ h_i$.
    Thus, if $S=S^\prime$, then $v_i \circ h_i = \mathrm{id}_1$ for $i \in S$ and $e$ for $i \not\in S$, so that $P_{n,S^\prime} \circ E_{n,S} = \mathrm{id}_{|S|}$.
    If $S \backslash S^\prime \neq \emptyset$, then we choose $j \in S \backslash S^\prime$.
    Then
    $\prod^{m}_{i=1} v_i \circ h_i = \prod^{j-1}_{i=1} v_i \circ h_i \times \epsilon \times \prod^{m}_{i=j+1} v_i \circ h_i \in \mathtt{I}^{\mathsf{pr}}_{\mathcal{C}} (r^\prime, r)$.
    Indeed, this follows from Lemma \ref{202601051527} and the fact $\epsilon \in \mathtt{I}^{\mathsf{pr}}_{\mathcal{C}} (0, 1)$.
\end{proof}

\begin{Lemma} \label{202603161242}
    $\delta_{\mathcal{C}} \circ \alpha_{\mathcal{C}} = \mathrm{id}$.
\end{Lemma}
\begin{proof}
    By the definitions of $\alpha_\mathcal{C}$ and $\delta_\mathcal{C}$, it is sufficient to show that, for $(n_1 , \cdots ,n_m )$ an $m$-partition of $n$ and $g_k \in \mathcal{C}(n_k,1), ~k\in {\bf m}$,
    \begin{align} \label{202603161229}
    \left( \bigotimes^{m}_{k=1} [g_k]  \right) \otimes e = \sum_\rho F_\rho \left( \prod^{m}_{k=1} g_k \right) \otimes \sigma_\rho .
    \end{align}
    Let $\rho_0$ be the map associated with the partition $(n_1 , \cdots ,n_m )$ (see Definition \ref{202604121701}).
    It is immediate from Definition \ref{202603161239} that $\sigma_{\rho_0} =e \in \mathfrak{S}_n$.
    Recall from (\ref{202606182102}) the equality $\prod^{m}_{k=1} g_k = (f_1,\cdots,f_m)$ with $f_k = g_k \circ P_{n,\rho_0^{-1} (k)}, ~k\in {\bf m}$.
    We then obtain
    $$F_\rho \left( \prod^{m}_{k=1} g_k \right)  = \bigotimes^{m}_{k=1} [f_k \circ E_{n,\rho^{-1}(k)} ] = \bigotimes^{m}_{k=1} [g_k \circ P_{n,\rho_0^{-1} (k)} \circ E_{n,\rho^{-1}(k)} ] .$$
    In particular, $F_{\rho_0} \left( \prod^{m}_{k=1} g_k \right)  = \bigotimes^{m}_{k=1} [g_k ]$ by Lemma \ref{202603141952}.
    We claim that if $\rho \neq \rho_0$, then $F_\rho \left( \prod^{m}_{k=1} g_k \right)  = 0$.
    Based on this claim, we can complete the proof of (\ref{202603161229}).
    Indeed, $$\sum_\rho F_\rho \left( \prod^{m}_{k=1} g_k \right)  \otimes \sigma_\rho = F_{\rho_0} \left( \prod^{m}_{k=1} g_k \right)  \otimes \sigma_{\rho_0} = \bigotimes^{m}_{k=1} [g_k ] \otimes e.$$

    We now prove the above claim.
    If there exists $k \in {\bf m}$ such that $|\rho_0^{-1} (k)| < |\rho^{-1}(k)|$, then Lemma \ref{202603141952} leads to $[g_k \circ P_{n,\rho_0^{-1} (k)} \circ E_{n,\rho^{-1}(k)} ] = 0$, so $F_\rho (g_1 \times \cdots \times g_m)= 0$.
    If $|\rho_0^{-1} (k)| \geq |\rho^{-1}(k)|$ for any $k\in {\bf m}$, then $|\rho_0^{-1} (k)| = |\rho^{-1}(k)|$ for any $k \in {\bf m}$, since $\sum^{m}_{k=1} |\rho_0^{-1} (k)| = n = \sum^{m}_{k=1} |\rho^{-1} (k)|$.
    By the hypothesis that $\rho \neq \rho_0$, there exists $k \in {\bf m}$ such that $\rho_0^{-1} (k) \neq \rho^{-1}(k)$.
    By Lemma \ref{202603141952}, for such $k$, $[g_k \circ P_{n,\rho_0^{-1} (k)} \circ E_{n,\rho^{-1}(k)} ] = 0$, so $F_\rho (g_1 \times \cdots \times g_m)= 0$.
\end{proof}

\subsubsection{Proof of $\alpha_\mathcal{C} \circ \delta_\mathcal{C} = \mathrm{id}$}

In order to complete the proof of Theorem \ref{202603161027}, we now show that $\alpha_\mathcal{C} \circ \delta_\mathcal{C}$ is the identity.
This argument heavily relies on a congruence description of the iterative comulitplication $\Delta^{(m)}_n$ modulo $\mathtt{I}^{\mathsf{pr}}_{\mathcal{C}}$ using $E_{\tilde{\rho}}$'s where $\tilde{\rho}$ is introduced in Definition \ref{202606101643}.
Since $\tilde{\rho}$ is injective, $E_{\tilde{\rho}}$ in Definition \ref{202603161044} is well-defined.
To formulate the results, we begin with some lemmas:

\begin{Defn}
    Let $a,m \in \mathds{N}$ and $r = (a+1)m$.
    Let $\kappa = \kappa_{a,m} : {\bf r} \to \mathds{N}$ be the map defined as
    \begin{align*}
        \kappa (i) {:=}
        \begin{cases}
            i -1 + \lceil i/ a\rceil & \mathrm{if~} i \leq am , \\
            (i-am) (a +1) & \mathrm{otherwise}  ,
        \end{cases}
    \end{align*}
    where $\lceil l \rceil$ denotes the smallest integer greater than or equal to $l$.
\end{Defn}

\begin{Lemma}
    $\kappa$ is a bijection onto the set ${\bf r}$, i.e. $\kappa \in \mathfrak{S}_r$.
    Therefore, $E_{\kappa}$ is defined as in Definition \ref{202603161044}.
\end{Lemma}
\begin{proof}
    For $i \in {\bf r}$, if $i \leq am$, then $\kappa (i) = i -1 + \lceil i/ a\rceil \leq i-1 + m <(a+1)m=r$.
    If $am < i \leq (a+1)m$, then $\kappa (i) = (i-am) (a +1) \leq m(a+1) = r$.
    Hence, $\kappa$ is a map onto ${\bf r}$.
    To show that $\kappa$ is a bijection onto ${\bf r}$, it suffices to prove that this map is injective.
    Let $i \in {\bf r}$.
    Note that, for $i \in \mathds{N}$ with $q = \lceil i/ a \rceil$, there exists a unique $i^\prime \in {\bf a}$ such that $i= (q-1)a + i^\prime$.
    By the definition of $\kappa$, we have $\kappa (i) = (q-1)(a+1) + i^\prime$ if $i \leq am$, and $\kappa (i) = ( ( q-m-1)a + i^\prime)(a+1)$ otherwise.
    These expressions, with their residues modulo $a+1$, determine $q$ and $i^\prime$ uniquely.
    Hence, $\kappa (i) = \kappa (j)$ implies $i= j$.
\end{proof}

The permutation $\kappa_{a,m}$ naturally appears when we compare the morphisms $\Delta^{(m)}_{a+1}$ and $\Delta^{(m)}_{a} \times \Delta^{(m)}_{1}$:

\begin{Lemma} \label{202603141912}
    For $a,m \in \mathds{N}$, we have $\Delta^{(m)}_{a+1} = E_{\kappa} \circ ( \Delta^{(m)}_{a} \times \Delta^{(m)}_{1})$.
\end{Lemma}
\begin{proof}
    Let $r= (a+1)m$ as before.
    By the universality of products in $\mathcal{C}$, it suffices to show $p_{r,s} \circ \Delta^{(m)}_{a+1} = p_{r,s} \circ E_{\kappa} \circ ( \Delta^{(m)}_{a} \times \Delta^{(m)}_{1})$ for $s \in {\bf r}$.
    To this end, it is useful to recall Definition \ref{202603161044}.
    Let $q = \lceil s / (a+1) \rceil \in {\bf m}$.
    There exists a unique $s^\prime \in {\bf a+1}$ such that $s = (a+1)(q-1) + s^\prime$.
    If $s^\prime \leq a$, then 
    $$
    p_{r,s} \circ E_{\kappa} \circ ( \Delta^{(m)}_{a} \times \Delta^{(m)}_{1}) = p_{r,a (q-1) + s^\prime} \circ ( \Delta^{(m)}_{a} \times \Delta^{(m)}_{1}) = p_{a+1,s^\prime} = p_{r,s} \circ \Delta^{(m)}_{a+1} .
    $$
    If $s^\prime = a +1$, then 
    $$
    p_{r,s} \circ E_{\kappa} \circ ( \Delta^{(m)}_{a} \times \Delta^{(m)}_{1}) = p_{r,am + q} \circ ( \Delta^{(m)}_{a} \times \Delta^{(m)}_{1}) = p_{a+1,a+1} = p_{r,s} \circ \Delta^{(m)}_{a+1} .
    $$
\end{proof}

\begin{Lemma} \label{202606191206}
    Let $\rho_1 : {\bf a} \to {\bf m}$ and $\rho_2 : {\bf 1} \to {\bf m}$ be maps.
    Let $g_{\rho_1,\rho_2} = \nabla \circ (\rho_1 \amalg \rho_2)$ where $\nabla : {\bf 2m} = {\bf m+m} \to {\bf m}$ denotes the folding map.
    Then $\kappa \circ (\tilde{\rho}_1 \amalg \tilde{\rho}_2) = \widetilde{g_{\rho_1,\rho_2}}$.
\end{Lemma}
\begin{proof}
    We freely use the result in Lemma \ref{202603170954} below.
    For $i \in {\bf a +1}$, we shall show $(\kappa \circ (\tilde{\rho}_1 \amalg \tilde{\rho}_2)) (i) = \widetilde{g_{\rho_1,\rho_2}} (i)$.
    If $i \leq a$, then $(\kappa \circ (\tilde{\rho}_1 \amalg \tilde{\rho}_2)) (i)$ is computed as
    $$
    i \stackrel{\tilde{\rho}_1 \amalg \tilde{\rho}_2}{\mapsto} a ( \rho_1 ( i) -1) + i
    \stackrel{\kappa_{a,m}}{\mapsto} a ( \rho_1 (i) -1) + i-1 + \lceil \rho_1 (i) -1 + i/a \rceil = (a+1) ( \rho_1 (i) -1) + i , $$
    and we have
    $$
    \widetilde{g_{\rho_1,\rho_2}} (i) = (a+1) (g_{\rho_1,\rho_2}(i)-1) + i = (a+1) ( \rho_1 (i) -1) + i .
    $$
    If $i = a+1$, then the assertion follows from 
    $$(\kappa \circ (\tilde{\rho}_1 \amalg \tilde{\rho}_2)) (a+1) = \kappa ( am + \tilde{\rho}_2 (1)) = \kappa ( am + \rho_2 (1)) = \rho_2 (1) ( a+1), $$ 
    and $$\widetilde{g_{\rho_1,\rho_2}}  (a+1) = (a+1) ( g_{\rho_1,\rho_2}(a+1) -1 ) + a+1 = (a+1) ( \rho_2 (1) -1) + a+1=\rho_2 (1)  (a+1).$$
\end{proof}

Based on these preliminaries, we obtain the following.

\begin{Lemma} \label{202603142003}
    For $n,m \in \mathds{N}$, we have
    \begin{align*}
        \Delta_n^{(m)} \equiv \sum_{\rho}  E_{\tilde{\rho}} \mod{\mathtt{I}^{\mathsf{pr}}_{\mathcal{C}}(nm,n)}  ,
    \end{align*}
    where $\rho$ runs over maps ${\bf n}\to {\bf m}$.
\end{Lemma}
\begin{proof}
    The statement is trivial for $m=0$ or $n=0$.
    The remaining proof is based on the induction on $n$.
    We first prove the case $n=1$.   
    In this case, the statement is equivalent to saying that
    \begin{align} \label{202606101536}
        \Delta^{(m)} \equiv \sum^{m}_{k=1} E_{m,k},
    \end{align}
    since $\tilde{\rho} = \rho$ for any map $\rho : {\bf 1} \to {\bf m}$.
    This assertion is proved using the induction on $m$ as follows.
    This is true when $m=1$, since $\Delta^{(1)} = \mathrm{id}_1 = E_{1,1}$.
    For $m=2$, this follows from $\Delta \equiv \eta \times \mathrm{id}_1 + \mathrm{id}_1 \times \eta$.
    If $m \geq 3$, then the inductive hypothesis implies that
    \begin{align*}
        &\Delta^{(m)} = (\Delta \times \mathrm{id}_1) \circ \Delta^{(m-1)} \equiv (\Delta \times \mathrm{id}_1) \circ \sum^{m-1}_{k} E_{m-1,k} = (\Delta \times \mathrm{id}_1) \circ E_{m-1,1} + \sum^{m-1}_{k=2} (\Delta \times \mathrm{id}_1) \circ E_{m-1,k} .
    \end{align*}
    By the definition of $E_{m-1,1}$, the first term $(\Delta \times \mathrm{id}_1) \circ E_{m-1,1}$ is equal to $ \Delta \times \eta^{\times m-2} = (\mathrm{id}_2 \times \eta^{\times m-2} ) \circ \Delta$.
    This coincides with the following modulo $\mathtt{I}^{\mathsf{pr}}_{\mathcal{C}} (m,1)$:
    $$(\mathrm{id}_2 \times \eta^{\times m-2} ) \circ (\mathrm{id}_1 \times \eta + \eta \times \mathrm{id}_1) = \mathrm{id}_1 \times \eta^{\times m-1} + \eta \times \mathrm{id}_1 \times \eta^{\times m-2} = E_{m,1} + E_{m,2} .$$
    Furthermore, for $2 \leq k \leq m-1$, we have $(\Delta \times \mathrm{id}_1) \circ E_{m-1,k} = E_{m,k+1}.$
    So, $\sum^{m-1}_{k=2} (\Delta \times \mathrm{id}_1) \circ E_{m-1,k} = \sum^{m}_{k=3} E_{m,k}$.
    These results prove (\ref{202606101536}) by induction on $m$.

    We now prove the claim for the case that $n = n_1+1$, assuming that it is true for $n = n_1$.
    By Lemma \ref{202603141912} and the assumption,
    \begin{align*}
        \Delta^{(m)}_{n_1+1} =& E_{\kappa} \circ ( \Delta^{(m)}_{n_1} \times \Delta^{(m)}_{1}) = \sum E_\kappa \circ ( E_{\tilde{\rho}_1} \times E_{\tilde{\rho}_2}) = \sum E_{\kappa \circ (\tilde{\rho}_1 \amalg \tilde{\rho}_2)} .
    \end{align*}
    where $\kappa = \kappa_{n_1,m}$, and $\rho_1$ and $\rho_2$ run over maps ${\bf n_1} \to {\bf m}$ and ${\bf 1} \to {\bf m}$.
    Recall the result in Lemma \ref{202606191206}.
    Since all the maps $\rho : {\bf n_1 +1} \to {\bf m}$ arise in the form of $g_{\rho_1,\rho_2}$ for unique $\rho_1: {\bf n_1} \to {\bf m}$ and $\rho_2: {\bf 1} \to {\bf m}$, the above sum equals $\sum_{\rho} E_{\tilde{\rho}}$ where $\rho : {\bf n_1+1} \to {\bf m}$ runs over all maps.
\end{proof}

\begin{Lemma} \label{202603161244}
    $\alpha_{\mathcal{C}} \circ \delta_{\mathcal{C}} = \mathrm{id}$.
\end{Lemma}
\begin{proof}
    By the definitions of $\alpha_{\mathcal{C}}$ and $ \delta_{\mathcal{C}}$, it suffices to prove that, for $f = (f_1,\cdots,f_m) \in \mathcal{C}(n,m)$,
    $$
    (f_1,\cdots,f_m)  = \sum_{\rho} \left( \prod^{m}_{k=1} f_k \circ E_{n,\rho^{-1}(k)} \right) \circ E_{\sigma_\rho} \mod{\mathtt{I}^{\mathsf{pr}}_{\mathcal{C}}(m,n)} ,
    $$
    where the sum is taken over maps $\rho : {\bf n} \to {\bf m}$.
    Recall the notation in Definition \ref{202606101643}.
    By Lemma \ref{202606191144}, we obtain
    $$\left( \prod^{m}_{k=1} f_k \circ E_{n,\rho^{-1}(k)} \right) \circ E_{\sigma_\rho} = \prod^{m}_{k=1} f_k \circ \prod^{m}_{k=1} E_{n,\rho^{-1}(k)} \circ E_{\sigma_\rho}= \prod^{m}_{k=1} f_k \circ E_{\coprod^{m}_{k=1} \iota_k}  \circ E_{\sigma_\rho} = \prod^{m}_{k=1} f_k \circ E_{\tilde{\rho}} .$$
    So, the above equality follows from Lemma \ref{202603142003}, since $(\prod^{m}_{k=1} f_k) \circ \Delta_n^{(m)} = (f_1, \cdots, f_m)$.    
\end{proof}

\section{The PROP associated with the primitivity ideal}
\label{202603211037}

Let $\mathcal{C}$ be a Lawvere theory with a zero object.
Recall the eigenmonad $\Phi_\mathcal{C}$ and the right $\mathtt{L}_{\mathfrak{S}}$-module $\Psi_\mathcal{C}$.
The goal of this section is to establish structural properties of the eigenmonad $\Phi_\mathcal{C}$.
We show that the $\mathds{k}$-linear category induced by $\Phi_\mathcal{C}$ is a ($\mathds{k}$-linear) PROP.
To better understand this PROP, we endow $\Psi_\mathcal{C}$ with a natural coaugmented comonoid species structure.
By Lemma \ref{202603211811} and Theorem \ref{202603161027}, there is an embedding
$$\Phi_\mathcal{C} (m,-) \hookrightarrow \mathtt{L}_{\mathcal{C}}/\mathtt{I}^{\mathsf{pr}}_{\mathcal{C}} (m,-) \cong \mu \Psi_\mathcal{C} (m,-) \cong \Psi_\mathcal{C}^{\odot m}.$$
We then describe $\Phi_\mathcal{C}$ in terms of the coaugmented comonoid species structure on $\Psi_\mathcal{C}$.

\subsection{The PROP structure and operadicity}

In this section, we show that the Lawvere theory structure on $\mathcal{C}$ induces a natural PROP structure on the $\mathds{k}$-linear category corresponding to $\Phi_\mathcal{C}$.

\begin{Lemma} \label{202603170852}
    The idealizer monad $\mathrm{Iz} (\mathtt{I}^{\mathsf{pr}}_{\mathcal{C}})\subset \mathtt{L}_{\mathcal{C}}$ is closed under the product of the Lawvere theory $\mathcal{C}$.
\end{Lemma}
\begin{proof}
    Let $m_1,m_2,n_1,n_2 \in \mathds{N}$ and $u_i \in (\mathrm{Iz} (\mathtt{I}^{\mathsf{pr}}_{\mathcal{C}})) (m_i,n_i),~ i \in {\bf 2}$.
    By definition, we have $\bar{\Delta}^{m_i,k} \circ u_i \in \mathtt{I}^{\mathsf{pr}}_{\mathcal{C}} (m_i+1,n_i)$ for any $k \in {\bf m_i}$.
    It suffices to show that, $\bar{\Delta}^{m_1+m_2,k}\circ (u_1 \times u_2)$ for $k \in {\bf m_1 + m_2}$.
    If $1 \leq k \leq m_1$, then $\bar{\Delta}^{m_1+m_2,k}\circ (u_1 \times u_2) = (\bar{\Delta}^{m_1,k} \circ u_1) \times u_2$, which lies in the idealizer monad by Lemma \ref{202601051527}.
    Similarly, one may treat the case that $m_1 < k \leq m_1+m_2$.
\end{proof}

\begin{theorem} \label{202604071720}
     Let $\tilde{\Phi}_{\mathcal{C}}$ be the $\mathds{k}$-linear category with object set $\mathds{N}$ induced by the monad $\Phi_{\mathcal{C}}$.
    This inherits a $\mathds{k}$-linear PROP structure from the product in the Lawvere theory $\mathcal{C}$.
\end{theorem}
\begin{proof}
    For $u_i \in (\mathrm{Iz} (\mathtt{I}^{\mathsf{pr}}_{\mathcal{C}})) (m_i,n_i),~ i \in \{1,2\}$, we define the monoidal operation $[u_1] \boxtimes [u_2] {:=} [u_1 \times u_2] \in \Phi_{\mathcal{C}}(m_1+m_2 , n_1+n_2)$ where $[-]$ denotes the equivalence class modulo $\mathtt{I}^{\mathsf{pr}}_{\mathcal{C}}$.
    This is well-defined by Lemma \ref{202603170852}.
    Furthermore, Lemma \ref{202603211051} proves that the symmetry in the Lawvere theory $\mathcal{C}$ induces a symmetry with respect to this monoidal operation.
\end{proof}

Let $\mathcal{P}$ be a PROP, and define the underlying operad $\mathfrak{O}$ by $\mathfrak{O}(n) := \mathcal{P}(n,1)$.
The operadic composition is induced by the composition and the monoidal operation in the PROP $\mathcal{P}$.
Then there is a PROP $\mathsf{Cat}\mathfrak{O}$ canonically associated to $\mathfrak{O}$, together with a natural morphism of PROPs $\mathsf{Cat}\mathfrak{O} \to \mathcal{P}$.
We say that $\mathcal{P}$ is {\it operadic} if this morphism gives an isomorphism.

\begin{Defn} \label{202603241701}
    Let $\mathfrak{O}_{\mathcal{C}}$ denote the underlying operad of the PROP $\tilde{\Phi}_{\mathcal{C}}$.
    We define a natural monad map $\mola_{\mathcal{C}}:\mu \mathfrak{O}_{\mathcal{C}} \to \Phi_{\mathcal{C}}$ by $\mathds{k}$-linearly extending the map
    $$
    \bigotimes^{m}_{k=1} \mathfrak{O}_{\mathcal{C}} (n_k) \otimes_{\mathfrak{S}_{n_1} \times \cdots \times \mathfrak{S}_{n_m}} \mathds{k}[\mathfrak{S}_n] \longrightarrow  \Phi_{\mathcal{C}} (m,n),
    $$
    given by
    $$
    \bigotimes^m_{k=1} [f_k ] \otimes \sigma \mapsto [\left(\prod^m_{k=1} f_k \right)\circ E_{\sigma}] ,
    $$
    where $(n_1,\cdots ,n_m)$ is an $m$-partition of $n$.
\end{Defn}

\begin{Lemma} \label{202604051240}
    The PROP $\tilde{\Phi}_{\mathcal{C}}$ is operadic if and only if the monad map $\mola_\mathcal{C}$ is an isomorphism.
\end{Lemma}
\begin{proof}
    Under the equivalence between the category of monads in $\mathsf{Mat}_{\mathds{k}}$ on $\mathds{N}$ and the category of $\mathds{k}$-linear categories with object set $\mathds{N}$, the functor $\mathsf{Cat}\mathfrak{O}_{\mathcal{C}} \to \tilde{\Phi}_\mathcal{C}$ corresponds to the map $\mola_\mathcal{C}$ (see also Remark \ref{202604051201}).
    This proves the claim.
\end{proof}

\subsection{Comonoid species associated to Lawvere theories}
\label{202603201642}

In this section, we introduce a natural coaugmented comonoid structure on the $\mathds{k}$-linear species $\Psi_\mathcal{C}$.
We identify $\mathfrak{O}_\mathcal{C}$ with primitive elements of $\Psi_\mathcal{C}$.
More generally, we characterize the whole eigenmonad $\Phi_\mathcal{C}$ from the coaugmented comonoid $\Psi_\mathcal{C}$,
and derive some sufficient condition under which the map $\mola_\mathcal{C} : \mu\mathfrak{O}_{\mathcal{C}} \to \Phi_{\mathcal{C}}$ is an isomorphism.

Recall that $\mathtt{L}_\mathcal{C}/\mathtt{I}^{\mathsf{pr}}_{\mathcal{C}}$ is an $(\mathtt{L}_\mathcal{C}, \Phi_{\mathcal{C}})$-bimodule.
The underlying $(\mathtt{L}_\mathcal{C},\mathtt{L}_\mathfrak{S})$-bimodule of $\mathtt{L}_\mathcal{C}/\mathtt{I}^{\mathsf{pr}}_{\mathcal{C}}$ induces a functor $\mathcal{C} \to \mathsf{Mod}\mbox{-}\mathtt{L}_{\mathfrak{S}}$ which assigns the right $\mathtt{L}_{\mathfrak{S}}$-module $\mathtt{L}_\mathcal{C}/\mathtt{I}^{\mathsf{pr}}_{\mathcal{C}} (m,-)$ to $m \in \mathds{N}$.
Via the equivalence of right $\mathtt{L}_{\mathfrak{S}}$-modules and $\mathds{k}$-linear species, this yields a functor $\mathtt{Q} : \mathcal{C} \to\mathsf{Sp}_{\mathds{k}}$ such that
$$\mathtt{Q}(m) = \mathtt{L}_\mathcal{C}/\mathtt{I}^{\mathsf{pr}}_{\mathcal{C}} (m,-) .$$

\begin{prop} \label{202603201557}
    The functor $\mathtt{Q}$ enhances to a symmetric (strong) monoidal functor $(\mathcal{C} , \times, 0) \to (\mathsf{Sp}_{\mathds{k}}, \odot, \mathds{k})$ with the structure morphisms $\phi$ and $\xi_{m,l}$ defined below.
\end{prop}
The structure morphisms are described as follows.
By the results in Example \ref{202604041521}, there is an isomorphism $\phi : \mathds{k} \to \mathtt{Q} (0)$ of $\mathds{k}$-linear species.
For $m,l \in \mathds{N}$, define $$\xi_{m,l} : \bigoplus_{(n_1,n_2)} \left( \mathtt{L}_\mathcal{C}/\mathtt{I}^{\mathsf{pr}}_{\mathcal{C}} (m,n_1) \otimes \mathtt{L}_\mathcal{C}/\mathtt{I}^{\mathsf{pr}}_{\mathcal{C}} (l,n_2) \right) \otimes_{\mathfrak{S}_{n_1}\times\mathfrak{S}_{n_2}} \mathds{k}[\mathfrak{S}_n] \to \mathtt{L}_\mathcal{C}/\mathtt{I}^{\mathsf{pr}}_{\mathcal{C}} (m+l,n)$$ by $\xi_{m,l}  ([ f_1] \otimes [f_2] \otimes \sigma) {:=} [(f_1 \times f_2) \circ E_\sigma]$ where the direct sum is taken over 2-partitions $(n_1,n_2)$ of $n \in \mathds{N}$.
This yields a right $\mathtt{L}_{\mathfrak{S}}$-module map $\xi_{m,l} : \mathtt{Q} (m) \odot \mathtt{Q} (l) \to \mathtt{Q} (m+l)$, which naturally extends to a map between the induced $\mathds{k}$-linear species.
Lemmas \ref{202601051527} and \ref{202509071205} imply that this is well-defined.

The proof of Proposition \ref{202603201557} is deferred to Section \ref{202604041854}.
Note that $\mathtt{Q}(1) = \Psi_\mathcal{C}$.
By this proposition, the morphisms $\Delta \in \mathcal{C}(1,2)$, $\epsilon \in \mathcal{C}  (1,0)$ and $\eta \in \mathcal{C}(0,1)$ induce the following maps, where we slightly abuse notation:
\begin{align*}
    \Delta: \Psi_\mathcal{C} \to \Psi_\mathcal{C} \odot \Psi_\mathcal{C}, \quad \epsilon: \Psi_\mathcal{C} \to \mathds{k}, \quad \eta: \mathds{k} \to \Psi_\mathcal{C}.
\end{align*}
Furthermore, if $\mathcal{C}$ satisfies the condition (M), then the multiplication $\nabla \in \mathcal{C} (2,1)$ induces a map $\nabla : \Psi_\mathcal{C} \odot \Psi_\mathcal{C} \to \Psi_\mathcal{C}$.

\begin{Corollary} \label{202603241657}
    \begin{enumerate}
        \item $(\Psi_\mathcal{C}, \Delta,\epsilon,\eta)$ is a connected coaugmented cocommutative comonoid species.
        \item If $\mathcal{C}$ satisfies (M), then $(\Psi_\mathcal{C}, \Delta, \nabla, \epsilon,\eta)$ is further endowed with a bimonoid species structure.  
        Furthermore, if $1 \in \mathcal{C}$ is a Hopf monoid or a commutative bimonoid, then $(\Psi_\mathcal{C}, \Delta, \nabla, \epsilon,\eta)$ is a Hopf monoid species or a commutative bimonoid species, respectively.
    \end{enumerate}
\end{Corollary}
\begin{proof}
    Since $\mathcal{C}$ is a Lawvere theory, the object $1 \in \mathrm{Obj} (\mathcal{C})$ is naturally a coaugmented cocommutative comonoid in $\mathcal{C}$.
    Hence, the symmetric monoidal functor in Proposition \ref{202603201557} sends this to a coaugmented cocommutative comonoid species.
    Furthermore, by Example \ref{202603221114}, this is connected.
    (2) is proved analogously.
\end{proof}

We now present an alternative description of the eigenmonad $\Phi_\mathcal{C}$ in terms of the coaugmented comonoid $\Psi_\mathcal{C}$.
Recall $\mathrm{Pr}^m$ from Definition \ref{202603241659}.
\begin{theorem} \label{202603221137}
    Let $m \in \mathds{N}$.
    The map $\Phi_\mathcal{C} (m,-) \hookrightarrow \Psi_\mathcal{C}^{\odot m}$ induces an isomorphism
    $$
    \mathrm{Pr}^m(\Psi_\mathcal{C}) \cong \Phi_\mathcal{C} (m,-) .
    $$
    In particular, $\mathrm{Pr}(\Psi_\mathcal{C}) \cong \mathfrak{O}_{\mathcal{C}}$.
\end{theorem}
\begin{proof}
    The case $m=0$ is trivial.
    Assume that $m \in \mathds{N}^\ast$.
    By (2) of Proposition \ref{202401241330}, $\Phi_\mathcal{C} (m,n)$ is isomorphic to
    $$
    \{ v \in (\mathtt{L}_\mathcal{C}/\mathtt{I}^{\mathsf{pr}}_\mathcal{C}) (m,n) \mid \overline{\Delta}^{m,k} \rhd v = 0 , ~ k \in {\bf m} \} .
    $$
    Hence, by the definition of the coaugmented comonoid $\Psi_\mathcal{C}$, we obtain
    $$
    \Phi_\mathcal{C} (m,-) \cong \bigcap^{m}_{k=1} \mathrm{Ker} \left(\mathrm{id}_{\Psi_\mathcal{C}}^{\odot k-1} \odot \bar{\Delta} \odot \mathrm{id}_{\Psi_\mathcal{C}}^{\odot m-k}: \Psi_\mathcal{C}^{\odot m} \to \Psi_\mathcal{C}^{\odot m+1} \right) .
    $$
    This intersection coincides with $\mathrm{Pr}^m(\Psi_\mathcal{C})$ by the definition of the latter.
    If we apply $m=1$ to the above result, then we also obtain $\mathrm{Pr}(\Psi_\mathcal{C}) \cong \Phi_\mathcal{C} (1,-) = \mathfrak{O}_{\mathcal{C}}$.
\end{proof}

\begin{Corollary} \label{202603201735}
    For $n,m\in\mathds{N}$, the map $\mola_\mathcal{C} : \mu\mathfrak{O}_\mathcal{C} (m,n) \to \Phi_\mathcal{C} (m,n)$ is an isomorphism if and only if the map $\mathrm{Pr} (\Psi_\mathcal{C})^{\odot m} (n) \to \mathrm{Pr}^{m} (\Psi_\mathcal{C}) (n)$ is an isomorphism.
    Therefore, the PROP $\tilde{\Phi}_\mathcal{C}$ is operadic if and only if the map $\mathrm{Pr} (\Psi_\mathcal{C})^{\odot m} \to \mathrm{Pr}^{m} (\Psi_\mathcal{C})$ is an isomorphism for any $m \in \mathds{N}$.
\end{Corollary}
\begin{proof}
    By Lemma \ref{202604051240}, it suffices to prove the first statement.
    By Theorem \ref{202603221137}, $\Phi_\mathcal{C} (m,-) \cong \mathrm{Pr}^m (\Psi_\mathcal{C})$.
    So, $\mathfrak{O}_\mathcal{C} = \Phi_\mathcal{C} (1,-) \cong  \mathrm{Pr} (\Psi_\mathcal{C})$.
    In particular, $\mu\mathfrak{O}_\mathcal{C} (m,-) \cong \mathfrak{O}_\mathcal{C}^{\odot m} =\mathrm{Pr} ( \Psi_\mathcal{C})^{\odot m}$.
    These isomorphisms identify the map $\mola_{\mathcal{C}}$ with the map $\mathrm{Pr} (\Psi_\mathcal{C})^{\odot m} \to \mathrm{Pr}^{m} (\Psi_\mathcal{C})$.
    This proves the desired claim.
\end{proof}

\begin{Corollary} \label{202603221141}
    For $n\in\mathds{N}$, the map $\mola_\mathcal{C}$ gives a $\mathds{k}$-algebra isomorphism
    $\mu\mathfrak{O}_\mathcal{C} (n,n) \to \Phi_\mathcal{C} (n,n)$.
    In particular, we have $\Phi_\mathcal{C} (n,n) \cong \mathfrak{O}_\mathcal{C} (1) \thicksim \mathfrak{S}_n$.
\end{Corollary}
\begin{proof}   
    Applying Lemma \ref{202603221113} to $C = \Psi_\mathcal{C}$, we obtain $\mathrm{Pr}^n (\Psi_\mathcal{C}) (n) \cong \mathrm{Pr} (\Psi_\mathcal{C})^{\odot n} (n)$.
    Thus, Corollary \ref{202603201735} leads to $\mu\mathfrak{O}_\mathcal{C} (n,n) \cong \Phi_\mathcal{C} (n,n)$.
    This is a $\mathds{k}$-algebra isomorphism since $\mola_\mathcal{C}$ is a monad map.
    The last result follows from the fact that $\mu\mathfrak{O}_\mathcal{C} (n,n) = \mathfrak{O}_\mathcal{C} (1) \thicksim \mathfrak{S}_n$.
\end{proof}

\begin{Corollary} \label{202604051251}
    If the ground ring $\mathds{k}$ is a field, then the monad map $\mola_\mathcal{C} : \mu\mathfrak{O}_\mathcal{C} \to \Phi_\mathcal{C}$ is an isomorphism.
    In particular, the PROP $\tilde{\Phi}_\mathcal{C}$ is operadic.
\end{Corollary}
\begin{proof}
    The result follows from Lemma \ref{202604011630} and Corollary \ref{202603201735}.
\end{proof}

\subsection{The case $\mathcal{C}=\G^{\mathsf{o}}$}
\label{202512051104}

In this section, we prove that the eigenmonad $\Phi_{\G^{\mathsf{o}}}$ is isomorphic to the monad $\mu\mathfrak{Lie}$ associated with the Lie operad.
The proof is based on the previous general results.
As a consequence, we show that the eigenmonad adjunction associated with $(\mathtt{L}_{\G^{\mathsf{o}}}, \mathtt{I}^{\mathsf{pr}}_{\G^{\mathsf{o}}})$ recovers Powell's adjunction \cite{powell2024analytic}.

The object $1 \in \G^{\mathsf{o}}$ is a cocommutative Hopf monoid, so by Corollary \ref{202603241657}, $\Psi_{\G^{\mathsf{o}}}$ is a cocommutative Hopf monoid species.
Recall from Example \ref{202604101536} that $\mathsf{Lin}$ is also a cocommutative Hopf monoid species.
\begin{Lemma} \label{202604101649}
    The isomorphism $\mathcal{F} : \mathsf{Lin} \stackrel{\cong}{\to} \Psi_{\G^{\mathsf{o}}}$ in Theorem \ref{202509091931} preserves the Hopf monoid structures.
\end{Lemma}
\begin{proof}
    It suffices to show that $\mathcal{F}$ preserves the bimonoid structures.
    It is obvious that $\mathcal{F}$ preserves the counit, the unit and the multiplication.
    To show that the comultiplication is preserved, we compute the comultiplication of $\mathcal{F} ( \theta_1 \cdots \theta_n) = [x_1\cdots x_n] \in \Psi_{\G^{\mathsf{o}}} (n)$ which is the generator of the right $\mathfrak{S}_n$-module $\Psi_{\G^{\mathsf{o}}} (n)$.
    Note that the diagonal map $\Delta \in \G^{\mathsf{o}}(1,2)$ is given by $(x_1,x_1) \in \mathsf{F}_1^{\times 2} \cong \G^{\mathsf{o}}(1,2)$.
    An argument similar to Lemma \ref{202603171154} gives
    \begin{align*}
        (x_1,x_1) \circ (x_1\cdots x_n) = (x_1\cdots x_n,x_1\cdots x_n) \equiv 
        \sum (x_{i_1}\cdots x_{i_l}, x_{i_{l+1}} \cdots x_{i_n}) \mod{\mathtt{I}^{\mathsf{pr}}_{\G^{\mathsf{o}}}(2,n)}
    \end{align*}
    where the sum is taken over all partitions $\{i_1 < \cdots < i_l \} \amalg \{i_{l+1} < \cdots < i_n \} = {\bf n}$.
    By the definition of the comultiplication on $\mathsf{Lin}$, this shows that $\mathcal{F}$ preserves the comlutiplication.
\end{proof}

Recall $\tilde{\mu}\mathsf{Lin}$ from Notation \ref{202604101645}.
\begin{theorem} \label{202603121700}
    We have a monad isomorphism $\Phi_{\G^{\mathsf{o}}} \cong\mu \mathfrak{Lie}$.
    Hence, $\tilde{\mu}\mathsf{Lin}$ is a $(\mathtt{L}_{\G^{\mathsf{o}}} , \mu\mathfrak{Lie})$-bimodule.
    The eigenmonad adjunction associated with $(\mathtt{L}_{\G^{\mathsf{o}}}, \mathtt{I}^{\mathsf{pr}}_{\G^{\mathsf{o}}})$ induces the following adjunction:
    $$
    \begin{tikzcd}
              \tilde{\mu}\mathsf{Lin} \otimes_{\mu\mathfrak{Lie}} (-) : \mu\mathfrak{Lie}\mbox{-}\mathsf{Mod} \arrow[r, shift right=1ex, ""{name=G}] & \mathtt{L}_{\G^{\mathsf{o}}}\mbox{-}\mathsf{Mod} : \mathrm{Hom}_{\mathtt{L}_{\G^{\mathsf{o}}}} ( \tilde{\mu}\mathsf{Lin} , - ) \arrow[l, shift right=1ex, ""{name=F}]
            \arrow[phantom, from=G, to=F, , "\scriptscriptstyle\boldsymbol{\top}"] .
    \end{tikzcd}
    $$
\end{theorem}

\begin{remark} 
    The $(\mathtt{L}_{\G^{\mathsf{o}}} , \mu\mathfrak{Lie})$-bimodule $\tilde{\mu}\mathsf{Lin}$ and the above adjunction recovers the $\mathds{k}$-bilinear functor ${}_{\triangle} \mathsf{Cat}{\mathfrak{Ass}}^{u}$ and Powell's adjunction \cite{powell2024analytic,kim2024analytic}, via the identification $\mathsf{Lin} = \mathfrak{Ass}^{u}$.
\end{remark}

\begin{proof}[Proof of Theorem \ref{202603121700}]
    We begin by giving an operad isomorphism $\mathfrak{Lie} \cong  \mathfrak{O}_{\G^{\mathsf{o}}}$.
    Proposition \ref{202604021602} and Lemma \ref{202604011516} imply that $ \mathrm{Pr} ( \mathsf{Lin}) = \mathfrak{Lie}$.
    Hence, Lemma \ref{202604101649} yields a right $\mathtt{L}_{\mathfrak{S}}$-module isomorphism $\mathfrak{Lie} \cong \mathrm{Pr} ( \Psi_{\G^{\mathsf{o}}}) = \mathfrak{O}_{\G^{\mathsf{o}}}$.
    We denote this map by $\xi : \mathfrak{Lie} \to \mathfrak{O}_{\G^{\mathsf{o}}}$, and show that this is an operad map.
    Let $\iota : \mathsf{Lin} \to \mathtt{L}_{\G^{\mathsf{o}}} (1,-)$ be the map which assigns $x_{1}\cdots x_{n} \in \mathsf{F}_n \cong \G^{\mathsf{o}} (n,1)$ to $\theta_1\cdots\theta_n \in \mathsf{Lin}(n)$.
    By the definition of the map $\xi$, we have $$\xi (f) = \iota ( f) \mod \mathtt{I}^{\mathsf{pr}}_{\G^{\mathsf{o}}} (1,n)$$ for $n \in \mathds{N}$ and $f \in \mathfrak{Lie} (n)$.
    Since $\iota : \mathsf{Lin} \to \mathtt{L}_{\G^{\mathsf{o}}} (1,-)$ is an operad map where $\mathtt{L}_{\G^{\mathsf{o}}} (1,-)$ is induced by the PROP structure of $\G^{\mathsf{o}}$, the map $\xi$ preserves the operad structures.

    By Proposition \ref{202604021602} and Corollary \ref{202603201735}, the map $\mola_{\G^{\mathsf{o}}} : \mu \mathfrak{O}_{\G^{\mathsf{o}}} \to \Phi_{\G^{\mathsf{o}}}$ is an isomorphism.
    Thus, the above result proves $\Phi_{\G^{\mathsf{o}}} \cong\mu \mathfrak{Lie}$.
    The assertion for the adjunction is immediate from Corollary \ref{202512081855}.
\end{proof}

\subsection{Proof of Proposition \ref{202603201557}}
\label{202604041854}

In this section, we prove Proposition \ref{202603201557}.

\begin{Lemma} 
    The maps $\xi_{m,l} : \mathtt{Q} (m) \odot \mathtt{Q} (l) \to \mathtt{Q} (m+l)$ form a natural transformation.
\end{Lemma}
\begin{proof}
    Let $g \in \mathcal{C}(m,m^\prime)$ and $h \in \mathcal{C} (l,l^\prime)$.
    It suffices to show that the following diagram commutes:
    $$
    \begin{tikzcd}
        \mathtt{Q} (m) \odot \mathtt{Q} (l) \ar[r, "\xi_{m,l}"] \ar[d, "g_\ast \odot h_\ast"] & \mathtt{Q} (m+l) \ar[d, "(g\times h)_\ast"] \\
        \mathtt{Q} (m^\prime) \odot \mathtt{Q} (l^\prime) \ar[r, "\xi_{m^\prime,l^\prime}"] & \mathtt{Q} (m^\prime+l^\prime)
    \end{tikzcd}
    $$
    Note that the right $\mathfrak{S}_n$-module $(\mathtt{Q} (m) \odot \mathtt{Q} (l))(n)$ is generated by $[f_1] \otimes [f_2] \in (\mathtt{Q} (m))(n_1) \otimes (\mathtt{Q} (l))(n_2)$ for partitions $(n_1,n_2)$ of $n$.
    For such generators, we have 
    \begin{align*}
        &(g\times h)_\ast ( \xi_{m,l} ([ f_1] \otimes [f_2] )) = (g\times h)_\ast ([f_1 \times f_2]) = [(g\times h) \circ (f_1\times f_2) ]  , \\
        =& [g\circ f_1 \times h\circ f_2] = \xi_{m^\prime,l^\prime} ( [g\circ f_1]\otimes [h\circ f_2]) = \xi_{m^\prime,l^\prime} ( (g_\ast \odot h_\ast) ( [ f_1] \otimes [f_2] )) .
    \end{align*}
\end{proof}

\begin{proof}[Proof of Proposition \ref{202603201557}]
    We begin by giving some preliminaries.
    By Theorem \ref{202603161027}, there is an isomorphism $\alpha_m : \Psi_{\mathcal{C}}^{\odot m} \to \mathtt{Q} (m)$ of species preserving the right $\mathfrak{S}_m$-actions.
    Moreover, $\alpha_0 = \phi$.
    Then the following diagram commutes:
    $$
    \begin{tikzcd}
        \Psi_{\mathcal{C}}^{\odot m} \odot \Psi_{\mathcal{C}}^{\odot l} \ar[r, "\cong"] \ar[d, "\alpha_m \odot \alpha_l"] & \Psi_{\mathcal{C}}^{\odot m+l} \ar[d, "\alpha_{m+l}"] \\
        \mathtt{Q} (m) \odot \mathtt{Q} (l) \ar[r, "\xi_{m,l}"] & \mathtt{Q} (m+l)
    \end{tikzcd}
    $$
    Here, the upper isomorphism is induced by the associator of the monoidal category $(\mathsf{Sp}_{\mathds{k}},\odot, \mathds{k})$.
    To prove this, consider $n \in \mathds{N}$, $u = \bigotimes^m_{i=1} [f_i] \otimes \sigma_1 \in \Psi^{\odot m} (n_1)$ and $v =\bigotimes^l_{j=1} [g_j] \otimes \sigma_2 \in \Psi^{\odot l} (n_2)$ with $n_1 +n_2=n$.
    We regard $u \otimes v$ an element of $ (\Psi_{\mathcal{C}}^{\odot m} \odot \Psi_{\mathcal{C}}^{\odot l}) (n)$ and compare the two compositions on this element.
    Then $$\xi_{m,l} ( ( \alpha_m \odot \alpha_l ) (u \otimes v ) )= \xi_{m,l} ( [\prod^m_{i=1} f_i] \otimes [\prod^l_{j=1} g_j] ) = [\prod^m_{i=1} f_i \times \prod^l_{j=1} g_j].$$
    Moreover, $\alpha_{m+l} ( \bigotimes^m_{i=1} [f_i] \otimes \bigotimes^l_{j=1} [g_j]) = [\prod^m_{i=1} f_i \times \prod^l_{j=1} g_j]$.
    Since the right $\mathfrak{S}_n$-module $ (\Psi_{\mathcal{C}}^{\odot m} \odot \Psi_{\mathcal{C}}^{\odot l}) (n)$ is generated by elements of the form $u \otimes v$, the two compositions coincide. Hence, the diagram commutes.

    The commutativity implies that the functor $\mathtt{Q}$ together with $\phi$ and $\xi_{m,l}$ enhances to a strong monoidal functor.
    In fact, we can prove that the functor $\mathtt{Q}$ with $\xi_{m,l}$'s preserve the associator using the above commutative diagram and the associator axiom of $(\mathsf{Sp}_{\mathds{k}},\odot, \mathds{k})$.
    Likewise, the functor $\mathtt{Q}$ with $\xi_{m,l}$'s and $\phi$ preserve the unitor, since $\alpha_0 = \phi$.
    Furthermore, $\xi_{m,l}$'s form a natural isomorphism since $\alpha_m$'s are isomorphisms.
    
    We now show that the functor $\mathtt{Q}$ with $\xi_{m,l}$ preserves the symmetry.
    For $m,l\in\mathds{N}$, let $\sigma_{m,l} \in \mathfrak{S}_{m+l}$ denote the bijection such that $\sigma (k) = k+l$ if $1 \leq k \leq m$ and $\sigma (k) = k-m$ for $m+1 \leq k \leq m+l$.
    Then $E_{\sigma_{m,l}} \in \mathcal{C} (m+l, l+m)$ is the symmetry in $\mathcal{C}$.
    We shall show the following commutative diagram where $s$ is the symmetry of linear species:
    $$
    \begin{tikzcd}
        \mathtt{Q}(m) \odot \mathtt{Q}(l) \ar[r, "s"] \ar[d, "\xi_{m\mbox{,}l}"] & \mathtt{Q}(l) \odot \mathtt{Q} (m) \ar[d, "\xi_{l\mbox{,}m}"] \\
        \mathtt{Q}(m+l) \ar[r, "(E_{\sigma_{m,l}})_\ast"] & \mathtt{Q}(l+m)
    \end{tikzcd}
    $$
    Let $[f] \otimes [g] \otimes e \in (\mathtt{Q}(m)(n_1) \otimes \mathtt{Q}(l)(n_2)) \otimes_{\mathfrak{S}_{n_1}\times\mathfrak{S}_{n_2}} \mathds{k}[\mathfrak{S}_n]$.
    We have
    \begin{align*}
        (E_{\sigma_{m,l}})_\ast ( \xi_{m,l} ( [f] \otimes [g] \otimes e )) =& (E_{\sigma_{m,l}})_\ast ([f \times g]) = [E_{\sigma_{m,l}} \circ (f \times g)]. 
    \end{align*}
    By the naturality of products in $\mathcal{C}$, this equals $$[(g\times f) \circ E_{\sigma_{n_1,n_2}}] = \xi_{l,m} ( [g] \otimes [f] \otimes \sigma_{n_1,n_2}) = \xi_{l,m} ( s ( [f] \otimes [g] \otimes e)).$$
    Since the right $\mathfrak{S}_n$-module $(\mathtt{Q}(m) \odot \mathtt{Q}(l))(n)$ is generated by elements of the form $[f] \otimes [g] \otimes e $ as above, the diagram commutes.
\end{proof}

\section{Generation of the quotient $\mathtt{L}_{\mathcal{C}}/\mathtt{I}^{\mathsf{pr}}_{\mathcal{C}}$ by $\mathfrak{O}_{\mathcal{C}}(1)$}
\label{202509201216}

Recall the conditions on Lawvere theories introduced in Section \ref{202606181932}.
Let $\mathcal{C}$ be a Lawvere theory satisfying (ZM).
In this section, we introduce a natural $(\mathtt{L}_{\mathbf{W}^{\mathsf{o}}}, \mathtt{L}_\mathfrak{S})$-bimodule map $\E_\mathcal{C} : \mu\mathsf{Lin} \otimes_{\mathrm{H}} \mathfrak{O}_{\mathcal{C}}(1)^{\otimes} \to \mathtt{L}_{\mathcal{C}}/\mathtt{I}^{\mathsf{pr}}_{\mathcal{C}}$, where the generalized Hadamard product $\otimes_{\mathrm{H}}$ is introduced in this section.
We investigate several useful properties of this map and, as a consequence, obtain a criterion for it to be an epimorphism.
As a biproduct, we also construct a bimonoid map $\mathsf{Lin} \otimes_{\mathrm{H}} \mathfrak{O}_{\mathcal{C}}(1)^{\otimes} \to \Psi_\mathcal{C}$ and obtain a similar criterion.

\subsection{Hadamard product}

In this section, we introduce a generalization of the Hadamard product of symmetric modules.

\begin{Defn} \label{202603301023}
    Let $\mathcal{X},\mathcal{Y}$ be sets.
    Let $\mathtt{M}$ and $\mathtt{N}$ be right $\mathtt{L}_{\mathfrak{S}}$-modules with codomains $\mathcal{X}$ and $\mathcal{Y}$ respectively.
    We define the {\it generalized Hadamard} product $\mathtt{M} \otimes_{\mathrm{H}} \mathtt{N}$ to be the right $\mathtt{L}_{\mathfrak{S}}$-module with codomain $\mathcal{X}\times\mathcal{Y}$ such that $(\mathtt{M} \otimes_{\mathrm{H}} \mathtt{N}) ((X,Y),n) {:=} \mathtt{M}(X,n) \otimes \mathtt{N}(Y,n)$ equipped with the diagonal right $\mathfrak{S}_n$-action.
\end{Defn}

\begin{Lemma} \label{202606121427}
    Let $\mathtt{T}$ be a monad on $\mathcal{X}$ and $\mathtt{M}$ be a $(\mathtt{T}, \mathtt{L}_\mathfrak{S})$-bimodule.
    Let $\mathtt{N}$ be a right $\mathtt{L}_\mathfrak{S}$-module.
    Then $\mathtt{M} \otimes_{\mathrm{H}} \mathtt{N}$ is a $(\mathtt{T},\mathtt{L}_\mathfrak{S})$-bimodule.
\end{Lemma}
\begin{proof}
    For $f \in \mathtt{T} (Y,X)$, the action of $f$ on $\mathtt{M}$ induces a map $\mathtt{M} (X, n) \otimes \mathtt{N} (n) \to \mathtt{M} (Y, n) \otimes \mathtt{N} (n)$ which is compatible with the diagonal $\mathfrak{S}_n$-action.
\end{proof}

Recall the right $\mathtt{L}_{\mathfrak{S}}$-module $V^{\otimes}$ from Example \ref{202604062227}.

\begin{Lemma} \label{202603240328}
    Let $V$ be a $\mathds{k}$-module.
    The endofunctor on $\mathsf{Mod}\mbox{-}\mathtt{L}_{\mathfrak{S}}$ which assigns $\mathtt{M} \otimes_{\mathrm{H}} V^{\otimes}$ to $\mathtt{M}$ is a {\it symmetric} monoidal functor.
\end{Lemma}
\begin{proof}
    Via the equivalence in (\ref{202603251315}), the functor $(-)\otimes_{\mathrm{H}} V^{\otimes}$ can be regarded as an endofunctor on $\mathds{k}$-linear species.
    Note that $V^{\otimes} (X_1 \amalg X_2) \cong V^{\otimes}(X_1) \otimes V^{\otimes} (X_2)$ for finite sets $X_1,X_2$.
    For the unit $\mathds{k}$, we have $\mathds{k} \otimes_{\mathrm{H}} V^{\otimes} \cong \mathds{k}$.
    We also have a natural isomorphism $(\mathtt{M} \odot \mathtt{N} ) \otimes_{\mathrm{H}} V^{\otimes} \cong (\mathtt{M} \otimes_{\mathrm{H}} V^{\otimes} ) \odot (\mathtt{N} \otimes_{\mathrm{H}} V^{\otimes} )$ which is compatible with the associator, the unitors and the symmetry.
    Indeed, this is given by the composition of following switching isomorphisms, where $X$ is a finite set $X$:
    \begin{align*}
        &\bigoplus_{X=X_1\amalg X_2} \mathtt{M} (X_1) \otimes \mathtt{N} (X_2)  \otimes V^{\otimes} (X) \cong \bigoplus_{X=X_1\amalg X_2} \mathtt{M} (X_1) \otimes \mathtt{N} (X_2)  \otimes V^{\otimes} (X_1) \otimes V^{\otimes} (X_2) , \\
        \cong& \bigoplus_{X=X_1\amalg X_2} \mathtt{M} (X_1) \otimes V^{\otimes} (X_1) \otimes \mathtt{N} (X_2)  \otimes V^{\otimes} (X_2) \cong \bigoplus_{X=X_1\amalg X_2} (\mathtt{M} \otimes_{\mathrm{H}} V^{\otimes} ) (X_1) \otimes (\mathtt{N} \otimes_{\mathrm{H}} V^{\otimes} ) (X_2) .
    \end{align*}
\end{proof}

\subsection{Construction of the map $\E_\mathcal{C}$}
\label{202603241131}

Let $\mathcal{C}$ be a Lawvere theory with a zero object.
In this section, we construct the map $\E_\mathcal{C} : \mu\mathsf{Lin} \otimes_{\mathrm{H}} \mathfrak{O}_{\mathcal{C}}(1)^{\otimes} \to \mathtt{L}_{\mathcal{C}}/\mathtt{I}^{\mathsf{pr}}_{\mathcal{C}}$.

\begin{Defn}
    We define a right $\mathtt{L}_{\mathfrak{S}}$-module $\Phi_\mathcal{C}^{\Delta}$ by $\Phi_\mathcal{C}^{\Delta} (n) = \Phi_\mathcal{C} (n,n), ~n \in \mathds{N}$ where we regard $\Phi_\mathcal{C} (n,n)$ as a right $\mathfrak{S}_n$-module using the conjugate action.
\end{Defn}

By Lemma \ref{202606121427}, the underlying $(\mathtt{L}_\mathcal{C},\mathtt{L}_\mathfrak{S})$-bimodule of $\mathtt{L}_\mathcal{C}/\mathtt{I}^{\mathsf{pr}}_{\mathcal{C}}$ induces an $(\mathtt{L}_\mathcal{C}, \mathtt{L}_\mathfrak{S})$-bimodule $\mathtt{L}_\mathcal{C}/\mathtt{I}^{\mathsf{pr}}_{\mathcal{C}} \otimes_{\mathrm{H}} \Phi_\mathcal{C}^{\Delta}$.

\begin{Lemma} 
    The right $\Phi_\mathcal{C}$-action map $\mathtt{L}_\mathcal{C}/\mathtt{I}^{\mathsf{pr}}_{\mathcal{C}} \otimes \Phi_\mathcal{C} \to \mathtt{L}_\mathcal{C}/\mathtt{I}^{\mathsf{pr}}_{\mathcal{C}}$ induces, by restriction, an $(\mathtt{L}_\mathcal{C}, \mathtt{L}_\mathfrak{S})$-bimodule map $\mathtt{L}_\mathcal{C}/\mathtt{I}^{\mathsf{pr}}_{\mathcal{C}} \otimes_{\mathrm{H}} \Phi_\mathcal{C}^{\Delta} \to \mathtt{L}_\mathcal{C}/\mathtt{I}^{\mathsf{pr}}_{\mathcal{C}}$.
\end{Lemma}
\begin{proof}
    Let $m,n\in \mathds{N}$.
    For $[f] \in \mathtt{L}_\mathcal{C}/\mathtt{I}^{\mathsf{pr}}_{\mathcal{C}} (m,n)$ and $[g] \in \Phi_\mathcal{C} (n,n)$ and $\sigma \in \mathfrak{S}_n$, the map gives
    $$[f] \lhd \sigma \otimes [g] \lhd \sigma = [f \circ E_\sigma] \otimes [ E_{\sigma^{-1}} \circ g \circ E_{\sigma}] \mapsto [(f \circ E_\sigma) \circ (E_{\sigma^{-1}} \circ g \circ E_{\sigma})] = [f \circ g \circ E_{\sigma}] = [f\circ g] \lhd \sigma .$$
\end{proof}

\begin{Lemma}
    The map $\mathfrak{O}_\mathcal{C} (1)^{\otimes n} \to \Phi_\mathcal{C} (n,n); \bigotimes^n_{i=1} [f_i] \mapsto [\prod^n_{i=1} f_i]$ induces a right $\mathtt{L}_\mathfrak{S}$-module map $\mathfrak{O}_\mathcal{C} (1)^{\otimes} \to \Phi_\mathcal{C}^{\Delta}$.
\end{Lemma}
\begin{proof}
    Applying Lemma \ref{202606182054} to $n=m$ and $n_k = 1$, we obtain $E_\sigma \circ \prod^n_{i=1} f_i = \prod^n_{i=1} f_{\sigma^{-1}(i)} \circ E_\sigma$ for $\sigma \in \mathfrak{S}_n$ and $f_i \in \mathcal{C} (1,1)$.
    This proves the statement.
\end{proof}

\begin{Defn} \label{202603261506}
    Let $\tilde{\E}$ be the $(\mathtt{L}_\mathcal{C},\mathtt{L}_\mathfrak{S})$-bimodule map obtained from the following composition:    $$\mathtt{L}_\mathcal{C}/\mathtt{I}^{\mathsf{pr}}_{\mathcal{C}} \otimes_{\mathrm{H}} \mathfrak{O}_\mathcal{C} (1)^{\otimes} \to \mathtt{L}_\mathcal{C}/\mathtt{I}^{\mathsf{pr}}_{\mathcal{C}} \otimes_{\mathrm{H}} \Phi_\mathcal{C}^{\Delta} \to \mathtt{L}_\mathcal{C}/\mathtt{I}^{\mathsf{pr}}_{\mathcal{C}}.$$
    We also denote its $(m,-)$-component by $$\tilde{\E}_m : \mathtt{L}_\mathcal{C}/\mathtt{I}^{\mathsf{pr}}_{\mathcal{C}} (m,-) \otimes_{\mathrm{H}} \mathfrak{O}_\mathcal{C} (1)^{\otimes} \to \mathtt{L}_\mathcal{C}/\mathtt{I}^{\mathsf{pr}}_{\mathcal{C}} (m,-) .$$
    In particular, $\tilde{\E}_1$ is a map $\Psi_\mathcal{C} \otimes_{\mathrm{H}} \mathfrak{O}_\mathcal{C} (1)^{\otimes} \to \Psi_\mathcal{C}$.
\end{Defn}

Let $\mathcal{C}$ be a Lawvere theory subject to the condition (ZM).
By the universality of $\mathbf{W}^{\mathsf{o}}$, we have a functor $\mathbf{W}^{\mathsf{o}} \to \mathcal{C}$.
Via this, we regard $\mathtt{L}_\mathcal{C}/\mathtt{I}^{\mathsf{pr}}_{\mathcal{C}}$ as an $(\mathtt{L}_{\mathbf{W}^{\mathsf{o}}},\mathtt{L}_{\mathfrak{S}})$-bimodule.
Then there is a $(\mathtt{L}_{\mathbf{W}^{\mathsf{o}}},\mathtt{L}_{\mathfrak{S}})$-bimodule map $\mathtt{L}_{\mathbf{W}^{\mathsf{o}}}/\mathtt{I}^{\mathsf{pr}}_{\mathbf{W}^{\mathsf{o}}} \to \mathtt{L}_\mathcal{C}/\mathtt{I}^{\mathsf{pr}}_{\mathcal{C}}$.

Using the map $\tilde{\E}$, we obtain an $(\mathtt{L}_{\mathbf{W}^{\mathsf{o}}},\mathtt{L}_{\mathfrak{S}})$-bimodule map by the following composition:
\begin{align} \label{202603262250}
    \mathtt{L}_{\mathbf{W}^{\mathsf{o}}}/\mathtt{I}^{\mathsf{pr}}_{\mathbf{W}^{\mathsf{o}}} \otimes_{\mathrm{H}} \mathfrak{O}_\mathcal{C}(1)^{\otimes} \to \mathtt{L}_\mathcal{C}/\mathtt{I}^{\mathsf{pr}}_{\mathcal{C}} \otimes_{\mathrm{H}} \mathfrak{O}_\mathcal{C}(1)^{\otimes} \stackrel{\tilde{\E}}{\to} \mathtt{L}_\mathcal{C}/\mathtt{I}^{\mathsf{pr}}_{\mathcal{C}}.
\end{align}

Consider the right $\mathtt{L}_{\mathfrak{S}}$-module map $\mathsf{Lin} \to \Psi_{\mathbf{W}^{\mathsf{o}}}$ which assigns $[x_{1} \cdots x_{n}] \in \Psi_{\mathbf{W}^{\mathsf{o}}} (n)$ to $\theta_{1} \cdots \theta_{n} \in \mathsf{Lin}(n)$.
Here, $x_{1} \cdots x_{n} \in \mathsf{W}_n = \mathbf{W}^{\mathsf{o}} (n,1)$.
Using the results in Theorem \ref{202603161027}, we obtain an $(\mathtt{L}_\mathfrak{S},\mathtt{L}_\mathfrak{S})$-bimodule map $\mu\mathsf{Lin} \to \mu\Psi_{\mathbf{W}^{\mathsf{o}}} \cong \mathtt{L}_{\mathbf{W}^{\mathsf{o}}}/\mathtt{I}^{\mathsf{pr}}_{\mathbf{W}^{\mathsf{o}}}$.
    
\begin{Defn} \label{202603301040}
    We define a natural $(\mathtt{L}_\mathfrak{S},\mathtt{L}_\mathfrak{S})$-bimodule map, $$\E = \E_\mathcal{C} : \mu\mathsf{Lin} \otimes_{\mathrm{H}} \mathfrak{O}_{\mathcal{C}}(1)^{\otimes} \to \mathtt{L}_\mathcal{C}/\mathtt{I}^{\mathsf{pr}}_{\mathcal{C}}$$ by the composition, $$\mu\mathsf{Lin} \otimes_{\mathrm{H}} \mathfrak{O}_\mathcal{C} (1)^{\otimes} \to \mathtt{L}_{\mathbf{W}^{\mathsf{o}}}/\mathtt{I}^{\mathsf{pr}}_{\mathbf{W}^{\mathsf{o}}} \otimes_{\mathrm{H}} \mathfrak{O}_\mathcal{C} (1)^{\otimes} \to \mathtt{L}_\mathcal{C}/\mathtt{I}^{\mathsf{pr}}_{\mathcal{C}} .$$
    For $m\in\mathds{N}$, we denote by $\E_m:\mu\mathsf{Lin} (m,-)\otimes_{\mathrm{H}} \mathfrak{O}_{\mathcal{C}}(1)^{\otimes} \to \mathtt{L}_\mathcal{C}/\mathtt{I}^{\mathsf{pr}}_{\mathcal{C}} (m,-)$ the evaluation at the $(m,-)$-components.    
    Evaluated at the components $(1,-)$, we obtain the following map:
    \begin{align} \label{202603262216}
        \E_1 :  \mathsf{Lin} \otimes_{\mathrm{H}}\mathfrak{O}_\mathcal{C}(1)^{\otimes} \to \Psi_\mathcal{C} .
    \end{align}
\end{Defn}

\begin{Example} \label{202604091518}
    Under $\mathfrak{O}_{\G^{\mathsf{o}}} (1) \cong \mathds{k}$ (see Proposition \ref{202509091604}), the map $\E_{\G^{\mathsf{o}}}$ coincides with the isomorphism $\mu\mathsf{Lin} \to \mathtt{L}_{\G^{\mathsf{o}}} / \mathtt{I}^{\mathsf{pr}}_{\G^{\mathsf{o}}}$ in Theorem \ref{202509091931}.
    In particular, the map in (\ref{202603262216}) applied to $\mathcal{C} = \G^{\mathsf{o}}$ agrees with $\mathcal{F}$ in (\ref{202603171311}).
\end{Example}

The following gives an explicit description of $\E_1$ in (\ref{202603262216}).
\begin{Lemma} \label{202512042132}
    For $f_i \in \mathcal{C}(1,1), ~ i \in {\bf n}$ and $\sigma \in \mathfrak{S}_n$, we have
    $$
    \E_1 \left( (\theta_{\sigma(1)} \theta_{\sigma(2)} \cdots \theta_{\sigma(n)}   ) \otimes \bigotimes^{n}_{i=1} [f_i] \right) = [p^\ast_{n,\sigma(1)} (f_{\sigma(1)}) \star p^\ast_{n,\sigma(2)} (f_{\sigma(2)}) \star \cdots \star p^\ast_{n,\sigma(n)} (f_{\sigma(n)})] .
    $$
    Here, $p_{n,k} \in \mathcal{C}(n,1)$ denotes the $k$-th projection.
\end{Lemma}
\begin{proof}
    From the definitions, we have
    \begin{align*}
        \E_1 \left( \theta_{\sigma(1)} \theta_{\sigma(2)} \cdots \theta_{\sigma(n)} \otimes \bigotimes^{n}_{i=1} [f_i] \right) = \left(p_{n,\sigma(1)} \star p_{n,\sigma(2)} \star \cdots \star p_{n,\sigma(n)} \right) \circ (f_1 \times \cdots \times f_n) \mod{\mathtt{I}^{\mathsf{pr}}_{\mathcal{C}}} .
    \end{align*}
    This coincides with the right hand side of the equation in the statement of the lemma, since $(f_1 \times \cdots \times f_n)^\ast : \mathcal{C}_n \to \mathcal{C}_n$ is a monoid map and $(f_1 \times \cdots \times f_n)^\ast (p_{n,k}) = f_k \circ p_{n,k} = p_{n,k}^\ast (f_k)$.
\end{proof}

\subsection{An epimorphism criterion for $\E$}

The following is the main theorem of this section:
\begin{theorem} \label{202603192120}
    If $\mathcal{C}$ satisfies (ZM*), then the map $\E_{\mathcal{C}} : \mu\mathsf{Lin} \otimes_{\mathrm{H}} \mathfrak{O}_{\mathcal{C}} (1)^{\otimes} \to \mathtt{L}_{\mathcal{C}} / \mathtt{I}^{\mathsf{pr}}_{\mathcal{C}}$ is an epimorphism.    
\end{theorem}

The theorem provides a general method for studying the quotient $\mathtt{L}_\mathcal{C} /\mathtt{I}^{\mathsf{pr}}_{\mathcal{C}}$.
Before proving the theorem, we give several useful consequences of the theorem.

\begin{Corollary} \label{202511270956}
    If $\mathcal{C}$ satisfies the condition (ZM*), then we have $\mathtt{I}^{\mathsf{pr}}_{\mathcal{C}} = \mathtt{L}_{\mathcal{C}}$ if and only if $\mathfrak{O}_{\mathcal{C}}(1) \cong 0$.
\end{Corollary}
\begin{proof}
    Since $\mathfrak{O}_\mathcal{C}(1) = \mathtt{L}_{\mathcal{C}} / \mathtt{I}^{\mathsf{pr}}_{\mathcal{C}} (1,1)$ by Lemma \ref{202601081632}, the result is immediate from Theorem \ref{202603192120}.
\end{proof}

We give an application to $\mathcal{C} = \mathbf{W}^{\mathsf{o}}$.
Note that, by Lemma \ref{202601081632} and Proposition \ref{202509091604}, $\E_{\mathbf{W}^{\mathsf{o}}}$ can be viewed as a map $\mu\mathsf{Lin} \to \mathtt{L}_{\mathbf{W}^{\mathsf{o}}} / \mathtt{I}^{\mathsf{pr}}_{\mathbf{W}^{\mathsf{o}}}$.
Recall the $(\mathtt{L}_{\G^{\mathsf{o}}},\mu \mathfrak{Lie})$-bimodule $\tilde{\mu}\mathsf{Lin}$ from Theorem \ref{202603121700}.
\begin{Corollary} \label{202603251841}
    The composition $\Phi_{\mathbf{W}^{\mathsf{o}}} \to \Phi_{\G^{\mathsf{o}}} \cong \mu \mathfrak{Lie}$ gives a monad isomorphism $\Phi_{\mathbf{W}^{\mathsf{o}}} \cong\mu \mathfrak{Lie}$ for which the map $\E_{\mathbf{W}^{\mathsf{o}}} : \tilde{\mu} \mathsf{Lin} \to \mathtt{L}_{\mathbf{W}^{\mathsf{o}}} / \mathtt{I}^{\mathsf{pr}}_{\mathbf{W}^{\mathsf{o}}}$ gives an $(\mathtt{L}_{\mathbf{W}^{\mathsf{o}}},\mu \mathfrak{Lie})$-bimodule isomorphism.
    These fit into the following commutative diagram:
    $$
    \begin{tikzcd}[row sep=scriptsize, column sep=scriptsize]
    & \mu\mathfrak{Lie} \arrow[dl, hookrightarrow] \arrow[rr, "\cong"] \arrow[dd, equal] & & \Phi_{\mathbf{W}^{\mathsf{o}}} \arrow[dl, hookrightarrow] \arrow[dd, "\cong"] \\
    \tilde{\mu} \mathsf{Lin} \arrow[rr, crossing over, 
    "\E_{\mathbf{W}^{\mathsf{o}}}"] \arrow[dd,equal] & & \mathtt{L}_{\mathbf{W}^{\mathsf{o}}} / \mathtt{I}^{\mathsf{pr}}_{\mathbf{W}^{\mathsf{o}}} \\
    & \mu\mathfrak{Lie} \arrow[dl, hookrightarrow] \arrow[rr, "\cong"] & & \Phi_{\G^{\mathsf{o}}} \arrow[dl, hookrightarrow] \\
    \tilde{\mu} \mathsf{Lin} \arrow[rr, "\E_{\G^{\mathsf{o}}}"] & & \mathtt{L}_{\G^{\mathsf{o}}} / \mathtt{I}^{\mathsf{pr}}_{\G^{\mathsf{o}}} \arrow[from=uu, crossing over, "\cong"]\\
    \end{tikzcd}
    $$
    In particular, the eigenmonad adjunction associated with $(\mathtt{L}_{\mathbf{W}^{\mathsf{o}}}, \mathtt{I}^{\mathsf{pr}}_{\mathbf{W}^{\mathsf{o}}})$ coincides with the adjunction in Theorem \ref{202603121700}.
\end{Corollary}
\begin{proof}   
    We first show that $\E_{\mathbf{W}^{\mathsf{o}}} : \mu \mathsf{Lin} \to \mathtt{L}_{\mathbf{W}^{\mathsf{o}}} / \mathtt{I}^{\mathsf{pr}}_{\mathbf{W}^{\mathsf{o}}}$ is an isomorphism of $\mathds{N}\times\mathds{N}$-indexed modules.
    By Theorem \ref{202603192120}, the map $\E_{\mathbf{W}^{\mathsf{o}}}$ is an epimorphism.
    The composition $\mu\mathsf{Lin} \stackrel{\E_{\mathbf{W}^{\mathsf{o}}}}{\to} \mathtt{L}_{\mathbf{W}^{\mathsf{o}}} / \mathtt{I}^{\mathsf{pr}}_{\mathbf{W}^{\mathsf{o}}} \to \mathtt{L}_{\G^{\mathsf{o}}} / \mathtt{I}^{\mathsf{pr}}_{\G^{\mathsf{o}}}$ coincides with $\E_{\G^{\mathsf{o}}}$, which is an isomorphism by the observation in Example \ref{202604091518}.
    Thus, $\E_{\mathbf{W}^{\mathsf{o}}}$ is a monomorphism, hence an isomorphism.

    The above discussion also implies that the natural map $\mathtt{L}_{\mathbf{W}^{\mathsf{o}}} / \mathtt{I}^{\mathsf{pr}}_{\mathbf{W}^{\mathsf{o}}} \to \mathtt{L}_{\G^{\mathsf{o}}} / \mathtt{I}^{\mathsf{pr}}_{\G^{\mathsf{o}}}$ is an isomorphism.
    Since this preserves the $(\mathtt{L}_{\mathbf{W}^{\mathsf{o}}},\Phi_{\mathbf{W}^{\mathsf{o}}})$-bimodule structure by the naturality, the map $\E_{\mathbf{W}^{\mathsf{o}}} : \tilde{\mu}\mathsf{Lin} \to \mathtt{L}_{\mathbf{W}^{\mathsf{o}}} / \mathtt{I}^{\mathsf{pr}}_{\mathbf{W}^{\mathsf{o}}}$ gives an $(\mathtt{L}_{\mathbf{W}^{\mathsf{o}}},\Phi_{\mathbf{W}^{\mathsf{o}}})$-bimodule isomorphism.

    By the above results, the $\mathtt{I}^{\mathsf{pr}}_{\mathbf{W}^{\mathsf{o}}}$-vanishing modules of $\tilde{\mu}\mathsf{Lin}$ and $\mathtt{L}_{\mathbf{W}^{\mathsf{o}}} / \mathtt{I}^{\mathsf{pr}}_{\mathbf{W}^{\mathsf{o}}}$ are identical. 
    By Proposition \ref{202401241330}, Theorem \ref{202603121700} leads to a monad isomorphism $\Phi_{\mathbf{W}^{\mathsf{o}}} \cong\mu \mathfrak{Lie}$.
    The commutative diagram follows from the naturality of $\E_\mathcal{C}$, the canonical bimodules and the eigenmonads.
\end{proof}

\begin{Corollary} \label{202604131147}
    Let $\mathcal{C}$ be a Lawvere theory satisfying (ZM).
    Then the map $\E_{\mathcal{C}} : \tilde{\mu}\mathsf{Lin} \otimes_{\mathrm{H}} \mathfrak{O}_{\mathcal{C}} (1)^{\otimes} \to \mathtt{L}_{\mathcal{C}} / \mathtt{I}^{\mathsf{pr}}_{\mathcal{C}}$ is a $(\mathtt{L}_{\mathbf{W}^{\mathsf{o}}},\mathtt{L}_{\mathfrak{S}})$-bimodule map.
\end{Corollary}
\begin{proof}
    Recall that the map $\E_{\mathcal{C}}$ is a composition of the isomorphism $\tilde{\mu}\mathsf{Lin}\otimes_{\mathrm{H}} \mathfrak{O}_{\mathcal{C}} (1)^{\otimes} \cong \mathtt{L}_{\mathbf{W}^{\mathsf{o}}}/\mathtt{I}^{\mathsf{pr}}_{\mathbf{W}^{\mathsf{o}}} \otimes_{\mathrm{H}} \mathfrak{O}_{\mathcal{C}} (1)^{\otimes}$ obtained from Corollary \ref{202603251841} and the map in 
    (\ref{202603262250}), which are $(\mathtt{L}_{\mathbf{W}^{\mathsf{o}}},\mathtt{L}_{\mathfrak{S}})$-bimodule maps.
\end{proof}

\subsection{Proof of Theorem \ref{202603192120}}

In this section, we prove Theorem \ref{202603192120}.
To this end, we first show that the map $\E_1$ in (\ref{202603262216}) is an epimorphism.
Recall the notation $\mathcal{C}_n^{(i)}$ from Definition \ref{202512090939}:
\begin{Lemma} \label{202509041909}
    Let $N \in \mathds{N}^\ast$, $n\in \mathds{N}$ and $i \in {\bf n}$.
    For $w_1, \cdots, w_N \in \mathcal{C}_n^{(i)}$ and $f_1, \cdots, f_N \in \mathcal{C}_1$, we have
    \begin{align*}
        & w_1 \star p_{n,i}^\ast (f_1) \star w_2 \star p_{n,i}^\ast (f_2) \star \cdots \star w_N \star p_{n,i}^\ast (f_N)  \\
        \equiv& \sum^{N}_{r=1} w_1 \star w_2  \cdots \star w_r \star  p_{n,i}^\ast (f_r) \star w_{r+1} \cdots  \star w_N  \mod{\mathtt{I}^{\mathsf{pr}}_{\mathcal{C}}(1,n)} .
    \end{align*}
\end{Lemma}
\begin{proof}
    We give a sketch of the proof for $i=n$ based on the induction on $N$.
    The same proof applies to the other cases.
    The assertion for $N=1$ is trivial.
    Assume $N \geq 2$.
    Let
    $$
    \rho = w_1 \star p_{n,n}^\ast (f_1) \star w_2 \star p_{n,n}^\ast (f_2) \star \cdots \star w_N \star p_{n,n}^\ast (f_N) \in \mathcal{C}_n .
    $$
    Using $q = (p_{n+1,1},p_{n+1,2},\cdots,p_{n+1,n}) \in \mathcal{C}(n+1,n)$, we set
    $$
    \rho^\prime =  q^\ast (w_1) \star p_{n+1,n}^\ast (f_1) \star q^\ast (w_2) \star p_{n+1,n+1}^\ast (f_2) \star \cdots \star q^\ast (w_N) \star p_{n+1,n+1}^\ast (f_N) \in \mathcal{C}_{n+1} .
    $$
    Using the fact that the precomposition $(\mathrm{id}_{n-1} \times \Delta)^\ast : \mathcal{C}_{n+1} \to \mathcal{C}_n$ is a monoid map, we can deduce 
    $$(\mathrm{id}_{n-1} \times \Delta)^\ast ( \rho^\prime ) = \rho .$$
    Hence, by the definition of $\mathtt{I}^{\mathsf{pr}}_{\mathcal{C}}$, we have
    \begin{align*}
        \rho = (\mathrm{id}_{n-1} \times \Delta)^\ast ( \rho^\prime ) \equiv (\mathrm{id}_n \times \eta )^\ast (\rho^\prime)  + (\mathrm{id}_{n-1} \times \eta \times \mathrm{id}_1 )^\ast (\rho^\prime) \mod{\mathtt{I}^{\mathsf{pr}}_{\mathcal{C}}(1,n)} ,
    \end{align*}
    which equals, by the definitions,
    \begin{align*}
        w_1 \star p_{n,n}^\ast (f_1) \star w_2 \star \cdots \star w_N  + w_1 \star w_2 \star p_{n,n}^\ast (f_2) \star \cdots \star w_N \star p_{n,n}^\ast (f_N)  .
    \end{align*}
    By putting $\tilde{w}_1 = w_1 \star w_2$ and $\tilde{w}_{j} = w_{j+1}$, $\tilde{f}_{j-1} = f_{j}$ for $2 \leq j \leq N$, the final term coincides with
    $\tilde{w}_1 \star p_{n,n}^\ast (\tilde{f}_2) \star \cdots \star \tilde{w}_{N-1} \star p_{n,n}^\ast (\tilde{f}_N)$.
    By the inductive hypothesis, we can complete the proof.
\end{proof}

\begin{Lemma} \label{202509091615}
    For a Lawvere theory $\mathcal{C}$ satisfying the condition (ZM*), the map $\E_{1} : \mathsf{Lin} \otimes_{\mathrm{H}} \mathfrak{O}_{\mathcal{C}}(1)^{\otimes} \to \Psi_\mathcal{C}$ is an epimorphism.
\end{Lemma}
\begin{proof}
    Let $n \in \mathds{N}$.
    We will show that any $w \in \mathcal{C}_n \subset \mathtt{L}_{\mathcal{C}}(1,n)$, modulo $\mathtt{I}^{\mathsf{pr}}_{\mathcal{C}}$, is in the image of the map $\E_{1}: \mathsf{Lin} (n) \otimes \mathfrak{O}_{\mathcal{C}} (1)^{\otimes n} \to \Psi_{\mathcal{C}} (n)$.
    Condition (M*) implies that each component $w_j$ can be expressed as a $\star$-product of elements in $\bigcup^{n}_{i=1} p_{n,i}^\ast (\mathcal{C}_1)$.
    Hence, by iteratively applying Lemma \ref{202509041909}, such $w$, modulo $\mathtt{I}^{\mathsf{pr}}_{\mathcal{C}}$, is identified with a linear combination of elements of the form
    \begin{align} \label{202512042147}
      p_{n,i_1}^\ast (f_1) \star \cdots 
        \star p_{n,i_2}^\ast (f_2) \star \cdots \star  p_{n,i_n}^\ast (f_n) ,
    \end{align}
    where $\{i_1,i_2,\cdots,i_n\} \subset \mathds{N}^\ast$, $\{ f_1, \cdots f_n \} \subset \mathcal{C}_1$, and $l \geq m$.
    If $i_1,i_2,\cdots,i_n$ are not pairwise distinct, then we choose $k \in {\bf n}$ such that $\{i_1,\cdots,i_n\} \subset {\bf n} \backslash \{ k\}$.
    Then the term (\ref{202512042147}) is contained in $\mathcal{C}_n^{(k)}$ (see Definition \ref{202512090939}), since, for each $l\in {\bf n}$, there exists $j$ such that $p_{n,i_l} = p_{n-1,j} \circ P_{n,{\bf n}\backslash \{k\}}$.
    By Proposition \ref{202509042136}, such a term vanishes in $\Psi_{\mathcal{C}}(n)$.
    Therefore, it suffices to consider the case where $i_1,i_2,\cdots,i_n$ are pairwise distinct.
    Now, consider the element $g = \theta_{\sigma(1)} \theta_{\sigma(2)} \cdots \theta_{\sigma(n)} \in \mathsf{Lin}(n)$ where $\sigma\in\mathfrak{S}_n$ is defined by $\sigma(k) = i_k$.
    By Lemma \ref{202512042132}, $\E_{1} (g \otimes \bigotimes^n_{i=1} [f_{\sigma^{-1}(i)}] )$ coincides with (\ref{202512042147}).
    This proves the desired statement.
\end{proof}

Recall the map $\alpha_\mathcal{C}$ in Definition \ref{202603131134}.
\begin{Lemma} \label{202603271041}
    The following diagram commutes:
    $$
    \begin{tikzcd}
        \mathtt{L}_\mathcal{C}/\mathtt{I}^{\mathsf{pr}}_{\mathcal{C}} (m,-) \otimes_{\mathrm{H}} \mathfrak{O}_\mathcal{C} (1)^{\otimes} \ar[r, "\tilde{\E}_{m}"]  & \mathtt{L}_\mathcal{C}/\mathtt{I}^{\mathsf{pr}}_{\mathcal{C}} (m,-) \\
        \mathtt{L}_\mathcal{C}/\mathtt{I}^{\mathsf{pr}}_{\mathcal{C}} (1,-)^{\odot m} \otimes_{\mathrm{H}} \mathfrak{O}_\mathcal{C} (1)^{\otimes} \ar[u, "\alpha_\mathcal{C} \otimes_{\mathrm{H}} \mathrm{id}"] & \\
        (\mathtt{L}_\mathcal{C}/\mathtt{I}^{\mathsf{pr}}_{\mathcal{C}} (1,-) \otimes_{\mathrm{H}} \mathfrak{O}_\mathcal{C} (1)^{\otimes})^{\odot m} \ar[r, "\tilde{\E}_{1}^{\odot m}"] \ar[u, "\cong"] & \mathtt{L}_\mathcal{C}/\mathtt{I}^{\mathsf{pr}}_{\mathcal{C}} (1,-)^{\odot m}  \ar[uu, "\alpha_\mathcal{C}"]
    \end{tikzcd} 
    $$
    where the lower left isomorphism is obtained from the monoidal functor $(-)\otimes_{\mathrm{H}} \mathfrak{O}_{\mathcal{C}}(1)^{\otimes}$.
\end{Lemma}
\begin{proof}
    Let $u \in (\mathtt{L}_\mathcal{C}/\mathtt{I}^{\mathsf{pr}}_{\mathcal{C}} (1,-) \otimes_{\mathrm{H}} \mathfrak{O}_\mathcal{C} (1)^{\otimes})^{\odot m}(n)$ be an element of the form $$u = \bigotimes^m_{k=1} \left([f_k] \otimes \bigotimes^{n_k}_{i=1} [v_{d_{k-1}+i}] \right) \otimes \sigma$$ 
    where $(n_1,\cdots ,n_m)$ is an $m$-partition of $n$.
    Here, $[f_k] \in \mathtt{L}_\mathcal{C}/\mathtt{I}^{\mathsf{pr}}_{\mathcal{C}} (1,n_k) ~ (k \in {\bf m})$, $v_i \in \mathcal{C}_1~(i \in {\bf n})$ and $\sigma \in \mathfrak{S}_n$.
    Following the anticlockwise path, we obtain
    \begin{align*}
        \bigotimes^m_{k=1} \left([f_k] \otimes \bigotimes^{n_k}_{i=1} [v_{d_{k-1}+i}] \right) \otimes \sigma \stackrel{\tilde{\E}_1^{\odot m}}{\mapsto} \bigotimes^m_{k=1} [f_k \circ  \prod^{n_k}_{i=1} v_{d_{k-1}+i}] \otimes \sigma \stackrel{\alpha_\mathcal{C}}{\mapsto} [\prod^m_{k=1} \left( f_k \circ  \prod^{n_k}_{i=1} v_{d_{k-1}+i} \right) \circ E_\sigma] .
    \end{align*}
    The result equals $[\prod^m_{k=1} f_k \circ \prod^{n}_{i=1} v_{i} \circ E_\sigma]$.
    On the other hand, following the clockwise path, we obtain
    \begin{align*}
        &\bigotimes^m_{k=1} \left([f_k] \otimes \bigotimes^{n_k}_{i=1} [v_{d_{k-1}+i}] \right) \otimes \sigma \mapsto \left( \bigotimes^m_{k=1} [f_k] \otimes \sigma \right) \otimes \bigotimes^{n}_{i=1} v_{\sigma^{-1}(i)}  ,\\
        \stackrel{\alpha_\mathcal{C}\otimes_{\mathrm{H}}\mathrm{id}}{\mapsto}& [\prod^m_{k=1} f_k \circ E_\sigma] \otimes \bigotimes^{n}_{i=1} v_{\sigma^{-1}(i)}
        \stackrel{\tilde{\E}_m}{\mapsto} [\prod^m_{k=1} f_k \circ E_\sigma \circ \prod^{n}_{i=1} v_{\sigma^{-1}(i)}] = [\prod^m_{k=1} f_k \circ \prod^{n}_{i=1} v_{i} \circ E_\sigma] .
    \end{align*}
    where the last equality follows from Lemma \ref{202606182054}.
\end{proof}

\begin{Lemma} \label{202603271102}
    Under the isomorphism $\mu\mathsf{Lin} (m,-) \cong \mathsf{Lin}^{\odot m}$ (see Lemma \ref{202603211811}), the following diagram commutes:
    $$
    \begin{tikzcd}
        \mu\mathsf{Lin} (m,-) \otimes_{\mathrm{H}} \mathfrak{O}_{\mathcal{C}} (1)^{\otimes} \ar[r, "\E_m"] & \mathtt{L}_{\mathcal{C}} / \mathtt{I}^{\mathsf{pr}}_{\mathcal{C}} (m,-) \\
        (\mathsf{Lin} \otimes_{\mathrm{H}} \mathfrak{O}_{\mathcal{C}} (1)^{\otimes})^{\odot m} \ar[r, "\E_1^{\odot m}"] \ar[u, "\cong"] & \Psi_\mathcal{C}^{\odot m} \ar[u, "\alpha_\mathcal{C}"]
    \end{tikzcd}
    $$
\end{Lemma}
\begin{proof}
    By the definition of $\E$, it suffices to show that the following diagram commutes:
    $$
    \begin{tikzcd}
        \mathtt{L}_{\mathbf{W}^{\mathsf{o}}}/\mathtt{I}^{\mathsf{pr}}_{\mathbf{W}^{\mathsf{o}}} (m,-) \otimes_{\mathrm{H}} \mathfrak{O}_\mathcal{C}(1)^{\otimes} \ar[r] & \mathtt{L}_\mathcal{C}/\mathtt{I}^{\mathsf{pr}}_{\mathcal{C}} (m,-) \otimes_{\mathrm{H}} \mathfrak{O}_\mathcal{C}(1)^{\otimes} \ar[r, "\tilde{\E}_m"] & \mathtt{L}_\mathcal{C}/\mathtt{I}^{\mathsf{pr}}_{\mathcal{C}} (m,-) \\
        (\mathtt{L}_{\mathbf{W}^{\mathsf{o}}}/\mathtt{I}^{\mathsf{pr}}_{\mathbf{W}^{\mathsf{o}}} (1,-) \otimes_{\mathrm{H}} \mathfrak{O}_\mathcal{C}(1)^{\otimes})^{\odot m} \ar[r] \ar[u, "\cong"]  & (\mathtt{L}_\mathcal{C}/\mathtt{I}^{\mathsf{pr}}_{\mathcal{C}} (1,-) \otimes_{\mathrm{H}} \mathfrak{O}_\mathcal{C}(1)^{\otimes})^{\odot m} \ar[r, "\tilde{\E}_1^{\odot m}"] \ar[u, "\cong"] & \mathtt{L}_\mathcal{C}/\mathtt{I}^{\mathsf{pr}}_{\mathcal{C}} (1,-)^{\odot m} \ar[u, "\alpha_\mathcal{C}"]
    \end{tikzcd}
    $$
    where the sequence in the upper row is nothing but that in (\ref{202603262250}).
    The right square commutes by Lemma \ref{202603271041}.
    The left square commutes since $(-)\otimes_{\mathrm{H}} \mathfrak{O}_\mathcal{C}(1)^{\otimes}$ is a monoidal functor.
\end{proof}

\begin{proof}[Proof of Theorem \ref{202603192120}]
    We now complete the proof of Theorem \ref{202603192120}.
    The result follows from Lemmas \ref{202509091615} and \ref{202603271102}.
\end{proof}

\subsection{Commutative case} \label{202603251934}
Thus far, we have studied Lawvere theories satisfying (ZM).
We now consider those satisfying (ZCM).
Recall the exponential species $\mathsf{Exp}$ and the canonical map $\mathsf{Lin} \to \mathsf{Exp}$.

\begin{Lemma} \label{202604031840}
    For a Lawvere theory $\mathcal{C}$ satisfying (ZCM), the map $\E$ factors as follows:
    $$
    \begin{tikzcd}
        \mu\mathsf{Lin} \otimes_{\mathrm{H}} \mathfrak{O}_\mathcal{C} (1)^{\otimes} \ar[r, "\E"] \ar[d] & \mathtt{L}_{\mathcal{C}} /\mathtt{I}^{\mathsf{pr}}_{\mathcal{C}} \\
        \mu\mathsf{Exp} \otimes_{\mathrm{H}} \mathfrak{O}_\mathcal{C} (1)^{\otimes} \ar[ur, "\exists"'] &
    \end{tikzcd}
    $$
\end{Lemma}
\begin{proof}
    We first consider the map $\E$ for the $(1,-)$-components.
    Lemma \ref{202512042132} implies that the map $\E : \mathsf{Lin} \otimes_{\mathrm{H}} \mathfrak{O}_\mathcal{C} (1)^{\otimes} \to \Psi_{\mathcal{C}}$ factors through the quotient map $\mathsf{Lin} \otimes_{\mathrm{H}} \mathfrak{O}_\mathcal{C} (1)^{\otimes} \to \mathsf{Exp} \otimes_{\mathrm{H}} \mathfrak{O}_\mathcal{C} (1)^{\otimes}$.
    So, we obtain a map $ \mathsf{Exp} \otimes_{\mathrm{H}} \mathfrak{O}_\mathcal{C} (1)^{\otimes} \to \Psi_{\mathcal{C}}$.

    Let $m \in \mathds{N}$ and consider the $(m,-)$-components of the map $\E$.
    The following diagram is fundamental:
    $$
    \begin{tikzcd}
        \mu\mathsf{Lin} (m,-) \otimes_{\mathrm{H}} \mathfrak{O}_{\mathcal{C}} (1)^{\otimes} \ar[r, twoheadrightarrow] \ar[rr, bend left = 15, "\E_m"] & \mu\mathsf{Exp} (m,-) \otimes_{\mathrm{H}} \mathfrak{O}_{\mathcal{C}} (1)^{\otimes} & (\mathtt{L}_\mathcal{C}  /\mathtt{I}^{\mathsf{pr}}_{\mathcal{C}} )(m,-) \\
        (\mathsf{Lin} \otimes_{\mathrm{H}} \mathfrak{O}_{\mathcal{C}} (1)^{\otimes})^{\odot m} \ar[r, twoheadrightarrow] \ar[u, "\cong"] \ar[rr, bend right =15, "\E_1^{\odot m}"'] & (\mathsf{Exp} \otimes_{\mathrm{H}} \mathfrak{O}_{\mathcal{C}} (1)^{\otimes})^{\odot m} \ar[u, "\cong"] \ar[r] & \Psi_{\mathcal{C}}^{\odot m} \ar[u, "\cong"]
    \end{tikzcd}
    $$
    The diagram related to $\E_m$ and $\E_1^{\odot m}$ commutes by Lemma \ref{202603271102}, and the left square commutes by the monoidality of the functor $(-)\otimes_{\mathrm{H}} \mathfrak{O}_{\mathcal{C}} (1)^{\otimes}$.
    Thus, the map $\mu\mathsf{Exp} (m,-) \otimes_{\mathrm{H}} \mathfrak{O}_{\mathcal{C}} (1)^{\otimes} \to (\mathtt{L}_\mathcal{C}  /\mathtt{I}^{\mathsf{pr}}_{\mathcal{C}} )(m,-)$ obtained from the diagram gives the desired map.
\end{proof}

\begin{notation} \label{202604031739}
    Let $\E : \mu\mathsf{Exp} \otimes_{\mathrm{H}} \mathfrak{O}_\mathcal{C} (1)^{\otimes} \to \mathtt{L}_{\mathcal{C}} /\mathtt{I}^{\mathsf{pr}}_{\mathcal{C}}$ denote the map obtained from the lemma with a slight abuse of notations.
    As given in the proof, the $(1,-)$-components give rise to a map $\E_{1} :\mathfrak{O}_\mathcal{C} (1)^{\otimes} \cong \mathsf{Exp} \otimes_{\mathrm{H}} \mathfrak{O}_\mathcal{C} (1)^{\otimes} \to \Psi_{\mathcal{C}}$.
\end{notation}

\begin{Corollary} \label{202509191320}
    For $\mathcal{C}$ a Lawvere theory satisfying the condition (ZCM*), the map $\E : \mu\mathsf{Exp} \otimes_{\mathrm{H}} \mathfrak{O}_\mathcal{C} (1)^{\otimes} \to \mathtt{L}_{\mathcal{C}} /\mathtt{I}^{\mathsf{pr}}_{\mathcal{C}}$ is an epimorphism.
\end{Corollary}
\begin{proof}
    This follows from Theorem \ref{202603192120}.
\end{proof}

\subsection{A natural bimonoid map to $\Psi_\mathcal{C}$}

In this section, we record a biproduct of Lemma \ref{202603271041}.
We show that the map in (\ref{202603262216}) is a bimonoid map.

Let $\mathcal{C}$ be a Lawvere theory with a zero object.
By Lemma \ref{202603240328}, the functor $(-)\otimes_{\mathrm{H}} \mathfrak{O}_\mathcal{C}(1)^{\otimes}$ is monoidal.
Therefore, the coaugmented comonoid species $\Psi_\mathcal{C}$ of Corollary \ref{202603241657} naturally induces a coaugmented comonoid structure on $\Psi_\mathcal{C} \otimes_{\mathrm{H}} \mathfrak{O}_\mathcal{C}(1)^{\otimes}$.
Moreover, if $\mathcal{C}$ satisfies (M), then $\Psi_\mathcal{C} \otimes_{\mathrm{H}} \mathfrak{O}_\mathcal{C}(1)^{\otimes}$ is a bimonoid species.
\begin{prop} \label{202603251911}
    The map $\tilde{\E}_1 : \Psi_\mathcal{C} \otimes_{\mathrm{H}} \mathfrak{O}_\mathcal{C} (1)^{\otimes} \to \Psi_\mathcal{C}$ is a coaugmented comonoid map.
    If $\mathcal{C}$ further satisfies the condition (M), then it is a bimonoid map.
\end{prop}
\begin{proof}
    The map obviously preserves the counit and the coaugmentation.
    It preserves the comultiplication since the following diagram commutes:
    $$
    \begin{tikzcd}
        \Psi_{\mathcal{C}} \otimes_{\mathrm{H}} \mathfrak{O}_\mathcal{C} (1)^{\otimes} \ar[r, "\Delta \otimes_{\mathrm{H}} \mathrm{id}"] \ar[d, "\tilde{\E}_1"] & 
        \Psi_{\mathcal{C}}^{\odot 2} \otimes_{\mathrm{H}} \mathfrak{O}_\mathcal{C} (1)^{\otimes} \ar[r, "\cong"] & \mathtt{L}_\mathcal{C}/\mathtt{I}^{\mathsf{pr}}_{\mathcal{C}} (2,-) \otimes_{\mathrm{H}} \mathfrak{O}_\mathcal{C} (1)^{\otimes} \ar[d, "\tilde{\E}_2"]  & 
        (\Psi_{\mathcal{C}} \otimes_{\mathrm{H}} \mathfrak{O}_\mathcal{C} (1)^{\otimes})^{\odot 2} \ar[d, "\tilde{\E}_1^{\odot 2}"] \ar[l, "\cong"] \\
        \Psi_{\mathcal{C}} \ar[rr, "\Delta \rhd (-)"] & & 
        \mathtt{L}_\mathcal{C}/\mathtt{I}^{\mathsf{pr}}_{\mathcal{C}} (2,-)  \ar[r, "\cong"] &  
        \Psi_{\mathcal{C}}^{\odot 2}
    \end{tikzcd}
    $$
    Here, $\Delta \rhd (-)$ denotes the left $\mathtt{L}_{\mathcal{C}}$-action of the diagonal map $\Delta \in \mathcal{C} (1,2)$.
    The left square commutes by the functoriality of $(-)\otimes_{\mathrm{H}} \mathfrak{O}_\mathcal{C}(1)^{\otimes}$, and the right one commutes by Lemma \ref{202603271041}.
    Dually, if $\mathcal{C}$ satisfies (M), then one may deduce that $\tilde{\E}_1$ is compatible with the multiplication.
\end{proof}

\begin{Corollary} \label{202512051235}
    Let $\mathcal{C}$ be a Lawvere theory satisfying (ZM).
    The map $\E_1 : \mathsf{Lin} \otimes_{\mathrm{H}} \mathfrak{O}_{\mathcal{C}}(1)^{\otimes} \to \Psi_\mathcal{C}$ is a bimonoid map.
    If $\mathcal{C}$ further satisfies (ZCM), then the map $\E_1 : \mathfrak{O}_{\mathcal{C}}(1)^{\otimes} \to \Psi_\mathcal{C}$ is a bimonoid map.
\end{Corollary}
\begin{proof}
    Recall that $\E_1$ is the composition $\mathsf{Lin} \otimes_{\mathrm{H}} \mathfrak{O}_{\mathcal{C}}(1)^{\otimes} \to \Psi_\mathcal{C} \otimes_{\mathrm{H}}\mathfrak{O}_{\mathcal{C}}(1)^{\otimes} \stackrel{\tilde{\E}_1}{\to} \Psi_\mathcal{C}$.
    The first map is a bimonoid map, since this is induced by the bimonoid map $\mathsf{Lin} \cong \Psi_{\mathbf{W}^{\mathsf{o}}} \to \Psi_\mathcal{C}$.
    The second map is a bimonoid map by Proposition \ref{202603251911}.
    The last claim is immediate form the construction of $\E_1 : \mathfrak{O}_{\mathcal{C}}(1)^{\otimes} \to \Psi_\mathcal{C}$.
\end{proof}
\section{Relations with the polynomiality ideal} \label{202512141933}

The purpose of this section is to make precise and prove part of Theorems \ref{202604091210} and \ref{202512111518}. 
To this end, we clarify the relationship between the present paper and certain results of our previous paper \cite{kim2025poly}.
More precisely, we shall investigate several interactions between the primitivity ideal $\mathtt{I}^{\mathsf{pr}}_{\mathcal{C}}$ and the polynomiality ideal $\mathtt{I}^{(d)}_{\mathcal{C}}$ which we will recall in Section \ref{202601071103}.
Throughout this section, let $\mathcal{C}$ be a Lawvere theory satisfying the condition (ZM*) introduced in Section \ref{202606181932}.

\subsection{Recollection of polynomial functor theory} \label{202601071103}

In this section, we give a brief review of Eilenberg-MacLane's polynomial functor theory \cite{EML,hartl2015polynomial} based on our framework.
For later use in this section, we also recall some notation and results from \cite{kim2025poly}, not in full generality.

Let $d \in \mathds{N}^\ast$ and $(n_1,\cdots,n_d)$ be a $d$-partition of $n$.
Note that $$\rho_r = \prod^{r-1}_{k=1} \mathrm{id}_{n_k} \times e_{n_r} \times \prod^d_{k=r+1} \mathrm{id}_{n_k} \in \mathcal{C}(n,n), ~r \in {\bf d}$$ are commuting idempotents, where $e_{k} \in \mathcal{C}(k,k)$ denotes the zero morphism.
Thus, they induce an idempotent $\rho_S {:=} \rho_{i_1} \circ \cdots \circ \rho_{i_l}$, where $S = \{i_1,\cdots,i_l\} \subset {\bf d}$, and an idempotent in the algebra $\mathtt{L}_\mathcal{C} (n,n)$,
\begin{align} \label{202604091600}
    (\mathrm{id}_{n_1}-\rho_1) \circ \cdots \circ (\mathrm{id}_{n_d}-\rho_d) = \sum_{S \subset {\bf d}} (-1)^{|S|} \rho_S = \prod^d_{k=1} ( \mathrm{id}_{n_k} - e_{n_k})  .
\end{align}
Let $\mathtt{M}$ be a left $\mathtt{L}_{\mathcal{C}}$-module.
We define the $\mathds{k}$-module $\mathrm{cr}_d (\mathtt{M}) (n_1,\cdots,n_d)$ to be the joint kernel of the actions $\rho_r \rhd (-) : \mathtt{M} (n) \to \mathtt{M} (n)$ for $r \in {\bf d}$.
Equivalently, this is the image of the action of the idempotent (\ref{202604091600}) on $\mathtt{M}$.
The above construction defines a functor $\mathrm{cr}_d (\mathtt{M}) : \mathcal{C}^{\times d} \to \mathds{k}\mbox{-}\mathsf{Mod}$ called the {\it $d$-th cross effect} of $\mathtt{M}$.

For $d \in \mathds{N}$, a left $\mathtt{L}_{\mathcal{C}}$-module $\mathtt{M}$ is said to be of {\it polynomial degree at most $d$}, and we denote $\deg \mathtt{M} \leq d$, if $\mathrm{cr}_{d+1} (\mathtt{M}) \cong 0$.
Let $\mathtt{L}_{\mathcal{C}}\mbox{-}\mathsf{Mod}^{\leq d}$ denote the full subcategory of $\mathtt{L}_{\mathcal{C}}\mbox{-}\mathsf{Mod}$ consisting of modules of polynomial degree at most $d$.

Let $\mathrm{P}_d (\mathtt{M}) \subset \mathtt{M}$ denote the maximal $\mathtt{L}_{\mathcal{C}}$-submodule of polynomial degree at most $d$.
This yields an increasing filtration $\mathrm{P}_d (\mathtt{M}) \subset \mathrm{P}_{d+1} (\mathtt{M})$.
Then $\mathtt{M}$ is said to be {\it analytic} if $\bigcup_{d \in \mathds{N}} \mathrm{P}_d (\mathtt{M}) = \mathtt{M}$.

We now recall from \cite{kim2025poly} the $d$-th polynomiality ideal $\mathtt{I}^{(d)}_{\mathcal{C}} \subset \mathtt{L}_{\mathcal{C}}$ for comparison with the present work. 
As we will not use its definition explicitly, we refer to \cite{kim2025poly} for details, while the following theorem gives a characterization:

\begin{theorem}[\mbox{\cite[Theorem 4.6]{kim2025poly}}] \label{202604091540}
    The ideal $\mathtt{I}^{(d)}_{\mathcal{C}}$ is a two-sided ideal of $\mathtt{L}_{\mathcal{C}}$ characterized by the property that $\mathtt{M}$ is of polynomial degree at most $d$ if and only if $\mathtt{I}^{(d)}_{\mathcal{C}} \rhd \mathtt{M} \cong 0$.
    Moreover, $\mathrm{P}_d (\mathtt{M}) = \mathrm{V} ( \mathtt{M} ; \mathtt{I}^{(d)}_{\mathcal{C}} )$ (see Definition \ref{202311281540}).
\end{theorem}

By the condition (ZM) on $\mathcal{C}$, $\mathcal{C}_n^{\times m}$ is a monoid.
Let $\mathfrak{I}(\mathcal{C}_{n}^{\times m}) \subset \mathds{k}[\mathcal{C}_{n}^{\times m}]$ denote the augmentation ideal over $\mathds{k}$.
In \cite{kim2025poly}, we establish a structural result on the ideal $\mathtt{I}^{(d)}_{\mathcal{C}}$ in terms of the augmentation ideal $\mathfrak{I}(\mathcal{C}_{n}^{\times m})$:
\begin{theorem}[\mbox{\cite[Theorem 8]{kim2025poly}}] \label{202604042007}
    For a Lawvere theory $\mathcal{C}$ satisfying (ZM*), under the isomorphism $\mathtt{L}_{\mathcal{C}}  (m,n)\cong \mathds{k}[\mathcal{C}_{n}^{\times m}]$, we have $\mathtt{I}^{(d)}_{\mathcal{C}}(m,n) = \mathfrak{I}(\mathcal{C}_{n}^{\times m})^{d+1}$.
\end{theorem}

\subsection{The $\nu$-analyticity ideal}
\label{202403091338}

In this section, using the polynomiality ideals $\mathtt{I}^{(d)}_{\mathcal{C}}$, we introduce some left ideals whose vanishingly generated $\mathtt{L}_{\mathcal{C}}$-modules are analytic.
These are applied, in Section \ref{202511261310}, to the comparison of primitively generated modules and analytic modules.

\begin{Defn} \label{202403042256}
    Let $\nu : \mathds{N} \to \mathds{N}$ be a map.
    We define a submodule $\mathtt{I}^{\nu}_{\mathcal{C}} \subset \mathtt{L}_{\mathcal{C}}$ by
    $$
    \mathtt{I}^{\nu}_{\mathcal{C}} (m,n) {:=} \mathtt{I}^{(\nu (n))}_{\mathcal{C}} (m,n), \quad m,n \in \mathds{N} .
    $$
    This is proved to be a left ideal of $\mathtt{L}_{\mathcal{C}}$ below, so this is referred to as the {\it $\nu$-analyticity ideal}.
\end{Defn}

\begin{Lemma} \label{202410161124}
    The submodule $\mathtt{I}^{\nu}_{\mathcal{C}} \subset \mathtt{L}_{\mathcal{C}}$ is a left ideal.
\end{Lemma}
\begin{proof}
    Let $n \in \mathds{N}$ and consider $d = \nu (n)$.
    For $m,l \in \mathds{N}$, if $g \in \mathtt{L}_{\mathcal{C}} (l,m)$ and $f \in \mathtt{I}^{\nu}_{\mathcal{C}} (m,n)$, then we have $g \circ f \in \mathtt{I}^{\nu}_{\mathcal{C}} (l,n)$ since $\mathtt{I}^{(d)}_{\mathcal{C}} \subset \mathtt{L}_{\mathcal{C}}$ is a left ideal.
\end{proof}

\begin{Lemma} \label{202403041313}
    For a left $\mathtt{L}_{\mathcal{C}}$-module $\mathtt{M}$, the image of the action $\mathtt{L}_{\mathcal{C}} \otimes \mathrm{V} (\mathtt{M} ; \mathtt{I}^{\nu}_{\mathcal{C}}) \to \mathtt{M}$ is contained in $\bigcup_{n \in \mathds{N}} \mathrm{V} ( \mathtt{M} ; \mathtt{I}^{(\nu (n))}_{\mathcal{C}} )$.
\end{Lemma}
\begin{proof}
    Let $m,n \in \mathds{N}$.
    By definition, the map
    $\mathtt{L}_{\mathcal{C}} (m,n) \otimes (\mathrm{V} (\mathtt{M};\mathtt{I}^{\nu}_{\mathcal{C}} )) (n) \to \mathtt{M} (n)$ coincides with
    \begin{align*}
        \mathtt{L}_{\mathcal{C}} (m,n) \otimes (\mathrm{V}  (\mathtt{M} ; \mathtt{I}^{(\nu (n))}_{\mathcal{C}})) (n) \to \mathtt{M} (m) .
    \end{align*}
    Since $\mathtt{I}^{(\nu (n))}_{\mathcal{C}} \subset \mathtt{L}_{\mathcal{C}}$ is a two-sided ideal, the image of this map is contained in $(\mathrm{V}  (\mathtt{M} ; \mathtt{I}^{(\nu (n))}_{\mathcal{C}})) (m)$.
    This argument holds for any $n \in \mathds{N}$, so we obtain the inclusion in the statement.
\end{proof}

\begin{notation}
    For a map $\nu : \mathds{N} \to \mathds{N}$, let $\sup (\nu) \in \mathds{N} \cup \{ \infty \}$ denote the supremum of the image of $\nu$.
\end{notation}

\begin{prop} \label{202403041319}
    Let $\nu : \mathds{N} \to \mathds{N}$ be a map.
    An $\mathtt{I}^{\nu}_{\mathcal{C}}$-vanishingly generated left $\mathtt{L}_{\mathcal{C}}$-module is analytic.
    Furthermore, if $\sup (\nu) < \infty$, then an $\mathtt{I}^{\nu}_{\mathcal{C}}$-vanishingly generated left $\mathtt{L}_{\mathcal{C}}$-module is of polynomial degree at most $\sup(\nu)$.
\end{prop}
\begin{proof}
    Since $\mathtt{M}$ is $\mathtt{I}^{\nu}_{\mathcal{C}}$-vanishingly generated, by Lemma \ref{202403041313}, we have $\mathtt{M} = \bigcup_{n \in \mathds{N}} \mathrm{V} ( \mathtt{M} ; \mathtt{I}^{(\nu (n))}_{\mathcal{C}} )$.
    If $\sup (\nu)$ is finite, then the union equals $\mathrm{V} ( \mathtt{M} ; \mathtt{I}^{(\sup(\nu ))}_{\mathcal{C}} )$; and $\varinjlim_{d} \mathrm{V} ( \mathtt{M} ; \mathtt{I}^{(d)}_{\mathcal{C}})$ otherwise.
    These lead to the desired statement.
\end{proof}

\subsection{An inclusion relation}
\label{202511261310}

In this section, we present a (partial) inclusion relation between the primitivity ideal and the polynomiality ideal.
Furthermore, we also give some applications comparing primitively generated $\mathtt{L}_{\mathcal{C}}$-modules and analytic $\mathtt{L}_{\mathcal{C}}$-modules.

\begin{Lemma} \label{202509201154}
    Let $\mathcal{C}$ be a Lawvere theory satisfying the condition (ZM*).
    For $d \in \mathds{N}$, the left $\mathtt{L}_{\mathcal{C}}$-module $(\mathtt{L}_{\mathcal{C}} / \mathtt{I}^{\mathsf{pr}}_{\mathcal{C}}) (-,d)$ is of polynomial degree $\leq d$.
\end{Lemma}
\begin{proof}
    If $\mathds{k} \neq 0$, then the polynomial degree of the left $\mathtt{L}_{\mathbf{W}^{\mathsf{o}}}$-module $\tilde{\mu}\mathsf{Lin} (-,d)$ is $d$.
    For instance, this follows from \cite[Proposition A.8]{powell2024analytic} or \cite[Proposition 6.1]{kim2024analytic}.
    By Theorem \ref{202603192120} and Corollary \ref{202604131147}, $\E : \tilde{\mu}\mathsf{Lin} (-,d) \otimes \mathfrak{O}_\mathcal{C} (1)^{\otimes d} \to \mathtt{L}_{\mathcal{C}} / \mathtt{I}^{\mathsf{pr}}_{\mathcal{C}} (-,d)$ is an epimorphism of left $\mathtt{L}_{\mathbf{W}^{\mathsf{o}}}$-modules.
    Since polynomial $\mathtt{L}_{\mathbf{W}^{\mathsf{o}}}$-modules are closed under quotients, the $\mathtt{L}_{\mathbf{W}^{\mathsf{o}}}$-module $\mathtt{L}_{\mathcal{C}} / \mathtt{I}^{\mathsf{pr}}_{\mathcal{C}}$, hence the left $\mathtt{L}_{\mathcal{C}}$-module $\mathtt{L}_{\mathcal{C}} / \mathtt{I}^{\mathsf{pr}}_{\mathcal{C}}$, is of polynomial degree at most $d$.
\end{proof}

\begin{theorem}[Inclusion relation] \label{202509201240}
    For $d\in \mathds{N}$, we have $$\mathtt{I}^{(d)}_{\mathcal{C}} (-,d) \subset \mathtt{I}^{\mathsf{pr}}_{\mathcal{C}} (-,d).$$
\end{theorem}
\begin{proof}
    By Lemma \ref{202509201154}, we have $$\mathrm{V} \left( \mathtt{L}_{\mathcal{C}} / \mathtt{I}^{\mathsf{pr}}_{\mathcal{C}}   ; \mathtt{I}^{(d)}_{\mathcal{C}} \right) (-,d) = (\mathtt{L}_{\mathcal{C}} / \mathtt{I}^{\mathsf{pr}}_{\mathcal{C}} ) (-,d) .$$
    This gives $\mathtt{I}^{(d)}_{\mathcal{C}} (-,d) = \mathtt{I}^{(d)}_{\mathcal{C}} \circ \mathtt{L}_{\mathcal{C}}(-,d) \subset \mathtt{I}^{\mathsf{pr}}_{\mathcal{C}} (-,d)$.
\end{proof}

\begin{Corollary} \label{202403041744}
    Let $\nu : \mathds{N} \to \mathds{N}$ be a map such that $\nu (n) \geq n$.
    The primitivity ideal contains the $\nu$-analyticity ideal, i.e. $\mathtt{I}^{\nu}_{\mathcal{C}} \subset \mathtt{I}^{\mathsf{pr}}_{\mathcal{C}}$.
\end{Corollary}
\begin{proof}
    Let $d \in \mathds{N}$.
    By Theorem \ref{202509201240}, we have $\mathtt{I}^{\nu}_{\mathcal{C}} (-,d) = \mathtt{I}^{(\nu (d))}_{\mathcal{C}} (-,d)  \subset \mathtt{I}^{(d)}_{\mathcal{C}} (-,d)\subset \mathtt{I}^{\mathsf{pr}}_{\mathcal{C}}(-,d)$ where we use the fact that $\mathtt{I}^{(d_2)}_{\mathcal{C}}   \subset \mathtt{I}^{(d_1)}_{\mathcal{C}}$ for $d_1 \leq d_2$.
\end{proof}

\begin{Corollary} \label{202403060133}
    Let $\nu : \mathds{N} \to \mathds{N}$ be a map such that $\nu (n) \geq n$.
    For $\mathtt{L}_{\mathcal{C}}$-modules, we have the following implications:
    \begin{align*}
        \mathrm{primitively~generated} \Rightarrow \mathtt{I}^{\nu}_{\mathcal{C}}\mathrm{-vanishingly~generated} \Rightarrow \mathrm{analytic} .
    \end{align*}
\end{Corollary}
\begin{proof}
    By Lemma \ref{202405291701} and Corollary \ref{202403041744}, if a left $\mathtt{L}_{\mathcal{C}}$-module is primitively generated, i.e. $\mathtt{I}^{\mathsf{pr}}_{\mathcal{C}}$-vanishingly generated, then it is $\mathtt{I}^{\nu}_{\mathcal{C}}$-vanishingly generated.
    The last one follows from Proposition \ref{202403041319}.
\end{proof}

\begin{remark}
    This corollary, applied to $\mathcal{C} = \G^{\mathsf{o}}$, recovers \cite[Proposition 9.9]{kim2024analytic}.
\end{remark}

\subsection{A supplement relation} \label{202512051208}

In this section, we give another relation between the primitivity ideal and the polynomiality ideal, which we call a supplement relation.
This provides a key step for proving the statements  under the assumption (ZM*) in Theorems \ref{202604091210} and \ref{202512111518}.
In particular, we show that the eigenmonad associated with the primitivity ideal is upper triangular in the sense of Definition \ref{202512071707}.

\begin{theorem}[Supplement relation] \label{202509201139}
    For $d \in \mathds{N}$, we have $$(\mathtt{I}^{\mathsf{pr}}_{\mathcal{C}} + \mathtt{I}^{(d)}_{\mathcal{C}})(-,r) = \mathtt{L}_{\mathcal{C}}(-,r), \quad r > d.$$
\end{theorem}

For the following, we recall from Notation \ref{202512162232} the monoid operation $\star$ on $\mathcal{C}_n^{\times m}$.

\begin{proof}
    Since $\mathtt{I}^{\mathsf{pr}}_{\mathcal{C}} + \mathtt{I}^{(d)}_{\mathcal{C}}$ is a left ideal of $\mathtt{L}_{\mathcal{C}}$, it suffices to prove that
    $$
    \mathrm{id}_r \in (\mathtt{I}^{\mathsf{pr}}_{\mathcal{C}} + \mathtt{I}^{(d)}_{\mathcal{C}}) (r,r), \quad r > d
    $$
    where $\mathrm{id}_r \in \mathcal{C}(r,r)\subset \mathtt{L}_{\mathcal{C}} (r,r)$ is the identity on $r$.
    Let $g_k \in \mathcal{C}_{r}^{\times r},~k\in {\bf d+1} = \{1,2,\cdots,d+1\}$ be
    \begin{align*}
        g_k = 
        \begin{cases}
            ( e , \cdots, e , \overbrace{p_{r,k}}^{k\mathrm{-th}} ,e , \cdots, e) & \mathrm{if~} k \leq d , \\
            ( e , \cdots, e ,\overbrace{p_{r,d+1}, p_{r,d+2}, \cdots, p_{r,r}}^{\text{from the $(d+1)$-th}}) & \mathrm{if~} k = d+1 .
        \end{cases}
    \end{align*}
    Recall $\mathcal{C}_{r}^{(i)}$ from Definition \ref{202512090939}.
    By the result in Example \ref{202606191434}, for distinct $i, j \in {\bf d+1}$, we have $g_j \in (\mathcal{C}_{r}^{(i)})^{\times r}$.
    Since $\mathcal{C}^{(i)}_r \subset \mathcal{C}_r$ is a submonoid, we have
    $g_{i_1} \star \cdots \star g_{i_l} \in (\mathcal{C}_{r}^{(i)})^{\times r}$, 
    where $i \in {\bf d+1}$ and $\{i_1,\cdots,i_l\} \subset {\bf d+1} \backslash \{ i\}$, so that $$g_{i_1} \star \cdots \star g_{i_l} \in \mathtt{I}^{\mathsf{pr}}_{\mathcal{C}} (r,r)$$ follows from Proposition \ref{202509042136}.
    Hence, we obtain
    $$(g_1 -e) \star (g_2 -e)\star \cdots \star(g_{d+1}-e) \equiv g_1 \star g_2 \star \cdots \star  g_{d+1} \mod{\mathtt{I}^{\mathsf{pr}}_{\mathcal{C}} (r,r)}. $$
    On the other hand, Theorem \ref{202604042007} implies that  
    $$
    (g_1 -e) \star (g_2 -e)\star \cdots \star(g_{d+1}-e) \equiv 0 \mod{\mathtt{I}^{(d)}_{\mathcal{C}} (r,r)} .
    $$
    By definition, $g_1 \star g_2 \star \cdots \star  g_{d+1}$ coincides with $( p_{r,1}, p_{r,2} , \cdots, p_{r,r})$, which corresponds to $\mathrm{id}_r$ as in (\ref{202606190854}).
    Therefore, we obtain $\mathrm{id}_r \in (\mathtt{I}^{\mathsf{pr}}_{\mathcal{C}} + \mathtt{I}^{(d)}_{\mathcal{C}}) (r,r)$.
\end{proof}

\begin{Corollary} \label{202509201331}
    Let $\mathtt{M}$ be a left $\mathtt{L}_{\mathcal{C}}$-module of polynomial degree at most $d$.
    Then we have 
    $$
    \mathrm{V} ( \mathtt{M} ; \mathtt{I}^{\mathsf{pr}}_{\mathcal{C}}) (r) \cong 0 , \quad r > d .
    $$
\end{Corollary}
\begin{proof}
    Let $v \in \mathrm{V} ( \mathtt{M} ; \mathtt{I}^{\mathsf{pr}}_{\mathcal{C}}) (r)$.
    By the definition of $\mathrm{V} ( \mathtt{M} ; \mathtt{I}^{\mathsf{pr}}_{\mathcal{C}})$, we have $\mathtt{I}^{\mathsf{pr}} \rhd v \cong 0$.
    Since $\deg \mathtt{M} \leq d$, we obtain $\mathtt{I}^{(d)} \rhd v \cong 0$ from Theorem \ref{202604091540}.
    Thus, Theorem \ref{202509201139} implies that $\mathtt{L}_{\mathcal{C}} \rhd v = 0$, so $v = 0$.
\end{proof}

\begin{Corollary} \label{202512090929}
    The monad $\Phi_{\mathcal{C}}$ on $\mathds{N}$ is upper triangular.
\end{Corollary}
\begin{proof}
    By Proposition \ref{202401241330}, we have $\Phi_{\mathcal{C}} (m,n) \cong \mathrm{V}(\mathtt{L}_{\mathcal{C}} /\mathtt{I}^{\mathsf{pr}}_{\mathcal{C}} ;\mathtt{I}^{\mathsf{pr}}_{\mathcal{C}}) (m,n)$.
    So, it suffices to show that the latter is trivial for $m > n$.
    We obtain the desired statement by applying Corollary \ref{202509201331} to $\mathtt{M} = (\mathtt{L}_{\mathcal{C}} /\mathtt{I}^{\mathsf{pr}}_{\mathcal{C}}) (-,n)$ which is of polynomial degree at most $n$ by Lemma \ref{202509201154}.
\end{proof}

\subsection{Polynomial $\mathtt{L}_{\mathcal{C}}$-modules and truncated $\Phi_{\mathcal{C}}$-modules}

In this section, we improve Corollary \ref{202512081855} by considering polynomial degrees and the notion introduced in Definition \ref{202512090927}.

\begin{Lemma} \label{202509201334}
    For a $d$-truncated left $\Phi_{\mathcal{C}}$-module $\mathtt{N}$, the induced left $\mathtt{L}_{\mathcal{C}}$-module $\mathtt{L}_{\mathcal{C}}/ \mathtt{I}^{\mathsf{pr}}_{\mathcal{C}}\otimes_{\Phi_{\mathcal{C}}} \mathtt{N}$ is of polynomial degree at most $d$.
\end{Lemma}
\begin{proof}
    By the hypothesis, we have
    $$(\mathtt{L}_{\mathcal{C}}/ \mathtt{I}^{\mathsf{pr}}_{\mathcal{C}}\otimes \mathtt{N}) (n) = \bigoplus_{k\in\mathds{N}} (\mathtt{L}_{\mathcal{C}}/ \mathtt{I}^{\mathsf{pr}}_{\mathcal{C}}) (n,k) \otimes \mathtt{N} (k) = \bigoplus_{k\leq d} (\mathtt{L}_{\mathcal{C}}/ \mathtt{I}^{\mathsf{pr}}_{\mathcal{C}} )(n,k) \otimes \mathtt{N} (k).$$
    Hence, by Lemma \ref{202509201154}, the left $\mathtt{L}_{\mathcal{C}}$-module $\mathtt{L}_{\mathcal{C}}/ \mathtt{I}^{\mathsf{pr}}_{\mathcal{C}}\otimes \mathtt{N}$ is of polynomial degree at most $d$, since the polynomial degree of $(\mathtt{L}_{\mathcal{C}}/ \mathtt{I}^{\mathsf{pr}}_{\mathcal{C}} )(-,k)$ does not exceed $k$, so $d$.
    Since $\mathtt{L}_{\mathcal{C}}/ \mathtt{I}^{\mathsf{pr}}_{\mathcal{C}}\otimes_{\Phi_{\mathcal{C}}} \mathtt{N}$ is a quotient of $\mathtt{L}_{\mathcal{C}}/ \mathtt{I}^{\mathsf{pr}}_{\mathcal{C}}\otimes \mathtt{N}$, we obtain the claim.
\end{proof}

\begin{Lemma} \label{202509201612}
    For $k \in {\bf d}$, 
    \begin{align*} 
        \bar{\Delta}^{d,k} \circ (\mathrm{id}_1 - e)^{\times d} \in \mathtt{I}^{(d)}_{\mathcal{C}} (d+1,d) .
    \end{align*}
\end{Lemma}
\begin{proof}
    A direct verification shows that
    $$\bar{\Delta} \circ (\mathrm{id}_1 - e) = (\mathrm{id}_1  - e)^{\times 2} \circ \Delta.$$
    Using this, we obtain
    \begin{align*}
        &\bar{\Delta}^{d,k} \circ (\mathrm{id}_1 - e)^{\times d} = (\mathrm{id}_1 - e)^{\times (k-1)} \times \left( \bar{\Delta} \circ (\mathrm{id}_1 - e) \right) \times (\mathrm{id}_1 - e)^{\times (d-k)},\\
        =& (\mathrm{id}_1 -e)^{\times (d+1)} \circ (\mathrm{id}_{k-1} \times \Delta \times \mathrm{id}_{d-k}) .
    \end{align*}
    Recall from (\ref{202606201739}) that $(\mathrm{id}_{k-1} \times \Delta \times \mathrm{id}_{d-k}) =  \left( \prod^{d+1}_{i=1} g_i \right) \circ \Delta^{(d+1)}_d$, where $g_i = p_{d,i}$ if $i \leq k$, and $g_i = p_{d,i-1}$ if $i > k$.
    This leads to
    \begin{align*}
        (\mathrm{id}_1 - e)^{\times (d+1)} \circ ( \mathrm{id}_{k-1} \times \Delta \times \mathrm{id}_{d-k}) = (\mathrm{id}_1 - e)^{\times (d+1)} \circ  \left( \prod^{d+1}_{i=1} g_i \right) \circ \Delta^{(d+1)}_d = \left( \prod^{d+1}_{i=1} ( g_i - e) \right) \circ \Delta^{(d+1)}_d .
    \end{align*}
    Let $h_i = ( e , \cdots, e , g_i ,e , \cdots, e)$ with $g_i$ in the $i$-th component.
    Then the above result proves
    $$(\mathrm{id}_1 - e)^{\times (d+1)} \circ (\mathrm{id}_{k-1} \times \Delta \times \mathrm{id}_{d-k}) = (h_1 -e ) \star \cdots \star (h_{d+1} -e) \in \mathfrak{I}(\mathcal{C}_{d}^{\times (d+1)})^{d+1} .$$
    Therefore, the desired statement follows from Theorem \ref{202604042007}.
\end{proof}

\begin{Lemma} \label{202509201609}
    Let $\mathtt{M}$ be a left $\mathtt{L}_{\mathcal{C}}$-module.
    For $d \in \mathds{N}$, we have $\left( \mathrm{V} (\mathtt{M} ; \mathtt{I}^{\mathsf{pr}}_{\mathcal{C}}) \right) (d) \subset (\mathrm{cr}_d (\mathtt{M})) (1,\cdots,1)$.
    If $\mathtt{M}$ is of polynomial degree at most $d$, then they coincide with each other.
\end{Lemma}
\begin{proof}
    Let $v \in \left( \mathrm{V} (\mathtt{M} ; \mathtt{I}^{\mathsf{pr}}_{\mathcal{C}}) \right) (d)$.
    By the relation $(\epsilon \times \epsilon) \circ \bar{\Delta} = - \epsilon$, for any $j \in {\bf d}$, we obtain 
    $$(\mathrm{id}_{j-1} \times \epsilon \times \mathrm{id}_{d-j} ) \rhd v = - \left( h \circ \bar{\Delta}^{d,j} \right) \rhd v = - h \rhd (\bar{\Delta}^{d,j} \rhd v) = 0$$ 
    where $h = \mathrm{id}_{j-1} \times \epsilon \times \epsilon \times \mathrm{id}_{d-j}$.
    Hence, $v \in (\mathrm{cr}_d (\mathrm{M})) (1,\cdots,1)$.

    We now assume that $\deg \mathtt{M} \leq d$.
    Let $v \in (\mathrm{cr}_d (\mathtt{M})) (1,\cdots,1)$.
    By the definition of the cross-effect, we have $(\mathrm{id}_1 - e)^{\times d} \rhd v = v$.
    Hence, Theorem \ref{202604091540} and Lemma \ref{202509201612} lead to 
    $\bar{\Delta}^{d,k} \rhd v = \left( \bar{\Delta}^{d,k} \circ (\mathrm{id}_1 - e)^{\times d} \right) \rhd v = 0$.
    This implies $v \in \left( \mathrm{V} (\mathtt{M} ; \mathtt{I}^{\mathsf{pr}}_{\mathcal{C}}) \right) (d)$.
\end{proof}

\begin{theorem} \label{202509201658}
    For $d \in \mathds{N}$, the eigenmonad adjunction given in Corollary \ref{202512081855} induces 
    $$
    \begin{tikzcd}
             \mathbb{L} : \Phi_{\mathcal{C}}\mbox{-}\mathsf{Mod}^{d-\mathsf{trun}} \arrow[r, shift right=1ex, ""{name=G}] & \mathtt{L}_{\mathcal{C}}\mbox{-}\mathsf{Mod}^{\leq d} :  \mathbb{R} \arrow[l, shift right=1ex, ""{name=F}]
            \arrow[phantom, from=G, to=F, , "\scriptscriptstyle\boldsymbol{\top}"].
    \end{tikzcd}
    $$
    Furthermore, for a $\mathtt{L}_{\mathcal{C}}$-module $\mathtt{M}$ of polynomial degree at most $d$, we have $(\mathbb{R}(\mathtt{M})) (d) = \mathrm{cr}_{d} (\mathtt{M}) (1, \cdots, 1)$.
\end{theorem}
\begin{proof}
    The result follows from the restrictions of the eigenmonad adjunction to the subcategories.
    By Corollary \ref{202509201331} and Lemma \ref{202509201334}, these restrictions induce the well-defined functors in the statement.
    The last assertion is immediate from Lemma \ref{202509201609}.
\end{proof}

\begin{remark}
    The adjunction does not give an equivalence in general, as shown in Example \ref{202512091138}.
    On the other hand, we will see an example in Theorem \ref{202509201759} where this yields an equivalence of categories.
\end{remark}

{\Large \part{Further applications} \label{202604131423}}

In this part, we apply the results of Part \ref{202512151908} to two examples of Lawvere theories arising from:
\begin{enumerate}
    \item finitely generated free $\mathcal{R}$-semisimple groups, where $\mathcal{R}$ is a radical functor for groups; and
    \item free $R$-modules of finite rank, where $R$ is a unital ring.
\end{enumerate}
As an important special case of the former, we also study the opposite category of finitely generated free nilpotent groups of class $\leq c$, where $c \in \mathds{N}^\ast$.
We compute the associated canonical bimodules and the eigenmonads.
As a consequence, we obtain proofs of Theorems \ref{202512111051}, \ref{202604111753} and \ref{202512111638}.

\vspace{4mm}
\section{The primitivity ideal for free $\mathcal{R}$-semisimple groups} \label{202512171132}

Let $\mathcal{R}$ be a radical functor for groups.
By Corollary \ref{202512081855} and Theorem \ref{202603161027}, there is an adjunction between the category of left $\mathtt{L}_{\G_{\mathcal{R}}^{\mathsf{o}}}$-modules and the category of left $\Phi_{\G_{\mathcal{R}}^{\mathsf{o}}}$-modules:
\begin{equation} \label{202604051252}
\begin{tikzcd}
         \tilde{\mu}\Psi_{\G_{\mathcal{R}}^{\mathsf{o}}} \otimes_{\Phi_{\G_{\mathcal{R}}^{\mathsf{o}}}} (-) : \Phi_{\G_{\mathcal{R}}^{\mathsf{o}}}\mbox{-}\mathsf{Mod} \arrow[r, shift right=1ex, ""{name=G}] & \mathtt{L}_{\G_{\mathcal{R}}^{\mathsf{o}}}\mbox{-}\mathsf{Mod} : \mathrm{Hom}_{\G_{\mathcal{R}}^{\mathsf{o}}} (\tilde{\mu}\Psi_{\G_{\mathcal{R}}^{\mathsf{o}}} , -) \arrow[l, shift right=1ex, ""{name=F}]
        \arrow[phantom, from=G, to=F, , "\scriptscriptstyle\boldsymbol{\top}"].
\end{tikzcd}
\end{equation}
Moreover, as in Theorem \ref{202509201658}, this gives a correspondence between polynomial $\mathtt{L}_{\G_{\mathcal{R}}^{\mathsf{o}}}$-modules and truncated $\Phi_{\G_{\mathcal{R}}^{\mathsf{o}}}$-modules.
In the case where $\mathcal{R}$ is trivial, this recovers Theorem \ref{202603121700}, in particular Powell's adjunction \cite{powell2024analytic}.
In this section, we examine this adjunction for general $\mathcal{R}$ in more detail.

\subsection{The map $\Psi_{\G^{\mathsf{o}}} \to \Psi_{\G_{\mathcal{R}}^{\mathsf{o}}}$ and its kernel}

The free $\mathcal{R}$-semisimple group on $n$ generators is naturally identified with $\mathsf{F}_n / \mathcal{R}(\mathsf{F}_n)$.
The quotient maps $\mathsf{F}_n \to \mathsf{F}_n / \mathcal{R}(\mathsf{F}_n)$ assemble into a functor $\G^{\mathsf{o}} \to \G_{\mathcal{R}}^{\mathsf{o}}$, which induces an epimorphism $\Psi_{\G^{\mathsf{o}}} \to \Psi_{\G_{\mathcal{R}}^{\mathsf{o}}}$.
In this section, we study $\Psi_{\G_{\mathcal{R}}^{\mathsf{o}}}$ by computing the kernel of this map.
We begin with an alternative description of the map using Theorem \ref{202509091931}:
\begin{Defn} \label{202604021150}
    We define the right $\mathtt{L}_\mathfrak{S}$-module $\K_\mathcal{R}$ to be the kernel of the map $\mathsf{Lin} \cong \Psi_{\G^{\mathsf{o}}} \to \Psi_{\G_{\mathcal{R}}^{\mathsf{o}}}$.
\end{Defn}
By the definition, we have an isomorphism $\Psi_{\G_{\mathcal{R}}^{\mathsf{o}}} \cong \mathsf{Lin}/ \K_{\mathcal{R}}$.

\begin{Example}
    If $\mathcal{R} (G) \equiv 1$, then $\K_{\mathcal{R}}  \cong 0$.
    Further examples arise in Section \ref{202512011517}.
\end{Example}

Recall the map $\T_{n}$ from Definition \ref{202509110944}.
\begin{prop} \label{202603292144}
    For $n \in \mathds{N}^\ast$, let $q_n : \mathds{k}[\mathsf{F}_n] \to \mathds{k}[\mathsf{F}_n/ \mathcal{R}(\mathsf{F}_n)]$ be the quotient map.  
    We have
    $$\K_{\mathcal{R}} (n) = \T_{n} ( \mathrm{Ker} ( q_{n})) \subset \mathsf{Lin}(n) .$$
\end{prop}

This proposition allows us to compute $\K_\mathcal{R}$ more effectively.
Before proving the proposition, we give some easy consequences.
Recall $|\mathcal{R}|$ from Definition \ref{202512311144}, and put $r = | \mathcal{R}|$
\begin{Corollary} \label{202512102003}
    $\K_{\mathcal{R}} (1)
    \cong r \mathds{k}$.
\end{Corollary}
\begin{proof}
    It is clear that the $\mathds{k}$-submodule $\mathrm{Ker} (q_1) \subset \mathds{k}[\mathsf{F}_1]$ is generated by $x_1^{k+r} - x_1^{k}$ where $k \in \mathds{Z}$.
    Since $\T_1 (x_1^{k+r} - x_1^{k}) = r\theta_1$, by Proposition \ref{202603292144}, we obtain the result.
\end{proof}

This result yields the following corollaries:

\begin{Corollary} \label{202512102031}
    We have an isomorphism of $\mathds{k}$-algebras $\Phi_{\G_{\mathcal{R}}^{\mathsf{o}}} (n,n) \cong (\mathds{k}/r\mathds{k}) [\mathfrak{S}_n]$.
\end{Corollary}
\begin{proof}
    By Corollary \ref{202603221141}, we have $\Phi_{\G_{\mathcal{R}}^{\mathsf{o}}} (n,n) \cong \Phi_{\G_{\mathcal{R}}^{\mathsf{o}}} (1,1) \thicksim \mathfrak{S}_n$.
    It follows from Lemma \ref{202601081632} that $\Phi_{\G_{\mathcal{R}}^{\mathsf{o}}} (1,1) \cong \Psi_{\G_{\mathcal{R}}^{\mathsf{o}}} (1)$.
    By Corollary \ref{202512102003}, this is isomorphic to $\mathsf{Lin}(1)/\K_{\mathcal{R}}(1) \cong \mathds{k} / r \mathds{k}$.
    Thus, the claim follows from $\mathds{k}/r\mathds{k} \thicksim \mathfrak{S}_n \cong  (\mathds{k}/r\mathds{k}) [\mathfrak{S}_n]$.
\end{proof}

\begin{Corollary} 
    We have $\mathtt{I}^{\mathsf{pr}}_{\G_{\mathcal{R}}^{\mathsf{o}}} = \mathtt{L}_{\G_{\mathcal{R}}^{\mathsf{o}}}$ if and only if $r \cdot 1_{\mathds{k}}$ is invertible in $\mathds{k}$.
\end{Corollary}
\begin{proof}
    Note that $r \cdot 1_{\mathds{k}}$ is invertible in $\mathds{k}$ if and only if $\mathds{k}/r\mathds{k} =0$.
    By Corollary \ref{202512102031}, the latter condition is equivalent to $\mathfrak{O}_{\G_{\mathcal{R}}^{\mathsf{o}}} (1) = \Phi_{\G_{\mathcal{R}}^{\mathsf{o}}}(1,1)\cong 0$.
    Hence, the result follows from Corollary \ref{202511270956}.
\end{proof}

We conclude this section with a proof of Proposition \ref{202603292144}.
\begin{proof}[Proof of Proposition \ref{202603292144}]
    First, note that the composition $\mathds{k}[\mathsf{F}_n] \stackrel{\T_n}{\to} \mathsf{Lin}(n) \cong \Psi_{\G^{\mathsf{o}}} \twoheadrightarrow \Psi_{\G_{\mathcal{R}}^{\mathsf{o}}}$ coincides with $\mathds{k}[\mathsf{F}_n] \stackrel{q_n}{\to} \mathds{k}[\mathsf{F}_n/ \mathcal{R}(\mathsf{F}_n)] \twoheadrightarrow \Psi_{\G_{\mathcal{R}}^{\mathsf{o}}}$.
    Hence, the map $\T_n$ induces $\mathrm{Ker}(q_n) \to \K_{\mathcal{R}}(n)$.
    For the proof of Proposition \ref{202603292144}, it suffices to show that this map is surjective.
    Let $W \subset \mathds{k}[\mathsf{F}_n]$ be the intersection of $\mathrm{Ker} (q_n)$ and $\mathtt{I}^{\mathsf{pr}}_{\G^{\mathsf{o}}} (1,n)$.
    Then we have a commutative diagram as follows.
    The third column, the second and third rows are exact by the definitions.    
    By Theorem \ref{202509091931}, the second column is also exact.
    $$
    \begin{tikzcd}
        & 0 \ar[d]& 0 \ar[d]& 0 \ar[d]& & \\
        0 \ar[r]& W\ar[d] \ar[r] & \mathtt{I}^{\mathsf{pr}}_{\G^{\mathsf{o}}} (1,n) \ar[d] \ar[r] & \mathtt{I}^{\mathsf{pr}}_{\G_{\mathcal{R}}^{\mathsf{o}}} (1,n) \ar[d] \ar[r] & 0 \\
        0 \ar[r] & \mathrm{Ker} (q_n) \ar[d] \ar[r] & \mathds{k}[\mathsf{F}_n] \ar[d, "\T_{n}"] \ar[r, "q_n"] &   \mathds{k}[\mathsf{F}_n/ \mathcal{R}(\mathsf{F}_n)] \ar[d] \ar[r] & 0 \\
        0 \ar[r] & \K_\mathcal{R}(n) \ar[d] \ar[r] & \mathsf{Lin} (n) \ar[d] \ar[r] &  \Psi_{\G_{\mathcal{R}}^{\mathsf{o}}} \ar[d] \ar[r] & 0\\
        & 0& 0& 0& &
    \end{tikzcd}
    $$
    
    Next, we show that the first row is exact.
    Note that, since the canonical functor $\G^{\mathsf{o}} \to \G_{\mathcal{R}}^{\mathsf{o}}$ is full, the map $\mathtt{I}^{\mathsf{pr}}_{\G^{\mathsf{o}}} (1,n) \to \mathtt{I}^{\mathsf{pr}}_{\G_{\mathcal{R}}^{\mathsf{o}}} (1,n)$ is surjective by the definition of the primitivity ideals.
    If this map sends $v \in \mathtt{I}^{\mathsf{pr}}_{\G^{\mathsf{o}}} (1,n)$ to zero, then by the commutative diagram, there exists $w \in \mathrm{Ker}(q_n)$ such that $v = w \in \mathds{k}[\mathsf{F}_n]$, which implies $v \in W$.
    It is trivial that $W \to \mathtt{I}^{\mathsf{pr}}_{\G^{\mathsf{o}}} (1,n)$ is injective.
    
    As a consequence, by the nine-lemma, the first column is exact, which proves the statement.
\end{proof}

\subsection{The case of a field of characteristic zero}

In this section, we improve the adjunction given in (\ref{202604051252}) by using some results of Powell \cite{powell2024analytic}.
For this purpose, we assume, only in this section, that the ground ring $\mathds{k}$ is a field of characteristic zero.

The following theorem is a consequence of Powell's equivalence theorem:
\begin{theorem} \label{202512082006}
    We assume that the ground ring $\mathds{k}$ is a field of characteristic zero.
    The refined eigenmonad adjunction associated with $(\mathtt{L}_{\G^{\mathsf{o}}}, \mathtt{I}^{\mathsf{pr}}_{\G^{\mathsf{o}}})$, given in Proposition \ref{202402011440}, induces the following adjoint equivalence:
    $$
    \begin{tikzcd}
             \tilde{\mu}\mathsf{Lin} \otimes_{\mu\mathfrak{Lie}} (-) : \mu\mathfrak{Lie}\mbox{-}\mathsf{Mod} \arrow[r, shift right=1ex, ""{name=G}] & \mathtt{L}_{\G^{\mathsf{o}}}\mbox{-}\mathsf{Mod}^{\mathsf{prim}} : \mathrm{Hom}_{\mathtt{L}_{\G^{\mathsf{o}}}} ( \tilde{\mu}\mathsf{Lin} , - ) \arrow[l, shift right=1ex, ""{name=F}]
            \arrow[phantom, from=G, to=F, , "\scriptscriptstyle\boldsymbol{\top}"] .
    \end{tikzcd}
    $$
    Furthermore, the following statements hold:
    \begin{enumerate}
        \item A left $\mathtt{L}_{\G^{\mathsf{o}}}$-module is primitively generated if and only if it is analytic.
        \item Every left $\mu\mathfrak{Lie}$-module is $\mathtt{I}^{\mathsf{pr}}_{\G^{\mathsf{o}}}$-vanishingly extensible.
    \end{enumerate}
\end{theorem}
\begin{proof}
    Recall that by Theorem \ref{202603121700}, the eigenmonad adjunction associated with $(\mathtt{L}_{\G^{\mathsf{o}}}, \mathtt{I}^{\mathsf{pr}}_{\G^{\mathsf{o}}})$ yields Powell's adjunction.
    In \cite[Theorem 9.19]{powell2024analytic}, it is established that, via the adjunction, the category of $\mu\mathfrak{Lie}$-modules is equivalent to that of analytic $\mathtt{L}_{\G^{\mathsf{o}}}$-modules.
    By Theorem \ref{202603121700}, the image of the left adjoint consists of primitively generated $\mathtt{L}_{\G^{\mathsf{o}}}$-modules.
    (1) follows from these arguments.
    (2) is proved from the natural isomorphism $\mathbb{R} \circ \mathbb{L} \cong \mathrm{Id}$ that is deduced from the Powell's equivalence result.
\end{proof}

\begin{remark}
    The following results heavily depend on Theorem \ref{202512082006}.
    Powell used some arguments related to Koszul duality to prove the equivalence.
    In our sequel paper, we will prove a general statement using a structure that arises in the monad $\mathtt{L}_{\G^{\mathsf{o}}}$ without invoking Koszul duality.
\end{remark}

Let $\mathcal{R}$ be a radical functor for groups.
In the following, we establish that, for the Lawvere theory $\G_{\mathcal{R}}^{\mathsf{o}}$, the operad $\mathfrak{O}_{\G_{\mathcal{R}}^{\mathsf{o}}}$, introduced in Definition \ref{202603241701}, plays the role of $\mathfrak{Lie}$ in Theorem \ref{202512082006}.

\begin{Lemma} \label{202512091030}
    We have an isomorphism of $(\mathtt{L}_{\G^{\mathsf{o}}}, \Phi_{\G_{\mathcal{R}}^{\mathsf{o}}})$-bimodules:
    $$
    \tilde{\mu}\mathsf{Lin}\otimes_{\mu\mathfrak{Lie}} \Phi_{\G_{\mathcal{R}}^{\mathsf{o}}} \cong \tilde{\mu}\Psi_{\G_{\mathcal{R}}^{\mathsf{o}}}.
    $$
\end{Lemma}
\begin{proof}
    By Lemma \ref{202509201154}, for $d \in \mathds{N}$, the $\mathtt{L}_{\G^{\mathsf{o}}}$-module $(\mathtt{L}_{\G_{\mathcal{R}}^{\mathsf{o}}} /\mathtt{I}^{\mathsf{pr}}_{\G_{\mathcal{R}}^{\mathsf{o}}})(-,d)$ is polynomial, so analytic.
    By definition, the $\mathtt{L}_{\G^{\mathsf{o}}}$-module $\tilde{\mu}\Psi_{\G_{\mathcal{R}}^{\mathsf{o}}}(-,d)$ satisfies the same property.
    Hence, the counit of the adjunction given in Theorem \ref{202512082006} induces the isomorphism
    $$
     \tilde{\mu}\mathsf{Lin}\otimes_{\mu\mathfrak{Lie}} \mathrm{V}(\tilde{\mu}\Psi_{\G_{\mathcal{R}}^{\mathsf{o}}} (-,d) ; \mathtt{I}^{\mathsf{pr}}_{\G_{\mathcal{R}}^{\mathsf{o}}}) \stackrel{\cong}{\longrightarrow} \tilde{\mu}\Psi_{\G_{\mathcal{R}}^{\mathsf{o}}} (-,d).
    $$
    Therefore, the lemma follows from the equivalent descriptions of the eigenmonad $\Phi_{\G_{\mathcal{R}}^{\mathsf{o}}}$, given in Proposition \ref{202401241330}.
\end{proof}

\begin{theorem} \label{202509201759}
    We assume that the ground ring $\mathds{k}$ is a field of characteristic zero.
    Then the refined eigenmonad adjunction associated with $(\mathtt{L}_{\G_{\mathcal{R}}^{\mathsf{o}}}, \mathtt{I}^{\mathsf{pr}}_{\G_{\mathcal{R}}^{\mathsf{o}}})$, given in Proposition \ref{202402011440}, leads to the following adjoint equivalence:
    $$
    \begin{tikzcd}
            \tilde{\mu}\Psi_{\G_{\mathcal{R}}^{\mathsf{o}}} \otimes_{\mu\mathfrak{O}_{\G_{\mathcal{R}}^{\mathsf{o}}}} (-) : \mu\mathfrak{O}_{\G_{\mathcal{R}}^{\mathsf{o}}}\mbox{-}\mathsf{Mod} \arrow[r, shift right=1ex, ""{name=G}] & \mathtt{L}_{\G_{\mathcal{R}}^{\mathsf{o}}}\mbox{-}\mathsf{Mod}^{\mathsf{prim}} :  \mathrm{Hom}_{\mathtt{L}_{\G_{\mathcal{R}}^{\mathsf{o}}}} (\tilde{\mu}\Psi_{\G_{\mathcal{R}}^{\mathsf{o}}} , -) \arrow[l, shift right=1ex, ""{name=F}]
            \arrow[phantom, from=G, to=F, , "\scriptscriptstyle\boldsymbol{\top}"].
    \end{tikzcd}
    $$
    Furthermore, the following statements hold:
    \begin{enumerate}
        \item A left $\mathtt{L}_{\G_{\mathcal{R}}^{\mathsf{o}}}$-module is primitively generated if and only if it is analytic.
        \item Every left $\mu\mathfrak{O}_{\G_{\mathcal{R}}^{\mathsf{o}}}$-module is $\mathtt{I}^{\mathsf{pr}}_{\G_{\mathcal{R}}^{\mathsf{o}}}$-vanishingly extensible.
        \item The equivalence induces $\mu\mathfrak{O}_{\G_{\mathcal{R}}^{\mathsf{o}}}\mbox{-}\mathsf{Mod}^{d-\mathsf{trun}} \simeq \mathtt{L}_{\G_{\mathcal{R}}^{\mathsf{o}}}\mbox{-}\mathsf{Mod}^{\leq d}$.
    \end{enumerate}
\end{theorem}
\begin{proof}
    Since $\mathds{k}$ is a field, by Corollary \ref{202604051251}, $\mu\mathfrak{O}_{\G_{\mathcal{R}}^{\mathsf{o}}} \cong \Phi_{\G_{\mathcal{R}}^{\mathsf{o}}}$.
    Under this isomorphism, the adjunction of the statement arises from the adjunction given in (\ref{202604051252}).
    Let $\mathbb{L},\mathbb{R}$ be the left and right adjoint functors.
    We prove that these functors yield an equivalence of categories using Theorem \ref{202512082006}.
    We put $\mathtt{S} =\Phi_{\G_{\mathcal{R}}^{\mathsf{o}}}$.
    Let $\pi : \mathtt{L}_{\G^{\mathsf{o}}} \to \mathtt{L}_{\G_{\mathcal{R}}^{\mathsf{o}}}$ be the canonical epimorphism, and let $\phi: \mu\mathfrak{Lie} \to \mathtt{S}$ be the monad map obtained from Proposition \ref{202512071739}.
    Consider a primitively generated left $\mathtt{L}_{\G_{\mathcal{R}}^{\mathsf{o}}}$-module $\mathtt{M}$.
    By definition, we have $\pi (\mathtt{I}^{\mathsf{pr}}_{\G^{\mathsf{o}}}) = \mathtt{I}^{\mathsf{pr}}_{\G_{\mathcal{R}}^{\mathsf{o}}}$, so we have $\mathrm{V}( \mathtt{M}; \mathtt{I}^{\mathsf{pr}}_{\G_{\mathcal{R}}^{\mathsf{o}}}) \cong \mathrm{V}( \pi^\ast (\mathtt{M}); \mathtt{I}^{\mathsf{pr}}_{\G^{\mathsf{o}}})$.
    In particular, $\pi^\ast (\mathtt{M})$ is primitively generated, since the following diagram commutes and the lower row is assumed to be an epimorphism:
    $$
    \begin{tikzcd}
        \mathtt{L}_{\G^{\mathsf{o}}} \otimes \mathrm{V}( \pi^\ast (\mathtt{M}); \mathtt{I}^{\mathsf{pr}}_{\G^{\mathsf{o}}}) \ar[r] \ar[d, twoheadrightarrow] & \pi^\ast (\mathtt{M})  \ar[d, equal] \\
        \mathtt{L}_{\G_{\mathcal{R}}^{\mathsf{o}}} \otimes \mathrm{V}( \mathtt{M};   \mathtt{I}^{\mathsf{pr}}_{\G_{\mathcal{R}}^{\mathsf{o}}}) \ar[r] & \mathtt{M} .
    \end{tikzcd}
    $$
    Using Lemma \ref{202512091030}, we obtain 
    \begin{align*}
        \mathbb{L}(\mathbb{R}(\mathtt{M})) \cong \tilde{\mu}\Psi_{\G_{\mathcal{R}}^{\mathsf{o}}} \otimes_{\mathtt{S}} \mathrm{V}(\mathtt{M}; \mathtt{I}^{\mathsf{pr}}_{\G_{\mathcal{R}}^{\mathsf{o}}}) \cong (\tilde{\mu}\mathsf{Lin} \otimes_{\mu\mathfrak{Lie}} \mathtt{S})\otimes_{\mathtt{S}} \mathrm{V}(\mathtt{M}; \mathtt{I}^{\mathsf{pr}}_{\G_{\mathcal{R}}^{\mathsf{o}}}) \cong \tilde{\mu}\mathsf{Lin} \otimes_{\mu\mathfrak{Lie}} \mathrm{V}(\pi^\ast\mathtt{M}; \mathtt{I}^{\mathsf{pr}}_{\G^{\mathsf{o}}}) \stackrel{\mathrm{Thm}\ref{202512082006}}{\cong} \mathtt{M}.
    \end{align*}
    
    Let $\mathtt{N}$ be a left  $\mathtt{S}$-module.
    We then have $\mathbb{R}(\mathbb{L}(\mathtt{N})) \cong \mathrm{V} ( \tilde{\mu}\Psi_{\G_{\mathcal{R}}^{\mathsf{o}}} \otimes_{\mathtt{S}} \mathtt{N} ; \mathtt{I}^{\mathsf{pr}}_{\G_{\mathcal{R}}^{\mathsf{o}}})$, and
    \begin{align*}
        \mathrm{V} ( \tilde{\mu}\Psi_{\G_{\mathcal{R}}^{\mathsf{o}}} \otimes_{\mathtt{S}} \mathtt{N} ; \mathtt{I}^{\mathsf{pr}}_{\G_{\mathcal{R}}^{\mathsf{o}}}) \stackrel{\mathrm{Lem}\ref{202512091030}}{\cong} \mathrm{V} ((\tilde{\mu}\mathsf{Lin} \otimes_{\mu\mathfrak{Lie}} \mathtt{S})\otimes_{\mathtt{S}} \mathtt{N} ; \mathtt{I}^{\mathsf{pr}}_{\G^{\mathsf{o}}}) \cong \mathrm{V} (\tilde{\mu}\mathsf{Lin} \otimes_{\mu\mathfrak{Lie}} \phi^\ast (\mathtt{N}) ; \mathtt{I}^{\mathsf{pr}}_{\G^{\mathsf{o}}}) \stackrel{\mathrm{Thm}\ref{202512082006}}{\cong} \mathtt{N} .
    \end{align*}
    
    The statements (1) and (2) are proved using the method analogous to Theorem \ref{202512082006}.
    The statement (3) follows from Theorem \ref{202509201658} together with the above results.
\end{proof}

\begin{Corollary} \label{202512181438}
    Let $d \in \mathds{N}$.
    We have an equivalence of categories:
    \begin{align*}
        \mathtt{L}_{\G_{\mathcal{R}}^{\mathsf{o}}}\mbox{-}\mathsf{Mod}^{\leq d} / \mathtt{L}_{\G_{\mathcal{R}}^{\mathsf{o}}}\mbox{-}\mathsf{Mod}^{\leq (d-1)} \simeq
        \begin{cases}
            \mathds{k}[\mathfrak{S}_d]\mbox{-}\mathsf{Mod} & \mathrm{if~}|\mathcal{R}| = 0 , \\
            \ast & \mathrm{otherwise}.
        \end{cases}
    \end{align*}
\end{Corollary}
\begin{proof}
    By (3) of Theorem \ref{202509201759}, the quotient category in the statement is equivalent to $$\Phi_{\G_{\mathcal{R}}^{\mathsf{o}}}\mbox{-}\mathsf{Mod}^{d-\mathsf{trun}}/\Phi_{\G_{\mathcal{R}}^{\mathsf{o}}}\mbox{-}\mathsf{Mod}^{(d-1)-\mathsf{trun}}.$$
    Thus, the corollary follows from Proposition \ref{202512081028}, as we have $\Phi_{\G_{\mathcal{R}}^{\mathsf{o}}} (d,d) \cong (\mathds{k}/r\mathds{k})[\mathfrak{S}_d]$ by Corollary \ref{202512102031}
\end{proof}

\begin{remark}
    The corollary refines the result of \cite{kim2025poly} stating that, for each $d \in \mathds{N}$, there {\it exists} a $\G_{\mathcal{R}}^{\mathsf{o}}$-module of polynomial degree $d$ when $|\mathcal{R}| = 0$.
\end{remark}

\subsection{The primitivity ideal for free nilpotent groups of class $\leq c$} \label{202512011517}

In this section, we specialize to $\mathcal{R} = \gamma_{c+1}$ the $(c+1)$-th component of the lower central series.
Recall that $\G_{\gamma_{c+1}}$ coincides with the category $\mathbf{N}_c$ of finitely generated free nilpotent groups of class $\leq c$.
We study the associated eigenmonad and the eigenmonad adjunction.

The ground ring $\mathds{k}$ is assumed to be an arbitrary unital commutative ring unless otherwise specified.

\begin{theorem} \label{202604021256}
    Let $c \in \mathds{N}^\ast$.
    \begin{enumerate}
        \item We have a monad isomorphism $\Phi_{\mathbf{N}_{c}^{\mathsf{o}}} \cong \mu\mathfrak{Lie}_c$.
        In particular, the operad $\mathfrak{O}_{\mathbf{N}_{c}^{\mathsf{o}}}$ is $\mathfrak{Lie}_c$.
        Hence, we can regard $\tilde{\mu}\Psi_{\mathbf{N}_c^{\mathsf{o}}}$(see Notation \ref{202603201036}) as an $(\mathtt{L}_{\mathbf{N}_{c}^{\mathsf{o}}}, \mu\mathfrak{Lie}_{c})$-bimodule.
        \item  There is an isomorphism of $(\mathtt{L}_{\G^{\mathsf{o}}}, \mu\mathfrak{Lie}_{c})$-bimodules:
        $$
        \tilde{\mu}\Psi_{\mathbf{N}_c^{\mathsf{o}}} \cong \tilde{\mu}\mathsf{Lin}\otimes_{\mu\mathfrak{Lie}} \mu\mathfrak{Lie}_{c} .
        $$
    \end{enumerate}
\end{theorem}

The proof is deferred to the end of this section.

\begin{remark}
    (2) of the theorem gives $$\Psi_{\mathbf{N}_c^{\mathsf{o}}} \cong \mathsf{Lin}\otimes_{\mu\mathfrak{Lie}} \mu\mathfrak{Lie}_{c} \cong \mathsf{Lin} / (\mathsf{Lin} \lhd \mu\mathfrak{Lie}_{c}),$$
    via which $\Psi_{\mathbf{N}_c^{\mathsf{o}}}$ can be viewed as a nilpotent truncation of $\mathsf{Lin}$.
\end{remark}

\begin{Corollary} \label{202512030930}
    The eigenmonad adjunction associated with $(\mathtt{L}_{\mathbf{N}_{c}^{\mathsf{o}}}, \mathtt{I}^{\mathsf{pr}}_{\mathbf{N}_{c}^{\mathsf{o}}})$ induces the following adjunction:
    \begin{equation} \label{202606081554}
    \begin{tikzcd}
             \tilde{\mu}\mathsf{Lin}\otimes_{\mu\mathfrak{Lie}} U(-)  : \mu\mathfrak{Lie}_{c}\mbox{-}\mathsf{Mod} \arrow[r, shift right=1ex, ""{name=G}] & \mathtt{L}_{\mathbf{N}_{c}^{\mathsf{o}}}\mbox{-}\mathsf{Mod} : \mathrm{Hom}_{\mathtt{L}_{\mathbf{N}_{c}^{\mathsf{o}}}} ( \tilde{\mu}\Psi_{\mathbf{N}_c^{\mathsf{o}}} , - ) \arrow[l, shift right=1ex, ""{name=F}]
            \arrow[phantom, from=G, to=F, , "\scriptscriptstyle\boldsymbol{\top}"],
    \end{tikzcd}
    \end{equation}
    where $U : \mu\mathfrak{Lie}_{c}\mbox{-}\mathsf{Mod} \to \mu\mathfrak{Lie}\mbox{-}\mathsf{Mod}$ denotes the functor induced by the canonical operad map $\mathfrak{Lie} \to \mathfrak{Lie}_{c}$.
    Furthermore, it fulfills the following conditions where $\mathbb{L}$ and $\mathbb{R}$ denote the left and right adjoint functors:
    \begin{enumerate}
        \item The image of $\mathbb{L}$ is contained in primitively generated $\mathtt{L}_{\mathbf{N}_{c}^{\mathsf{o}}}$-modules, in particular analytic $\mathtt{L}_{\mathbf{N}_{c}^{\mathsf{o}}}$-modules.
        \item For $d\in\mathds{N}$, the functors $\mathbb{L}$ and $\mathbb{R}$ give a correspondence between $d$-truncated $\mu\mathfrak{Lie}_{c}$-modules and $\mathtt{L}_{\mathbf{N}_{c}^{\mathsf{o}}}$-modules of polynomial degree at most $d$.
        \item For a $\mathtt{L}_{\mathbf{N}_{c}^{\mathsf{o}}}$-module $\mathtt{M}$ of polynomial degree at most $d$, we have $\mathbb{R}(\mathtt{M})(d) = (\mathrm{cr}_{d} (\mathtt{M})) (1,\cdots, 1)$.
    \end{enumerate}
\end{Corollary}
\begin{proof}
    The adjunction follows from Corollary \ref{202512081855} and Theorem \ref{202604021256}, where the left adjoint functor is calculated as follows:
    $$\tilde{\mu}\Psi_{\mathbf{N}_c^{\mathsf{o}}} \otimes_{\mu\mathfrak{Lie}_{c}} (-) \cong \left( \tilde{\mu}\mathsf{Lin}\otimes_{\mu\mathfrak{Lie}} \mu\mathfrak{Lie}_{c} \right) \otimes_{\mu\mathfrak{Lie}_{c}} (-) \cong \tilde{\mu}\mathsf{Lin}\otimes_{\mu\mathfrak{Lie}}U (-) .$$
    The part (1) follows from Proposition \ref{202402011440}.
    The assertions (2) and (3) follow from Theorem \ref{202509201658}.
\end{proof}

\begin{Corollary} \label{202604051720}
    If $\mathds{k}$ is a field of characteristic zero, then the adjunction given in (\ref{202606081554}) yields an adjoint equivalence between $\mu\mathfrak{Lie}_{c}$-modules and analytic $\mathtt{L}_{\mathbf{N}_{c}^{\mathsf{o}}}$-modules.
\end{Corollary}
\begin{proof}
    This follows from Theorem \ref{202509201759}.
\end{proof}

\begin{Example} \label{202509192025}
    For instance, consider $c=1$.
    Clearly, $\mathbf{N}_{1}$ equals the category $\mathbf{M}_{\mathds{Z}}$ of finitely generated free abelian groups.
    By Example \ref{202603100930}, $\mu\mathfrak{Lie}_1 \cong \mathtt{L}_{\mathfrak{S}}$.
    Furthermore, the $(\mathtt{L}_{\mathfrak{S}},\mathtt{L}_{\mathfrak{S}})$-bimodule isomorphism $\mu\mathsf{Lin}\otimes_{\mu\mathfrak{Lie}} \mathtt{L}_{\mathfrak{S}} \cong \mu\mathsf{Exp}$, with (2) of Theorem \ref{202604021256}, induces an $(\mathtt{L}_{\mathbf{M}_{\mathds{Z}}^{\mathsf{o}}} , \mathtt{L}_{\mathfrak{S}})$-bimodule structure on $\mu\mathsf{Exp}$, which we denote by $\tilde{\mu}\mathsf{Exp}$.
    Then the adjunction (\ref{202606081554}) yields the following.
    In Section \ref{202410161800}, we shall present a generalization of this construction.
    $$
    \begin{tikzcd}
             \tilde{\mu}\mathsf{Exp} \otimes_{\mathtt{L}_{\mathfrak{S}}} (-) : \mathtt{L}_{\mathfrak{S}}\mbox{-}\mathsf{Mod} \arrow[r, shift right=1ex, ""{name=G}] & \mathtt{L}_{\mathbf{M}_{\mathds{Z}}^{\mathsf{o}}}\mbox{-}\mathsf{Mod} : \mathrm{Hom}_{\mathtt{L}_{\mathbf{M}_{\mathds{Z}}^{\mathsf{o}}}} ( \tilde{\mu}\mathsf{Exp} , - ) \arrow[l, shift right=1ex, ""{name=F}]
            \arrow[phantom, from=G, to=F, , "\scriptscriptstyle\boldsymbol{\top}"].
    \end{tikzcd}
    $$
    When $\mathds{k}$ is a field of characteristic zero, Corollary \ref{202604051720} applied to this case recovers some well-known results in polynomial functor theory.
    The related arguments appear, for instance, in \cite[section 2.4]{vespa2022functors} or \cite[Remark 8.6]{powell2024analytic}.
\end{Example}

In the remaining part of this section, we prove Theorem \ref{202604021256}.
To this end, we give a refinement of Lemma \ref{202604011516}.
Recall $\K_{\gamma_{c+1}}$ from Definition \ref{202604021150}.
\begin{Lemma} \label{202604021122}
    There exists a PBW datum $(\{ \mathfrak{B}_X \mid X\subset \mathds{N}^\ast \}, \preceq )$ for the bimonoid $\mathsf{Lin}$, as in Lemma \ref{202604011516}, satisfying the following additional condition: for $c \in \mathds{N}^\ast$ and $X \subset \mathds{N}^\ast$, the $\mathds{k}$-submodule $\K_{\gamma_{c+1}} (X) \subset \mathsf{Lin}(X)$ is freely generated by $h_1 \cdots h_r,~r\in\mathds{N}^\ast$ where $h_j \in \mathfrak{B}_{X_j}$ and $X = \coprod^r_{j=1} X_j$ such that $h_1 \succ \cdots \succ h_r$ and $|X_j| > c$ for some $j$.
\end{Lemma}
The proof heavily depends on classical structure theorems for free Lie algebras.
We defer the proof to Appendix \ref{202512031342}, as they require some notions that are not used elsewhere in this paper.

\begin{Lemma} \label{202604111808}
    There exists a PBW datum $(\{ \mathfrak{B}^c_X \mid X\subset \mathds{N}^\ast \}, \preceq )$ for $\Psi_{\mathbf{N}_c^{\mathsf{o}}}$ such that the inclusion $\mathfrak{O}_{\mathbf{N}_c^{\mathsf{o}}} \hookrightarrow \Psi_{\mathbf{N}_c^{\mathsf{o}}}$ induces $\mathfrak{O}_{\mathbf{N}_c^{\mathsf{o}}} (X) \cong \mathds{k} [\mathfrak{B}^c_X]$ for $X \subset \mathds{N}^\ast$.
\end{Lemma}
\begin{proof}
    Let $X \subset \mathds{N}^\ast$.
    Consider the following short exact sequence in the category of $\mathds{k}$-modules:
    $$ 0 \to \K_{\gamma_{c+1}} (X) \to \mathsf{Lin}(X) \to \Psi_{\mathbf{N}_c^{\mathsf{o}}} (X) \to 0 .$$
    Note that this splits by Lemma \ref{202604021122}.
    In fact, we can choose a specific section $s_X : \Psi_{\mathbf{N}_c^{\mathsf{o}}} (X) \to \mathsf{Lin}(X)$ which assigns $h_1 \cdots h_r \in \mathsf{Lin}(X)$ to the class $[h_1 \cdots h_r] \in \mathsf{Lin} (X) /  \K_{\gamma_{c+1}} (X) \cong \Psi_{\mathbf{N}_c^{\mathsf{o}}} (X)$ where $X = \coprod^r_{j=1} X_j$ is a decomposition such that $|X_j| \leq c$, $h_j \in \mathfrak{B}_{X_j}$, and $h_1 \succ \cdots \succ h_r$. 
    Using the map $s_X$, we regard $\Psi_{\mathbf{N}_c^{\mathsf{o}}} (X)$ as a $\mathds{k}$-submodule of $\mathsf{Lin}(X)$.
    We now define 
    \begin{align*}
        \mathfrak{B}^c_{X} {:=}
        \begin{cases}
            \mathfrak{B}_{X} & \mathrm{if~} |X| \leq c, \\
            \emptyset & \mathrm{otherwise}.
        \end{cases}
    \end{align*}
    By the choice of $s_X$, the sets $\mathfrak{B}^c_{X}$ satisfy part of the conditions of a PBW datum, namely (2) of Definition \ref{202604021629}.
    Furthermore, $h \in \mathfrak{B}^c_X$ is primitive in $\Psi_{\mathbf{N}_c^{\mathsf{o}}}$, since it is already primitive in $\mathsf{Lin}$.
    Hence, $(\{ \mathfrak{B}^c_X \mid X \in \mathcal{P}_f (\mathds{N}^{\ast}) \}, \preceq )$ gives a PBW datum for $\Psi_{\mathbf{N}_c^{\mathsf{o}}}$.
    By the definition of $\mathfrak{B}^c_{X}$ and Proposition \ref{202604021602}, we obtain $\mathfrak{O}_{\mathbf{N}_c^{\mathsf{o}}} (X) = \mathrm{Pr} (\Psi_{\mathbf{N}_c^{\mathsf{o}}}) (X) \cong  \mathds{k} [\mathfrak{B}^c_X]$ for $X \subset \mathds{N}^\ast$.
\end{proof}

Consider the composition $\mathfrak{Lie} \cong \mathfrak{O}_{\G^{\mathsf{o}}} \to \mathfrak{O}_{\mathbf{N}_c^{\mathsf{o}}}$ where the isomorphism appears in Theorem \ref{202603121700} and the other map is obtained from Proposition \ref{202512071739}.

\begin{Lemma} \label{202604021826}
    The map $\mathfrak{Lie}\to \mathfrak{O}_{\mathbf{N}_c^{\mathsf{o}}}$ factors through $\mathfrak{Lie} \twoheadrightarrow \mathfrak{Lie}_c$ and gives an operad isomorphism $\mathfrak{Lie}_c \cong \mathfrak{O}_{\mathbf{N}_c^{\mathsf{o}}}$.
\end{Lemma}
\begin{proof}
    We use the result in Lemma \ref{202604111808}.
    For $n \in \mathds{N}$, using the PBW data for $\Psi_{\mathbf{N}_c^{\mathsf{o}}}$ and $\mathsf{Lin}$, we obtain the following commutative diagram:
    $$
    \begin{tikzcd}
         \mathds{k} [\mathfrak{B}_{\bf n}]\ar[d, "\cong"]  \ar[r] & \mathds{k} [\mathfrak{B}^c_{\bf n}] \ar[d, "\cong"]\\
        \mathfrak{Lie}(n) \ar[r] & \mathfrak{O}_{\mathbf{N}_c^{\mathsf{o}}} (n) .
    \end{tikzcd}
    $$
    Therefore, by the definition of $\mathfrak{B}^c_{\bf n}$, the kernel of $\mathfrak{Lie}(n) \to \mathfrak{O}_{\mathbf{N}_c^{\mathsf{o}}} (n)$ equals $\mathfrak{Lie}(n)$ if $n > c$; and vanishes otherwise.
    This yields the desired result.
\end{proof}

\begin{proof}[Proof of Theorem \ref{202604021256}]
    We first prove (1) of Theorem \ref{202604021256}.
    By Proposition \ref{202604021602} and Corollary \ref{202603201735}, the map $\mola_{\mathbf{N}_c^{\mathsf{o}}} : \mu\mathfrak{O}_{\mathbf{N}_c^{\mathsf{o}}} \stackrel{\cong}{\to} \Phi_{\mathbf{N}_c^{\mathsf{o}}}$ is an isomorphism.
    Hence, $\Phi_{\mathbf{N}_c^{\mathsf{o}}} \cong \mu \mathfrak{Lie}_{c}$ by Lemma \ref{202604021826}.

    Next, we prove (2) of Theorem \ref{202604021256}.
    Let $\pi: \tilde{\mu}\mathsf{Lin} \twoheadrightarrow \tilde{\mu}\Psi_{\mathbf{N}_c^{\mathsf{o}}}$ be the quotient map.
    Considering the $(\mathtt{L}_{\G^{\mathsf{o}}}, \mu\mathfrak{Lie})$-bimodule $\tilde{\mu}\Psi_{\mathbf{N}_c^{\mathsf{o}}}$, the map $\pi$ is a $(\mathtt{L}_{\G^{\mathsf{o}}}, \mu\mathfrak{Lie})$-bimodule map.
    Hence, the composition of $\tilde{\mu}\mathsf{Lin}\otimes \mu\mathfrak{Lie}_{c} \stackrel{\pi\otimes\mathrm{id}}{\to}   \tilde{\mu}\Psi_{\mathbf{N}_c^{\mathsf{o}}} \otimes \mu\mathfrak{Lie}_{c}\stackrel{\lhd}{\to} \tilde{\mu}\Psi_{\mathbf{N}_c^{\mathsf{o}}}$ factors through the identification map $\tilde{\mu}\mathsf{Lin}\otimes \mu\mathfrak{Lie}_{c} \twoheadrightarrow \tilde{\mu}\mathsf{Lin}\otimes_{\mu\mathfrak{Lie}} \mu\mathfrak{Lie}_{c}$.
    So, we obtain a canonical $(\mathtt{L}_{\G^{\mathsf{o}}},\mu\mathfrak{Lie}_{c} )$-bimodule map
    $$\alpha : \tilde{\mu}\mathsf{Lin}\otimes_{\mu\mathfrak{Lie}} \mu\mathfrak{Lie}_{c} \to \tilde{\mu}\Psi_{\mathbf{N}_c^{\mathsf{o}}}.$$
    
    We construct the inverse of $\alpha$ as follows.
    Let $\beta^\prime : \tilde{\mu}\mathsf{Lin} \to \tilde{\mu}\mathsf{Lin}\otimes_{\mu\mathfrak{Lie}} \mu\mathfrak{Lie}_{c}$ be the following composition:
    $$
    \tilde{\mu}\mathsf{Lin} \cong \tilde{\mu}\mathsf{Lin} \otimes \mathds{I}_{\mathds{N}} \to \tilde{\mu}\mathsf{Lin} \otimes \mu\mathfrak{Lie}_c \twoheadrightarrow \tilde{\mu}\mathsf{Lin}\otimes_{\mu\mathfrak{Lie}} \mu\mathfrak{Lie}_{c} .
    $$
    We first show that this map factors through $\pi: \tilde{\mu}\mathsf{Lin} \twoheadrightarrow \tilde{\mu}\Psi_{\mathbf{N}_c^{\mathsf{o}}}$:
    $$
    \begin{tikzcd}
         \tilde{\mu}\mathsf{Lin} \ar[d, twoheadrightarrow, "\pi"] \ar[r, "\beta^\prime"] & \tilde{\mu}\mathsf{Lin}\otimes_{\mu\mathfrak{Lie}} \mu\mathfrak{Lie}_{c} \\
         \tilde{\mu}\Psi_{\mathbf{N}_c^{\mathsf{o}}} \ar[ur, "\exists \beta"'] & 
    \end{tikzcd}
    $$
    Let $\mathtt{K}$ be the kernel of the epimorphism $\tilde{\mu}\mathsf{Lin} \twoheadrightarrow \tilde{\mu}\Psi_{\mathbf{N}_c^{\mathsf{o}}}$.
    For $m,n\in\mathds{N}$, we can identify $\mathtt{K}(m,n)$ with the kernel of $\mathsf{Lin}^{\odot m} (n) \to (\Psi_{\mathbf{N}_c^{\mathsf{o}}})^{\odot m} (n)$ (see Lemma \ref{202603211811}).
    Recall that $\K_{\gamma_{c+1}}$ is defined to be the kernel of the map $\mathsf{Lin} \twoheadrightarrow \Psi_{\mathbf{N}_c^{\mathsf{o}}}$.
    So, the $\mathds{k}$-module $\mathtt{K} (m,n)$ is generated by $\bigotimes^m_{i=1} w_i \in \bigotimes^{m}_{i=1} \mathsf{Lin}(X_i) \subset  \mathsf{Lin}^{\odot m} (n)$ where ${\bf n} = \coprod^m_{i=1} X_i$, $w_i \in \mathsf{Lin}(X_i)$ and there exists $j_0\in {\bf m}$ such that $w_{j_0} \in \K_{\gamma_{c+1}}(X_{j_0})$.
    It suffices to show that the map $\beta^\prime$ sends such generators to zero.
    We can assume that $j_0 = 1$, since $\beta^\prime$ is compatible with the left $\mathtt{L}_\mathfrak{S}$-action.
    By Lemma \ref{202604021122}, we can also assume that $w_1 = h_1 \cdots h_r$ where $X_1 = \coprod^r_{i=1} Y_i$ is a decomposition such that $|Y_{j_1}| > c$ for some $j_1$, $h_i \in \mathfrak{B}_{Y_i}$ and $h_1 \succ \cdots \succ h_r$. 
    Then we have
    \begin{align*}
        &\bigotimes^m_{i=1} w_i = (h_1 \cdots h_r) \otimes  \bigotimes^m_{i=2} w_i ,\\
        =& \left( (h_1 \cdots h_{j_1-1} \theta_{n+1} h_{j_1+1} \cdots h_r) \otimes  \bigotimes^m_{i=2} w_i \right) \lhd ( \theta_1 \otimes \cdots \otimes \theta_n \otimes h_{j_1}).
    \end{align*}
    where $\lhd$ denotes the right $\mu\mathfrak{Lie}$-action on $\tilde{\mu}\mathsf{Lin}$.
    Hence, by the condition $|Y_{j_1}| > c$, we obtain
    \begin{align*}
        & \beta (\bigotimes^m_{i=1} w_i ) =[\bigotimes^m_{i=1} w_i  \otimes 1_n] ,\\
        =& [( h_1 \cdots h_{j_1-1} \theta_{n+1} h_{j_1+1} \cdots h_r \otimes  \bigotimes^m_{i=2} w_i ) \otimes ( \theta_1 \otimes \cdots \otimes \theta_n \otimes h_{j_1})]  = 0 .
    \end{align*}
    Here, $[-]$ denotes the class of $(\tilde{\mu}\mathsf{Lin}\otimes_{\mu\mathfrak{Lie}} \mu\mathfrak{Lie}_{c}) (m,n)$ induced by an element of $(\tilde{\mu}\mathsf{Lin}\otimes \mu\mathfrak{Lie}_{c} )(m,n)$.
    
    Let $\beta : \tilde{\mu}\Psi_{\mathbf{N}_c^{\mathsf{o}}}\to (\tilde{\mu}\mathsf{Lin}\otimes_{\mu\mathfrak{Lie}} \mu\mathfrak{Lie}_{c})$ the map obtained from $\beta^\prime$, as discussed above.
    It is immediate from the definitions that the following diagram commutes, which leads to the proof of $\alpha \circ \beta = \mathrm{id}$:
     $$
     \begin{tikzcd}
        \tilde{\mu}\mathsf{Lin} \ar[r, "\beta^\prime"] \ar[d, twoheadrightarrow, "\pi"] & \tilde{\mu}\mathsf{Lin}\otimes_{\mu\mathfrak{Lie}} \mu\mathfrak{Lie}_{c} \ar[d, "\alpha"] \\ 
         \tilde{\mu}\Psi_{\mathbf{N}_c^{\mathsf{o}}} \ar[ur, "\beta"] \ar[r, equal] & \tilde{\mu}\Psi_{\mathbf{N}_c^{\mathsf{o}}} .
     \end{tikzcd}
     $$
     
     To prove $\beta \circ \alpha = \mathrm{id}$, consider $$u \otimes v\in \tilde{\mu}\mathsf{Lin} (m,l) \otimes \mu\mathfrak{Lie}_{c}(l,n) \subset (\tilde{\mu}\mathsf{Lin}\otimes \mu\mathfrak{Lie}_{c}) (m,n) ,$$ where $u \in \tilde{\mu}\mathsf{Lin} (m,l)$ and $v \in \mu\mathfrak{Lie}_{c}(l,n)$.
     It suffices to prove $(\beta \circ \alpha) ( [u \otimes v]) = [u \otimes v]$.
     Let $q : \mu\mathfrak{Lie} \to \mu\mathfrak{Lie}_c$ be the epimorphism and take $v^\prime \in \mu\mathfrak{Lie} (l,n)$ such that $q (v^\prime) = v$.
     By the definitions of $\alpha$ and $\beta$, we obtain $\beta ( \alpha ( [u \otimes v])) = \beta (\pi (u) \lhd v) = \beta ( \pi (u) \lhd q (v^\prime))$.
     Using $\pi (u) \lhd q (v^\prime) = \pi (u \lhd v^\prime)$ (see Proposition \ref{202512071739}), we obtain $$\beta ( \alpha ( [u \otimes v])) = \beta ( \pi (u \lhd v^\prime)) = [(u \lhd v^\prime) \otimes 1_n] = [u \otimes q(v^\prime)] = [u \otimes v].$$
\end{proof}

\section{The primitivity ideal for free $R$-modules}
\label{202410161800}

Let $R$ be a unital ring which is not necessarily commutative.
Recall from Example \ref{202606091733} the Lawvere theory $\mathbf{M}_{R}$ of finitely generated free right $R$-modules.
In this section, we investigate the primitivity ideal $\mathtt{I}^{\mathsf{pr}}_{\mathcal{C}}$ for $\mathcal{C} = \mathbf{M}_{R}$, and the associated canonical bimodule and the eigenmonad, based on a discussion analogous to the previous sections.
As a consequence, we obtain Theorem \ref{202512111638}.

For convenience in what follows, we adopt the following (nonstandard) notation:
\begin{notation}\label{202511271410}
    We write $\mathfrak{e}^{X} \in \mathbf{M}_{R} (n,m)$ for the morphism corresponding to $X \in \mathrm{M}_{m,n} (R)$ where $\mathfrak{e}$ is merely a formal symbol.
    In spite of the familiarity, we attempt to avoid the matrix expression in order that we distinguish sums in the $\mathds{k}$-module $\mathtt{L}_{\mathbf{M}_{R}} (m,n) = \mathds{k} [\mathbf{M}_{R} (n,m)]$ from that of matrices.
\end{notation}

Using this notation, we describe several well-known structures on $\mathbf{M}_{R}$.
The category $\mathbf{M}_{R}$ has a binary product given by the direct sum of $R$-modules.
For matrices $X,Y$, we have $\mathfrak{e}^{X} \oplus \mathfrak{e}^{Y} = \mathfrak{e}^{X \oplus Y}$.
Furthermore, the composition of $\mathbf{M}_{R}$ is described by 
\begin{align} \label{202408081447}
    \mathfrak{e}^{Y} \circ \mathfrak{e}^{X} = \mathfrak{e}^{YX} , \quad X \in \mathrm{M}_{m,n}(R) ,  Y \in \mathrm{M}_{l,m} (R).
\end{align}

The category $\mathbf{M}_{R}$ is a Lawvere theory satisfying the condition (ZCM*).
See Section \ref{202509031755} for notations.
More precisely, for $\mathcal{C} = \mathbf{M}_{R}$, we regard $\mathcal{C}_n = \mathbf{M}_{R} (n,1)$ as the commutative monoid by
\begin{align} \label{202509221012}
    \mathfrak{e}^{X} \star \mathfrak{e}^{Y} = \mathfrak{e}^{X+Y}, \quad X,Y \in \mathrm{M}_{1,n}(R) .
\end{align}

\subsection{The algebra $\mathfrak{O}_{\mathbf{M}_{R}} (1)$}

We recall from Definition \ref{202603241701} the operad $\mathfrak{O}_{\mathbf{M}_{R}}$.
In this section, we investigate the $\mathds{k}$-algebra $\mathfrak{O}_{\mathbf{M}_{R}} (1)$.

\begin{notation}
    Let $R_{\mathds{k}} {:=} \mathds{k} \otimes_{\mathds{Z}} R$.
    We also use the notation $\lambda r = \lambda \otimes r$ for $r \in R$ and $\lambda \in \mathds{k}$.
\end{notation}

\begin{prop} \label{202509221014}
    We have an isomorphism $\mathfrak{O}_{\mathbf{M}_{R}} (1) \cong R_{\mathds{k}}$ of unital $\mathds{k}$-algebras.
\end{prop}
\begin{proof}
    By Lemma \ref{202601081632}, we have $\mathfrak{O}_{\mathbf{M}_{R}} (1) = \Phi_{\mathbf{M}_{R}}(1,1) \cong \mathtt{L}_{\mathbf{M}_{R}}/\mathtt{I}^{\mathsf{pr}}_{\mathbf{M}_{R}}(1,1)$.
    Clearly, the assignment $\mathfrak{e}^{r}$ to $r \in R$ gives $\mathtt{L}_{\mathbf{M}_{R}}(1,1) \cong \mathds{k}[R]$.
    Recall that $\mathtt{I}^{\mathsf{pr}}_{\mathbf{M}_{R}}(1,1)$ is generated by $f \circ \bar{\Delta}$ for $f \in \mathtt{L}_{\mathbf{M}_{R}} (1,2)$ where $\bar{\Delta} \in \mathtt{L}_{\mathbf{M}_{R}} (2,1)$ is the reduced comultiplication.
    Using Notation \ref{202511271410}, one may verify that $\bar{\Delta} = \mathfrak{e}^{[1,1]^{\mathrm{t}}} -\mathfrak{e}^{[1,0]^{\mathrm{t}}} -\mathfrak{e}^{[0,1]^{\mathrm{t}}}$, where $\mathrm{t}$ denotes the transposition of matrices.
    Note that $f$ is expressed as a linear combination of $\mathfrak{e}^{[r,r^\prime]}$ for $r,r^\prime \in R$.
    By the formula in (\ref{202408081447}), we obtain 
    $$
    \mathfrak{e}^{[r,r^\prime]} \circ \bar{\Delta} = \mathfrak{e}^{r+r^\prime} - \mathfrak{e}^r - \mathfrak{e}^{r^\prime} , \quad \{r,r^\prime\} \subset R.
    $$
    So, $\mathfrak{e}^{r+r^\prime} \equiv \mathfrak{e}^r + \mathfrak{e}^{r^\prime} \mod{\mathtt{I}^{\mathsf{pr}}_{\mathbf{M}_{R}}}$.
    By this relation, we obtain a $\mathds{k}$-linear map $R_{\mathds{k}} \to \Phi_{\mathbf{M}_{R}}  (1,1) ; r \mapsto \mathfrak{e}^r$.
    The inverse is given by $\Phi_{\mathbf{M}_{R}}  (1,1) \to R_{\mathds{k}} ; \mathfrak{e}^r \mapsto r$.
    Moreover, by (\ref{202408081447}), we also have $\mathfrak{e}^{r^\prime}\circ \mathfrak{e}^{r} = \mathfrak{e}^{r^\prime r}$, so the previous construction gives a $\mathds{k}$-algebra isomorphism.
\end{proof}

We give a simple application:
\begin{Corollary}
    We have $\mathtt{I}^{\mathsf{pr}}_{\mathbf{M}_{R}} = \mathtt{L}_{\mathbf{M}_{R}}$ if and only if $R_{\mathds{k}} \cong 0$.
\end{Corollary}
\begin{proof}
    This follows from Corollary \ref{202511270956}.
\end{proof}

\subsection{The canonical bimodule and the eigenmonad for $\mathbf{M}_{R}$}
\label{202410041316}

In this section, we give some structural results for the canonical bimodule and the eigenmonad associated with the primitivity ideal $\mathtt{I}^{\mathsf{pr}}_{\mathbf{M}_{R}}$.

Recall from Notation \ref{202604031739} the map $\E$, applied to $\mathcal{C} = \mathbf{M}_{R}$.
Since the Lawvere theory $\mathbf{M}_{R}$ satisfies the condition (ZCM*), by Proposition \ref{202509221014} and Corollaries \ref{202509191320} and \ref{202512051235}, we obtain a bimonoid epimorphism:
\begin{align} \label{202604061308}
    \E= \E_1 :  R_{\mathds{k}}^{\otimes} \to \Psi_{\mathbf{M}_{R}} .
\end{align}

\begin{Lemma} \label{202403191752}
    Let $n\in \mathds{N}$.
    For $r_1,\cdots, r_n \in R$,
    we have $\E \left( r_1 \otimes \cdots \otimes r_n  \right) = 
    [ \mathfrak{e}^{[r_1, \cdots, r_n]}]$.
\end{Lemma}
\begin{proof}
    By the definition of $\E :  R_{\mathds{k}}^{\otimes} \to \Psi_{\mathbf{M}_{R}}$, we have $\E ( r_1 \otimes \cdots \otimes r_n ) =  [\prod^{m}_{k=1} \mathfrak{e}^{r_ke_k} ]$
    where $\prod$ denotes the successive product with respect to the monoid structure $\star$ on $\mathbf{M}_{R}(n,1)$ (explained in Notation \ref{202511271410}), and $e_k {:=} [0, \cdots, 0,1,0,\cdots,0] \in \mathrm{M}_{1,n}(R)$ has $1 \in R$ in the $k$-th entry.
    Using (\ref{202509221012}), we obtain $\E ( r_1 \otimes \cdots r_n )= [\mathfrak{e}^{\sum^m_{k=1} r_ke_k} ] = [ \mathfrak{e}^{[r_1, \cdots, r_n]}]$.
\end{proof}

\begin{remark}
    Note that the assignment of $[r_1, \cdots, r_n]$ to $(r_1, \cdots, r_n) \in R^{\times n}$ itself is not bilinear; but the above lemma implies that, embedding matrices into $\mathtt{L}_{\mathbf{M}_{R}}$ by using $\mathfrak{e}^{(-)}$, it becomes bilinear modulo the primitivity ideal.
\end{remark}

\begin{Defn}
    We define $\E^{-1} : \Psi_{\mathbf{M}_{R}} \to  R_{\mathds{k}}^{\otimes}$ a map of right $\mathtt{L}_\mathfrak{S}$-modules.
    For $n \in \mathds{N}$, we define $\E^{-1} : \Psi_{\mathbf{M}_{R}} (n) \to  R_{\mathds{k}}^{\otimes} (n) =  R_{\mathds{k}}^{\otimes n}$ by the $\mathds{k}$-linear extension of
    $$
    \E^{-1} ( \mathfrak{e}^{[r_1, \cdots, r_n]}) {:=} r_1 \otimes \cdots \otimes r_n .
    $$
\end{Defn}

\begin{Lemma} \label{202408081425}
    The map $\E^{-1}$ is well-defined, and it gives the inverse of $\E :  R_{\mathds{k}}^{\otimes} \to \Psi_{\mathbf{M}_{R}}$.
\end{Lemma}
\begin{proof}
    To prove the well-definedness, we will show that $\E^{-1} ( \mathfrak{e}^{X} \circ  \bar{\Delta}_{n,i} ) = 0$ for $i \in {\bf n}$ and $X \in \mathrm{M}_{1,n+1}(R)$.
    By applying (\ref{202408081447}), it is equivalent to
    $$
    \E^{-1} \left(\mathfrak{e}^{(I_{i-1}\oplus [1~1]^{\mathrm{t}} \oplus I_{n-i})X} - \mathfrak{e}^{(I_{i-1}\oplus [0~1]^{\mathrm{t}} \oplus I_{n-i})X} - \mathfrak{e}^{(I_{i-1}\oplus [1~0]^{\mathrm{t}} \oplus I_{n-i})X} \right) = 0.
    $$
    If $X = [r_1, \cdots, r_{i-1},r_i,r_i^\prime,r_{i+1},\cdots,r_n]$, then this condition means that
    $$\E^{-1} \left( \mathfrak{e}^{[r_1, \cdots, r_{i-1},r_i+r_i^\prime,r_{i+1},\cdots,r_n]} - \mathfrak{e}^{[r_1, \cdots, r_{i-1},r_i^\prime,r_{i+1},\cdots,r_n]} -\mathfrak{e}^{[r_1, \cdots, r_{i-1},r_i,r_{i+1},\cdots,r_n]}\right)  = 0 .$$
    This follows from the definition of $\E^{-1}$.
    Applying Lemma \ref{202403191752}, one may check $\E^{-1} \circ \E = \mathrm{id}$ and $\E \circ \E^{-1} = \mathrm{id}$.
    This completes the proof.
\end{proof}

Recall the operad $I_{R_{\mathds{k}}}$ from Example \ref{202606121312}.
\begin{theorem} \label{202603241057}
    \begin{enumerate}
        \item The map in (\ref{202604061308}) gives an isomorphism of bimonoid species. 
        \item We have an operad isomorphism $\mathfrak{O}_{\mathbf{M}_{R}} \cong I_{R_{\mathds{k}}}$ which fits into the following commutative diagram:
        $$
         \begin{tikzcd}
              R_{\mathds{k}}^{\otimes} \ar[r, "\E_1"] & \Psi_{\mathbf{M}_{R}}  \\
             I_{R_{\mathds{k}}} \ar[r, "\cong"] \ar[u, hookrightarrow] & \mathfrak{O}_{\mathbf{M}_{R}} . \ar[u, hookrightarrow]
         \end{tikzcd}
        $$
    \end{enumerate}
\end{theorem}
\begin{proof}
    (1) is immediate from Lemma \ref{202408081425}
    We now prove (2).
    By (1), we have $\mathfrak{O}_{\mathbf{M}_{R}} = \mathrm{Pr} ( \Psi_{\mathbf{M}_{R}}) \cong \mathrm{Pr} ( R_{\mathds{k}}^{\otimes})$.
    So, we can use Proposition \ref{202603252001} to obtain $I_{R_{\mathds{k}}} \stackrel{\cong}{\to} \mathfrak{O}_{\mathbf{M}_{R}}$ which is obviously an operad map.
\end{proof}

\begin{Corollary}\label{202604051802}
    \begin{enumerate}
        \item We have an $(\mathtt{L}_{\mathfrak{S}},\mathtt{L}_\mathfrak{S})$-bimodule isomorphism $\mathtt{L}_{\mathbf{M}_{R}} / \mathtt{I}^{\mathsf{pr}}_{\mathbf{M}_{R}} \cong \mu R_{\mathds{k}}^{\otimes}$.
        \item The monad map $\mola_{\mathbf{M}_{R}} : \mu I_{R_{\mathds{k}}} \cong \mu\mathfrak{O}_{\mathbf{M}_{R}} \to \Phi_{\mathbf{M}_{R}}$ gives an isomorphism.
        In particular, the PROP $\tilde{\Phi}_{\mathbf{M}_{R}}$ is operadic.
    \end{enumerate}
\end{Corollary}
\begin{proof}
    The first claim follows from Theorems \ref{202603161027} and \ref{202603241057}.
    The second claim follows from (2) of Corollary \ref{202603201735} and Proposition \ref{202603252001}.
\end{proof}

\begin{notation} \label{202604061405}
    By Corollary \ref{202604051802}, the $(\mathtt{L}_{\mathbf{M}_{R}} , \mu I_{R_{\mathds{k}}})$-bimodule structure on $\mathtt{L}_{\mathbf{M}_{R}} / \mathtt{I}^{\mathsf{pr}}_{\mathbf{M}_{R}}$ is transported to $\mu  R_{\mathds{k}}^{\otimes}$.
    We denote the resulting bimodule by $\tilde{\mu} R_{\mathds{k}}^{\otimes}$.
\end{notation}

\begin{Example}
    Note that $\mathds{Z}_{\mathds{k}} =\mathds{k}$.
    The theorem applied to $R = \mathds{Z}$ yields $\mathtt{L}_{\mathbf{M}_{\mathds{Z}}}/\mathtt{I}^{\mathsf{pr}}_{\mathbf{M}_{\mathds{Z}}} \cong \tilde{\mu}\mathds{k}^{\otimes} \cong \tilde{\mu}\mathsf{Exp}$ and $\Phi_{\mathbf{M}_{\mathds{Z}}} \cong \mu I_{\mathds{k}} \cong \mathtt{L}_{\mathfrak{S}}$.
    This agrees with the results in Example \ref{202509192025} under the isomorphisms $\mathrm{M}_{\mathds{Z}} \cong \mathrm{M}_{\mathds{Z}}^{\mathsf{o}} \cong \mathbf{N}_1^{\mathsf{o}}$.
\end{Example}

\begin{remark} \label{202604051807}
    In this remark, we briefly describe the bimodule $\tilde{\mu} R_{\mathds{k}}^{\otimes}$ introduced in Notation \ref{202604061405}.
    \begin{enumerate}
        \item We begin with the left $\mathtt{L}_{\mathbf{M}_{R}}$-module structure.
        We may interpret this structure as the functor $\mathbf{M}_{R} \to \mathsf{Sp}_{\mathds{k}} ; m \mapsto \mathtt{L}_{\mathbf{M}_{R}}/ \mathtt{I}^{\mathsf{pr}}_{\mathbf{M}_{R}} (m,-) \cong ( R_{\mathds{k}}^{\otimes})^{\odot m}$.
        By Proposition \ref{202603201557}, this enhances to a symmetric monoidal functor, whose structure is encoded in $ R_{\mathds{k}}^{\otimes}$ as the bicommutative Hopf monoid species equipped with a left $R$-action.
        A similar, though not identical, construction appears in \cite{Touze2021}.
        Using Notation \ref{202511271410}, the left $R$-action is described as $$\mathfrak{e}^{a} \rhd \bigotimes^n_{i=1} r_i = \bigotimes^n_{i=1} ar_i, \quad \mathrm{where~}  \{r_1,\cdots,r_n\} \subset R_{\mathds{k}},~a\in R .$$
        \item  
        For the right $\mu I_{R_{\mathds{k}}}$-module structure on $\tilde{\mu} R_{\mathds{k}}^{\otimes}$, it suffices to describe the action $$\tilde{\mu} R_{\mathds{k}}^{\otimes} (m,n) \otimes \mu I_{R_{\mathds{k}}}(n,n) \to \mu I_{R_{\mathds{k}}} (m,n) \quad \mathrm{for~} m,n \in \mathds{N} ,$$ 
        since $\mu I_{R_{\mathds{k}}}$ is diagonal.
        Recall the isomorphism $\mu I_{R_{\mathds{k}}} (n,n) \cong R_{\mathds{k}} \thicksim \mathfrak{S}_n$ from Proposition \ref{202604061428}.
        For an $m$-partition $(n_1,\cdots,n_m)$ of $n$, consider an element $\bigotimes^m_{k=1} v_k \otimes \tau \in \tilde{\mu} R_{\mathds{k}}^{\otimes} (m,n)$ where $v_k \in  R_{\mathds{k}}^{\otimes n_k}$ and $\tau \in \mathfrak{S}_n$ (see Definition \ref{202603131024}).
        For $u \otimes \sigma \in  R_{\mathds{k}}^{\otimes n} \otimes \mathds{k} [\mathfrak{S}_n] = R_{\mathds{k}} \thicksim \mathfrak{S}_n$, we have
        $$
        (\bigotimes^m_{k=1} v_k \otimes \tau ) \lhd (u \otimes \sigma) = (\bigotimes^m_{k=1} v_k) (\tau \rhd u) \otimes \tau\sigma .
        $$
    \end{enumerate}
\end{remark}

\begin{remark}
    We record an equivalent description of (1) of Corollary \ref{202604051802}, although we will not make use of the following observation later on.
    The following map, given in Notation \ref{202604031739}, is also an isomorphism $\E : \mu\mathsf{Exp} \otimes_{\mathrm{H}}  R_{\mathds{k}}^{\otimes} \to \mathtt{L}_{\mathbf{M}_{R}} / \mathtt{I}^{\mathsf{pr}}_{\mathbf{M}_{R}}$.
    This follows from the construction of this map in the proof of Lemma \ref{202604031840}.
\end{remark}

\subsection{The refined eigenmonad adjunction}
\label{202404171454}

In this section, using the previous results, we study the refined eigenmonad adjunction associated with the primitivity ideal $\mathtt{I}^{\mathsf{pr}}_{\mathbf{M}_{R}}$.

\begin{Lemma} \label{202404161135}
    There exists a retract of the map $\mu I_{R_{\mathds{k}}} \to \tilde{\mu}  R_{\mathds{k}}^{\otimes}$ preserving the right $\mu I_{R_{\mathds{k}}}$-action.
\end{Lemma}
\begin{proof}
    The identity on $R_{\mathds{k}}$ induces a monomorphism $I_{R_{\mathds{k}}} \hookrightarrow  R_{\mathds{k}}^{\otimes}$.
    Note that the map $I_{R_{\mathds{k}}} \hookrightarrow  R_{\mathds{k}}^{\otimes}$ has a unique retract in the category of right $\mathtt{L}_\mathfrak{S}$-modules.
    Through the $\mu$ construction, the retract extends to a map $\varphi : \tilde{\mu}  R_{\mathds{k}}^{\otimes} \to \mu I_{R_{\mathds{k}}}$.
    For an $m$-partition $(n_1,\cdots,n_m)$ of $n$, consider an element $\bigotimes^m_{k=1} v_k \otimes \tau \in \tilde{\mu} R_{\mathds{k}}^{\otimes} (m,n)$ where $v_k \in  R_{\mathds{k}}^{\otimes n_k}$ and $\tau \in \mathfrak{S}_n$ (see Definition \ref{202603131024}).
    By the construction of $\varphi$, we have $\varphi ( \bigotimes^m_{k=1} v_k \otimes \tau ) = \bigotimes^m_{k=1} v_k \otimes \tau$ if $m=n$ and $n_1=\cdots =n_m = 1$; and $\varphi ( \bigotimes^m_{k=1} v_k \otimes \tau ) =0$ otherwise.   
    Using the description in (2) of Remark \ref{202604051807}, we can verify that $\varphi$ preserves the right $\mu I_{R_{\mathds{k}}}$-action.
\end{proof}

\begin{Lemma} \label{202512151805}
    Every left $\Phi_{\mathbf{M}_{R}}$-module is $\mathtt{I}^{\mathsf{pr}}_{\mathbf{M}_{R}}$-vanishingly extensible.
\end{Lemma}
\begin{proof}
    This follows from Lemmas \ref{202512081840} and \ref{202404161135}.
\end{proof}

\begin{theorem} \label{202404161429}
    The refined eigenmonad adjunction, given in Proposition \ref{202402011440}, induces the following adjunction:
    $$
    \begin{tikzcd}
             \mathbb{L} : \mu I_{R_{\mathds{k}}} \mbox{-} \mathsf{Mod} \arrow[r, shift right=1ex, ""{name=G}] & \mathtt{L}_{\mathbf{M}_{R}}\mbox{-}\mathsf{Mod}^{\mathsf{prim}} : \mathbb{R} , \arrow[l, shift right=1ex, ""{name=F}]
            \arrow[phantom, from=G, to=F, , "\scriptscriptstyle\boldsymbol{\top}"]
    \end{tikzcd}
    $$
    with faithful adjoint functors $\mathbb{L} = \tilde{\mu}  R_{\mathds{k}}^{\otimes} \otimes_{\mu I_{R_{\mathds{k}}}} (-)$ and $\mathbb{R} = \mathrm{Hom}_{\mathtt{L}_{\mathbf{M}_{R}}} (\tilde{\mu}  R_{\mathds{k}}^{\otimes}, -)$.
    \begin{enumerate}
        \item For $d\in\mathds{N}$, the functors $\mathbb{L}$ and $\mathbb{R}$ give a correspondence between $d$-truncated $\mu I_{R_{\mathds{k}}}$-modules and $\mathtt{L}_{\mathbf{M}_{R}}$-modules of polynomial degree at most $d$.
        \item For a $\mathtt{L}_{\mathbf{M}_{R}}$-module $\mathtt{M}$ of polynomial degree at most $d$, we have $(\mathbb{R}(\mathtt{M}))(d) = (\mathrm{cr}_{d} (\mathtt{M}) )(1,\cdots, 1)$.
    \end{enumerate}
\end{theorem}
\begin{proof}
    The refined eigenmonad adjunction gives the adjunction in the statement by Lemma \ref{202512151805} and the definition of primitive $\mathtt{L}_{\mathbf{M}_{R}}$-modules.
    The statements (1) and (2) follow from Theorem \ref{202509201658}.
\end{proof}

\begin{remark}
    To the best of our knowledge, this result does not appear explicitly in the literature, while its certain aspects are classically known.
    Indeed, by (2) of Theorem \ref{202404161429}, applied to $\mathds{k} = \mathds{Z}$, the adjunction restricts to the functors appearing in Pirashvili's equivalence \cite[Theorem 2.1]{PolyPira}: 
    $$\mathtt{L}_{\mathbf{M}_{R}}\mbox{-}\mathsf{Mod}^{\leq d}/ \mathtt{L}_{\mathbf{M}_{R}}\mbox{-}\mathsf{Mod}^{\leq (d-1)} \simeq \mu I_{R_{\mathds{k}}}\mbox{-}\mathsf{Mod}^{d-\mathsf{trun}} / \mu I_{R_{\mathds{k}}} \mbox{-}\mathsf{Mod}^{(d-1)-\mathsf{trun}} \simeq (R_{\mathds{k}} \thicksim \mathfrak{S}_d)\mbox{-}\mathsf{Mod}, $$
    where the final equivalence follows from Propositions \ref{202512081028} and \ref{202604061428}.
    Similarly, if $\mathds{k} = \mathds{Q}$, then the adjunction restricts to the functors appearing in Macdonald's equivalences \cite[Theorem 4]{wreathMac} \cite[Appendix I.A]{Mac1995}.
\end{remark}

\begin{Example} \label{202512091138}
    We present an example of Theorem \ref{202404161429} which does not yield an equivalence of categories.
    Suppose that $\mathds{k}$ is a nontrivial ring such that $2 \cdot 1_{\mathds{k}} = 0$ and $R = \mathds{Z}$.
    We have $R_{\mathds{k}} = \mathds{k}$, so $I_{ R_{\mathds{k}}^{\otimes}} \cong L_{\mathfrak{S}}$ and $\tilde{\mu}  R_{\mathds{k}}^{\otimes} \cong \tilde{\mu}\mathds{k}^{\otimes}$.
    Consider a nontrivial $\mathds{k}$-module $W$ endowed with the trivial $\mathfrak{S}_2$-action.
    Let $\mathtt{N}$ be the left $\mathtt{L}_{\mathfrak{S}}$-module concentrated on $2$-component as $\mathtt{N}(2) = W$.
    Then $\mathbb{L}( \mathtt{N}) (1) \cong \tilde{\mu}\mathds{k}^{\otimes} (1,2) \otimes_{\mathfrak{S}_2} W \cong \mathds{k} \otimes_{\mathfrak{S}_2} W \cong W$.
    The reduced comultiplication $\bar{\Delta} \in \mathtt{L}_{\mathbf{M}_{\mathds{Z}}} (2,1)$ trivially acts on $\mathbb{L}( \mathtt{N})$, as proved below.
    Therefore, we obtain $$\mathbb{R} ( \mathbb{L}( \mathtt{N} )) (1)) = \mathrm{V} ( \mathbb{L}( \mathtt{N}) ; \mathtt{I}^{\mathsf{pr}}_{\mathbf{M}_{\mathds{Z}}}) (1) = \mathrm{Ker} (\bar{\Delta} \rhd (-) : \mathbb{L}( \mathtt{N}) (1) \to \mathbb{L}( \mathtt{N}) (2) ) = W ,$$
    which implies $\mathtt{N}\not\cong \mathbb{R} ( \mathbb{L}( \mathtt{N} ))$.
\end{Example}

\begin{proof}[Proof of the claim in Example \ref{202512091138}]
    Under the isomorphism $\tilde{\mu}\mathds{k}^{\otimes} (m,n) \cong (\mathds{k}^{\otimes})^{\odot m} (n)$, the left $\mathtt{L}_{\mathbf{M}_{\mathds{Z}}}$-module structure on $\tilde{\mu}\mathds{k}^{\otimes}$ is equivalent to the bicommutative Hopf monoid structure on $\mathds{k}^{\otimes}$ considered in Example \ref{202604091327}.
    Using this argument, we can directly verify that the action
    $\bar{\Delta} \rhd (-) : \mathbb{L}(\mathtt{N}) (1) \to \mathbb{L}(\mathtt{N}) (2)$ gives a map 
    $W = \mathds{k} \otimes_{\mathfrak{S}_2} W \to \bigoplus_{{\bf 2}\to {\bf 2}} \mathds{k} \otimes_{\mathfrak{S}_2} W = \mathbb{L}(\mathtt{N})(2)$ assigns to $w \in W$ the element $(\mathrm{id}_{\mathbf{2}} + \tau) \otimes w = 2 \mathrm{id}_2 \otimes w = 0$, where $\tau : {\bf 2} \to {\bf 2}$ denotes the switching map.
\end{proof}

\appendix

\section{Proofs of Lemmas \ref{202604011630} and \ref{202603321744}} \label{202604101548}

In this appendix, we prove Lemmas \ref{202604011630} and \ref{202603321744}.
It is useful to recall the notion of a chain complex in the category $\mathsf{Sp}_{\mathds{k}}$ of $\mathds{k}$-linear species. 
A chain complex in $\mathsf{Sp}_{\mathds{k}}$ is a sequence of $\mathds{k}$-linear species, $\{ C_q\}_{q \in \mathds{Z}}$, with species maps $\partial_q : C_q \to C_{q-1}$ such that $\partial_{q-1} \circ \partial_q = 0$.
We denote by $H_q (C_\bullet)$ the $q$-th homology which is taken in the category $\mathsf{Sp}_{\mathds{k}}$.
A chain complex in $\mathsf{Sp}_{\mathds{k}}$ can be regarded as a functor from $\mathsf{Bij}^{\mathsf{o}}$ to the category of chain complexes of $\mathds{k}$-modules (see Section \ref{202604091335} for our convention on species).
In particular, a chain complex $C_\bullet$ in $\mathsf{Sp}_{\mathds{k}}$ induces a chain complex of $\mathds{k}$-modules $C_\bullet (X) : \cdots \to C_q (X) \to C_{q-1}(X) \to \cdots $ for any finite set $X$.
It is clear that $H_q (C_\bullet) (X) = H_q ( C_\bullet (X))$.

For chain complexes $C_\bullet,D_\bullet$ in $\mathsf{Sp}_\mathds{k}$, we define a chain complex $C_\bullet \odot D_\bullet$ in $\mathsf{Sp}_\mathds{k}$ to be $$(C_\bullet  \odot D_\bullet) (X) {:=} \bigoplus_{X_1\amalg X_2 = X} C_\bullet (X_1) \otimes D_\bullet (X_2),$$ where the tensor product denotes that for chain complexes.
As usual, we have a map of linear species, $\bigoplus_{n=p+q} H_p (C_\bullet ) \odot H_q (D_\bullet) \to H_{n} (C_\bullet \odot D_\bullet)$.
By definition, evaluated at $X$, this coincides with the following map:
$$
\bigoplus_{n=p+q} \bigoplus_{X_1 \amalg X_2 = X} H_p( C_\bullet (X_1)) \otimes H_q( D_\bullet (X_2)) \to \bigoplus_{X_1 \amalg X_2 = X} H_{n} (C_\bullet (X_1) \otimes D_\bullet (X_2)) .
$$

We now specialize to a coaugmented comonoid $\mathtt{C}$ in $\mathsf{Sp}_{\mathds{k}}$.
We regard $\bar{\mathtt{C}} \stackrel{\bar{\Delta}}{\to} \bar{\mathtt{C}}^{\odot 2}$ as a chain complex $K_\bullet$ in $\mathsf{Sp}_{\mathds{k}}$ by $K_0 = \bar{\mathtt{C}}$, $K_{-1} = \bar{\mathtt{C}}^{\odot 2}$ and $\partial_0 = \bar{\Delta}$.
By Lemma \ref{202603221101}, we have $$H_0 (K_\bullet^{\odot m}) = \mathrm{Pr}^m (\mathtt{C}).$$
Therefore, $\mathrm{Pr}(\mathtt{C})^{\odot m} \stackrel{\cong}{\to} \mathrm{Pr}^m (\mathtt{C})$ if and only if the canonical map $g : H_0 ( K_\bullet )^{\odot m} \to H_0 ( K_\bullet^{\odot m} )$ is an isomorphism.
\begin{proof}[Proof of Lemma \ref{202604011630}]
    By the previous discussion, it suffices to prove that
    $g: \bigotimes^m_{i=1} H_0 ( K_\bullet (X_i) ) \to H_0 ( \bigotimes^m_{i=1} K_\bullet (X_i))$ is an isomorphism for any decompositions $X= \amalg^{m}_{i=1} X_i$.
    Since $\mathds{k}$ is assumed to be a field, the K\"unneth formula ensures the statement.
\end{proof}

\begin{proof}[Proof of Lemma \ref{202603321744}]
    We continue to use the above notation.
    It suffices to show that for finite subsets $X_1,\cdots,X_m$, the canonical map $g : \bigotimes^m_{i=1} H_0 ( K_\bullet (X_i) ) \to H_0 ( \bigotimes^m_{i=1} K_\bullet (X_i))$ is an isomorphism.
    Let $|X_i| = n_i$.
    Fixing a bijection $X_i \cong {\bf n_i}$, we define a map $\varphi^\prime_{i} : K_{-1} (X_i) \to K_0 (X_i)$ by the commposition
    $$\varphi^\prime_{i} :  \bar{\mathtt{C}}^{\odot 2} (X_i) \cong \bar{\mathtt{C}}^{\odot 2} (n_i) \stackrel{\varphi_{n_i}}{\to} \bar{\mathtt{C}} (n_i) \cong \bar{\mathtt{C}} (X_i).$$
    By the hypothesis on $\varphi_n$'s, we have $\bar{\Delta} = \bar{\Delta} \circ \varphi^\prime_i \circ \bar{\Delta}$ which implies $\partial_0 \circ ( \mathrm{id}_{K_0(X_i)} - \varphi^\prime_i \circ \partial_0) = 0$.
    Since $K_{1} = 0$ by the definition of $K_\bullet$, the map $\varphi^\prime_i$ induces a map $r_i {:=} (\mathrm{id}_{K_0(X_i)} - \varphi^\prime_i \circ \partial_0) : K_0 (X_i) \to  H_0 ( K_\bullet (X_i) )$.
    Let $\iota : H_0 ( \bigotimes^m_{i=1} K_\bullet (X_i)) \hookrightarrow \bigotimes^m_{i=1} K_0 (X_i)$ denote the inclusion.
    We define a map $f$ by the composition,
    $$
    H_0 ( \bigotimes^m_{i=1} K_\bullet (X_i)) \stackrel{\iota}{\hookrightarrow} \bigotimes^m_{i=1} K_0 (X_i) \stackrel{\bigotimes^m_{i=1} r_{i}}{\longrightarrow} \bigotimes^m_{i=1} H_0 ( K_\bullet (X_i) ) .
    $$
    We will show that $f$ is the inverse of the canonical map $g$.
    Observe that $H_0 (K_\bullet (X_i)) \hookrightarrow K_0 (X_i) \stackrel{r_i}{\longrightarrow} H_0 (K_\bullet (X_i))$ is the identity.
    Thus, the following commutative diagram leads to $f \circ g = \mathrm{id}$:
    $$
    \begin{tikzcd}
        \bigotimes^m_{i=1} H_0 ( K_\bullet (X_i) ) \ar[r, "g"] \ar[d] & H_0 ( \bigotimes^m_{i=1} K_\bullet (X_i)) \ar[r, "f"] \ar[d, hookrightarrow, "\iota"] & \bigotimes^m_{i=1} H_0 ( K_\bullet (X_i) ) \ar[d, equal] \\
        \bigotimes^m_{i=1} K_0 (X_i) \ar[r, equal] & \bigotimes^m_{i=1} K_0 (X_i) \ar[r, "\bigotimes^m_{i=1} r_{i}"] & \bigotimes^m_{i=1} H_0 ( K_\bullet (X_i) ) .
    \end{tikzcd}
    $$
    Moreover, the following commutative diagram proves $g \circ f = \mathrm{id}$.
    $$
    \begin{tikzcd}
        H_0 ( \bigotimes^m_{i=1} K_\bullet (X_i)) \ar[r, "f"] \ar[d, hookrightarrow, "\iota"] & \bigotimes^m_{i=1} H_0 ( K_\bullet (X_i) ) \ar[r, "g"] \ar[d, equal] & H_0 ( \bigotimes^m_{i=1} K_\bullet (X_i)) \ar[d, hookrightarrow, "\iota"] \\
        \bigotimes^m_{i=1} K_0 (X_i) \ar[r, "\bigotimes^m_{i=1} r_{i}"] & \bigotimes^m_{i=1} H_0 ( K_\bullet (X_i) )  \ar[r] & \bigotimes^m_{i=1} K_0 (X_i) .
    \end{tikzcd}
    $$
    In fact, by the definition of $H_0 (K_\bullet^{\odot m})$, we have $(\mathrm{id}_{K_0(X_1)} \otimes \cdots \otimes \partial_0 \otimes \cdots \mathrm{id}_{K_0(X_m)}) \circ \iota= 0$, so the anticlockwise composition above coincides with $$\left( \bigotimes^m_{i=1} (\mathrm{id}_{K_0(X_i)} - \varphi^\prime_{i} \circ \partial_0) \right) \circ \iota = \left( \bigotimes^{m}_{i=1} \mathrm{id}_{K_0(X_i)} + \sum \pm \cdots \otimes (\varphi^\prime_{i} \circ \partial_0) \otimes \cdots \right) \circ \iota = \iota.$$

    We now prove the last claim.
    The map $\varphi_n \circ \bar{\Delta}$ is an idempotent since $(\varphi_n \circ \bar{\Delta})^2 = \varphi_n \circ \bar{\Delta} \circ \varphi_n \circ \bar{\Delta} = \varphi_n \circ \bar{\Delta}$, and hence, so is $\mathrm{id}_{\mathtt{C}(n)} - \varphi_n \circ \bar{\Delta}$.
    Furthermore, $\bar{\Delta} \circ ( \mathrm{id}_{\mathtt{C}(n)} - \varphi_n \circ \bar{\Delta}) = 0$ ensures that the image of the idempotent is contained in $\mathrm{Pr}(\mathtt{C})(n)$.
    For $v \in \mathrm{Pr}(\mathtt{C})(n)$, we also have $(\mathrm{id}_{\mathtt{C}(n)} - \varphi_n \circ \bar{\Delta})(v) = v$.
\end{proof}

\section{Existence of a PBW datum} \label{202512031342}

In this appendix, we give the deferred proofs of Lemmas \ref{202604011516} and \ref{202604021122}.
The strategy is based on a refinement of a basis theorem for free Lie algebras.

\subsection{Graded coaugmented coalgebras} \label{202601071503}

In this section, we give a general framework for treating graded coaugmented coalgebras.

Let $V$ be a commutative monoid with the monoid product $+$ and the unit $0 \in V$.
A $V$-graded ($\mathds{k}$-)module is a $\mathds{k}$-module $C$ endowed with a decomposition $C = \bigoplus_{\delta \in V} C_{\delta}$.
A morphism between $V$-graded modules is a $\mathds{k}$-linear map which respects the grading.
The category of $V$-graded modules is a monoidal category whose tensor product is given by $(C \otimes D)_{\delta} {:=} \bigoplus_{\substack{\delta_1,\delta_2\in V \\ \delta = \delta_1+\delta_2}} C_{\delta_1} \otimes C_{\delta_2}$.

A $V$-graded bialgebra is a bimonoid in the monoidal category of $V$-graded modules.
A $V$-graded coaugmented coalgebra is a coaugmented comonoid in the monoidal category of $V$-graded modules.

Let $H$ be a set with a $V$-graded decomposition $H = \coprod_{\delta \in V} H_{\delta}$.
Then the free $\mathds{k}$-algebra $\mathds{k} \langle H \rangle$ has a compatible $V$-graded bialgebra structure.
The comultiplication is characterized by
\begin{align*} 
    \Delta ( t_1 t_2 \cdots t_r )= \sum t_{\sigma(1)} \cdots t_{\sigma(k)} \otimes t_{\sigma(k+1)} \cdots t_{\sigma(n)}, \quad t_i \in H
\end{align*}
where the sum is taken over integers $0 \leq k \leq n$ and $\sigma$ is a $(k,n-k)$-shuffle, i.e. $\sigma \in \mathfrak{S}_n$ such that $\sigma(1) < \cdots < \sigma (k)$ and $\sigma (k+1) < \cdots \sigma (n)$.
There are some submodules of $\mathds{k} \langle H \rangle$ determined by a total order on $H$:
\begin{Defn} \label{202601071527}
    When $H$ is endowed with a total order $\preceq$, we define $\mathds{k} \langle H , \preceq \rangle$ to be the $\mathds{k}$-submodule of $\mathds{k} \langle H \rangle$ generated by the products $t_1 t_2 \cdots t_r$ for $r \in \mathds{N}$ and $\{t_i \in H \mid i \in {\bf r}\}$ such that $t_1 \succeq t_2 \succeq \cdots \succeq t_r$.
\end{Defn}
We note that $\mathds{k} \langle H , \preceq \rangle$ inherits the $V$-graded coaugmented coalgebra structures from $\mathds{k} \langle H \rangle$ while $\mathds{k} \langle H , \preceq \rangle \subset \mathds{k} \langle H \rangle$ is not closed under the multiplication.

\begin{Defn}
    Let $|\cdot |: V \to \mathds{N}$ be a monoid homomorphism.
    Using this, we also define a new grading on $\mathds{k} \langle H \rangle$: for a homogeneous generator $t_1t_2\cdots t_r$, its {\it maximum $V$-degree} is defined as the maximum of $|\delta_j|$ where $\delta_j \in V$ is the $V$-degree of $t_j$ for $j \in {\bf r}$.
\end{Defn}
\begin{Defn}
    Let $\mathds{k} \langle H , \preceq \rangle_{\max > d}$ be the $\mathds{k}$-submodule of $\mathds{k} \langle H \rangle$ generated by homogeneous elements with the maximum $V$-degree $>d$.
\end{Defn}

\subsection{Proof of Lemma \ref{202604011516}}
\label{202408191340}

This section is devoted to a proof of Lemma \ref{202604011516}.
To this end, using the framework of Appendix \ref{202601071503}, we present a graded refinement of the basis theorem for free Lie algebras \cite{reutenauer2003free}.

\begin{notation}
    Let $A$ be a $\mathds{k}$-algebra.
    For $a,b \in A$, let $[a,b] {:=} ab -ba$. 
    More generally, for $a_1,\cdots,a_m \in A$, we recursively define the (normalized) higher commutator by $[a_1, \cdots, a_m] {:=} [ [ a_1, \cdots, a_{m-1}] , a_m]$.
\end{notation}

We briefly recall several notions from \cite{reutenauer2003free}.
Let $S$ be a set and $M(S)$ be the free magma generated by $S$.
Let the parenthesis $(-,-)$ denote the magma operation on $M(S)$.
A free magma is basically a collection of binary trees whose leaves are labeled by elements of $S$, so we call an element of $M(S)$ a {\it tree}.
For $t \in M(S) \backslash S$, we denote by $t^\prime$ and $t^{\prime\prime}$ the first and second {\it subtrees} of $t$, i.e. $t = (t^\prime, t^{\prime\prime})$.
Given any total order $\preceq$ on $M(S)$ satisfying $t \prec t^{\prime\prime}$ whenever $t = (t^\prime, t^{\prime\prime}) \in M(S)$, we can determine the {\it Hall set} $H(S) \subset M(S)$ through the following procedure.
This is recursively defined by declaring that $S \subset H(S)$, and that $t = (t^\prime, t^{\prime\prime}) \in H(S)$ if and only if $t^\prime, t^{\prime\prime} \in H(S)$, $t^\prime \prec t^{\prime\prime}$ and either $t^\prime \in S$ or $(t^\prime)^{\prime\prime} \succeq t^{\prime\prime}$.
An element of $H(S)$ is called a {\it Hall tree}.

Let $(t_1, \cdots, t_n)$ with $n \in \mathds{N}^\ast$ be a sequence of Hall trees.
We say that $(t_1, \cdots, t_n)$  is {\it standard} if, for any $i$, either $t_i \in S$, or $t_i^{\prime\prime} \succeq t_{i+1}, \cdots, t_n$.
A {\it legal rise} of $(t_1, \cdots, t_n)$ is an index $i$ such that $t_i \prec t_{i+1} \succeq t_{i+2}, \cdots, t_n$.
An {\it inversion} of the sequence is a pair $(i,j)$ such that $i < j$ and $t_i \prec t_j$.

We define a map $g : M(S) \to \mathds{k} \langle S \rangle$ recursively by 
\begin{align*}
    g(x) {:=} x \mathrm{~for~}x \in S;\quad g(t) {:=} [g(t^\prime), g(t^{\prime\prime})] \quad \mathrm{if~} t = (t^\prime, t^{\prime\prime}) .
\end{align*}
It extends to a bialgebra map $g : \mathds{k} \langle M(S) \rangle \to \mathds{k} \langle S \rangle$.
By Definition \ref{202601071527} applied to the above order on $H(S)$, we take the $\mathds{k}$-submodule $\mathds{k}\langle H(S), \preceq \rangle \subset \mathds{k} \langle H(S) \rangle$.
Then the restriction of $g$ to $\mathds{k} \langle H(S) , \preceq \rangle \subset \mathds{k}\langle M(S) \rangle$ gives a coaugmented coalgebra homomorphism,
$$g : \mathds{k} \langle H(S) , \preceq \rangle \to \mathds{k} \langle S \rangle.$$

Let $\Lie (S) \subset \mathds{k}\langle S \rangle$ be the Lie algebra (over $\mathds{k}$) generated by $\mathds{k}[S] \subset \mathds{k}\langle S \rangle$.

Recall the inclusion $\mathds{k} [H(S)] \subset \mathds{k} \langle H(S), \preceq \rangle$.

\begin{theorem}[\mbox{\cite[Theorems 3 and 4]{reutenauer2003free}}] \label{202601071540}
    The map $g : \mathds{k} \langle H(S) , \preceq \rangle \to \mathds{k} \langle S \rangle$ gives an isomorphism of coaugmented coalgebras.
    Furthermore, the map $g$ induces an isomorphism between $\mathds{k} [H(S)]$ and $\Lie (S)$.
\end{theorem}

In what follows, we refine this result.
For a set $S$, let $V_S$ denote the abelian monoid of maps $S \to \mathds{N}$ with finite support.
For $\delta \in V_S$, we introduce the $\delta$-components of the algebra $\mathds{k} \langle S \rangle$, the set $H(S)$ and the Lie algebra $\Lie (S)$ as follows.

\begin{notation}
    For $\delta \in V_S$, let $$|\delta| {:=} \sum_{x \in S} \delta (x) \in \mathds{N} .$$
\end{notation}

\begin{Defn}
    Let $\mathds{k} \langle S \rangle_\delta {:=} \mathds{k}$ if $\delta = 0$, and $\mathds{k} \langle S \rangle_{\delta} {:=} \mathds{k} s$ if $|\delta| =1$, where $s \in S$ is the unique element of the support of $\delta$.
    For $|\delta| \geq 2$, we recursively define $\mathds{k} \langle S \rangle_\delta$ by
    \begin{align*}
        \mathds{k} \langle S \rangle_\delta {:=}
            \sum \mathds{k} \langle S \rangle_{\delta_1} \mathds{k} \langle S \rangle_{\delta_2} ,
    \end{align*}
    where the sum is taken over $\delta_1,\delta_2 \in V_S$ such that $\delta_1 + \delta_2=\delta$ and $|\delta_1|, |\delta_2| < |\delta|$.
\end{Defn}

In other words, $\mathds{k} \langle S \rangle_\delta$ is the $\mathds{k}$-submodule of $\mathds{k} \langle S \rangle$ spanned by monomials containing exactly $\delta(x)$ occurrences of each $x \in S$.

\begin{Defn}
    For $\delta \in V_S$, let $H_{\delta} (S) \subset H(S)$ be the subset consisting of Hall trees which have $\delta (x)$ copies of $x \in S$ as leaves. In particular, $H_0 (S) = \emptyset$.
\end{Defn}

We have a decomposition $H(S) = \coprod_{\delta \in V_S} H_{\delta} (S)$ which yields a $V_S$-graded bialgebra $\mathds{k} \langle H(S) \rangle$.
Furthermore, its underlying coaugmented coalgebra structure restricts to one on $\mathds{k} \langle H(S) , \preceq \rangle$.

\begin{Defn}
    Let $\Lie (S)_\delta {:=} 0$ if $\delta = 0$ and $\Lie (S)_{\delta} {:=} \mathds{k} s$ if $|\delta| =1$, where $s \in S$ is the unique element of the support of $\delta$.
    For $|\delta| \geq 2$, we recursively define
    \begin{align*}
        \Lie (S)_{\delta} {:=} \sum [\Lie (S)_{\delta_1}, \Lie (S)_{\delta_2}] ,
    \end{align*}
    where the sum is taken over $\delta_1,\delta_2 \in V_S$ such that $\delta_1 + \delta_2=\delta$ and $|\delta_1|, |\delta_2| < |\delta|$.
\end{Defn}

Similarly, $\Lie (S)_{\delta}$ is the $\mathds{k}$-submodule of $\Lie (S)$ generated by Lie brackets in which $x \in S$ appears exactly $\delta(x)$ times.
For $\delta \in V_S$, we have $\Lie (S)_{\delta} \subset \mathds{k} \langle S \rangle_{\delta}$, so that $\Lie (S) = \bigoplus_{\delta \in V_{S}} \Lie (S)_{\delta}$ is a $V_S$-graded module.

\begin{Corollary} \label{202410131830}
    The map $g : \mathds{k} \langle H(S),\preceq \rangle \to \mathds{k} \langle S \rangle$ gives an isomorphism of $V_{S}$-graded coaugmented coalgebras.
    Furthermore, it induces an isomorphism of $V_{S}$-graded $\mathds{k}$-modules between $\mathds{k}[H(S)]$ and $\Lie (S)$.
\end{Corollary}
\begin{proof}
    By Theorem \ref{202601071540}, it suffices to prove that $g ( \mathds{k} \langle H(S),\preceq \rangle_\delta) \subset \mathds{k} \langle S \rangle_\delta$ for $\delta \in V_S$.
    The case that $\delta = 0$ is trivial.
    Using the induction on $|\delta|$, we prove the rest.
    If $|\delta| =1$, then the support of $\delta$ is $\{ x\} \subset S$, so $H_{\delta} (S) = \{ x\}$ and $\mathds{k} \langle H(S),\preceq \rangle_\delta = \mathds{k} x$.
    By definition, we have $g(x) = x \in \mathds{k} \langle S \rangle_{\delta}$.
    Hence, it completes the proof for the case that $|\delta| =1$.
    We now suppose that the assertion holds for $|\delta| \leq k$ where $k \geq 1$.
    Let $\delta^\prime \in V_S$ such that $|\delta^\prime| = k+1$.
    Recall that $\mathds{k} \langle H(S),\preceq \rangle_{\delta^\prime}$ is generated by $t_1 \cdots t_r$ where $\delta_i \in V_S \backslash \{0\}$, $t_i \in H_{\delta_i} (S)$, $\delta^\prime = \sum^r_{i=1} \delta_i$ and $t_1 \succeq \cdots \succeq t_r$.
    If $r \geq 2$, then $|\delta_i| < |\delta^\prime| = k+1$, so by the induction hypothesis, we have $g(t_i) \in \mathds{k} \langle S \rangle_{\delta_i}$.
    Hence, we obtain $g(t_1 \cdots t_r) = g(t_1) \cdots g(t_r) \in \mathds{k} \langle S \rangle_{\delta^\prime}$.
    If $r =1$, then the assumption $|\delta^\prime|=k+1 \geq 2$ implies that $t_1 \in H(S) \backslash S$.
    So, we can take subtrees as $t_1 = (t_1^\prime,t_1^{\prime\prime})$.
    Then there exist $\delta_1,\delta_2\in V_S$ such that $\delta^\prime = \delta_1 + \delta_2$, $t_1^\prime \in H_{\delta_1} (S)$ and $t_1^{\prime\prime} \in H_{\delta_2} (S)$.
    Since $\max (|\delta_1|,|\delta_2|) < |\delta^\prime| = k+1$, applying the induction hypothesis to $\delta_1$ and $\delta_2$, we obtain $g ( t_1) = [ g(t_1^\prime), g(t_1^{\prime\prime})] \in [\mathds{k} \langle S \rangle_{\delta_1},\mathds{k} \langle S \rangle_{\delta_2}] \subset \mathds{k} \langle S \rangle_{\delta^\prime}$.
\end{proof}

We now prove Lemma \ref{202604011516} by applying the previous results to $S = \Theta$.
Recall that $\Theta$ is the set fixed in Notation \ref{202606111507}.
\begin{proof}[Proof of Lemma \ref{202604011516}] 
    For a finite set $X \subset \mathds{N}^\ast$, let $X^\prime {:=} \{ \theta_i \in \Theta \mid i \in X \}$.
    Then $\mathsf{Lin} (X)$ is naturally isomorphic to the $\mathds{1}_{X^\prime}$-component of $\mathds{k} \langle \Theta \rangle$.
    Note that, under this identification, the bimonoid species structure on $\mathsf{Lin}$ (see Example \ref{202604101536})  coincides with the structure induced from the bialgebra structure on $\mathds{k} \langle \Theta \rangle$.
    Moreover, this identification restricts to $ \mathfrak{Lie} (X) \cong \Lie (\Theta)_{\mathds{1}_{X^\prime}}$.
    We now define $$\mathfrak{B}_{X}\subset \mathsf{Lin}(X)$$ to be the subset corresponding to $g (H_{\mathds{1}_{X^\prime}} (\Theta)) \subset \mathds{k} \langle \Theta \rangle_{\mathds{1}_{X^\prime}}$.
    The union $\coprod_{X \in \mathcal{P}_f (\mathds{N}^{\ast})} \mathfrak{B}_X$ inherits a total order from the Hall set $H(\Theta)$. 
    This construction yields the desired PBW datum.
    Indeed, any element of $\mathfrak{B}_X$ is primitive by the preceding identification of bimonoid structures.
    By Corollary \ref{202410131830},
    \begin{align} \label{202606211055}
        \mathds{k} \langle H(\Theta), \preceq \rangle_{\mathds{1}_{X^\prime}} \stackrel{\cong}{\to} \mathds{k} \langle \Theta \rangle_{\mathds{1}_{X^\prime}} \cong \mathsf{Lin} (X) ,
    \end{align}
    from which (2) of Definition \ref{202604021629} follows, since $\mathds{k} \langle H(\Theta), \preceq \rangle_{\mathds{1}_{X^\prime}}$ is freely generated by words $h_1 \cdots h_r$, where $r \in \mathds{N}^\ast,~ X^\prime = \coprod^r_{j=1} X_j ,~ h_j \in H_{\mathds{1}_{X_j}} (\Theta)$ such that $h_1 \succ \cdots \succ h_r$.
    Furthermore, by Corollary \ref{202410131830}, the isomorphisms (\ref{202606211055}) induce  
    $\mathfrak{Lie} (X)  \cong \mathds{k} [ \mathfrak{B}_{X} ]$, which proves the last statement in Lemma \ref{202604011516}.
\end{proof}

\subsection{The map $\T_S$ and commutators}

This section provides preliminaries for the proof of Lemma \ref{202604021122}.
Fix $n \in \mathds{N}$.
Recall from Definition \ref{202509110944} the map $\T_S : \mathds{k} [\mathsf{F}_n] \to \alg (n)_{\mathds{1}_S}$, where $S \subset {\bf n}$.
We compute the value of $\T_S$ assigned to several commutators in the free group $\mathsf{F}_n$.

\begin{Defn}
    For $S \subset {\bf n}$, let $\mathsf{F} (S)$ denote the subgroup of $\mathsf{F}_n$ generated by $x_i$ for $i\in S$.
    We define a subset $Z(S) \subset \mathsf{F}_n$ by
    $$
    Z(S) {:=} \{ a \in \mathsf{F}(S) \mid  \T_{S^\prime} (a) = 0 , \quad \mathrm{for~all~} \emptyset \neq S^\prime \subsetneq S  \} .
    $$
\end{Defn}

\begin{Example}
    \begin{enumerate}
        \item $Z (\emptyset ) = e$ and $e \in Z(S)$ for any $S \subset {\bf n}$. 
        \item For a singleton $S = \{i\}$, $Z(S)$ is the subgroup generated by $x_i \in \mathsf{F}_n$.
        \item If $S = \{ i ,j\}$, then $x_i^{-1} x_j^{-1} x_i x_j \in Z(S)$.
    \end{enumerate}
\end{Example}

\begin{Lemma} \label{202604030409}
    Let $S_1 \subset S_2 \subset {\bf n}$.
    For $a \in Z(S_1)$, $\T_{S^\prime} (a) = 0$ for $\emptyset \neq S^\prime \subsetneq S_2$ such that $S^\prime \neq S_1$.
\end{Lemma}
\begin{proof}
    Let $\emptyset \neq S^\prime \subsetneq S_2$.
    Since $a$ lies in the subgroup generated by $x_i,~i \in S_1$, we have $\T_{S^\prime} (a) = 0$ if $S^\prime \not\subset S_1$.
    If $S^\prime \subsetneq S_1$, then the hypothesis $a \in Z(S_1)$ implies that $\T_{S^\prime} (a) = 0$.
\end{proof}

\begin{Lemma} \label{202604011801}
    Let $\emptyset \neq S \subset {\bf n}$.
    If $a \in Z(S)$, then $a^{-1} \in Z(S)$ and $\T_{S} (a^{-1}) = - \T_{S} (a)$.
\end{Lemma}
\begin{proof}
    For a nonempty subset $S^\prime$ of $S$, by Lemma \ref{202509081651}, we have
    $$
    0 = \T_{S^\prime} (e) = \T_{S^\prime} (aa^{-1}) = \sum \T_{S_1} (a) \T_{S_2} (a^{-1}) , 
    $$
    where the sum is taken over decompositions $S^\prime = S_1 \amalg S_2$.
    Since $a \in Z(S)$, the terms except the cases of $(S_1,S_2) = (S^\prime, \emptyset)$ and $(S_1,S_2) = (\emptyset , S^\prime)$ should vanish.
    So, $0 = \T_{S^\prime}(a) + \T_{S^\prime} (a^{-1})$.
    If $S^\prime \subsetneq S$, then it implies $\T_{S^\prime} (a^{-1}) = 0$, so $a^{-1} \in Z(S)$.
    If $S^\prime = S$, then we obtain $\T_{S} (a^{-1}) = - \T_{S} (a)$.
\end{proof}

\begin{Lemma} \label{202604030431}
    Let $U,V \subset {\bf n}$ and $S = U \cup V$.
    Let $a \in Z(U)$ and $b \in Z(V)$.
    Then $a^{-1}b^{-1}ab \in Z(S)$ and
    \begin{align*}
        \T_{S} ( a^{-1}b^{-1}ab ) =
        \begin{cases}
            [ \T_{U} (a), \T_{V} (b)] & \mathrm{if~} S \neq \emptyset = U\cap V , \\
            1 & \mathrm{if~} S = \emptyset , \\
            0 & \mathrm{otherwise}.
        \end{cases}
    \end{align*}
\end{Lemma}
\begin{proof}
    The case where $U = \emptyset$ or $V = \emptyset$ is trivial.
    We assume that $U,V$ are not empty.
    Let $\emptyset \neq S^\prime \subset S$.
    By Lemma \ref{202509081651}, we obtain
    \begin{align} \label{202604030418}
    \T_{S^\prime} ( a^{-1}b^{-1}ab ) = \sum \T_{S_1} (a^{-1} ) \T_{S_2} (b^{-1} ) \T_{S_3} (a) \T_{S_4} (b) ,
    \end{align}
    where the sum is taken over decompositions $S^\prime = \coprod^4_{j=1} S_j$.
    We investigate each term in the summation.
    By Lemma \ref{202604030409}, unless $\{ S_1 , S_3 \} \subset \{ U, \emptyset \}$, the term in (\ref{202604030418}) should vanish.
    Analogously, unless $\{ S_2 , S_4 \} \subset \{ V, \emptyset \}$, the term vanishes.
    It is sufficient to consider the above sum only for the cases that $\{ S_1 , S_3 \} \subset \{ U, \emptyset \}$ and $\{ S_2 , S_4 \} \subset \{ V, \emptyset \}$.
    
    Assume that $\emptyset \neq S^\prime \subsetneq S$.
    If $S^\prime \not\in \{ U, V\}$, then the sum in (\ref{202604030418}) vanishes by the condition $S^\prime = \coprod^4_{j=1} S_j$.
    If $S^\prime = U \neq V$, then $S_2,S_4$ are empty, so the sum reduces to $\T_{U} (a^{-1}) + \T_{U} (a)$, which vanishes by Lemma \ref{202604011801}.
    If $S^\prime = U = V$, then
    the sum equals $$\T_{U} (a^{-1}) + \T_{V} (b^{-1}) + \T_{U} (a) + \T_{V} (b) = 0.$$
    The case where $S^\prime =V \neq U$ is proved similarly.
    Therefore, $a^{-1}b^{-1}ab \in Z(S)$.
    
    To prove the last assertions, we now consider the sum (\ref{202604030418}) applied to $S^\prime = S$.
    If $U \cap V \neq \emptyset$, then the condition that $\{ S_1 , S_3 \} \subset \{ U, \emptyset \}$ and $\{ S_2 , S_4 \} \subset \{ V, \emptyset \}$ implies that the sum should vanish, since $S_j$'s do not intersect.
    If $U \cap V = \emptyset$, then, by Lemma \ref{202604011801}, the sum in (\ref{202604030418}) reduces to
    \begin{align*}
        &\T_U(a^{-1}) T_V (b^{-1} ) + T_U(a^{-1}) T_V (b) + T_V (b^{-1}) T_U (a) + T_U (a) T_V (b) , \\
        =& \T_U(a) T_V (b ) - T_U(a) T_V (b) - T_V (b) T_U (a) + T_U (a) T_V (b) = [ \T_{U} (a), \T_{V} (b)] .
    \end{align*}
\end{proof}

For a group $G$ and $a_1, \cdots, a_m \in G$, we define the higher commutator recursively by 
\begin{align*}
    \langle a_1,a_2 \rangle {:=} a_1^{-1}a_2^{-1}a_1a_2, \quad \mathrm{and} \quad \langle a_1, \cdots, a_m \rangle {:=} \langle \langle a_1, \cdots, a_{m-1} \rangle, a_m \rangle
\end{align*}

\begin{Lemma} \label{202509151210}
    For $S \subset {\bf n}$ and $i_1, \cdots, i_l \in {\bf n}$, we have 
    \begin{align*}
        \T_{S} ( \langle x_{i_1}, x_{i_2}, \cdots, x_{i_l} \rangle^{\pm 1} ) =
        \begin{cases}
            \pm [\theta_{i_1}, \theta_{i_2}, \cdots, \theta_{i_l}] & |S| = l \mathrm{~and~} S = \{ i_1, i_2, \cdots, i_l \} , \\
            1 & S = \emptyset , \\
            0 & \mathrm{otherwise}.
        \end{cases}
    \end{align*}
\end{Lemma}
\begin{proof}
    It suffices to prove the statement for the case of the exponent $+1$, since this implies the other one by Lemma \ref{202604011801}.
    Moreover, we can assume $S \subset \{i_1,\cdots, i_l\}$, since if not, by the definition of $\T_S$, $\T_{S} ( \langle x_{i_1}, x_{i_2}, \cdots, x_{i_l} \rangle ) = 0$.
    
    By iteratively applying Lemma \ref{202604030431}, one obtains $\langle x_{i_1}, x_{i_2}, \cdots, x_{i_l} \rangle \in Z ( \{ i_1, \cdots, i_l \})$.
    So, if $S \subsetneq \{ i_1, \cdots, i_l\}$, we have $\T_S (\langle x_{i_1}, x_{i_2}, \cdots, x_{i_l} \rangle) =0$.
    
    We now assume that $S= \{ i_1, \cdots, i_l \}$.
    Put $S_k = \{ i_1, \cdots, i_k \}$ for $k \in {\bf l}$.
    First, consider the case $|S|=l$.
    Clearly, $i_1,\cdots,i_l$ are pairwise distinct.
    The assertion is true if $l = 1$.
    Based on the following formula obtained from Lemma \ref{202604030431}, we can also inductively prove the assertion: 
    $$\T_{S_k} ( \langle x_{i_1}, x_{i_2}, \cdots, x_{i_k} \rangle ) = [ \T_{S_{k-1}} (\langle x_{i_1}, x_{i_2}, \cdots, x_{i_{k-1}} \rangle), \T_{\{i_k\}} (x_{i_k}) ] = [ [\theta_{i_1}, \theta_{i_2}, \cdots, \theta_{i_{k-1}}], \theta_{i_k} ] .$$
    
    All that remains is the case $|S| \neq l$.
    We can choose the minimum $k_0 \in {\bf l}$ such that $i_{k_0} = i_{k_0+1}$.
    Then by $S_{k_0+1} = S_{k_0} \cup \{ i_{k_0+1} \}$ and $S_{k_0} \cap \{ i_{k_0+1} \} \neq \emptyset$, we obtain from Lemma \ref{202604030431},
    $$\T_{S_{k_0+1}} ( \langle x_{i_1} , \cdots, x_{i_{k_0+1}} \rangle) = \T_{S_{k_0+1}} ( \langle \langle x_{i_1} , \cdots, x_{i_{k_0}} \rangle ,x_{i_{k_0+1}} \rangle) = 0.$$
    Hence, for $k \geq k_0+1$, we can prove $\T_{S_{k}} ( \langle x_{i_1} , \cdots, x_{i_{k}} \rangle) = 0$ by iteratively applying Lemma \ref{202604030431}.    
\end{proof}

\subsection{Proof of Lemma \ref{202604021122}}
\label{202604031207}

In this section, we prove Lemma \ref{202604021122}.
We begin by recalling $\K_\mathcal{R}$ from Definition \ref{202604021150}.
The following description of $\K_\mathcal{R}$ for a general radical functor for groups is useful.
\begin{Lemma} \label{202509101818}  
    Let $n \in \mathds{N}$ and $U \subset \mathcal{R}(\mathsf{F}_n)$.
    If $U$ generates $\mathcal{R}(\mathsf{F}_n)$, then $\K_{\mathcal{R}}(n)$ coincides with the $\mathds{k}$-submodule of $\mathsf{Lin}(n)$ generated by $u \T_{{\bf n} \backslash S} (h)$ where $S \subsetneq {\bf n}$, $u \in \mathsf{Lin} ( S)$ and $h \in U$.
\end{Lemma}
\begin{proof}
    The case where $n=0$ is clear, since $\K_\mathcal{R} (0) = 0$, so we assume that $n >0$.
    Let $V\subset \mathsf{Lin}(n)$ be the $\mathds{k}$-submodule appearing in the lemma.
    Recall that by Proposition \ref{202603292144}, we have $\K_{\mathcal{R}}(n) = \T_n(\mathrm{Ker} (q_n))$.
    We will prove that $V = \T_n(\mathrm{Ker} (q_n))$.
    The group $\mathcal{R}(\mathsf{F}_n)$ is a normal subgroup of $\mathsf{F}_n$, so the $\mathds{k}$-module $\mathrm{Ker} (q_{n})$ is generated by $g(h-e) \in \mathds{k}[\mathsf{F}_n]$ for $g \in \mathsf{F}_n$ and $h \in \mathcal{R}(\mathsf{F}_n)$.
    Since $U$ generates $\mathcal{R}(\mathsf{F}_n)$, we may assume that $h \in U$.
    Thus, the $\mathds{k}$-module $\K_{\mathcal{R}}(n)$ is generated by $\T_{n} ( g(h-e) )$ for $h \in U$ and $g \in \mathsf{F}_n$.
    Applying Lemma \ref{202509081651}, we obtain
    \begin{align*}
        \T_{n} ( g(h-e) ) = \T_{n} ( gh) - \T_{n}(g)  = \sum_{\substack{S=S_1 \amalg S_2 \\ S_1 \neq {\bf n}}} \T_{S_1} (g) \T_{S_2} (h ) \in V.
    \end{align*}
    
    We now prove that $V \subset \T_n (\mathrm{Ker}(q_n))$.
    To this end, we shall show $u \T_{{\bf n}\backslash S} (h) \in \T_n (\mathrm{Ker}(q_n))$ where $u,h,S$ are given as in the statement.
    We prove this using the induction on $|S|$.
    The case where $|S|=0$ is clear since $\T_n ( h) = \T_n (h-e) \in \T_n (\mathrm{Ker}(q_n))$.
    Assume that $|S| > 0$.
    We choose $v \in \mathds{k}[\mathsf{F} (S)]$ such that $\T_{S} (v) = u$ where $\mathsf{F} ( S)$ denotes the free group generated by $S$.
    Modulo $\T_n (\mathrm{Ker}(q_n))$, we obtain
    $$
    0 \equiv \T_n ( v (h-e)) = \sum_{\substack{S=S_1 \amalg S_2 \\ S_1 \neq {\bf n}}} \T_{S_1} (v) \T_{S_2} (h ) .
    $$
    Note that in the sum, if $S_1 \not\subset S$, then $\T_{S_1} (v) = 0$ since $v \in \mathds{k} [\mathsf{F} (S)]$.
    So, we obtain 
    $$0 \equiv \sum_{\substack{S=S_1 \amalg S_2 \\ S_1 \subset S}} \T_{S_1} (v) \T_{S_2} (h ).$$
    By the induction hypothesis, if $S_1 \subsetneq S$, then $\T_{S_1} (v) \T_{S_2} (h ) \equiv 0$.
    Thus, $\T_{S} (v) \T_{{\bf n} \backslash S} (h ) \equiv 0$.
\end{proof}

\begin{Defn} \label{202601071542}
    We define $\IL_{d} (S) \subset \mathds{k} \langle S \rangle$ to be the two-sided ideal generated by $\bigoplus_{|\delta|>d} \Lie (S)_{\delta}$.
    It inherits the $V_{S}$-grading from $\mathds{k} \langle S \rangle$:
    $$
    \IL_{d} (S)_{\delta} = \sum_{\substack{\delta = \delta_1 + \delta_2 + \delta_3 \\ |\delta_2| > d}} \mathds{k} \langle S \rangle_{\delta_1} \cdot \Lie (S)_{\delta_2} \cdot \mathds{k} \langle S \rangle_{\delta_3} .
    $$
\end{Defn}

We recall from Section \ref{202601071503} that $\mathds{k} \langle H(S) , \preceq \rangle_{\max > d}$ is the $\mathds{k}$-submodule of $\mathds{k} \langle H(S) , \preceq \rangle$ generated by homogeneous elements of the maximum $V_{S}$-degree $>d$.
By carefully refining the proof of \cite[Theorem 3]{reutenauer2003free}, we obtain the following:
\begin{Lemma} \label{202509161922}
    The isomorphism $g : \mathds{k} \langle H(S) , \preceq \rangle \to \mathds{k} \langle S \rangle$ maps $\mathds{k} \langle H(S) , \preceq \rangle_{\max > d}$ onto $\IL_{d} (S)$.
\end{Lemma}
\begin{proof}
    Clearly, the map $g : \mathds{k} \langle H(S) , \preceq \rangle \to \mathds{k} \langle S \rangle$ in Corollary \ref{202410131830} satisfies $g (\mathds{k} \langle H(S) , \preceq \rangle_{\max > d}) \subset \IL_{d} (S)$.
    It suffices to prove that the inclusion $\mathds{k} \langle H(S) , \preceq \rangle_{\max > d} \hookrightarrow \IL_{d} (S)$ is surjective.
    We start with normalizing generators of $\IL_{d} (S)$.
    Observe that, the identity $[s_1, \cdots, s_l] s_{l+1} = s_{l+1} [s_1, \cdots, s_l] + [s_1, \cdots, s_l, s_{l+1}]$ for $s_1,\cdots,s_{l+1} \in S$ moves the commutator term to the right in the first summand, while the second summand contains a commutator whose length is increased by one.
    By iteratively applying this observation, one may see that $\IL_{d} (S)$ coincides with the {\it left} ideal of $\mathds{k} \langle S \rangle$ generated by $\bigoplus_{\delta \in V_S,|\delta| > d} \Lie (S)_{\delta}$.
    Hence, by Corollary \ref{202410131830}, the $\mathds{k}$-module $\IL_{d} (S)$ is generated by $s_1 \cdots s_k g(t) = g (s_1 \cdots s_k t)$ where $s_j \in S$, $t \in H_{\delta} (S)$ and $\delta \in V_S$ such that $|\delta| > d$.
    All that remains is to prove that such $g (s_1 \cdots s_k t)$ lies in the image $g (\mathds{k} \langle H(S) , \preceq \rangle_{\max > d})$.
    
    We now use the proof of \cite[Theorem 3]{reutenauer2003free}.
    In that proof, one encounters an algorithm that, given a standard sequence $(t_1, \cdots, t_n)$, systematically rewrites $g(t_1 \cdots t_n)$ as a linear combination of $g(u)$ with $u \in \mathds{k} \langle H(S), \preceq \rangle$.
    Given a standard sequence $s$, at each step, we choose a legal rise $i$ of $s$, obtaining two standard sequences $\lambda_i (s)$ with shorter length and $\rho_i (s)$ with one inversion less such that $g(s) = g(\lambda_i (s)) + g(\rho_i(s))$.    
    One can observe that, if the maximal $V_{S}$-degree of $s$ exceeds $d$, then so do the maximal $V_{S}$-degrees of $\lambda_i (s)$ and $\rho_i (s)$; thus, it follows by induction that $g(s)$ is equal to a sum of $g(u)$ with $u \in \mathds{k} \langle H(S), \preceq \rangle_{\max >d}$.
    We apply the previous algorithm to the standard sequence of Hall trees $(s_1 , \cdots , s_k , t)$.
\end{proof}

\begin{Lemma} \label{202604031247}
    For $n \in \mathds{N}$,
    $\IL_c (S)_{\mathds{1}} = \K_{\gamma_{c+1}}(n)$ where $S = \{ \theta_i \in \Theta \mid i \in {\bf n} \}$ and $\mathds{1} \in V_S$ is the constant map with value $1$.
\end{Lemma}
\begin{proof}
    We begin by proving that $\IL_{c} (S)_{\mathds{1}}$ contains $\K_{\gamma_{c+1}} (n)$.
    Recall the fact that the group $\gamma_{c+1} (\mathsf{F}_n) \subset \mathsf{F}_n$ is the normal closure of the higher commutators $\langle x_{i_1},x_{i_2},\cdots, x_{i_{c+1}}\rangle$ where $i_l$'s range over ${\bf n}$.
    Hence, by applying Lemma \ref{202509101818} to $\mathcal{R} = \gamma_{c+1}$, the $\mathds{k}$-module $\K_{\gamma_{c+1}} (n)$ is generated by 
    $$\theta_{j_1}\cdots \theta_{j_k} \T_{S^\prime} ( g^{-1} \langle x_{i_1},x_{i_2},\cdots, x_{i_{c+1}}\rangle g ),$$
    where $\{i_1,i_2,\cdots,i_{c+1}\}\subset {\bf n}$, $g \in \mathsf{F}_n$, $\emptyset \neq S^\prime \subset {\bf n}$ and $j_1,\cdots,j_k$ are distinct elements of ${\bf n} \backslash S^\prime$.
    Let $v = \langle x_{i_1},x_{i_2},\cdots, x_{i_{c+1}}\rangle$.
    By Lemma \ref{202509081651}, we have:
    \begin{align*}
        \T_{S^\prime} ( g^{-1} v g ) = \sum_{S^\prime = S_1 \amalg S_2 \amalg S_3} \T_{S_1} (g^{-1} ) \T_{S_3} (v) \T_{S_2} (g ) .
    \end{align*}
    Applying Lemma \ref{202509151210} to $\T_{S_3} (v)$, it follows that, if $\{ i_1, \cdots, i_{c+1} \} \subset S^\prime$ and $i_1, \cdots, i_{c+1}$ are pairwise distinct, then we obtain
    $$
    \T_{S^\prime} ( g^{-1} v g ) = \sum_{S^{\prime\prime} = S_1 \amalg S_2} \T_{S_1} (g^{-1} ) [\theta_{i_1},\theta_{i_2},\cdots ,\theta_{i_{c+1}}] \T_{S_2} (g) ,
    $$
    where $S^{\prime\prime} = S^\prime \backslash \{ i_1, \cdots, i_{c+1} \}$; otherwise, $\T_{S^\prime} ( g^{-1} v g )$ vanishes.
    In any case, we obtain $\T_{S^\prime} ( g^{-1} v g ) \in \IL_{c} (S)_{\mathds{1}}$.
    Thus, $\K_{\gamma_{c+1}} (n) \subset \IL_{c} (S)_{\mathds{1}}$.
    
    We now prove that $\IL_{c} (S)_{\mathds{1}} \subset \K_{\gamma_{c+1}} (n)$.
    By definition, the $\mathds{k}$-module $\IL_{c} (S)_{\mathds{1}}$ is generated by $w [\theta_{i_1},\theta_{i_2},\cdots ,\theta_{i_{l}}] w^\prime$ for $l >c$.
    Here, $w \in \mathds{k} \langle S_1 \rangle_{\mathds{1}}$, $w^\prime \in \mathds{k} \langle S_3 \rangle_{\mathds{1}}$ and $[\theta_{i_1},\theta_{i_2},\cdots ,\theta_{i_{l}}] \in L (S_2)_{\mathds{1}}$, where $S = S_1 \amalg S_2 \amalg S_3$ and $i_1, \cdots, i_l$ are distinct elements of $S_2$.
    Hence, it suffices to prove $w [\theta_{i_1},\theta_{i_2},\cdots ,\theta_{i_{l}}] w^\prime \in \K_{\gamma_{c+1}} (n)$.
    By iteratively applying the formula $[\theta_{i},\cdots ,\theta_{j}] \theta_k = [\theta_{i},\cdots ,\theta_{j} , \theta_k]+ \theta_k [\theta_{i},\cdots ,\theta_{j}]$, we may start with the case where $w^\prime = 1$, in particular $S_3 = \emptyset$.
    By applying the general result in Lemma \ref{202509101818} to $\mathcal{R}= \gamma_{c+1}$, we obtain $w [\theta_{i_1},\theta_{i_2},\cdots ,\theta_{i_{l}}] \in \K_{\gamma_{c+1}} (n)$, since $\T_{S_2}(\langle x_{i_1},x_{i_2},\cdots ,x_{i_{l}} \rangle)= [\theta_{i_1},\theta_{i_2},\cdots ,\theta_{i_{l}}]$ by Lemma \ref{202509151210}.
\end{proof}

We now establish Lemma \ref{202604021122}.
Recall from Section \ref{202408191340} the proof of Lemma \ref{202604011516} together with notations.

\begin{proof}[Proof of Lemma \ref{202604021122}]
    We shall prove that the PBW datum given in the proof of Lemma \ref{202604011516} provides the desired result.
    First, it is convenient to extend the notation $\K_{\gamma_{c+1}}$ to finite subsets $X$ of $\mathds{N}^\ast$.
    For such $X$, let $\K_{\gamma_{c+1}} (X)$ be the kernel of $\mathsf{Lin} (X) \to \Psi_{\mathbf{N}_c^{\mathsf{o}}} (X)$.
    Lemma \ref{202604031247} implies $\IL_c (X^\prime)_{\mathds{1}_{X^\prime}} = \K_{\gamma_{c+1}}(X)$, where $X^\prime = \{ \theta_i \in \Theta \mid i \in X \}$.
    By Lemma \ref{202509161922}, we obtain an isomorphism $$(\mathds{k} \langle H(\Theta) , \preceq \rangle_{\max >c})_{\mathds{1}_{X^\prime}} \cong K_{\gamma_{c+1}}(X) .$$
    By definition, $(\mathds{k} \langle H(\Theta) , \preceq \rangle_{\max >c})_{\mathds{1}_{X^\prime}}$ is freely generated by the words $t_1 \cdots t_r, ~ r \in \mathds{N}^\ast$ where $t_i \in H_{\delta_i} (\Theta)$ and $\delta_i \in V_\Theta$ such that $t_1 \succ \cdots \succ t_r$, $\mathds{1}_{X^\prime} = \sum^r_{i=1} \delta_i$ and $|\delta_j| >c$ for some $j$.
    Under the above identification, such $t_i$ corresponds to $g(t_i) = h_i \in \mathfrak{B}_{X_i}$ where $X_i$ is the support of $\delta_i$.
    Therefore, the $\mathds{k}$-submodule $\K_{\gamma_{c+1}} (X) \subset \mathsf{Lin}(X)$ is freely generated by $h_1 \cdots h_r,~r \in \mathds{N}^\ast$ where $h_j \in \mathfrak{B}_{X_j}$ and $X = \coprod^r_{j=1} X_j$ such that $h_1 \succ \cdots \succ h_r$ and $|X_j| > c$ for some $j$.
\end{proof}

\section*{Acknowledgements}

The author is supported by a KIAS Individual Grant MG093701 at Korea Institute for Advanced Study.
The author would like to express his sincere gratitude to Geoffrey Powell, Christine Vespa and Aurélien Djament for carefully reading the earlier versions of this paper and providing him valuable comments.
In particular, he thanks Aurélien Djament for bringing the reference \cite{PolyPira} to his attention.
The author also thanks many people he encountered at the International Category Theory Conference CT2024.
In particular, he thanks Rui Prezado for suggesting the reference \cite{betti1983variation}.
Finally, the author is deeply grateful to the anonymous referee for a careful reading and many insightful suggestions.

\bibliography{reference}
\bibliographystyle{alpha}
\end{document}